\documentclass[pdflatex,sn-nature]{sn-jnl}

\usepackage[utf8]{inputenc}
\usepackage[T1]{fontenc}
\usepackage{graphicx}%
\usepackage{multirow}%
\usepackage{amsmath,amssymb,amsfonts}%
\usepackage{booktabs}
\usepackage{float}
\usepackage{amsthm}%
\usepackage{mathrsfs}%
\usepackage[title]{appendix}%
\usepackage{xcolor}%
\usepackage{textcomp}%
\usepackage{manyfoot}%
\usepackage{booktabs}%
\usepackage{algorithm}%
\usepackage{algorithmicx}%
\usepackage{algpseudocode}%
\usepackage{listings}%
\usepackage{array}%
\usepackage{float}%
\usepackage{tikz}%
\usetikzlibrary{}
\graphicspath{{Figure/}{media/}}
\definecolor{navy}{RGB}{0,45,90}
\begin{document}

\title[OptiXDE]{OptiXDE: A fast optical-inspired solver for differential equations}

\author*[1]{\fnm{Yang} \sur{Yang}}\email{yangyang\_kmy@powerchina.cn}
\author*[2]{\fnm{Mingjiao} \sur{Yan}}\email{mingjiaoyan@hhu.edu.cn}
\author[1]{\fnm{Zongliang} \sur{Zhang}}

\affil*[1]{\orgname{PowerChina Kunming Engineering Corporation Limited}, \orgaddress{\city{Kunming}, \postcode{650051}, \state{Yunnan}, \country{China}}}
\affil[2]{\orgdiv{College of Water Conservancy and Hydropower Engineering}, \orgname{Hohai University}, \orgaddress{\city{Nanjing}, \postcode{210098}, \state{Jiangsu}, \country{China}}}

\abstract{OptiXDE is a matrix-free spectral operator framework for differential equations on uniform grids and embedded domains. Inspired by angular-spectrum propagation in Fourier optics, it maps transform-diagonal spatial operators to analytical modal multipliers and composes them with physical-space operators for nonlinearities, geometry and boundary enforcement. A common transform--operator--inverse-transform backbone is demonstrated across transient diffusion, periodic and embedded-domain Poisson problems, the cubic nonlinear Schr"odinger equation, viscous Burgers dynamics, the two-dimensional Allen--Cahn equation and incompressible flows from the Taylor--Green vortex to embedded-cylinder vortex shedding. Transform-compatible linear problems are recovered near the floating-point limit, whereas errors on the singular L-shaped domain remain localized near the re-entrant corner and regularized interface. Nonlinear benchmarks recover second-order temporal convergence and the expected conservative or dissipative behavior, while incompressibility remains near round-off level during long-time vortex shedding. The matrix-free updates require \(\mathcal{O}(N\log N)\) work and \(\mathcal{O}(N)\) memory. Device-resident transform workloads reach \(94.9\times\) GPU acceleration, and the complete embedded-cylinder solver achieves a \(42.1\times\) CPU--GPU speedup under matched numerical settings. These results establish OptiXDE as a deterministic and extensible operator-centric framework for structured and embedded-domain differential equations.}


\maketitle

\section*{Introduction}
\label{sec:introduction}
When light passes through a Fourier lens, a complex wavefront is decomposed into elementary spatial-frequency components. Each component can be filtered, phase-shifted or amplified, and an inverse transform reconstructs the resulting field. This transform--modify--reconstruct principle provides an intuitive optical interpretation of spectral computation. In the proposed OptiXDE framework, differential operators that are diagonalizable under Fourier, sine or cosine transforms are represented by their modal symbols, so that transient propagation or steady inversion can be evaluated modewise in transform space. Nonlinear terms, heterogeneous coefficients, irregular geometries and non-periodic boundaries are incorporated through operator splitting, pseudo-spectral evaluation, physical-space enforcement or iterative correction. The resulting computation therefore combines an analytically tractable spectral core with physical-space operators for components that do not admit the same closed-form modal treatment.

Partial differential equations underpin a broad range of physical and engineering models, and conventional numerical methods such as the finite element method (FEM) \citep{zienkiewicz2005finite,hughes2012finite}, finite volume method (FVM) \citep{versteeg2007introduction} and finite difference method (FDM) \citep{leveque2007finite} remain the principal computational tools. Their generality is accompanied by substantial numerical infrastructure: implicit, steady-state and coupled formulations commonly lead to large algebraic systems whose assembly, solution and preconditioning may dominate high-resolution simulations \citep{bathe2006finite,saad2003iterative,benzi2002preconditioning}, while explicit formulations may instead be constrained by stability-limited time increments. Complex geometries can additionally require mesh generation, element-quality control or adaptive refinement \citep{thompson1998handbook,freitag1997tetrahedral}. These limitations do not diminish the versatility of local discretization methods, but motivate complementary formulations that reduce global matrix operations when the governing operators possess exploitable transform structure.

Spectral and pseudo-spectral methods provide such an alternative by representing differential operators through global basis functions \citep{boyd2001chebyshev,trefethen2000spectral,canuto2007spectral,shen2011spectral}. The fast Fourier transform reduces Fourier-transform cost to $\mathcal{O}(N\log N)$ \citep{cooley1965algorithm,frigo2005design}, enabling fast Poisson solvers \citep{adams1999fast,buzbee1970fast}, wave and quantum-dynamical propagation \citep{kosloff1983fourier,feit1982solution,fornberg1998practical}, periodic micromechanics \citep{moulinec1998numerical,michel2001computational} and, more recently, massively parallel Fourier pseudo-spectral simulation on accelerator-based systems \citep{yeung2025gpu}. Analytical or semi-analytical modal propagation is likewise established in exponential time-differencing and integrating-factor schemes \citep{cox2002exponential,kassam2005fourth}, with recent high-order exponential spectral formulations extending this principle to multidimensional nonlinear parabolic problems with non-periodic boundary conditions \citep{wang2025exponential}. Split-step Fourier methods similarly separate transform-diagonal linear evolution from nonlinear physical-space dynamics \citep{feit1982solution}, while related $k$-space formulations derive equation-specific spectral correction factors directly from dispersion relations \citep{treeby2018nonstandard}. From a physical perspective, the angular-spectrum method in Fourier optics provides a particularly transparent interpretation of this mathematics: a field is decomposed into spatial-frequency components, each component is propagated through a transfer function, and the field is reconstructed by an inverse transform \citep{goodman2005introduction,voelz2011computational}.

The principal challenge is therefore not the Fourier representation itself, nor the isolated treatment of non-periodic boundaries or complex geometries, for which a substantial body of spectral methodology already exists. Fourier embedded-boundary formulations have combined regular transform grids with transient spectral evolution to solve Poisson and Laplace problems on irregular domains \citep{sabetghadam2009fourier}, while Fourier-continuation methods have enabled high-order elliptic, parabolic, hyperbolic and nonlinear-flow calculations on general domains, including variable-coefficient problems \citep{bruno2010fcad,albin2011fc,bruno2014fcad}. Immersed Boundary Smooth Extension methods similarly retain Fourier spectral discretizations on embedding grids while recovering high-order accuracy for elliptic and parabolic equations and, subsequently, incompressible Navier--Stokes flow on arbitrary smooth domains \citep{stein2016immersed,stein2017ibse}. On Cartesian domains, sine and cosine transforms provide direct treatment of non-periodic boundaries \citep{wise2021pstd}, and recent work has extended fast trigonometric-transform formulations to Fourier interaction-picture propagation and heterogeneous transient diffusion under combinations of Dirichlet and Neumann conditions \citep{hatharasinghe2025nonperiodic,sanoko2025fft}. These developments demonstrate that transform-based PDE computation can accommodate substantially broader settings than idealized periodic problems. At the same time, they rely on distinct choices of continuation, basis construction, temporal integration, interface treatment, constraint enforcement and computational backend, while recent FFT-based interface-enrichment studies further show that local interface accuracy can remain limiting even when the bulk spectral solution is highly accurate \citep{gehrig2025xfft}. The remaining opportunity is therefore to organize these established spectral capabilities around a reusable operator abstraction in which transform-diagonal propagation or inversion, physical-space interactions, geometry and constraint enforcement, and computational backends can be composed within a common matrix-free architecture.

Data-driven approaches provide a different route to operator-level computation. Physics-informed neural networks \citep{raissi2019physics}, DeepXDE \citep{lu2021deepxde}, DeepONet \citep{lu2021learning} and Fourier neural operators \citep{li2021fourier,kovachki2023neural} enable forward, inverse and parametric modelling using neural representations or learned mappings between function spaces. Once trained, neural operators can provide rapid inference across families of parameterized problems, but their predictive reliability depends on training-data coverage, optimization, architecture and out-of-distribution behaviour \citep{karniadakis2021physics,guo2022review}. Deterministic solvers therefore remain important when reproducibility, convergence control, physical interpretability and instance-specific accuracy are required. This distinction also motivates operator-centric deterministic frameworks in which the governing differential operator, rather than a learned surrogate, remains the primary computational object.

Against this background, we present OptiXDE, an optical-inspired spectral operator framework that organizes transform-domain differentiation, analytical modal propagation or steady inversion, physical-space correction and computational backends within a common matrix-free architecture. The contribution is therefore not a new Fourier identity, but a unified computational organization of established spectral ingredients around reusable differential-equation operators. For linear constant-coefficient components that are diagonalizable under Fourier, sine or cosine transforms, modal evolution or inversion is evaluated analytically using cached spectral operators. Nonlinear interactions, heterogeneous coefficients, embedded geometries and boundary corrections are composed with this spectral core rather than being assumed to share the same closed-form representation. We assess the resulting framework across parabolic, elliptic, dispersive, nonlinear and incompressible-flow problems, from transform-compatible analytical benchmarks to singular embedded geometries and long-time vortex shedding. The study separates errors arising from spectral representation, temporal splitting, nonlinear aliasing, finite precision and embedded-boundary enforcement, and examines both algorithmic complexity and CPU/GPU performance. These results establish the regimes in which the transform--operator--inverse-transform architecture retains its analytical and computational advantages, while also identifying the cases in which nonlinear or geometric corrections become the dominant source of error and cost.

\section*{Results}
\label{sec:results}
\begin{figure}[h]
    \centering
    \includegraphics[width=\linewidth]{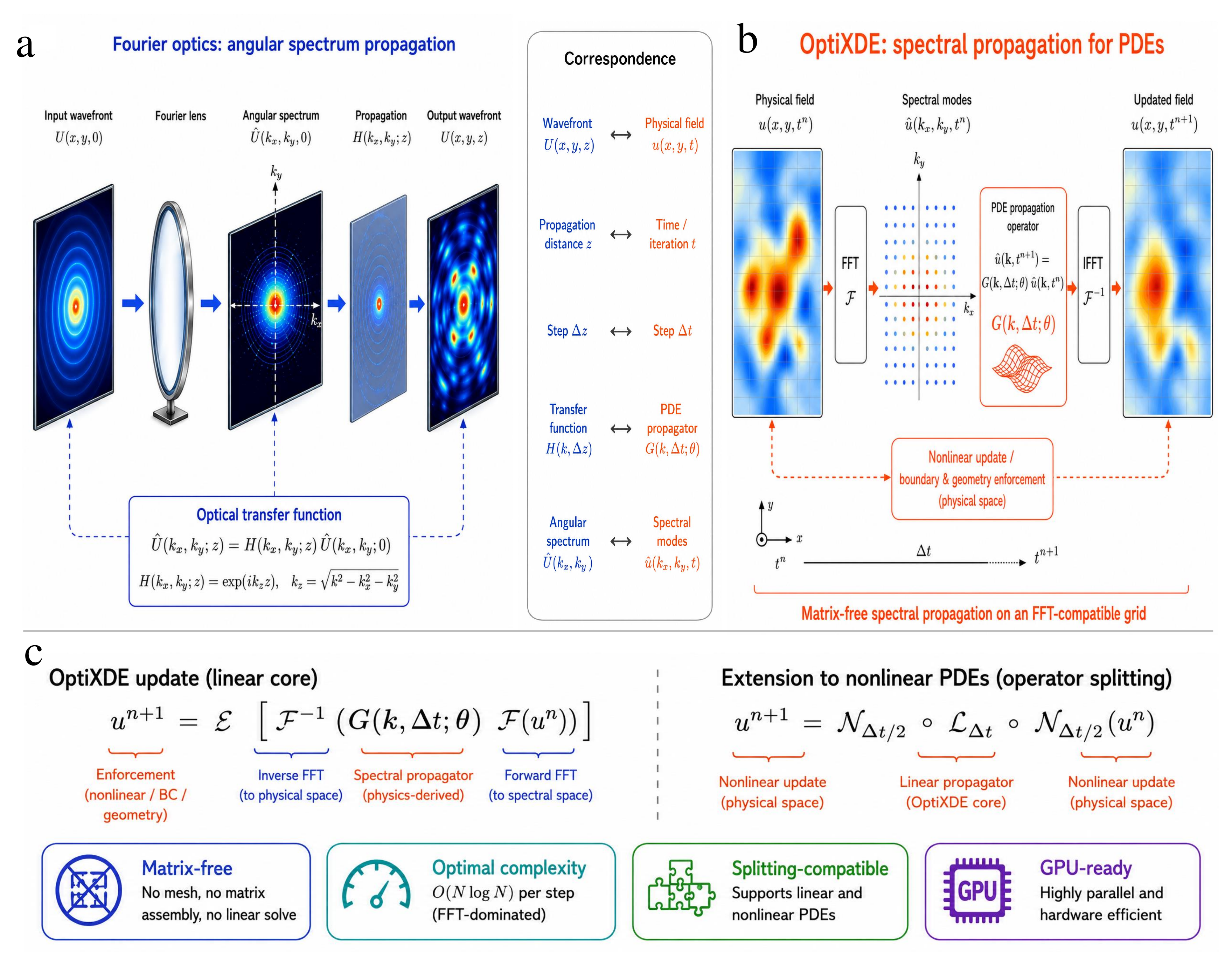}
    \caption{\textbf{Optical interpretation of the OptiXDE framework.} \textbf{a}, Fourier-optics angular-spectrum propagation decomposes an input wavefront into spectral components, propagates them through an optical transfer function, and reconstructs the output wavefront. \textbf{b}, OptiXDE maps the same sequence to differential-equation solving: a physical field is transformed into spectral modes, advanced by a PDE propagator, and transformed back. Boundary, geometry, and nonlinear effects are incorporated through physical-space enforcement or operator splitting. \textbf{c}, The resulting update combines a matrix-free linear spectral core with a splitting-compatible nonlinear extension.}
    \label{fig:optical_interpretation}
\end{figure}

As illustrated in Fig.~\ref{fig:optical_interpretation}, OptiXDE transfers the angular-spectrum viewpoint of Fourier optics to differential-equation solving. An optical wavefront is decomposed into spectral components, propagated by a transfer function, and reconstructed in physical space. In OptiXDE, the physical field is transformed into spectral modes, advanced by a physics-derived PDE propagator, and transformed back. Geometry, boundary conditions, and nonlinear terms are introduced as physical-space enforcement or splitting operators without changing the transform--propagate--inverse-transform core.

\subsection*{Round-off-limited accuracy of linear spectral operators}
\label{sec:results_linear}

\begin{figure}[t]
\centering
\includegraphics[width=\linewidth]{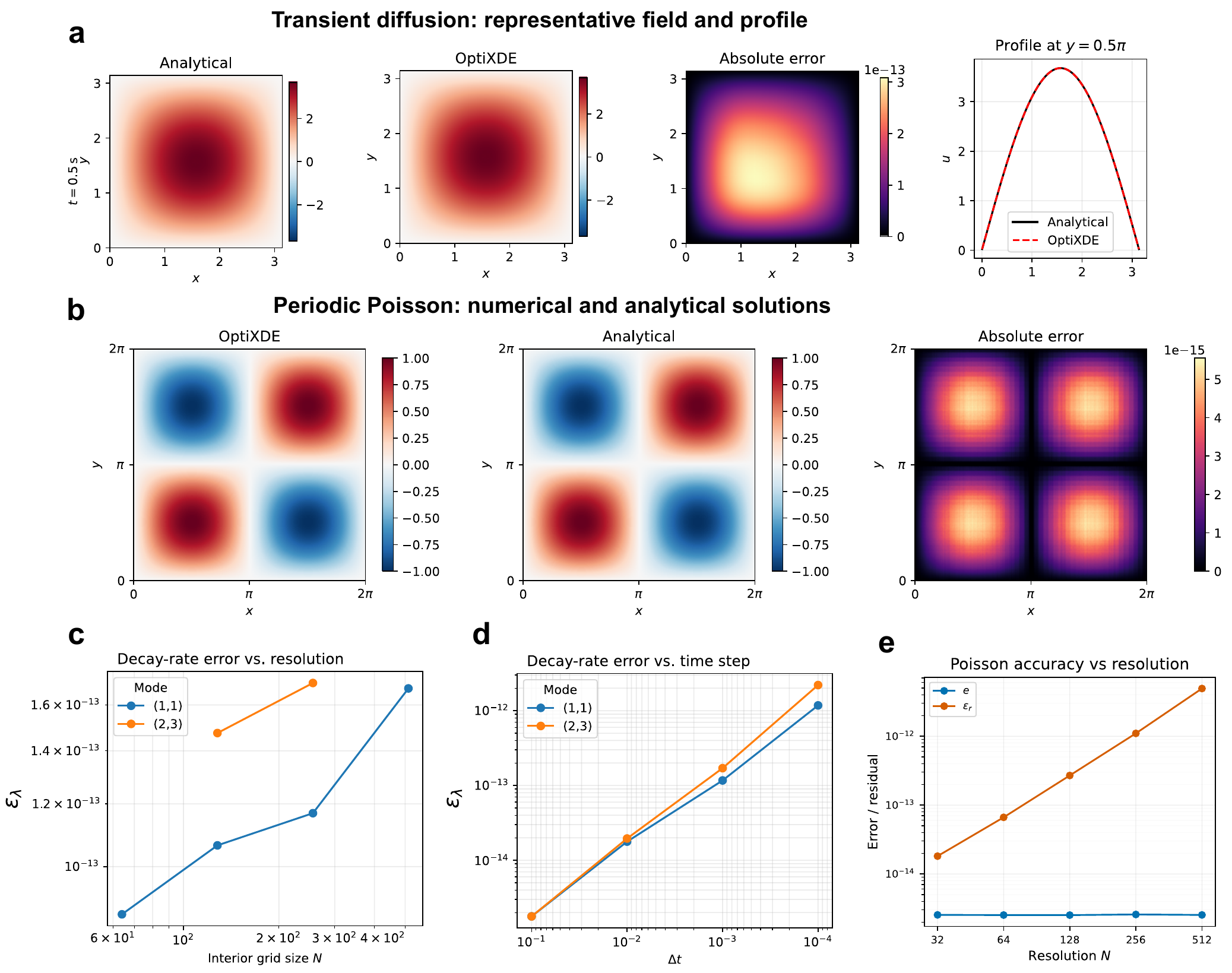}
\caption{\textbf{Round-off-limited accuracy of linear spectral operators.} \textbf{a}, Transient-diffusion benchmark showing the analytical field, OptiXDE solution, pointwise absolute error and central profile comparison. \textbf{b}, Periodic Poisson benchmark showing the OptiXDE solution, analytical solution and pointwise absolute error. \textbf{c}, Diffusion decay-rate error versus spatial resolution for the $(1,1)$ and $(2,3)$ Laplacian eigenmodes. \textbf{d}, Diffusion decay-rate error versus propagation interval $\Delta t$, showing the accumulation of finite-precision error under repeated transform--propagate cycles. \textbf{e}, Periodic Poisson solution error and reconstructed algebraic residual versus spatial resolution. The solution error remains at the floating-point floor, whereas the residual increases with transform size because high-wave-number differentiation amplifies weak round-off components.}
\label{fig:linear_spectral_accuracy}
\end{figure}

We first isolated the accuracy of the closed-form spectral propagator using source-free diffusion on $\Omega=[0,\pi]^2$ with homogeneous Dirichlet boundaries. Two Laplacian eigenmodes, $(k_x,k_y)=(1,1)$ and $(2,3)$, were selected to probe smooth low-frequency and more rapidly varying high-frequency dynamics, with analytical decay rates of $2$ and $13$, respectively. Because both modes are represented exactly by the discrete sine basis, this benchmark separates propagation accuracy from spatial-approximation error. The numerical fields reproduced the analytical modal structures without visible phase displacement or systematic amplitude distortion (Fig.~\ref{fig:linear_spectral_accuracy}a). At $\Delta t=10^{-3}$, the terminal error of the $(1,1)$ mode remained between $3.32\times10^{-14}$ and $6.31\times10^{-14}$ as the spatial resolution increased from $64^2$ to $512^2$, while the inferred decay-rate error remained $\mathcal{O}(10^{-13})$ (Fig.~\ref{fig:linear_spectral_accuracy}c). The more strongly attenuated $(2,3)$ mode exhibited comparable decay-rate accuracy, confirming consistent reproduction of mode-dependent exponential decay across distinct spatial frequencies (Supplementary Figs.~S1 and S2 and Supplementary Table~S3).

Reducing the propagation interval did not improve the diffusion solution because each linear modal update is already evaluated analytically. Decreasing $\Delta t$ from $10^{-1}$ to $10^{-4}$ increased the number of transform--propagate--inverse-transform cycles from $15$ to $15{,}000$, while the terminal error of the $(1,1)$ mode increased from $7.15\times10^{-16}$ to $4.42\times10^{-13}$. The corresponding decay-rate errors increase as the time increment decreases (Fig.~\ref{fig:linear_spectral_accuracy}d; Supplementary Fig.~S2 and Supplementary Table~S4). This behavior is opposite to conventional temporal-convergence trends and reflects the accumulation of finite-precision errors through repeated transforms rather than temporal truncation, indicating that transform-resolved diffusion modes are limited primarily by floating-point arithmetic.

We next tested the steady counterpart using a periodic Poisson problem on $\Omega=[0,2\pi)^2$ with the manufactured solution $u_{\mathrm{ex}}(x,y)=\sin(x)\sin(y)$. The solution is represented exactly by the Fourier basis, allowing the direct inverse-Laplacian operator to be assessed without embedded-boundary correction or outer iteration. Across spatial resolutions from $64^2$ to $512^2$, the discrete $L_2$ error remained between $2.56\times10^{-15}$ and $2.59\times10^{-15}$, with no systematic resolution dependence (Fig.~\ref{fig:linear_spectral_accuracy}e). At $512^2$, the pointwise absolute error remained below approximately $6\times10^{-15}$ throughout the domain, and the numerical and analytical fields were visually indistinguishable (Fig.~\ref{fig:linear_spectral_accuracy}b; Supplementary Fig.~S3 and Supplementary Table~S6). The reconstructed algebraic residual increased from $4.59\times10^{-14}$ to $4.40\times10^{-12}$ with increasing transform size (Fig.~\ref{fig:linear_spectral_accuracy}e), consistent with amplification of finite-precision noise in weak high-frequency coefficients during residual evaluation rather than deterioration of the resolved solution. Together, the diffusion and Poisson benchmarks show that the transient propagator and steady inverse operator recover transform-compatible linear problems at round-off-limited accuracy.

\subsection*{Embedded-domain accuracy on a singular L-shaped geometry}
\label{sec:results_lshape}

\begin{figure}[t]
\centering
\includegraphics[width=\linewidth]{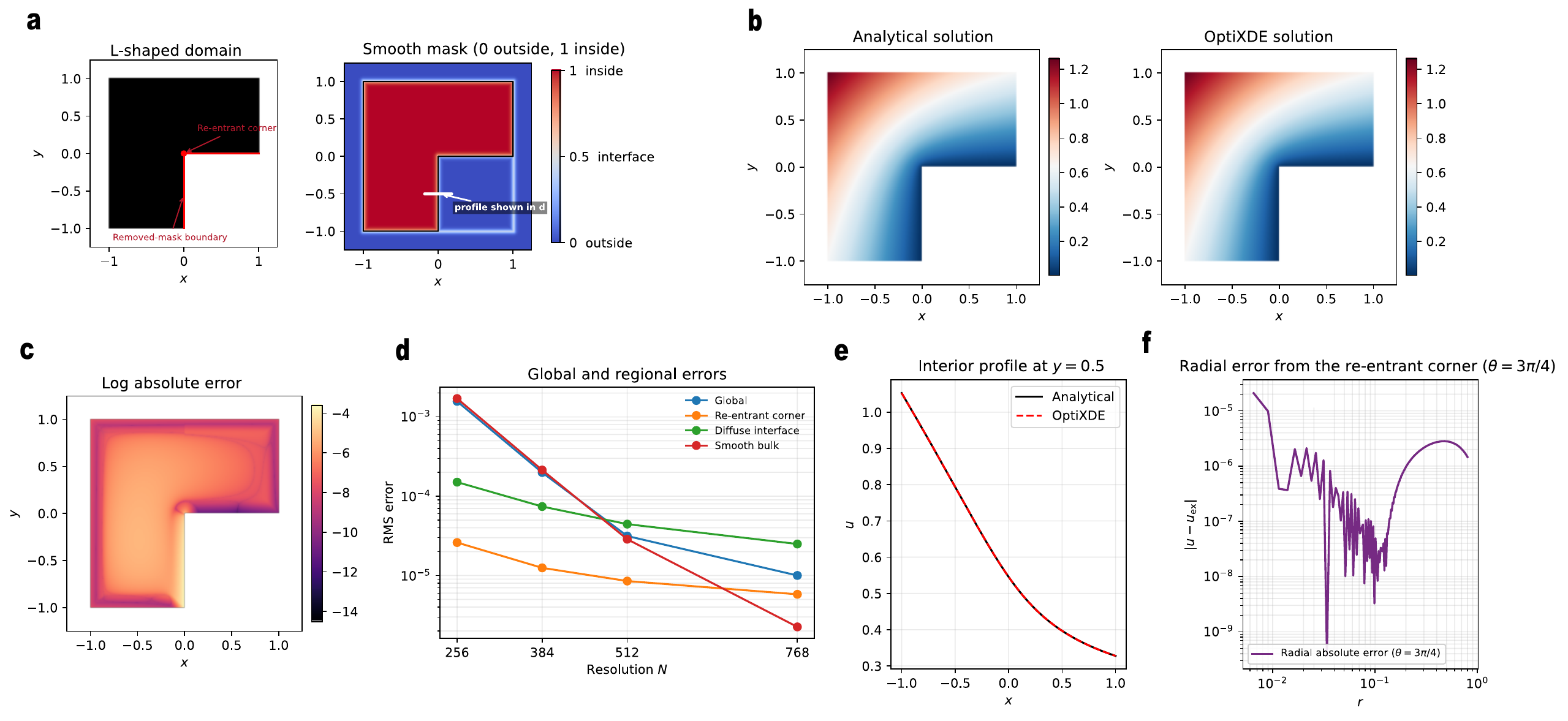}
\caption{\textbf{Embedded-domain accuracy on a singular L-shaped geometry.} \textbf{a}, L-shaped physical domain with the re-entrant corner and regularized mask used for embedded-domain representation. \textbf{b}, Analytical and OptiXDE solution fields at \(768^2\) resolution. \textbf{c}, Base-10 logarithm of the pointwise absolute error, showing localization near the embedded boundary and re-entrant corner. \textbf{d}, Global and regional root-mean-square errors as functions of spatial resolution, separating the re-entrant-corner, diffuse-interface and smooth-bulk contributions. \textbf{e}, Interior profile along \(y=0.5\), showing visual overlap between the analytical and OptiXDE solutions. \textbf{f}, Pointwise absolute error along the radial path from the re-entrant corner at \(\theta=3\pi/4\), revealing a multiscale and nonmonotonic spatial error structure.}
\label{fig:lshape_poisson}
\end{figure}

We next examined whether the spectral formulation remains accurate when the physical domain is not compatible with a global transform basis. A Poisson problem was solved on an L-shaped domain with a re-entrant corner of interior angle \(3\pi/2\) and singular analytical solution \(u_{\mathrm{ex}}(r,\theta)=r^{2/3}\sin(2\theta/3)\), for which \(|\nabla u_{\mathrm{ex}}|\sim r^{-1/3}\) near the corner. This benchmark simultaneously introduces nonperiodic Dirichlet boundaries, geometric embedding and reduced solution regularity. The physical domain was represented inside a rectangular transform domain through a regularized mask, while the constant-coefficient bulk operator remained diagonal in transform space (Fig.~\ref{fig:lshape_poisson}a). This separation retains the global spectral operator while enforcing irregular geometry and boundary data locally in physical space.

Despite the corner singularity and diffuse embedded interface, the solution improved systematically with spatial resolution. The global root-mean-square error decreased from \(1.5710\times10^{-3}\) at \(256^2\) to \(1.0007\times10^{-5}\) at \(768^2\), while the smooth-bulk error decreased from \(1.7056\times10^{-3}\) to \(2.2466\times10^{-6}\) (Fig.~\ref{fig:lshape_poisson}d). Over the same refinement range, the re-entrant-corner error decreased from \(2.6032\times10^{-5}\) to \(5.7914\times10^{-6}\), whereas the diffuse-interface error decreased from \(1.5051\times10^{-4}\) to \(2.4958\times10^{-5}\) and became the dominant localized contribution at the highest resolution. Because the physical interface thickness was fixed at \(\varepsilon=0.03\) while the penalty coefficient scaled as \(\eta=0.003h^2\), these results quantify systematic resolution improvement at fixed interface regularization rather than sharp-interface asymptotic convergence.

At \(768^2\), the OptiXDE and analytical fields are visually indistinguishable over most of the domain (Fig.~\ref{fig:lshape_poisson}b), whereas the logarithmic error field reveals strong localization near the embedded boundary and re-entrant corner (Fig.~\ref{fig:lshape_poisson}c). The maximum pointwise discrepancy is approximately \(2.4\times10^{-4}\) near the inner vertical boundary; away from these localized regions, the interior profile at \(y=0.5\) overlaps visually with the analytical solution, while the radial corner error spans several orders of magnitude and exhibits a multiscale, nonmonotonic structure (Fig.~\ref{fig:lshape_poisson}e,f). Controlled parameter variations further distinguish these errors from floating-point transform effects: reducing the penalty coefficient decreases the finite-penalty error, insufficient mask smoothing introduces high-frequency contamination into the bulk solution, and increasing the pseudo-time increment accelerates convergence only until spatial accuracy begins to deteriorate (Supplementary Fig.~S25). The production choices \(\eta/h^2=0.003\), \(\varepsilon=0.03\) and \(\Delta\tau/h^2=10\) therefore provide a stable accuracy--regularization--cost compromise rather than a narrowly tuned optimum. Together, these results show that once the smooth interior is resolved, the dominant accuracy limitation shifts from the bulk spectral representation to localized singularity and geometric enforcement, while the transform-based operator remains effective on a non-smooth physical domain that is not aligned with the global spectral basis.

\subsection*{Nonlinear dynamics through operator composition}
\label{subsec:nonlinear_dynamics}

\begin{figure}[h]
\centering
\includegraphics[width=\linewidth]{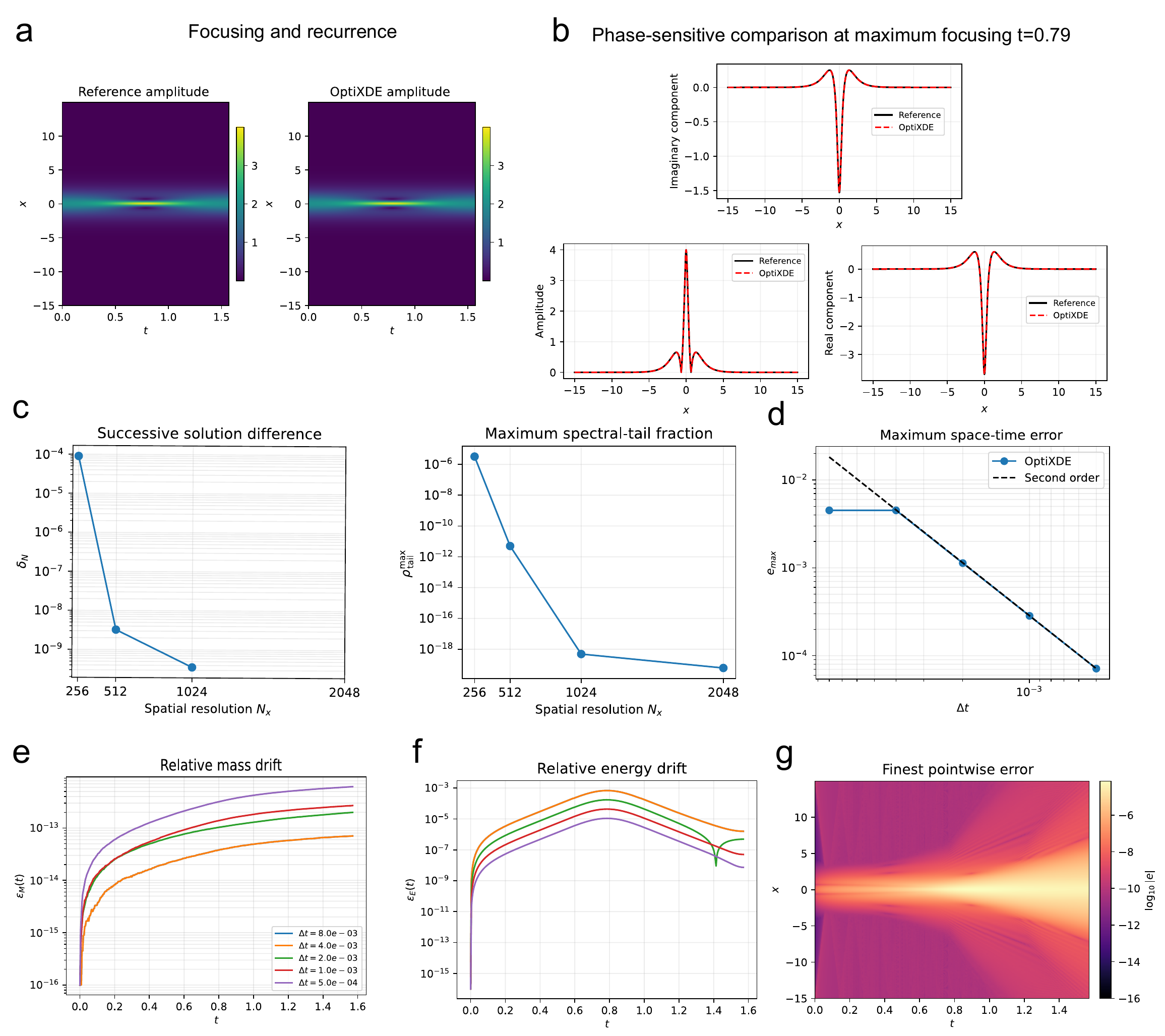}
\caption{\textbf{Conservative nonlinear dynamics in the focusing nonlinear Schr\"odinger equation.} \textbf{a}, Reference and OptiXDE amplitude fields over the complete focusing--recurrence cycle. \textbf{b}, Phase-sensitive comparison near maximum focusing (\(t\approx0.79\simeq\pi/4\)), showing the amplitude and the real and imaginary components of the reference and OptiXDE solutions. \textbf{c}, Spatial-resolution diagnostics showing the successive solution difference \(\delta_N\) and maximum high-wave-number spectral-tail fraction \(\rho_{\mathrm{tail}}^{\max}\). \textbf{d}, Maximum space--time complex-field error \(e_{\max}\) versus time increment \(\Delta t\), together with a second-order reference slope. \textbf{e}, Relative mass drift, demonstrating near-round-off preservation of the discrete mass. \textbf{f}, Relative Hamiltonian-energy drift, showing systematic reduction of the bounded splitting error under temporal refinement. \textbf{g}, Logarithmic pointwise complex-error distribution for the finest calculation, showing error localization through the focusing--recurrence cycle.}
\label{fig:nls_composition}
\end{figure}

Having isolated the accuracy of the linear spectral operators and embedded-domain treatment, we next examined whether spectral propagation remains accurate when repeatedly composed with nonlinear evolution. The focusing cubic nonlinear Schr\"odinger equation provides a conservative, complex-valued test in which dispersion competes with an amplitude-dependent phase rotation. OptiXDE reproduces the complete focusing--recurrence cycle together with the real and imaginary components through the maximum-focusing state (Fig.~\ref{fig:nls_composition}a,b). Spatial refinement reduces the successive solution difference from \(8.99\times10^{-5}\) at \(N_x=256\) to \(3.57\times10^{-10}\) at \(N_x=1024\), while the maximum high-wave-number spectral-tail fraction falls from \(3.24\times10^{-6}\) to below \(10^{-18}\) at the two finest resolutions (Fig.~\ref{fig:nls_composition}c). At \(N_x=2048\), successive halving of the time increment reduces the maximum space--time complex-field error from \(1.100375\times10^{-2}\) to \(4.319465\times10^{-5}\), with observed orders of \(1.995\), \(1.999\), \(2.000\) and \(2.000\) (Fig.~\ref{fig:nls_composition}d,g). The maximum relative mass drift remains between \(5.33\times10^{-14}\) and \(8.49\times10^{-13}\), whereas the bounded Hamiltonian excursion decreases from \(1.68\times10^{-3}\) to \(6.58\times10^{-6}\) under temporal refinement (Fig.~\ref{fig:nls_composition}e,f). Explicit removal of the upper spectral band changes the resolved trajectory by only \(2.46\times10^{-10}\), far below the discretization error (Supplementary Fig.~S26), confirming that the focusing dynamics are not sustained by unresolved spectral content.

\begin{figure}[h]
\centering
\includegraphics[width=\linewidth]{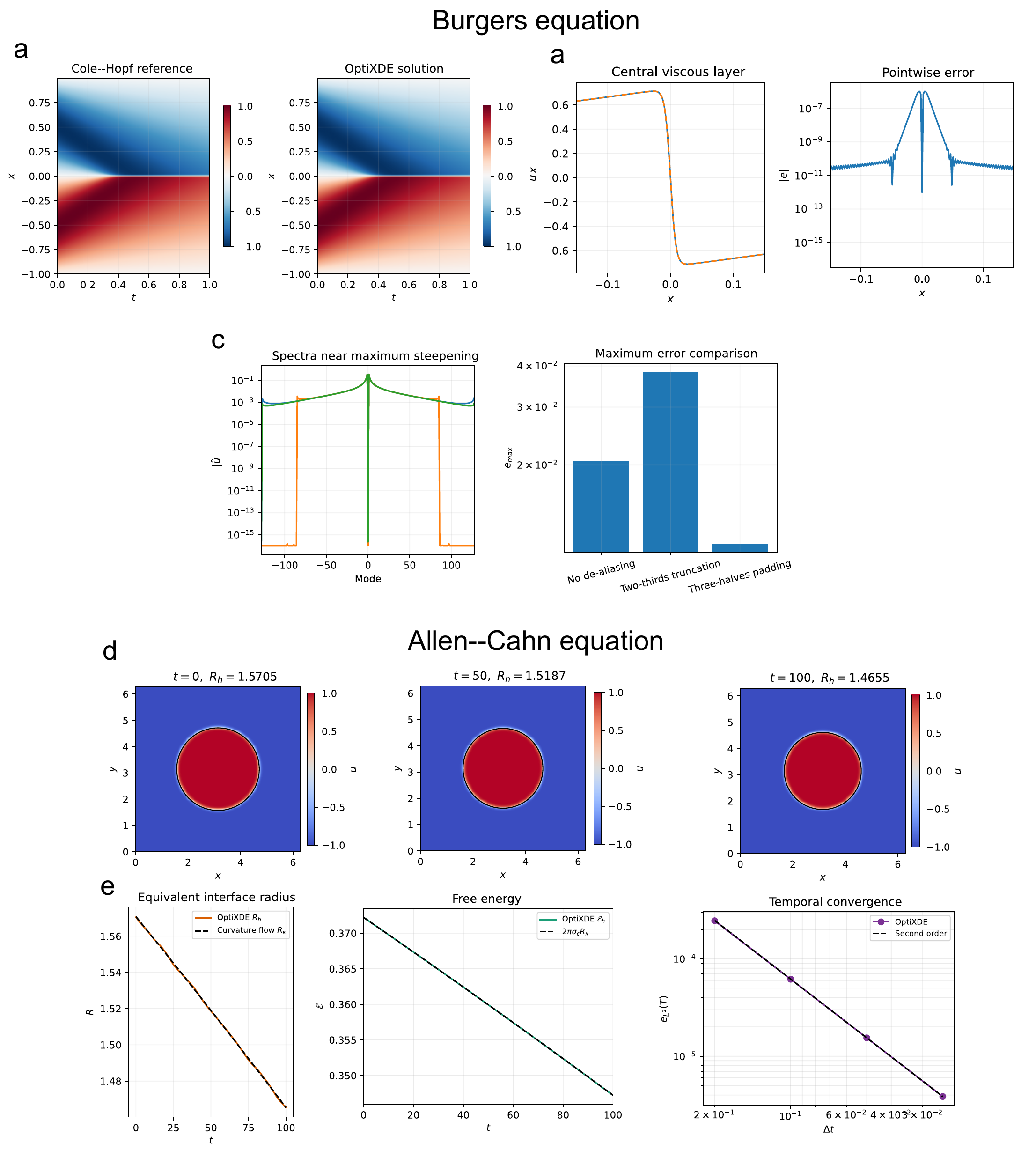}
\caption{\textbf{Nonlinear transport and dissipative phase-field dynamics.} \textbf{a}, Cole--Hopf reference and OptiXDE space--time solutions for the viscous Burgers equation, showing consistent reproduction of the global steepening and relaxation process. \textbf{b}, Enlarged comparison of the central viscous layer together with the corresponding pointwise error, showing that the remaining discrepancy is strongly localized near the narrow internal layer around \(x=0\). \textbf{c}, Effect of nonlinear de-aliasing in the Burgers calculation, showing the spectrum near maximum steepening and the corresponding maximum-error comparison for no de-aliasing, two-thirds truncation and three-halves padding. \textbf{d}, Allen--Cahn phase-field distributions at \(t=0\), \(50\) and \(100\), showing curvature-driven contraction of the diffuse circular interface. \textbf{e}, Equivalent interface radius, free-energy evolution and temporal-convergence results for the Allen--Cahn equation. The numerical radius closely follows the leading-order curvature-flow prediction, the discrete free energy decreases monotonically, and the terminal relative \(L_2\) error exhibits second-order decay under time-step refinement.}
\label{fig:nonlinear_transport_dissipation}
\end{figure}

We then considered the viscous Burgers equation, where nonlinear transport generates a narrow internal layer and transfers energy toward progressively higher spatial frequencies. The independently evaluated Cole--Hopf reference was verified to substantially higher accuracy than the OptiXDE solutions (Supplementary Fig.~S11 and Supplementary Table~S13). Increasing the spatial resolution from \(N_x=128\) to \(1024\) reduces the maximum space--time error from \(8.93\times10^{-2}\) to \(4.59\times10^{-6}\) and the terminal error from \(9.10\times10^{-3}\) to \(8.51\times10^{-8}\), with the remaining discrepancy localized around the central viscous layer (Fig.~\ref{fig:nonlinear_transport_dissipation}a,b; Supplementary Figs.~S12 and S14). At fixed \(N_x=1024\), temporal refinement recovers second-order convergence, with observed orders between \(1.992\) and \(1.998\) (Supplementary Fig.~S13). The nonlinear convolution treatment is also decisive: three-halves padding gives a maximum error of \(1.16\times10^{-2}\), compared with \(2.07\times10^{-2}\) without de-aliasing and \(3.84\times10^{-2}\) with two-thirds truncation (Fig.~\ref{fig:nonlinear_transport_dissipation}c). This distinction shows that suppressing high-wave-number content is not equivalent to removing aliasing, because the resolved high-frequency modes are required to represent the viscous layer. Kinetic-energy decay, homogeneous boundary values and odd symmetry are simultaneously preserved to numerical precision across the refinement study (Supplementary Fig.~S27), providing complementary physical-consistency checks.

Finally, we tested dissipative nonlinear dynamics using the two-dimensional Allen--Cahn equation, in which a nonconserved phase field relaxes toward \(u=\pm1\) while a diffuse interface contracts under curvature. An independently generated Fourier pseudo-spectral ETDRK4 solution provided the reference for the convergence studies (Supplementary Table~S17). Starting from a circular positive-phase region with \(\epsilon=0.04\), OptiXDE predicts a decrease in equivalent radius from \(1.57049\) to \(1.46551\) over \(T=100\), with a maximum relative discrepancy of \(6.25\times10^{-4}\) from the leading-order curvature-flow prediction (Fig.~\ref{fig:nonlinear_transport_dissipation}d,e). The discrete free energy decreases monotonically from \(0.372206\) to \(0.347221\), while the interface remains nearly circular and the phase bounds are preserved to within \(1.47\times10^{-7}\) (Fig.~\ref{fig:nonlinear_transport_dissipation}e; Supplementary Fig.~S16). At \(256^2\) resolution, reducing \(\Delta t\) from \(2.0\times10^{-1}\) to \(2.5\times10^{-2}\) decreases the terminal relative \(L_2\) error from \(2.45\times10^{-4}\) to \(3.86\times10^{-6}\), with approximately second-order convergence (Fig.~\ref{fig:nonlinear_transport_dissipation}e). Large-step tests further separate numerical boundedness from physical fidelity: monotonic free-energy dissipation is retained for \(\Delta t\le0.5\), whereas an energy increase first appears at \(\Delta t=1\) despite the field remaining bounded (Supplementary Fig.~S28). Together, the Schr\"odinger, Burgers and Allen--Cahn benchmarks show that the same operator-composition framework spans conservative dispersive waves, nonlinear transport and dissipative phase-field evolution while retaining the characteristic numerical and physical structure of each regime.

\subsection*{Incompressible flow and embedded vortex shedding}
\label{subsec:incompressible_flow}

\begin{figure}[t]
\centering
\includegraphics[width=\linewidth]{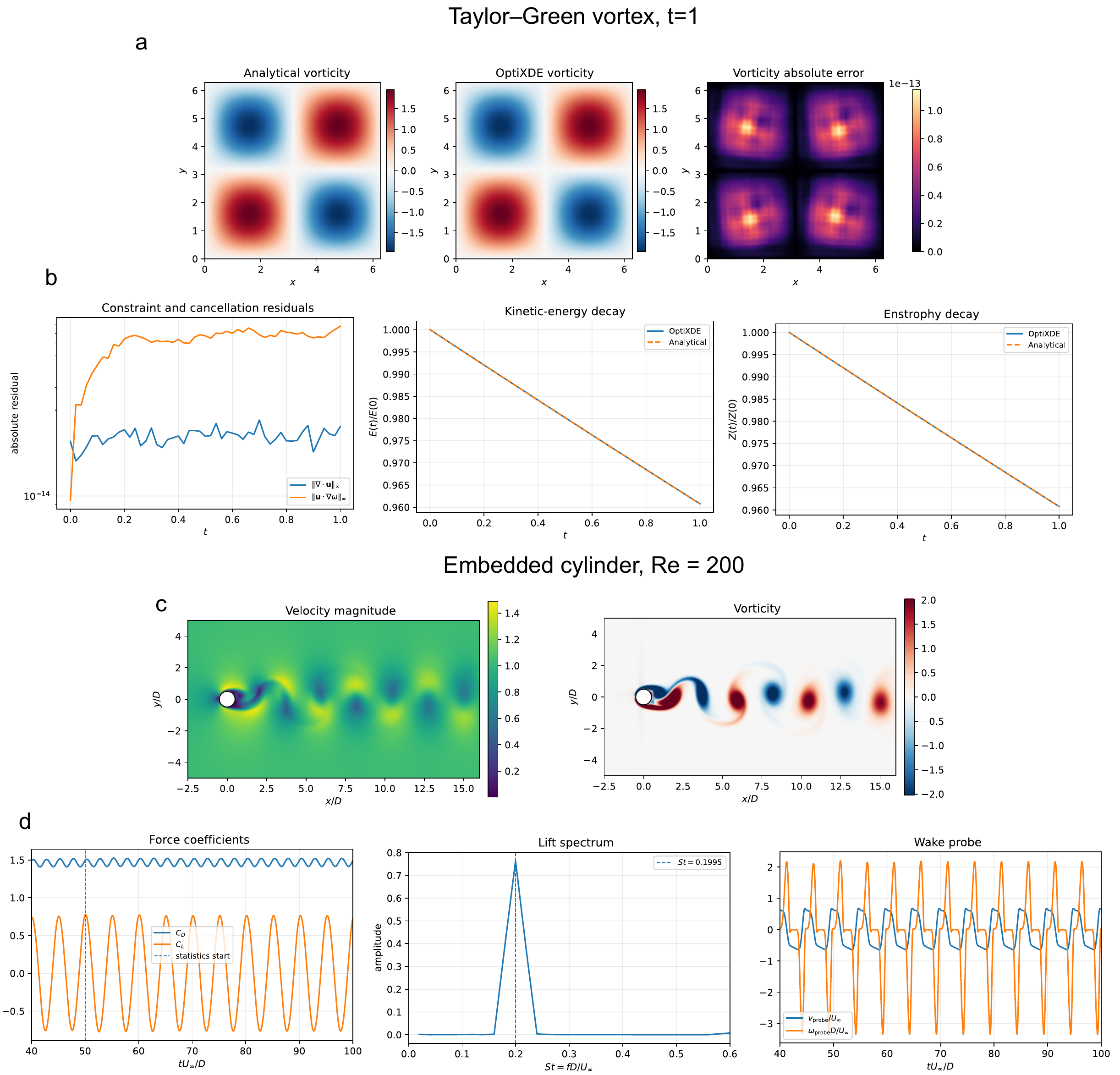}
\caption{\textbf{Incompressible flow from analytical vortices to embedded vortex shedding.} \textbf{a}, Taylor--Green vortex verification at \(t=1\), showing the analytical vorticity, OptiXDE solution and pointwise absolute vorticity error. \textbf{b}, Temporal physical diagnostics for the Taylor--Green vortex, including the incompressibility and nonlinear-cancellation residuals together with the normalized kinetic-energy and enstrophy decay. \textbf{c}, Instantaneous velocity-magnitude and vorticity fields for flow past an embedded circular cylinder at \(Re=200\) and \(t=100\), showing the downstream velocity deficit and sustained alternating vortex street. \textbf{d}, Force and shedding-frequency diagnostics for the cylinder calculation, including the drag and lift histories, lift-coefficient spectrum and downstream wake-probe signal. The independently extracted shedding frequencies give consistent Strouhal numbers, \(St=0.19955\) from the lift spectrum and \(St=0.19957\) from the wake probe.}
\label{fig:incompressible_flow}
\end{figure}

We next examined whether the spectral-operator framework extends from scalar nonlinear equations to constrained vector dynamics. The two-dimensional Taylor--Green vortex provides a strict analytical verification of the periodic incompressible-flow formulation, including vorticity transport, velocity recovery, pressure reconstruction and the divergence-free constraint. At \(128^2\) resolution with \(\Delta t=0.01\), the relative terminal errors in vorticity and the two velocity components are \(3.8340\times10^{-14}\), \(3.7865\times10^{-14}\) and \(3.7853\times10^{-14}\), respectively, while the independently reconstructed pressure has a relative error of \(7.3426\times10^{-14}\) (Fig.~\ref{fig:incompressible_flow}a). The corresponding incompressibility residual is \(2.4314\times10^{-14}\), and the relative kinetic-energy and enstrophy errors are \(7.5306\times10^{-14}\) and \(7.4744\times10^{-14}\). Spatial resolutions from \(32^2\) to \(256^2\) and time increments from \(0.04\) to \(0.005\) remain within or close to the same double-precision plateau (Supplementary Figs.~S18--S20), consistent with exact viscous propagation of the resolved Taylor--Green mode rather than conventional spatial or temporal truncation.

We then moved to a nonlinear flow without a closed-form solution by considering vortex shedding past an embedded circular cylinder at \(Re=200\). The cylinder is represented by a compact smooth Brinkman mask within a \(40D\times20D\) Fourier domain, while a downstream fringe region prevents the periodic wake from contaminating the inflow. The calculation develops a sustained alternating vortex street, accompanied by the expected downstream velocity deficit and unsteady pressure distribution (Fig.~\ref{fig:incompressible_flow}b). Increasing the resolution from \(N_D=24\) to \(40\) points per cylinder diameter reduces the successive changes in the principal wake statistics: the mean drag changes by \(1.47\%\) from \(N_D=24\) to \(32\) and by only \(0.55\%\) from \(32\) to \(40\), while the corresponding changes in the lift root-mean-square decrease from \(3.32\%\) to \(1.11\%\). At \(N_D=40\), the statistically stationary interval \(50\le t\le100\) gives \(\overline{C}_D=1.46716\) and \(C_{L,\mathrm{rms}}=0.54049\). The dominant shedding frequency obtained from the lift spectrum gives \(St=0.19955\), in close agreement with the independently extracted wake-probe value \(St=0.19957\) (Fig.~\ref{fig:incompressible_flow}c,d). Over the complete calculation, the Fourier projection maintains \(\max_t\|\nabla\cdot\mathbf{u}\|_\infty=9.79\times10^{-14}\), showing that the embedded-boundary treatment and long-time nonlinear wake evolution do not compromise the incompressibility constraint.

Controlled variations of the embedded-cylinder parameters further distinguish physical sensitivity from numerical instability. Strengthening the Brinkman penalty from \(\eta=0.01\) to \(0.0025\) reduces the solid-region velocity residual from \(7.17\times10^{-2}\) to \(2.44\times10^{-2}\), while the associated Strouhal number varies by less than \(0.45\%\). Reducing the mask transition half-width from \(h\) to \(0.75h\) changes the mean drag by only \(0.21\%\), whereas broadening it to \(1.5h\) increases the drag by \(4.30\%\), demonstrating that an excessively diffuse mask alters the effective hydrodynamic geometry. Varying the transverse domain height from \(16D\) to \(24D\) changes the mean drag by less than \(0.4\%\) relative to the \(20D\) production domain and the Strouhal number by less than \(0.2\%\) (Supplementary Fig.~S29). Together, the Taylor--Green and cylinder calculations show that the same spectral framework can preserve incompressibility at near-round-off accuracy while progressing from an analytically resolved periodic vortex to sustained nonlinear shedding around an embedded solid boundary.

\subsection*{Workload-dependent hardware acceleration}
\label{subsec:performance}

\begin{figure}[h]
\centering
\includegraphics[width=\linewidth]{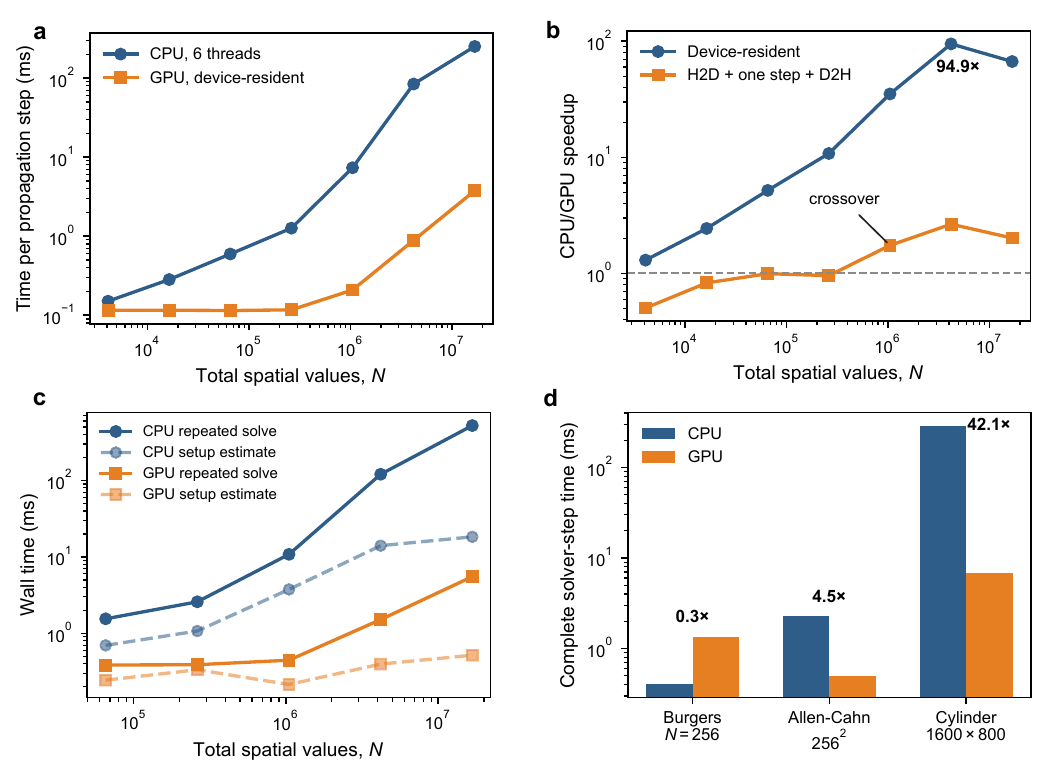}
\caption{\textbf{Workload-dependent computational acceleration of OptiXDE.} \textbf{a}, CPU and device-resident GPU execution time per propagation step for the two-dimensional diffusion problem over resolutions from \(64^2\) to \(4096^2\). \textbf{b}, Corresponding CPU/GPU speedup for device-resident execution and for a transfer-inclusive operation consisting of one host-to-device transfer, one propagation step and one device-to-host transfer. Device-resident acceleration reaches \(94.9\times\) at \(2048^2\), whereas transfer costs delay the practical crossover to larger grids. \textbf{c}, CPU and GPU repeated-solve times for the periodic Poisson problem with the inverse spectral operator cached; at \(4096^2\), the solve time decreases from \(522.46\) to \(5.522\) ms, corresponding to a \(94.6\times\) acceleration. \textbf{d}, Complete solver-step times for representative nonlinear applications using matched CPU and GPU numerical configurations. The small one-dimensional Burgers calculation remains accelerator-overhead dominated, whereas the \(256^2\) Allen--Cahn and \(1600\times800\) cylinder calculations achieve \(4.50\times\) and \(42.1\times\) GPU acceleration, respectively.}
\label{fig:computational_performance}
\end{figure}

We finally examined how the transform-based structure of OptiXDE translates into computational performance across problem scales. For the two-dimensional diffusion propagator, device-resident GPU acceleration increases rapidly as the transform workload grows, from only \(1.30\times\) at \(64^2\) to a maximum of \(94.9\times\) at \(2048^2\), where the median propagation time decreases from \(84.33\) ms on six CPU threads to \(0.889\) ms on the A100 GPU (Fig.~\ref{fig:computational_performance}a,b). At \(4096^2\), the GPU remains \(66.9\times\) faster despite the onset of stronger memory and backend effects. Including one host-to-device transfer and one device-to-host transfer substantially delays the practical crossover: CPU and GPU execution are approximately equal at \(256^2\), whereas a clear transfer-inclusive advantage emerges from \(1024^2\) and reaches \(2.66\times\) at \(2048^2\). The normalized execution cost approaches the transform-dominated regime at large resolutions, consistent with the \(\mathcal{O}(N\log N)\) complexity of the underlying spectral updates, while the small-grid GPU plateau reflects fixed launch and dispatch overheads rather than algorithmic scaling.

The same workload dependence is evident at the application level. A cached periodic Poisson solve at \(4096^2\) decreases from \(522.46\) ms on the CPU to \(5.522\) ms on the GPU, corresponding to a \(94.6\times\) acceleration. By contrast, the one-dimensional Burgers calculation at \(N=256\) is too small to amortize accelerator overhead and runs at only \(0.30\times\) the CPU performance. Increasing the workload to the two-dimensional Allen--Cahn problem at \(256^2\) produces a \(4.50\times\) acceleration, while the matched \(1600\times800\) \(Re=200\) cylinder calculation reduces the complete solver-step time from \(289.77\) to \(6.886\) ms, giving a \(42.1\times\) GPU speedup with a CPU--GPU terminal relative \(L_2\) difference of \(1.19\times10^{-13}\) (Fig.~\ref{fig:computational_performance}c,d). These results show that the computational advantage of OptiXDE is not associated with accelerator use alone, but emerges when sufficiently large spectral workloads expose the parallelism of the transform-based operators. Detailed timing, shared-memory scaling, de-aliasing overhead and peak-memory measurements are provided in Supplementary Figs.~S30 and S31 and Supplementary Table~S24.

\begin{table}[h]
\centering
\small
\caption{\textbf{Representative CPU--GPU performance of OptiXDE.} Timings correspond to matched numerical configurations. Diffusion reports device-resident propagation, Poisson reports cached repeated inversion, and the nonlinear cases report complete solver-step times.}
\label{tab:performance_summary}
\begin{tabular}{lcccc}
\toprule
Benchmark & Resolution & CPU & GPU & Speedup \\
\midrule
Diffusion & \(2048^2\) & \(84.328\) ms & \(0.889\) ms & \(94.9\times\) \\
Periodic Poisson & \(4096^2\) & \(522.459\) ms & \(5.522\) ms & \(94.6\times\) \\
Burgers & \(N=256\) & \(0.402\) ms & \(1.343\) ms & \(0.30\times\) \\
Allen--Cahn & \(256^2\) & \(2.267\) ms & \(0.504\) ms & \(4.50\times\) \\
Cylinder, \(Re=200\) & \(1600\times800\) & \(289.765\) ms & \(6.886\) ms & \(42.1\times\) \\
\bottomrule
\end{tabular}
\end{table}

\section*{Discussion}
\label{sec:discussion}

The results identify two distinct accuracy regimes within the OptiXDE framework. For linear constant-coefficient problems on transform-compatible domains, analytical modal propagation or inversion reduces the numerical error to a regime dominated by finite-precision arithmetic rather than conventional spatial or temporal truncation. The diffusion and periodic Poisson benchmarks demonstrate this behaviour directly. This round-off-limited regime is consistent with recent non-periodic Fourier propagation algorithms based on sine and cosine transforms, which likewise report machine-precision accuracy for a transform-compatible heat-equation benchmark \citep{hatharasinghe2025nonperiodic}. When irregular geometry is introduced, however, the dominant error shifts from the spectral bulk operator to the geometric representation and boundary treatment. The L-shaped problem shows that high accuracy can still be retained in the smooth interior, while the remaining discrepancy becomes localized near the re-entrant singularity and regularized embedded interface. Related interface-accuracy limitations have recently motivated enriched FFT formulations for thermal homogenization, where local solution fields obtained on regular grids lose accuracy near nonconforming material interfaces \citep{gehrig2025xfft}. Although the interface in that setting is a material interface rather than an external geometric boundary, both observations emphasize that once the bulk spectral operator is sufficiently resolved, interface representation can become the dominant source of local error. Geometry is therefore not eliminated from the numerical problem, but its influence is transferred from mesh construction and element quality to interface resolution, mask regularization and boundary enforcement.

The nonlinear and incompressible-flow benchmarks show that the analytical spectral propagator is best interpreted as a reusable computational backbone rather than a complete solver in isolation. Conservative and dissipative dynamics can be constructed by composing transform-space propagation with physical-space nonlinear evolution, while nonlinear transport additionally requires appropriate de-aliasing, conservative evaluation and time-step control. This composition follows the broader principles underlying split-step, integrating-factor and exponential spectral methods \citep{feit1982solution,cox2002exponential,kassam2005fourth}; the distinction in OptiXDE is their organization within a common operator architecture that also incorporates embedded geometry and constraint enforcement. The Taylor--Green vortex and embedded-cylinder calculations extend this operator composition to constrained vector dynamics, combining spectral Poisson inversion, differentiation, nonlinear transport, geometry enforcement and incompressibility projection within the same framework. These results suggest that the central reusable element of OptiXDE is not any individual propagation formula, but the decomposition of a differential equation into transform-diagonal propagation, physical-space interactions and constraint operators.

The computational results further show that the practical meaning of ``fast'' is workload dependent. Small problems may remain dominated by transform dispatch and accelerator overhead and can therefore be more efficient on the CPU, whereas sufficiently large two-dimensional workloads expose the parallel structure of the transform-based operators and benefit substantially from GPU execution. Device residency, transform size, the number of transforms required by each update and host--device communication all influence the realized acceleration. Recent exascale Fourier pseudo-spectral simulations similarly demonstrate that large-scale accelerator performance depends not only on FFT throughput but also on device-resident computation and communication-aware global transforms \citep{yeung2025gpu}. The observed OptiXDE performance therefore complements, rather than follows automatically from, the \(\mathcal{O}(N\log N)\) complexity of the spectral core. At still larger distributed-memory scales, communication associated with global transforms may become a limiting factor and will require dedicated decomposition and communication strategies.

The present results also define the current scope of OptiXDE. The framework is most effective when a dominant component of the governing operator is constant coefficient and diagonalizable under a Fourier, sine or cosine representation. Strongly heterogeneous coefficients, discontinuous material properties, shock-dominated solutions and highly localized geometric features introduce modal coupling or non-smooth structure that cannot generally be represented by a single fixed spectral multiplier, consistent with the established sensitivity of global spectral approximations to non-smooth fields \citep{boyd2001chebyshev,trefethen2000spectral}. Such problems require iterative correction, operator splitting, local regularization or hybrid global--local discretizations, and may not retain the near-round-off behaviour observed for transform-compatible linear problems. OptiXDE should therefore be regarded as a matrix-free spectral propagation architecture that complements, rather than universally replaces, finite-element, finite-volume and other locally adaptive methods.

Taken together, these results establish transform-diagonal spectral propagation as a reusable computational backbone for deterministic differential-equation solving, while clarifying where its accuracy and efficiency are controlled by nonlinear coupling, geometry and hardware scale. A natural next step is to extend the present operator abstraction from equation-specific implementations to a programmable PDE framework in which users specify governing equations, state variables, physical parameters, geometry and initial and boundary conditions at the mathematical level. The framework could then identify transform-compatible linear components automatically, construct the corresponding spectral symbols, propagators or steady inverse operators, and map nonlinear terms, sources and constraints to physical-space operators assembled through splitting or iterative composition. Automatic selection of transform bases, null-mode treatment and de-aliasing strategies would further reduce equation-specific implementation. Such a PDE-specification-to-operator pathway would extend OptiXDE from a set of implemented solvers toward a programmable spectral system that translates mathematical differential-equation descriptions into executable operator compositions, while retaining the deterministic, interpretable and matrix-free structure demonstrated here.

\section*{Methods}
\label{sec:methods}

OptiXDE separates a transform-diagonal bulk operator from physical-space operators that represent nonlinear interactions, geometry, boundary conditions and other constraints. The closed-form spectral propagators themselves follow established Fourier and spectral formulations; the framework contribution is to organize propagation, steady inversion, physical-space enforcement and constraint projection within a common matrix-free computational architecture. Transform space is used wherever the governing operator can be diagonalized analytically, whereas mode-coupling contributions are evaluated in physical space and composed with the spectral update. The resulting workflow is summarized in Fig.~\ref{fig:optical_interpretation}. Detailed transform conventions, benchmark-specific operators and reproducibility settings are given in Supplementary Section~1.

\subsection*{Spectral propagation and steady inversion}
\label{subsec:method_spectral}

Consider a linear evolution equation
\begin{equation}
\frac{\partial u}{\partial t}
=
\mathcal{L}u+s,
\label{eq:method_linear_evolution}
\end{equation}
where $\mathcal{L}$ is a constant-coefficient spatial differential operator and $s$ is a prescribed source. For a Fourier direction, spatial derivatives become analytical multipliers,
\begin{equation}
\widehat{\frac{\partial^p u}{\partial x_j^p}}
=
(\mathrm{i}k_j)^p\widehat{u},
\qquad
\widehat{\nabla^2u}
=
-|\mathbf{k}|^2\widehat{u},
\qquad
|\mathbf{k}|^2
=
\sum_{j=1}^{d}k_j^2.
\label{eq:method_spectral_derivatives}
\end{equation}
A transform basis that diagonalizes $\mathcal{L}$ therefore converts Eq.~\eqref{eq:method_linear_evolution} into independent modal equations,
\begin{equation}
\frac{\partial\widehat{u}(\mathbf{k},t)}{\partial t}
=
\lambda(\mathbf{k})\widehat{u}(\mathbf{k},t)
+
\widehat{s}(\mathbf{k},t),
\label{eq:method_modal_equation}
\end{equation}
where $\lambda(\mathbf{k})$ is the spectral symbol of $\mathcal{L}$.

The exact variation-of-constants representation over one propagation interval is
\begin{equation}
\widehat{u}^{\,n+1}
=
G(\mathbf{k},\Delta t)\widehat{u}^{\,n}
+
\int_{0}^{\Delta t}
\exp\!\left[
\lambda(\mathbf{k})(\Delta t-\tau)
\right]
\widehat{s}(\mathbf{k},t_n+\tau)\,
\mathrm{d}\tau,
\label{eq:method_variation_constants}
\end{equation}
with the homogeneous propagation operator
\begin{equation}
G(\mathbf{k},\Delta t)
=
\exp\!\left[
\lambda(\mathbf{k})\Delta t
\right].
\label{eq:method_modal_propagator}
\end{equation}
When the source is constant over the propagation interval, Eq.~\eqref{eq:method_variation_constants} reduces to
\begin{equation}
\widehat{u}^{\,n+1}
=
G(\mathbf{k},\Delta t)\widehat{u}^{\,n}
+
\Delta t\,
\varphi_1\!\left[
\lambda(\mathbf{k})\Delta t
\right]
\widehat{s}^{\,n},
\label{eq:method_source_update}
\end{equation}
where
\begin{equation}
\varphi_1(z)
=
\begin{cases}
\dfrac{\exp(z)-1}{z}, & z\neq0,\\[2mm]
1, & z=0.
\end{cases}
\label{eq:method_phi1}
\end{equation}
The limiting value in Eq.~\eqref{eq:method_phi1} provides the nonsingular update for zero-eigenvalue modes. For time-dependent sources, the temporal accuracy of the complete scheme additionally depends on the approximation used for the source integral in Eq.~\eqref{eq:method_variation_constants}.

For diffusion,
\begin{equation}
\lambda_{\mathrm{diff}}(\mathbf{k})
=
-D|\mathbf{k}|^2,
\qquad
G_{\mathrm{diff}}(\mathbf{k},\Delta t)
=
\exp\!\left(
-D|\mathbf{k}|^2\Delta t
\right).
\label{eq:method_diffusion_propagator}
\end{equation}
This multiplier has the same transform--multiply--inverse-transform structure as the transfer function in angular-spectrum propagation, although the diffusion propagator is dissipative rather than phase preserving \citep{goodman2005introduction,voelz2011computational}. Because the constant-coefficient linear subproblem is integrated analytically, its propagation is not subject to the explicit diffusion stability restriction. Restrictions may nevertheless arise from nonlinear transport, source approximation, splitting or physical-space enforcement.

Boundary-compatible transforms are selected independently in each separable Cartesian direction,
\begin{equation}
\mathcal{T}_{b}
=
\begin{cases}
\mathcal{F}, & b=\mathrm{periodic},\\
\mathcal{S}, & b=\mathrm{Dirichlet},\\
\mathcal{C}, & b=\mathrm{Neumann},
\end{cases}
\qquad
\mathcal{T}
=
\bigotimes_{j=1}^{d}
\mathcal{T}_{b_j},
\label{eq:method_transform_selection}
\end{equation}
where $\mathcal{F}$, $\mathcal{S}$ and $\mathcal{C}$ denote Fourier, sine and cosine transforms, respectively. Nonhomogeneous separable boundary data are reduced to homogeneous form through a lifting field. Mixed endpoint and null-mode conventions are described in Supplementary Sections~1.1--1.3.

A general linear OptiXDE propagation step can consequently be written as
\begin{equation}
\mathcal{P}_{\Delta t}[u]
=
\mathcal{T}^{-1}
\left[
G(\boldsymbol{\Lambda},\Delta t)
\odot
\mathcal{T}[u]
\right],
\label{eq:method_general_propagation}
\end{equation}
where $\boldsymbol{\Lambda}$ contains the spectral eigenvalues and $\odot$ denotes elementwise multiplication.

Steady constant-coefficient equations use the same transform infrastructure. For
\begin{equation}
\mathcal{L}u=f,
\label{eq:method_steady_equation}
\end{equation}
the transformed non-null modes satisfy
\begin{equation}
\widehat{u}(\mathbf{k})
=
\frac{\widehat{f}(\mathbf{k})}
{\lambda(\mathbf{k})},
\qquad
\lambda(\mathbf{k})\neq0.
\label{eq:method_steady_inversion}
\end{equation}
Null modes are handled through the corresponding solvability and normalization conditions. Thus, transient propagation and steady inversion are two instances of the same transform--operator--inverse-transform architecture.

\subsection*{Geometry and boundary enforcement}
\label{subsec:method_geometry}

OptiXDE distinguishes boundary conditions that can be incorporated directly through the transform basis from those requiring physical-space enforcement. Periodic, homogeneous Dirichlet and homogeneous Neumann conditions on separable Cartesian boundaries belong to the first category. Irregular geometries, nonseparable boundary partitions and embedded obstacles belong to the second and are treated without changing the transform-diagonal bulk operator.

An irregular physical domain $\Omega$ is embedded in a rectangular transform domain $\widetilde{\Omega}$,
\begin{equation}
\Omega
\subseteq
\widetilde{\Omega}
=
\prod_{j=1}^{d}[0,L_j),
\label{eq:method_embedded_domain}
\end{equation}
and represented through a signed-distance or level-set function,
\begin{equation}
\phi(\mathbf{x})<0
\quad\text{in }\Omega,
\qquad
\phi(\mathbf{x})=0
\quad\text{on }\partial\Omega,
\qquad
\phi(\mathbf{x})>0
\quad\text{in }\widetilde{\Omega}\setminus\Omega.
\label{eq:method_levelset}
\end{equation}
A regularized physical-domain mask and an interface-localization mask are defined by
\begin{equation}
m_{\Omega}(\mathbf{x})
=
H_{\varepsilon}\!\left[-\phi(\mathbf{x})\right],
\qquad
\chi_{\Gamma}(\mathbf{x})
=
4m_{\Omega}(\mathbf{x})
\left[
1-m_{\Omega}(\mathbf{x})
\right],
\label{eq:method_geometry_masks}
\end{equation}
where $H_{\varepsilon}$ is a regularized Heaviside function and $\varepsilon$ controls the interface width. The outward normal in the interface region is evaluated as
\begin{equation}
\mathbf{n}
=
\frac{\nabla\phi}
{\sqrt{\|\nabla\phi\|_2^2+\delta_n^2}},
\label{eq:method_boundary_normal}
\end{equation}
where $\delta_n>0$ regularizes the normalization. The detailed mask construction and geometric Boolean operations are given in Supplementary Section~1.4.

The central separation between spectral propagation and physical-space enforcement is expressed as
\begin{equation}
u^{n+1}
=
\mathcal{E}_{\Omega,\Gamma}
\left[
\mathcal{P}_{\Delta t}(u^n)
\right],
\label{eq:method_embedded_composition}
\end{equation}
where $\mathcal{P}_{\Delta t}$ advances the transform-compatible bulk equation and $\mathcal{E}_{\Omega,\Gamma}$ imposes geometry and boundary constraints. The spectral operator therefore remains independent of the detailed shape of $\Omega$.

Dirichlet, Neumann and Robin conditions can be represented by the general boundary operator
\begin{equation}
\mathcal{B}(u)
=
\alpha u
+
\beta\,
\mathbf{n}\cdot\mathbf{A}\nabla u
=
r
\qquad
\text{on }\Gamma_B,
\label{eq:method_general_boundary}
\end{equation}
where $\mathbf{A}$ is the coefficient tensor associated with the diffusive flux. Dirichlet, Neumann and Robin conditions correspond respectively to $\beta=0$, $\alpha=0$ and $\alpha\beta\neq0$.

For an intermediate field $v$, the corresponding boundary residual is
\begin{equation}
R_B(v)
=
\alpha v
+
\beta\,
\mathbf{n}\cdot\mathbf{A}\nabla v
-r.
\label{eq:method_boundary_residual}
\end{equation}
A localized residual correction is applied in the boundary band,
\begin{equation}
\mathcal{E}_{B}(v)
=
v
-
\omega_B\chi_B
\frac{R_B(v)}
{\sigma_B+\delta_B},
\label{eq:method_boundary_correction}
\end{equation}
where $0<\omega_B\leq1$ is a relaxation parameter, $\delta_B>0$ prevents division by zero and
\begin{equation}
\sigma_B
=
|\alpha|
+
|\beta|
\frac{a_n}{h},
\qquad
a_n
=
\mathbf{n}\cdot\mathbf{A}\mathbf{n},
\label{eq:method_boundary_scaling}
\end{equation}
provides a local scaling between value and flux residuals.

For embedded Dirichlet data $u=g$, the local penalty equation can instead be integrated exactly. The resulting relaxation is
\begin{equation}
\mathcal{E}_{D}[v]
=
g_{\mathrm{ext}}
+
\exp\!\left[
-\frac{\Delta\tau}{\eta_D}\chi_D
\right]
\left(
v-g_{\mathrm{ext}}
\right),
\label{eq:method_dirichlet_relaxation}
\end{equation}
where $g_{\mathrm{ext}}$ is an extension of the prescribed value, $\eta_D$ is the penalty parameter and $\Delta\tau$ is the enforcement interval. Exact local integration removes the explicit stability restriction associated with $\Delta\tau/\eta_D$, although $\eta_D$ and the interface resolution continue to influence boundary accuracy.

When the boundary contains several condition types or segments, the complete enforcement operation is composed from their individual maps,
\begin{equation}
\mathcal{E}_{\mathrm{mix}}
=
\mathcal{E}_{B_{N_b}}
\circ\cdots\circ
\mathcal{E}_{B_2}
\circ
\mathcal{E}_{B_1}.
\label{eq:method_mixed_boundary}
\end{equation}
For transient problems, this enforcement is composed with the physical time update. For steady embedded problems, the same propagation--enforcement composition is repeated as a pseudo-time or fixed-point iteration until convergence. Padding may additionally reduce interaction with periodic images but does not itself impose a physical boundary condition. Full penalty, mixed-boundary and convergence definitions are provided in Supplementary Section~1.5.

\subsection*{Nonlinear and incompressible equations}
\label{subsec:method_nonlinear}

For nonlinear problems, OptiXDE decomposes the governing equation into a transform-compatible linear part and a mode-coupling contribution,
\begin{equation}
\frac{\partial u}{\partial t}
=
\mathcal{L}u
+
\mathcal{N}
\left(
u,\nabla u,\mathbf{x},t
\right).
\label{eq:method_semilinear}
\end{equation}
Let $\mathcal{P}_{\tau}$ denote the closed-form flow generated by $\mathcal{L}$ and $\mathcal{Q}_{\tau}$ the nonlinear flow. The second-order Strang composition used in the nonlinear benchmarks is
\begin{equation}
u^{n+1}
=
\mathcal{P}_{\Delta t/2}
\circ
\mathcal{Q}_{\Delta t}
\circ
\mathcal{P}_{\Delta t/2}
\left(
u^n
\right).
\label{eq:method_strang}
\end{equation}

Nonlinear terms are evaluated pseudo-spectrally. The field and its derivatives are reconstructed from spectral coefficients,
\begin{equation}
u^n
=
\mathcal{F}^{-1}
\left[
\widehat{u}^{\,n}
\right],
\qquad
\nabla u^n
=
\mathcal{F}^{-1}
\left[
\mathrm{i}\mathbf{k}\widehat{u}^{\,n}
\right],
\label{eq:method_pseudospectral_derivatives}
\end{equation}
and nonlinear products are formed in physical space before transformation back,
\begin{equation}
\widehat{\mathcal{N}}^{\,n}
=
\mathcal{F}
\left[
\mathcal{N}
\left(
u^n,
\nabla u^n,
\mathbf{x},
t_n
\right)
\right].
\label{eq:method_pseudospectral_nonlinearity}
\end{equation}
Because physical-space multiplication corresponds to spectral convolution, nonlinear interactions can generate frequencies outside the resolved band. De-aliasing is therefore applied where required. The nonlinear Schr\"odinger and Allen--Cahn solvers use analytically solvable physical-space subflows, whereas the Burgers solver combines closed-form viscous propagation with conservative pseudo-spectral advection. Benchmark-specific subflows and de-aliasing rules are provided in Supplementary Sections~2.4--2.6.

The two-dimensional incompressible-flow formulation uses scalar vorticity $\omega$ and streamfunction $\psi$,
\begin{equation}
-\nabla^2\psi
=
\omega,
\qquad
\mathbf{u}
=
\overline{\mathbf{u}}
+
\left(
\frac{\partial\psi}{\partial y},
-\frac{\partial\psi}{\partial x}
\right),
\label{eq:method_streamfunction}
\end{equation}
where $\overline{\mathbf{u}}$ denotes an optional spatially uniform mean flow. The vorticity evolves according to
\begin{equation}
\frac{\partial\omega}{\partial t}
+
\mathbf{u}\cdot\nabla\omega
=
\nu\nabla^2\omega.
\label{eq:method_vorticity}
\end{equation}
The streamfunction is recovered directly in transform space,
\begin{equation}
\widehat{\psi}(\mathbf{k})
=
\frac{\widehat{\omega}(\mathbf{k})}
{|\mathbf{k}|^2},
\qquad
\mathbf{k}\neq\mathbf{0},
\label{eq:method_streamfunction_inversion}
\end{equation}
while spatial derivatives and viscous propagation remain diagonal spectral operations.

For embedded solid boundaries, the velocity is first relaxed locally toward the prescribed solid velocity using the same exponential enforcement principle as Eq.~\eqref{eq:method_dirichlet_relaxation}. Because physical-space enforcement can introduce a non-solenoidal component, the corrected velocity is subsequently projected onto the divergence-free Fourier subspace,
\begin{equation}
\widehat{\mathbf{u}}_{\perp}(\mathbf{k})
=
\left(
\mathbf{I}
-
\frac{
\mathbf{k}\mathbf{k}^{\mathsf T}
}{
|\mathbf{k}|^2
}
\right)
\widehat{\mathbf{u}}(\mathbf{k}),
\qquad
\mathbf{k}\neq\mathbf{0}.
\label{eq:method_incompressible_projection}
\end{equation}
The prescribed mean-flow component is retained separately at $\mathbf{k}=\mathbf{0}$. The incompressible solver therefore combines spectral inversion and propagation, physical-space nonlinear and geometric operations, and an exact transform-space constraint projection within the same operator-composition architecture. Full Taylor--Green and embedded-cylinder formulations are given in Supplementary Sections~2.7 and~2.8.

\subsection*{Implementation and computational cost}
\label{subsec:method_implementation}

OptiXDE is implemented as a matrix-free Python framework separating equation-specific solvers, spectral operators, geometry and boundary utilities, and transform backends. Periodic transforms are available through NumPy and optional PyTorch or CuPy backends, while compatible nonperiodic Cartesian directions use sine or cosine transforms. Solutions, wave-number arrays and cached operators remain on the selected device during repeated updates whenever the backend permits.

During preprocessing, OptiXDE constructs the spatial grid, transform basis, spectral eigenvalues and equation-specific propagation or inverse operators. Static embedded geometries additionally generate level-set fields, masks, normals and boundary partitions. Time-independent arrays are cached and reused. A subsequent update consists only of the transform operations, diagonal spectral multiplications, physical-space nonlinear or source evaluations, boundary enforcement and constraint projections required by the selected equation.

Let
\begin{equation}
N
=
\prod_{j=1}^{d}N_j
\label{eq:method_total_dofs}
\end{equation}
denote the total number of grid values and let $q$ be the fixed number of forward--inverse transform pairs required by one update. The principal computational costs are
\begin{equation}
C_{\mathrm{pre}}
=
\mathcal{O}(N),
\qquad
C_{\mathrm{step}}
=
\mathcal{O}\!\left(qN\log N\right)
+
\mathcal{O}(N),
\label{eq:method_step_complexity}
\end{equation}
where the second term represents pointwise propagation, nonlinear evaluation, masking and constraint operations. Consequently,
\begin{equation}
C_{\mathrm{transient}}
=
\mathcal{O}\!\left(
N_t qN\log N
\right),
\qquad
C_{\mathrm{embedded}}
=
\mathcal{O}\!\left(
n_{\mathrm{it}}qN\log N
\right),
\label{eq:method_total_complexity}
\end{equation}
for a transient calculation with $N_t$ updates and a steady embedded calculation requiring $n_{\mathrm{it}}$ outer iterations, respectively. A direct steady spectral inversion corresponds to $q=1$ and has $\mathcal{O}(N\log N)$ complexity.

The dominant memory requirement is
\begin{equation}
M
=
\mathcal{O}(N),
\label{eq:method_memory_complexity}
\end{equation}
because the framework stores field arrays, transform workspaces, spectral multipliers and geometric masks without assembling global stiffness, mass or differentiation matrices. The value of $q$ and the practical prefactor depend on the equation and integration scheme, but the asymptotic scaling remains $\mathcal{O}(N\log N)$ for a fixed-stage algorithm.

All calculations reported here use the transform conventions, numerical precision and benchmark-specific parameters documented in the Supplementary Information. Error measures, physical diagnostics, timing procedures, CPU/GPU synchronization, memory measurements and complete reproducibility settings are provided in Supplementary Sections~1.8--1.9 and~4.

\bmhead{Supplementary information}

Supplementary information accompanies this paper and contains detailed convergence tables, additional validation figures, Fourier-transform derivations, and operator-splitting formulas.

\bmhead{Funding}

This work was supported by the Xing Dian Talent Support Program of Yunnan Province (grant No. XDYC-QNRC-2022-0764) and the Yunnan Fundamental Research Projects (grant No. 202401CF070043).

\bmhead{Code availability}

The OptiXDE project website, including software documentation and usage information, is available at \url{https://yangylab.github.io/optixde-site/}. 








\bibliography{references}

\end{document}


\maketitle

\clearpage
\tableofcontents
\clearpage

\section{Numerical formulation, implementation and reproducibility}
\label{supp:numerical_implementation}

This section provides the numerical conventions, implementation details and reproducibility information used throughout the supplementary benchmark studies. Unless otherwise stated, all calculations are performed in double precision on uniform Cartesian discretizations. Spatial transforms, spectral multiplications and inverse transforms are evaluated without stiffness-matrix assembly, numerical quadrature or iterative linear-system solution. Spatial coordinates, wave-number arrays, propagation factors, geometric masks and boundary labels that remain unchanged during a simulation are constructed during preprocessing and cached for subsequent reuse.

\subsection{Discrete transforms and wave-number construction}
\label{supp:discrete_transforms}

Let the rectangular transform domain be defined as
\begin{equation}
\widetilde{\Omega}=\prod_{\ell=1}^{d}[0,L_{\ell}),
\label{eq:supp_computational_domain}
\end{equation}
where $d$ is the spatial dimension and $L_{\ell}$ is the domain length in the $\ell$-th coordinate direction. For a periodic direction discretized using $N_{\ell}$ uniformly spaced spatial points, the nodal coordinates and spatial spacing are defined as
\begin{equation}
x_{\ell,j}=j\Delta x_{\ell},\qquad \Delta x_{\ell}=\frac{L_{\ell}}{N_{\ell}},\qquad j=0,\ldots,N_{\ell}-1.
\label{eq:supp_periodic_coordinates}
\end{equation}

For a one-dimensional periodic discretization containing $N$ points, the forward and inverse discrete Fourier transforms are defined as
\begin{equation}
\widehat{u}_{m}=\sum_{j=0}^{N-1}u_j\exp\left(-\frac{2\pi\mathrm{i}jm}{N}\right),\qquad m=0,\ldots,N-1,
\label{eq:supp_forward_dft}
\end{equation}
and
\begin{equation}
u_j=\frac{1}{N}\sum_{m=0}^{N-1}\widehat{u}_{m}\exp\left(\frac{2\pi\mathrm{i}jm}{N}\right),\qquad j=0,\ldots,N-1.
\label{eq:supp_inverse_dft}
\end{equation}

In multiple dimensions, the transform is evaluated as the tensor product of the corresponding one-dimensional transforms. The inverse $d$-dimensional transform therefore contains the normalization factor $1/\prod_{\ell=1}^{d}N_{\ell}$.

For each periodic coordinate direction, the signed frequency index is arranged according to the standard fast-Fourier-transform ordering,
\begin{equation}
\mu_{\ell,m}=
\begin{cases}
m, & 0\leq m\leq \left\lfloor (N_{\ell}-1)/2\right\rfloor,\\
m-N_{\ell}, & \left\lfloor (N_{\ell}-1)/2\right\rfloor<m\leq N_{\ell}-1,
\end{cases}
\label{eq:supp_frequency_index}
\end{equation}
and the corresponding angular wave number is
\begin{equation}
k_{\ell,m}=\frac{2\pi}{L_{\ell}}\mu_{\ell,m}.
\label{eq:supp_periodic_wavenumber}
\end{equation}

For an even spatial resolution, this ordering stores the Nyquist mode at $\mu_{\ell,N_{\ell}/2}=-N_{\ell}/2$. When an odd-order derivative of a real-valued field is evaluated at an even spatial resolution, the multiplier associated with the Nyquist mode is set to zero so that the discrete derivative remains real valued to machine precision.

For a multidimensional frequency index $\boldsymbol{m}=(m_1,\ldots,m_d)$, the wave-number vector and its squared magnitude are defined as
\begin{equation}
\mathbf{k}_{\boldsymbol{m}}=\left(k_{1,m_1},\ldots,k_{d,m_d}\right),\qquad |\mathbf{k}_{\boldsymbol{m}}|^2=\sum_{\ell=1}^{d}k_{\ell,m_{\ell}}^2.
\label{eq:supp_multidimensional_wavenumber}
\end{equation}

Spatial derivatives are evaluated using analytical spectral multipliers,
\begin{equation}
\widehat{\frac{\partial^p u}{\partial x_{\ell}^p}}=\left(\mathrm{i}k_{\ell}\right)^p\widehat{u},\qquad \widehat{\nabla^2u}=-|\mathbf{k}|^2\widehat{u}.
\label{eq:supp_spectral_derivatives}
\end{equation}

Consequently, spatial differentiation requires only elementwise multiplication in transform space and does not involve finite-difference stencils, element-level integration or assembled differentiation matrices.

Homogeneous Dirichlet boundary conditions on a rectangular coordinate interval are represented using a sine basis. For $N_{\ell}$ interior spatial points, the coordinates are
\begin{equation}
x_{\ell,j}=\frac{(j+1)L_{\ell}}{N_{\ell}+1},\qquad j=0,\ldots,N_{\ell}-1,
\label{eq:supp_dirichlet_coordinates}
\end{equation}
and the associated wave numbers are
\begin{equation}
k_{\ell,m}^{(D)}=\frac{(m+1)\pi}{L_{\ell}},\qquad m=0,\ldots,N_{\ell}-1.
\label{eq:supp_dirichlet_wavenumber}
\end{equation}

The corresponding modal coefficients are evaluated using an orthonormally normalized type-I discrete sine transform.

Homogeneous Neumann boundary conditions are represented using a cosine basis on a spatial discretization containing both endpoints,
\begin{equation}
x_{\ell,j}=\frac{jL_{\ell}}{N_{\ell}-1},\qquad j=0,\ldots,N_{\ell}-1,
\label{eq:supp_neumann_coordinates}
\end{equation}
with the associated wave numbers
\begin{equation}
k_{\ell,m}^{(N)}=\frac{m\pi}{L_{\ell}},\qquad m=0,\ldots,N_{\ell}-1.
\label{eq:supp_neumann_wavenumber}
\end{equation}

The corresponding coefficients are evaluated using an orthonormally normalized type-I discrete cosine transform. Fourier, sine and cosine transforms may be combined dimension by dimension, allowing different Cartesian coordinate directions to use different boundary-condition representations.

The same transform normalization, frequency ordering and wave-number construction are retained throughout all spatial-resolution and time-step studies. Consequently, variations in the reported numerical errors arise from changes in spatial spacing, propagation interval, boundary representation or nonlinear approximation rather than from changes in the discrete-transform convention.

\subsection{Closed-form spectral propagation operators}
\label{supp:spectral_propagators}

Consider a linear evolution equation written as
\begin{equation}
\frac{\partial u}{\partial t}=\mathcal{L}u+s,
\label{eq:supp_linear_evolution}
\end{equation}
where $\mathcal{L}$ is a constant-coefficient spatial differential operator and $s$ is a prescribed source term. After applying a spatial transform whose basis diagonalizes $\mathcal{L}$, each spectral coefficient satisfies an independent scalar evolution equation,
\begin{equation}
\frac{\partial\widehat{u}}{\partial t}=\lambda(\mathbf{k})\widehat{u}+\widehat{s},
\label{eq:supp_modal_evolution}
\end{equation}
where $\lambda(\mathbf{k})$ is the spectral symbol of $\mathcal{L}$. The transformation therefore converts the original spatially coupled partial differential equation into a collection of mutually independent modal ordinary differential equations.

The exact variation-of-constants representation over the interval $[t_n,t_{n+1}]$, with $t_{n+1}=t_n+\Delta t$, is
\begin{equation}
\widehat{u}^{\,n+1}
=
\exp\left[\lambda(\mathbf{k})\Delta t\right]\widehat{u}^{\,n}
+
\int_{0}^{\Delta t}
\exp\left[\lambda(\mathbf{k})(\Delta t-\tau)\right]
\widehat{s}(t_n+\tau)\,\mathrm{d}\tau.
\label{eq:supp_variation_of_constants}
\end{equation}

When $\lambda(\mathbf{k})$ is time independent and the source is constant over the propagation interval, $\widehat{s}(t_n+\tau)=\widehat{s}^{\,n}$, the exact modal update becomes
\begin{equation}
\widehat{u}^{\,n+1}
=
G(\mathbf{k},\Delta t)\widehat{u}^{\,n}
+
\Delta t\,\varphi_1\!\left(\lambda(\mathbf{k})\Delta t\right)\widehat{s}^{\,n},
\label{eq:supp_exact_modal_update}
\end{equation}
where
\begin{equation}
G(\mathbf{k},\Delta t)
=
\exp\left[\lambda(\mathbf{k})\Delta t\right]
\label{eq:supp_general_propagator}
\end{equation}
is the homogeneous propagation operator and
\begin{equation}
\varphi_1(z)
=
\begin{cases}
\dfrac{\exp(z)-1}{z}, & z\neq0,\\[2mm]
1, & z=0
\end{cases}
\label{eq:supp_phi_function}
\end{equation}
is the first exponential-integrator function. The first term in Eq.~\eqref{eq:supp_exact_modal_update} propagates the existing modal state, whereas the second term accounts for the continuous accumulation and simultaneous propagation of the source contribution over the same interval. The update is exact for a time-independent modal eigenvalue and a source that is constant over the propagation interval. When the source varies in time, the temporal accuracy is determined by the approximation adopted for the source integral in Eq.~\eqref{eq:supp_variation_of_constants}.

For sufficiently small $|z|$, direct evaluation of $\exp(z)-1$ may suffer from subtractive cancellation. The function $\varphi_1(z)$ is therefore evaluated using a cancellation-free exponential-minus-one operation or, when necessary, the Taylor expansion
\begin{equation}
\varphi_1(z)
=
1+\frac{z}{2}+\frac{z^2}{6}+\frac{z^3}{24}+\mathcal{O}(z^4).
\label{eq:supp_phi_taylor}
\end{equation}

This treatment also provides the correct limiting update for modes satisfying $\lambda(\mathbf{k})=0$,
\begin{equation}
\widehat{u}^{\,n+1}
=
\widehat{u}^{\,n}
+
\Delta t\,\widehat{s}^{\,n}.
\label{eq:supp_zero_eigenvalue_update}
\end{equation}

For the diffusion equation
\begin{equation}
\frac{\partial u}{\partial t}
=
\kappa\nabla^2u,
\label{eq:supp_diffusion_equation}
\end{equation}
the transformed equation is
\begin{equation}
\frac{\partial\widehat{u}}{\partial t}
=
-\kappa|\mathbf{k}|^2\widehat{u}.
\label{eq:supp_diffusion_modal_equation}
\end{equation}

The corresponding modal eigenvalue and propagation factor are
\begin{equation}
\lambda_{\mathrm{diff}}(\mathbf{k})
=
-\kappa|\mathbf{k}|^2,
\qquad
G_{\mathrm{diff}}(\mathbf{k},\Delta t)
=
\exp\left(-\kappa|\mathbf{k}|^2\Delta t\right).
\label{eq:supp_diffusion_propagator}
\end{equation}

For $\kappa>0$, the propagation factor satisfies
\begin{equation}
0<G_{\mathrm{diff}}(\mathbf{k},\Delta t)\leq1,
\label{eq:supp_diffusion_contraction}
\end{equation}
so that nonzero spatial-frequency modes decay without numerical amplification. High-frequency modes decay more rapidly because their attenuation rate is proportional to $|\mathbf{k}|^2$, reproducing the physical smoothing mechanism of diffusion. The zero-frequency mode satisfies
\begin{equation}
G_{\mathrm{diff}}(\mathbf{0},\Delta t)=1,
\label{eq:supp_diffusion_zero_mode}
\end{equation}
which is consistent with conservation of the spatial mean in the absence of a source term. The homogeneous constant-coefficient diffusion subproblem is therefore integrated exactly in time and is not subject to the explicit diffusion stability restriction $\Delta t=\mathcal{O}(\Delta x^2)$. Nevertheless, the propagation interval may still affect the accuracy of time-dependent source approximations, nonlinear splitting, embedded-boundary corrections and time-resolved output.

For the linear Schr\"odinger equation
\begin{equation}
\mathrm{i}\frac{\partial u}{\partial t}
+
\frac{1}{2}\nabla^2u
=
0,
\label{eq:supp_linear_schrodinger_equation}
\end{equation}
the transformed modal equation is
\begin{equation}
\frac{\partial\widehat{u}}{\partial t}
=
-\frac{\mathrm{i}}{2}|\mathbf{k}|^2\widehat{u}.
\label{eq:supp_schrodinger_modal_equation}
\end{equation}

The corresponding modal eigenvalue and propagation factor are
\begin{equation}
\lambda_{\mathrm{Sch}}(\mathbf{k})
=
-\frac{\mathrm{i}}{2}|\mathbf{k}|^2,
\qquad
G_{\mathrm{Sch}}(\mathbf{k},\Delta t)
=
\exp\left(-\frac{\mathrm{i}}{2}|\mathbf{k}|^2\Delta t\right).
\label{eq:supp_schrodinger_propagator}
\end{equation}

Because
\begin{equation}
\left|G_{\mathrm{Sch}}(\mathbf{k},\Delta t)\right|=1,
\label{eq:supp_schrodinger_unit_modulus}
\end{equation}
the linear propagation step changes only the modal phase and introduces neither modal amplification nor numerical dissipation. In nonlinear Schr\"odinger problems, this closed-form propagator is applied to the linear dispersive subproblem, whereas the nonlinear contribution is advanced separately in physical space through operator splitting.

For the wave equation
\begin{equation}
\frac{\partial^2u}{\partial t^2}
=
c^2\nabla^2u,
\label{eq:supp_wave_equation}
\end{equation}
the auxiliary velocity variable
\begin{equation}
v=\frac{\partial u}{\partial t}
\label{eq:supp_wave_velocity}
\end{equation}
is introduced. Each transformed mode then satisfies
\begin{equation}
\frac{\mathrm{d}}{\mathrm{d}t}
\begin{bmatrix}
\widehat{u}\\
\widehat{v}
\end{bmatrix}
=
\begin{bmatrix}
0 & 1\\
-\omega^2 & 0
\end{bmatrix}
\begin{bmatrix}
\widehat{u}\\
\widehat{v}
\end{bmatrix},
\qquad
\omega=c|\mathbf{k}|.
\label{eq:supp_wave_first_order_system}
\end{equation}

The exact modal update is
\begin{equation}
\begin{bmatrix}
\widehat{u}^{\,n+1}\\
\widehat{v}^{\,n+1}
\end{bmatrix}
=
\begin{bmatrix}
\cos(\omega\Delta t) & \dfrac{\sin(\omega\Delta t)}{\omega}\\[3mm]
-\omega\sin(\omega\Delta t) & \cos(\omega\Delta t)
\end{bmatrix}
\begin{bmatrix}
\widehat{u}^{\,n}\\
\widehat{v}^{\,n}
\end{bmatrix}.
\label{eq:supp_wave_propagator}
\end{equation}

Equivalently, the displacement and velocity coefficients are updated as
\begin{equation}
\widehat{u}^{\,n+1}
=
\cos(\omega\Delta t)\widehat{u}^{\,n}
+
\frac{\sin(\omega\Delta t)}{\omega}\widehat{v}^{\,n},
\label{eq:supp_wave_displacement_update}
\end{equation}
and
\begin{equation}
\widehat{v}^{\,n+1}
=
-\omega\sin(\omega\Delta t)\widehat{u}^{\,n}
+
\cos(\omega\Delta t)\widehat{v}^{\,n}.
\label{eq:supp_wave_velocity_update}
\end{equation}

This propagation matrix is the exact solution operator of the corresponding modal harmonic oscillator and advances the displacement and velocity without artificial damping or time-discretization phase error. For the zero-frequency mode, the limiting values
\begin{equation}
\lim_{\omega\rightarrow0}
\frac{\sin(\omega\Delta t)}{\omega}
=
\Delta t,
\qquad
\lim_{\omega\rightarrow0}
\omega\sin(\omega\Delta t)
=
0
\label{eq:supp_wave_zero_mode}
\end{equation}
are used explicitly. The zero-frequency velocity therefore remains constant, while the corresponding mean displacement evolves linearly according to
\begin{equation}
\widehat{u}^{\,n+1}(\mathbf{0})
=
\widehat{u}^{\,n}(\mathbf{0})
+
\Delta t\,\widehat{v}^{\,n}(\mathbf{0}),
\qquad
\widehat{v}^{\,n+1}(\mathbf{0})
=
\widehat{v}^{\,n}(\mathbf{0}).
\label{eq:supp_wave_zero_mode_update}
\end{equation}

Steady linear equations are treated through direct modal inversion. For
\begin{equation}
\mathcal{L}u=f,
\label{eq:supp_steady_linear_equation}
\end{equation}
the transformed equation is
\begin{equation}
\lambda(\mathbf{k})\widehat{u}(\mathbf{k})
=
\widehat{f}(\mathbf{k}),
\label{eq:supp_steady_modal_equation}
\end{equation}
and the solution is obtained directly as
\begin{equation}
\widehat{u}(\mathbf{k})
=
\frac{\widehat{f}(\mathbf{k})}{\lambda(\mathbf{k})},
\qquad
\lambda(\mathbf{k})\neq0.
\label{eq:supp_direct_modal_inversion}
\end{equation}

The sign of the inversion factor is determined by the adopted form of the governing equation. For example, if
\begin{equation}
-\nabla^2u=f,
\label{eq:supp_negative_laplacian_poisson}
\end{equation}
then
\begin{equation}
\lambda(\mathbf{k})
=
|\mathbf{k}|^2,
\qquad
\widehat{u}(\mathbf{k})
=
\frac{\widehat{f}(\mathbf{k})}{|\mathbf{k}|^2},
\qquad
\mathbf{k}\neq\mathbf{0}.
\label{eq:supp_negative_laplacian_inversion}
\end{equation}

If the equation is instead written as
\begin{equation}
\nabla^2u=f,
\label{eq:supp_positive_laplacian_poisson}
\end{equation}
the spectral symbol and modal solution are
\begin{equation}
\lambda(\mathbf{k})
=
-|\mathbf{k}|^2,
\qquad
\widehat{u}(\mathbf{k})
=
-\frac{\widehat{f}(\mathbf{k})}{|\mathbf{k}|^2},
\qquad
\mathbf{k}\neq\mathbf{0}.
\label{eq:supp_positive_laplacian_inversion}
\end{equation}

For a periodic Poisson equation, the zero-frequency eigenvalue vanishes,
\begin{equation}
\lambda(\mathbf{0})=0.
\label{eq:supp_poisson_zero_eigenvalue}
\end{equation}

Solvability therefore requires the source to satisfy the compatibility condition
\begin{equation}
\widehat{f}(\mathbf{0})=0,
\label{eq:supp_poisson_compatibility_spectral}
\end{equation}
which is equivalent to
\begin{equation}
\int_{\widetilde{\Omega}}
f(\mathbf{x})\,\mathrm{d}\mathbf{x}
=
0.
\label{eq:supp_poisson_compatibility_physical}
\end{equation}

When necessary, the compatibility condition is enforced numerically by subtracting the discrete spatial mean,
\begin{equation}
f
\leftarrow
f-\overline{f}.
\label{eq:supp_poisson_mean_removal}
\end{equation}

After the nonzero modes have been inverted, the zero-frequency solution coefficient is set to
\begin{equation}
\widehat{u}(\mathbf{0})=0,
\label{eq:supp_poisson_zero_mean_selection}
\end{equation}
thereby selecting the unique zero-mean representative from the family of solutions that differ only by an additive constant.

The closed-form updates described above apply directly to linear constant-coefficient operators that are diagonal in the selected Fourier, sine or cosine basis. Spatially varying coefficients, nonlinear terms and embedded-boundary corrections generally couple spectral modes and cannot be represented by a single fixed multiplier. These contributions are therefore treated separately through operator splitting, pseudo-spectral evaluation, explicit or implicit correction steps, or iterative embedded-boundary updates, while the diagonal constant-coefficient component retains its closed-form spectral propagation.

All arrays that depend only on the spatial discretization, physical parameters and propagation interval are precomputed and cached. These arrays include the componentwise wave numbers $\mathbf{k}$, the squared wave-number magnitude $|\mathbf{k}|^2$, the modal eigenvalue $\lambda(\mathbf{k})$, the propagation factor $G(\mathbf{k},\Delta t)$ and, for the wave equation, the quantities
\begin{equation}
\cos(\omega\Delta t),
\qquad
\frac{\sin(\omega\Delta t)}{\omega},
\qquad
-\omega\sin(\omega\Delta t).
\label{eq:supp_cached_wave_factors}
\end{equation}

The limiting value $\sin(\omega\Delta t)/\omega=\Delta t$ is stored at $\omega=0$. The source-integration factor
\begin{equation}
\Delta t\,\varphi_1\!\left(\lambda(\mathbf{k})\Delta t\right)
\label{eq:supp_cached_source_factor}
\end{equation}
is similarly cached whenever the spatial discretization, operator parameters and propagation interval remain unchanged.

Consequently, a homogeneous linear propagation step reduces to
\begin{equation}
\widehat{u}^{\,n}
=
\mathcal{T}\!\left[u^n\right],
\qquad
\widehat{u}^{\,n+1}
=
G\odot\widehat{u}^{\,n},
\qquad
u^{n+1}
=
\mathcal{T}^{-1}\!\left[\widehat{u}^{\,n+1}\right],
\label{eq:supp_cached_linear_step}
\end{equation}
where $\mathcal{T}$ denotes the boundary-condition-compatible spatial transform and $\odot$ denotes elementwise multiplication. Each homogeneous propagation step therefore requires only one forward transform, one cached elementwise spectral multiplication and one inverse transform, with no stiffness-matrix assembly, numerical quadrature or iterative linear-system solution.

\subsection{Boundary-condition extensions}
\label{supp:boundary_extensions}

Boundary conditions determine the spatial basis used to diagonalize the constant-coefficient bulk operator. Periodic conditions are imposed directly through the discrete Fourier representation, so the field and the spatial derivatives required by the governing equation are periodic across opposite boundaries of $\widetilde{\Omega}$ and no additional boundary correction is required.

For homogeneous Dirichlet conditions on a Cartesian interval,
\begin{equation}
u(0,t)=u(L,t)=0,
\label{eq:supp_homogeneous_dirichlet}
\end{equation}
the field is represented by an odd extension,
\begin{equation}
u_{\mathrm{odd}}(-x,t)
=
-u_{\mathrm{odd}}(x,t),
\qquad
u_{\mathrm{odd}}(2L-x,t)
=
-u_{\mathrm{odd}}(x,t),
\qquad
x\in[0,L].
\label{eq:supp_odd_extension}
\end{equation}
Its $2L$-periodic continuation is equivalent, on $[0,L]$, to the sine-transform representation defined in Sec.~\ref{supp:discrete_transforms}.

For homogeneous Neumann conditions,
\begin{equation}
\frac{\partial u}{\partial x}(0,t)
=
\frac{\partial u}{\partial x}(L,t)
=
0,
\label{eq:supp_homogeneous_neumann}
\end{equation}
the field is represented by an even extension,
\begin{equation}
u_{\mathrm{even}}(-x,t)
=
u_{\mathrm{even}}(x,t),
\qquad
u_{\mathrm{even}}(2L-x,t)
=
u_{\mathrm{even}}(x,t),
\qquad
x\in[0,L].
\label{eq:supp_even_extension}
\end{equation}
Its $2L$-periodic continuation is equivalent to the cosine-transform representation. In practice, the appropriate sine or cosine transform is applied directly, and the doubled extension is not constructed explicitly.

Homogeneous mixed Dirichlet--Neumann conditions at opposite endpoints are represented by half-integer modes:
\begin{equation}
\begin{aligned}
u(0,t)=0,\quad
\frac{\partial u}{\partial x}(L,t)=0
&:
\quad
\phi_m(x)=\sin(k_m x),\\
\frac{\partial u}{\partial x}(0,t)=0,\quad
u(L,t)=0
&:
\quad
\phi_m(x)=\cos(k_m x),\\
k_m
&=
\frac{\left(m+\tfrac{1}{2}\right)\pi}{L},
\qquad
m=0,1,\ldots.
\end{aligned}
\label{eq:supp_mixed_endpoint_modes}
\end{equation}
The corresponding discrete sine or cosine transform type is selected consistently with the endpoint sampling convention.

Nonhomogeneous boundary data on separable Cartesian boundaries are treated by a lifting decomposition,
\begin{equation}
u(\mathbf{x},t)
=
w(\mathbf{x},t)
+
g_{\mathrm{ext}}(\mathbf{x},t),
\label{eq:supp_boundary_lifting}
\end{equation}
where $g_{\mathrm{ext}}$ is a sufficiently smooth extension of the prescribed boundary data and $w$ satisfies the associated homogeneous condition. For the linear boundary operator
\begin{equation}
\mathcal{B}u
=
\alpha u
+
\beta\frac{\partial u}{\partial n}
=
q,
\label{eq:supp_general_boundary_operator}
\end{equation}
$\alpha$ and $\beta$ are prescribed boundary coefficients, not both zero, $q$ is the prescribed boundary datum and $\partial u/\partial n=\mathbf{n}\cdot\nabla u$ is the derivative along the outward unit normal $\mathbf{n}$. Choosing $g_{\mathrm{ext}}$ such that $\mathcal{B}g_{\mathrm{ext}}=q$ gives $\mathcal{B}w=0$. Dirichlet, Neumann and Robin conditions correspond to $\beta=0$, $\alpha=0$ and $\alpha\beta\neq0$, respectively.

For the evolution equation in Eq.~\eqref{eq:supp_linear_evolution}, substitution of Eq.~\eqref{eq:supp_boundary_lifting} gives
\begin{equation}
\frac{\partial w}{\partial t}
=
\mathcal{L}w
+
\widetilde{s},
\qquad
\widetilde{s}
=
s
+
\mathcal{L}g_{\mathrm{ext}}
-
\frac{\partial g_{\mathrm{ext}}}{\partial t},
\label{eq:supp_lifted_evolution}
\end{equation}
where $\mathcal{L}$ is the constant-coefficient bulk operator, $s$ is the original source term and $\widetilde{s}$ is the modified source after lifting. The homogeneous variable $w$ is advanced using the corresponding spectral propagator, and the physical field is recovered from
\begin{equation}
u^{n+1}(\mathbf{x})
=
w^{n+1}(\mathbf{x})
+
g_{\mathrm{ext}}(\mathbf{x},t_{n+1}).
\label{eq:supp_lifting_recovery}
\end{equation}
The extension $g_{\mathrm{ext}}$ should be sufficiently smooth to avoid artificial high-frequency content and slow spectral convergence.

For a multidimensional rectangular domain, the transform is constructed independently in each coordinate direction:
\begin{equation}
\mathcal{T}
=
\mathcal{T}_{1}
\otimes
\mathcal{T}_{2}
\otimes
\cdots
\otimes
\mathcal{T}_{d},
\label{eq:supp_tensor_product_transform}
\end{equation}
where $\otimes$ denotes the tensor product and $\mathcal{T}_{\ell}$ is the Fourier, sine or cosine transform associated with the boundary conditions in the $\ell$-th coordinate direction. Periodic, Dirichlet, Neumann and mixed endpoint conditions can therefore be combined direction by direction whenever the domain and boundary partition are separable.

For the pure Neumann Poisson problem
\begin{equation}
-\nabla^{2}u=f
\quad
\text{in }\Omega,
\qquad
\frac{\partial u}{\partial n}=q
\quad
\text{on }\partial\Omega,
\label{eq:supp_neumann_poisson_problem}
\end{equation}
the constant mode belongs to the null space of the Laplacian. Solvability requires
\begin{equation}
\int_{\Omega}f\,\mathrm{d}\Omega
+
\int_{\partial\Omega}q\,\mathrm{d}\Gamma
=
0,
\label{eq:supp_neumann_compatibility}
\end{equation}
and an additional normalization, such as
\begin{equation}
\widehat{u}_{\mathbf{0}}
=
0,
\label{eq:supp_neumann_zero_mode_constraint}
\end{equation}
is imposed to fix the arbitrary additive constant. Here, $\widehat{u}_{\mathbf{0}}$ denotes the zero-frequency coefficient. This mode is excluded from division by the Laplacian eigenvalue. In transient Neumann problems, the zero mode is instead propagated according to the governing equation and the prescribed initial mean.

General homogeneous Robin conditions are not, in general, diagonalized by the standard Fourier, sine or cosine transforms. They require either a dedicated Robin eigenbasis or the physical-space boundary-enforcement procedure described in Sec.~\ref{supp:boundary_enforcement}. Curved interfaces, embedded geometries and nonseparable boundary segments are treated in the same way: the transform-space propagator advances the constant-coefficient bulk component, while localized physical-space corrections enforce the boundary and geometric constraints.

\subsection{Embedded geometries and mask construction}
\label{supp:geometry_masks}

An irregular physical domain $\Omega$ is embedded in the rectangular transform domain $\widetilde{\Omega}$. Its boundary is represented by a signed-distance or level-set function $\phi(\mathbf{x})$ with the sign convention
\begin{equation}
\begin{cases}
\phi(\mathbf{x})<0, & \mathbf{x}\in\Omega,\\
\phi(\mathbf{x})=0, & \mathbf{x}\in\partial\Omega,\\
\phi(\mathbf{x})>0, & \mathbf{x}\in\widetilde{\Omega}\setminus\overline{\Omega},
\end{cases}
\label{eq:supp_level_set_definition}
\end{equation}
where $\overline{\Omega}$ denotes the closure of $\Omega$. For an exact signed-distance function, $\lVert\nabla\phi\rVert_{2}=1$ almost everywhere near a smooth boundary, and $\nabla\phi$ points outward under the adopted sign convention.

Geometric primitives are combined through constructive solid geometry. For level-set functions $\phi_A$ and $\phi_B$ representing domains $A$ and $B$,
\begin{equation}
\phi_{A\cup B}
=
\min(\phi_A,\phi_B),
\qquad
\phi_{A\cap B}
=
\max(\phi_A,\phi_B),
\qquad
\phi_{A\setminus B}
=
\max(\phi_A,-\phi_B).
\label{eq:supp_csg_operations}
\end{equation}
These operations preserve the zero contour and the interior--exterior sign convention, but the resulting function is not generally an exact signed-distance function. Consequently, when $\phi$ is not re-distanced, the parameter $\epsilon$ introduced below specifies a level-set transition width rather than an exact Euclidean distance everywhere.

The binary indicator of the physical domain is
\begin{equation}
m_{\Omega}^{0}(\mathbf{x})
=
\begin{cases}
1, & \phi(\mathbf{x})\leq0,\\
0, & \phi(\mathbf{x})>0.
\end{cases}
\label{eq:supp_binary_mask}
\end{equation}
To avoid an abrupt transition at the spatial-sampling scale, the regularized Heaviside function
\begin{equation}
H_{\epsilon}(s)
=
\begin{cases}
0, & s\leq-\epsilon,\\[1mm]
\dfrac{1}{2}
\left[
1+\dfrac{s}{\epsilon}
+\dfrac{1}{\pi}
\sin\left(\dfrac{\pi s}{\epsilon}\right)
\right],
& |s|<\epsilon,\\[3mm]
1, & s\geq\epsilon
\end{cases}
\label{eq:supp_regularized_heaviside}
\end{equation}
is used to define the smoothed interior and exterior masks,
\begin{equation}
m_{\Omega}(\mathbf{x})
=
H_{\epsilon}\!\left[-\phi(\mathbf{x})\right],
\qquad
m_{\mathrm{ext}}(\mathbf{x})
=
1-m_{\Omega}(\mathbf{x}).
\label{eq:supp_smoothed_masks}
\end{equation}
Here, $\epsilon>0$ is the interface half-width. Thus, $m_{\Omega}\approx1$ inside $\Omega$, $m_{\Omega}\approx0$ outside $\Omega$, and $0<m_{\Omega}<1$ within a transition layer of total width approximately $2\epsilon$ when $\phi$ is a signed-distance function.

A dimensionless interface-localization mask is defined by
\begin{equation}
\chi_{\Gamma}(\mathbf{x})
=
4m_{\Omega}(\mathbf{x})
\left[
1-m_{\Omega}(\mathbf{x})
\right].
\label{eq:supp_interface_band}
\end{equation}
It vanishes away from the interface and reaches unity at $\phi=0$. The function $\chi_{\Gamma}$ is used as a smooth localization weight and is not a normalized approximation of the surface Dirac delta. When only a compact binary band is required, it may be replaced by the indicator of $|\phi|<\epsilon$.

The outward unit normal is evaluated in the interface band as
\begin{equation}
\mathbf{n}(\mathbf{x})
=
\frac{\nabla\phi(\mathbf{x})}
{\sqrt{\lVert\nabla\phi(\mathbf{x})\rVert_{2}^{2}+\delta_{n}^{2}}},
\label{eq:supp_boundary_normal}
\end{equation}
where $\delta_n>0$ prevents division by zero. Analytical gradients are used for elementary primitives when available. For composite or tabulated geometries, $\nabla\phi$ is evaluated using second-order centered differences, with one-sided differences at the outer boundary of $\widetilde{\Omega}$. At nonsmooth corners, the normal is interpreted segmentwise and is not uniquely defined at the corner point itself.

When $\partial\Omega$ is divided into $N_b$ Dirichlet, Neumann or Robin segments, nonnegative partition functions $\rho_b(\mathbf{x})$ are assigned such that
\begin{equation}
\rho_b(\mathbf{x})\geq0,
\qquad
\sum_{b=1}^{N_b}\rho_b(\mathbf{x})=1
\quad
\text{where }
\chi_{\Gamma}(\mathbf{x})>0.
\label{eq:supp_boundary_partition_unity}
\end{equation}
Here, $\rho_b$ localizes the $b$-th boundary segment and may overlap smoothly with neighbouring partitions near segment junctions. The corresponding enforcement mask is
\begin{equation}
\chi_b(\mathbf{x})
=
\chi_{\Gamma}(\mathbf{x})\rho_b(\mathbf{x}).
\label{eq:supp_segment_mask}
\end{equation}
This construction permits different boundary operators to be applied to different portions of the same embedded interface without modifying the interior spectral operator.

The smoothing width is selected relative to the Cartesian spacing:
\begin{equation}
\epsilon
=
\max\left(\epsilon_0,c_{\epsilon}h\right),
\qquad
h
=
\max_{1\leq\ell\leq d}\Delta x_{\ell},
\label{eq:supp_mask_width}
\end{equation}
where $\Delta x_{\ell}$ is the spacing in direction $\ell$, $c_{\epsilon}>0$ controls the resolved interface thickness and $\epsilon_0\geq0$ is an optional physical lower bound. For a signed-distance function, the transition layer spans approximately $2c_{\epsilon}$ intervals in the coarsest coordinate direction when $\epsilon=c_{\epsilon}h$. Benchmark-specific choices of $\epsilon_0$ and $c_{\epsilon}$ are reported with the corresponding embedded-domain calculations.

For static geometries, $\phi$, $m_{\Omega}^{0}$, $m_{\Omega}$, $m_{\mathrm{ext}}$, $\chi_{\Gamma}$, $\mathbf{n}$ and all segment masks are generated once during preprocessing and cached. For moving geometries, these quantities must be updated whenever the interface changes.

\subsection{Penalty and boundary enforcement}
\label{supp:boundary_enforcement}

The geometry-independent spectral propagator is combined with a physical-space enforcement map. For one evolution step,
\begin{equation}
u^{n+1}
=
\mathcal{E}_{\Delta t}
\left[
\mathcal{P}_{\Delta t}(u^{n})
\right],
\label{eq:supp_propagate_enforce}
\end{equation}
where $\mathcal{P}_{\Delta t}$ advances the constant-coefficient bulk equation over $\Delta t$, and $\mathcal{E}_{\Delta t}$ enforces the embedded boundary or geometric constraints. This composition separates the transform-diagonal bulk update from localized physical-space corrections.

For an embedded Dirichlet condition,
\begin{equation}
u=g
\qquad
\text{on }\Gamma_D,
\label{eq:supp_embedded_dirichlet_condition}
\end{equation}
the penalized evolution equation is written as
\begin{equation}
\frac{\partial u}{\partial t}
=
\mathcal{L}u+s
-
\frac{\chi_D(\mathbf{x})}{\eta_D}
\left[
u-g_{\mathrm{ext}}(\mathbf{x},t)
\right],
\label{eq:supp_dirichlet_penalized_evolution}
\end{equation}
where $\mathcal{L}$ is the bulk differential operator, $s$ is the source term, $\chi_D$ is a nonnegative Dirichlet enforcement mask, $g_{\mathrm{ext}}$ is an extension of the prescribed boundary value and $\eta_D>0$ is the penalty parameter. The mask is constructed from the exterior and interface masks defined in Sec.~\ref{supp:geometry_masks} and vanishes where no Dirichlet correction is required.

During a local enforcement interval $\Delta\tau$, the penalty subproblem is
\begin{equation}
\frac{\partial u}{\partial \tau}
=
-
\frac{\chi_D(\mathbf{x})}{\eta_D}
\left[
u-g_{\mathrm{ext}}(\mathbf{x})
\right].
\label{eq:supp_dirichlet_penalty_ode}
\end{equation}
If $\chi_D$ and $g_{\mathrm{ext}}$ are held fixed over $\Delta\tau$, its pointwise exact solution is
\begin{equation}
\mathcal{E}_{D,\Delta\tau}(v)
=
g_{\mathrm{ext}}
+
\exp
\left[
-\frac{\Delta\tau}{\eta_D}\chi_D
\right]
\left(
v-g_{\mathrm{ext}}
\right),
\label{eq:supp_exact_dirichlet_penalty}
\end{equation}
where $v$ is the field before enforcement. The update leaves $v$ unchanged where $\chi_D=0$ and relaxes it exponentially toward $g_{\mathrm{ext}}$ where $\chi_D>0$. Exact integration removes the explicit stability restriction associated with $\Delta\tau/\eta_D$, although the penalty magnitude still affects boundary accuracy and conditioning.

For a second-order diffusion-type operator, the penalty parameter is scaled as
\begin{equation}
\eta_D
=
\eta_{\mathrm{cell}}
\frac{h^{2}}{a_{\mathrm{ref}}},
\label{eq:supp_penalty_scaling}
\end{equation}
where $h$ is the representative Cartesian spacing defined in Eq.~\eqref{eq:supp_mask_width}, $a_{\mathrm{ref}}>0$ is a representative diffusion coefficient and $\eta_{\mathrm{cell}}>0$ is dimensionless. Benchmark-specific values are reported with the corresponding embedded-domain calculations, while the same scaling rule is retained within each spatial-resolution study.

For steady embedded problems, the bulk propagation and penalty map are applied within a pseudo-time or outer fixed-point iteration,
\begin{equation}
u^{(m+1)}
=
\mathcal{E}_{\Delta\tau}
\left[
\mathcal{P}_{\Delta\tau}
\left(
u^{(m)}
\right)
\right],
\qquad
m=0,1,\ldots,
\label{eq:supp_steady_enforcement_iteration}
\end{equation}
until the prescribed solution-change or residual tolerance is satisfied. Thus, the converged field solves the coupled bulk--boundary problem rather than a single post-processed spectral inversion.

Neumann and Robin conditions are represented by
\begin{equation}
\mathcal{B}(u)
=
\alpha u
+
\beta\,\mathbf{n}\cdot\mathbf{A}\nabla u
=
r
\qquad
\text{on }\Gamma_B,
\label{eq:supp_flux_boundary_operator}
\end{equation}
where $\alpha$ and $\beta$ are prescribed coefficients that are not simultaneously zero, $r$ is the boundary datum, $\mathbf{n}$ is the outward unit normal and $\mathbf{A}$ is the coefficient tensor associated with the diffusive flux. The sign convention uses $\mathbf{n}\cdot\mathbf{A}\nabla u$; an alternative flux sign may be absorbed into $\beta$ or $r$. Pure Neumann and Robin conditions correspond to $\alpha=0$ and $\alpha\beta\neq0$, respectively.

For an intermediate field $v$, the boundary residual is
\begin{equation}
R_B(v)
=
\alpha v
+
\beta\,\mathbf{n}\cdot\mathbf{A}\nabla v
-
r.
\label{eq:supp_boundary_residual}
\end{equation}
A localized fixed-point correction is applied in the boundary mask $\chi_B$:
\begin{equation}
\mathcal{E}_B(v)
=
v
-
\omega_B\chi_B
\frac{R_B(v)}
{\sigma_B+\delta_B},
\label{eq:supp_boundary_relaxation}
\end{equation}
where $0<\omega_B\leq1$ is the relaxation factor, $\delta_B>0$ prevents division by zero and $\sigma_B$ scales the residual to the field magnitude. For a second-order scalar diffusion operator,
\begin{equation}
\sigma_B
=
|\alpha|
+
|\beta|
\frac{a_n}{h},
\qquad
a_n
=
\mathbf{n}\cdot\mathbf{A}\mathbf{n},
\label{eq:supp_boundary_scaling}
\end{equation}
where $a_n>0$ is the diffusivity in the normal direction. If $\mathbf{A}=a\mathbf{I}$, then $a_n=a$. The gradient in Eq.~\eqref{eq:supp_boundary_residual} is evaluated using the same spatial differentiation convention as the governing equation.

For a boundary divided into $N_b$ segments, each segment has its own mask, data and enforcement map. One mixed-boundary sweep is written as
\begin{equation}
\mathcal{E}_{\mathrm{mix}}
=
\mathcal{E}_{B_{N_b}}
\circ
\cdots
\circ
\mathcal{E}_{B_2}
\circ
\mathcal{E}_{B_1},
\label{eq:supp_mixed_boundary_enforcement}
\end{equation}
where $\circ$ denotes operator composition and the rightmost map is applied first. The segment masks form the partition defined in Eq.~\eqref{eq:supp_boundary_partition_unity}, so interactions between different corrections are confined to segment junctions. The ordering and number of correction sweeps are held fixed within each convergence study and are reported with the corresponding benchmark settings.

After enforcement, only values in the physical domain are used for error measures and physical diagnostics. An embedded-domain representation of the restricted field is
\begin{equation}
u_{\Omega}(\mathbf{x})
=
m_{\Omega}^{0}(\mathbf{x})u(\mathbf{x}),
\label{eq:supp_solution_restriction}
\end{equation}
where $m_{\Omega}^{0}$ is the binary interior mask. The smoothed masks are used for enforcement and interface localization, whereas the binary mask is used to select the physical-domain values for diagnostics and visualization.

\subsection{Nonlinear splitting schemes}
\label{supp:nonlinear_splitting}

Consider a semilinear evolution equation
\begin{equation}
\frac{\partial u}{\partial t}
=
\mathcal{L}u
+
\mathcal{N}\!\left(u,\nabla u,\mathbf{x},t\right),
\label{eq:supp_semilinear_equation}
\end{equation}
where $\mathcal{L}$ is the linear constant-coefficient operator advanced in transform space and $\mathcal{N}$ contains the nonlinear and variable physical-space contributions. Let $\mathcal{P}_{\tau}$ denote the exact or closed-form flow generated by $\mathcal{L}$ over an interval $\tau$, and let $\mathcal{Q}_{\tau}$ denote the corresponding pointwise or pseudospectral nonlinear flow. The operator composition $\mathcal{A}\circ\mathcal{B}$ means that $\mathcal{B}$ is applied first.

The first-order Lie splitting is
\begin{equation}
u^{n+1}
=
\mathcal{P}_{\Delta t}
\circ
\mathcal{Q}_{\Delta t}
\left(u^n\right),
\label{eq:supp_lie_splitting}
\end{equation}
whereas the second-order Strang splitting used in the nonlinear benchmarks is
\begin{equation}
u^{n+1}
=
\mathcal{P}_{\Delta t/2}
\circ
\mathcal{Q}_{\Delta t}
\circ
\mathcal{P}_{\Delta t/2}
\left(u^n\right).
\label{eq:supp_strang_splitting}
\end{equation}
For sufficiently smooth autonomous subflows, Lie and Strang splitting have global temporal orders one and two, respectively. For nonperiodic or embedded problems, boundary enforcement is incorporated symmetrically into the split update so that the second-order composition is not destroyed by an asymmetric correction.

For the Allen--Cahn equation
\begin{equation}
\frac{\partial u}{\partial t}
=
\kappa_{\mathrm{AC}}\nabla^2u
+
u-u^3,
\label{eq:supp_allen_cahn_equation}
\end{equation}
where $\kappa_{\mathrm{AC}}>0$ is the diffusion coefficient, the linear flow is
\begin{equation}
\widehat{u}(\mathbf{k},t+\tau)
=
\exp
\left(
-\kappa_{\mathrm{AC}}
\lVert\mathbf{k}\rVert_2^2
\tau
\right)
\widehat{u}(\mathbf{k},t),
\label{eq:supp_allen_cahn_linear_flow}
\end{equation}
where $\mathbf{k}$ is the wave-number vector. The nonlinear subproblem
\begin{equation}
\frac{\partial u}{\partial t}
=
u-u^3
\label{eq:supp_allen_cahn_nonlinear_subproblem}
\end{equation}
is integrated exactly at each spatial location:
\begin{equation}
\mathcal{Q}_{\tau}^{\mathrm{AC}}(u)
=
\frac{u\exp(\tau)}
{\sqrt{
1+u^2\left[\exp(2\tau)-1\right]
}}.
\label{eq:supp_allen_cahn_nonlinear_flow}
\end{equation}
The complete update consists of two diffusion half-steps and one nonlinear full-step according to Eq.~\eqref{eq:supp_strang_splitting}.

For the cubic nonlinear Schr\"odinger equation
\begin{equation}
\mathrm{i}\frac{\partial h}{\partial t}
+
\frac{1}{2}\frac{\partial^2h}{\partial x^2}
+
|h|^2h
=
0,
\label{eq:supp_nls_equation}
\end{equation}
where $h$ is complex valued and $\mathrm{i}^2=-1$, the first linear half-step is
\begin{equation}
\widehat{h}^{\,*}(k)
=
\exp
\left(
-\frac{\mathrm{i}k^2\Delta t}{4}
\right)
\widehat{h}^{\,n}(k).
\label{eq:supp_nls_first_linear_half_step}
\end{equation}
After applying the inverse transform to obtain $h^*$, the nonlinear flow is evaluated exactly:
\begin{equation}
h^{**}
=
\exp
\left(
\mathrm{i}|h^*|^2\Delta t
\right)
h^*.
\label{eq:supp_nls_nonlinear_step}
\end{equation}
The second linear half-step is
\begin{equation}
\widehat{h}^{\,n+1}(k)
=
\exp
\left(
-\frac{\mathrm{i}k^2\Delta t}{4}
\right)
\widehat{h}^{\,**}(k).
\label{eq:supp_nls_second_linear_half_step}
\end{equation}
Both split subproblems are advanced in closed form, so no nonlinear iteration is required.

For the viscous Burgers equation
\begin{equation}
\frac{\partial u}{\partial t}
+
u\frac{\partial u}{\partial x}
=
\nu\frac{\partial^2u}{\partial x^2},
\label{eq:supp_burgers_equation}
\end{equation}
where $\nu>0$ is the kinematic viscosity, the diffusion flow is
\begin{equation}
\widehat{u}(k,t+\tau)
=
\exp
\left(
-\nu k^2\tau
\right)
\widehat{u}(k,t).
\label{eq:supp_burgers_diffusion_flow}
\end{equation}
The nonlinear subproblem is
\begin{equation}
\frac{\partial u}{\partial t}
=
-u\frac{\partial u}{\partial x}.
\label{eq:supp_burgers_nonlinear_subproblem}
\end{equation}
At every Runge--Kutta stage, $\partial u/\partial x$ is evaluated spectrally using the multiplier $\mathrm{i}k$, the product $u\,\partial u/\partial x$ is formed in physical space, and the nonlinear flow is advanced with the classical fourth-order Runge--Kutta method over the full nonlinear interval.

Quadratic pseudospectral products are de-aliased using either the two-thirds truncation rule or the three-halves padding rule. Let $\mu_{\ell,m_{\ell}}$ denote the signed integer Fourier-mode index in direction $\ell$, and define the two-thirds projector by
\begin{equation}
D_{2/3}(\boldsymbol{\mu})
=
\begin{cases}
1,
&
|\mu_{\ell,m_{\ell}}|
\leq
\left\lfloor N_{\ell}/3\right\rfloor
\quad
\text{for all }\ell,\\
0,
&
\text{otherwise}.
\end{cases}
\label{eq:supp_two_thirds_projector}
\end{equation}
For the Burgers nonlinearity, the filtered state and nonlinear term are evaluated as
\begin{equation}
u^{<}
=
\mathcal{F}^{-1}
\left[
D_{2/3}\widehat{u}
\right],
\qquad
\mathcal{N}_{2/3}(u)
=
\mathcal{F}^{-1}
\left\{
D_{2/3}
\mathcal{F}
\left[
-u^{<}
\frac{\partial u^{<}}{\partial x}
\right]
\right\}.
\label{eq:supp_two_thirds_rule}
\end{equation}
Filtering is applied at every nonlinear Runge--Kutta stage. This construction removes aliasing for quadratic products represented on the retained Fourier band.

Under three-halves padding, each spectral direction is enlarged to
\begin{equation}
N_{\ell}^{\mathrm{pad}}
=
\left\lceil
\frac{3N_{\ell}}{2}
\right\rceil,
\label{eq:supp_three_halves_size}
\end{equation}
the retained modes are zero padded, the nonlinear product is evaluated on the enlarged physical array, and the result is transformed back and truncated to the original spectral resolution. For nonlinear maps that are not purely quadratic, such as the exact Allen--Cahn and nonlinear Schr\"odinger pointwise flows, spectral truncation acts as a high-mode filter rather than an exact quadratic de-aliasing identity.

All nonlinear benchmark calculations use fixed time increments. The linear propagation factors are constructed once for the prescribed $\Delta t$ and reused throughout the simulation. When a complex FFT backend is used for a problem with a real-valued solution, the residual imaginary component after the inverse transform is discarded only after verifying that it remains at the level of double-precision round-off.

\subsection{Computational complexity and storage requirements}
\label{supp:computational_complexity}

This subsection derives the asymptotic computational complexity and storage requirements of the principal \OptiXDE{} operations. The analysis includes preprocessing, spectral propagation, nonlinear updates, boundary enforcement and physical diagnostics, but excludes file input, disk output and visualization. Constant prefactors depend on the spatial dimension, transform type, numerical backend, arithmetic precision, boundary treatment, splitting scheme and de-aliasing strategy.

For a $d$-dimensional Cartesian discretization containing $N_{\ell}$ spatial values in the $\ell$-th coordinate direction, the total number of spatial values is

\begin{equation}
N
=
\prod_{\ell=1}^{d}N_{\ell}.
\label{eq:supp_complexity_total_points}
\end{equation}

The Fourier, sine and cosine transforms used in \OptiXDE{} are separable tensor-product transforms. Along the $\ell$-th coordinate direction, the algorithm evaluates $N/N_{\ell}$ independent one-dimensional transforms of length $N_{\ell}$. The associated operation count is therefore

\begin{equation}
C_{\mathcal{T},\ell}
=
c_{\mathcal{T},\ell}
\frac{N}{N_{\ell}}
N_{\ell}\log_{2}N_{\ell}
=
c_{\mathcal{T},\ell}N\log_{2}N_{\ell},
\label{eq:supp_complexity_directional_transform}
\end{equation}

where $c_{\mathcal{T},\ell}$ is a constant determined by the transform type, arithmetic representation and numerical backend. Summing over all coordinate directions gives

\begin{equation}
C_{\mathcal{T}}
=
N\sum_{\ell=1}^{d}
c_{\mathcal{T},\ell}\log_{2}N_{\ell}.
\label{eq:supp_complexity_multidimensional_transform}
\end{equation}

Defining

\begin{equation}
c_{\mathcal{T},\max}
=
\max_{1\leq\ell\leq d}
c_{\mathcal{T},\ell},
\label{eq:supp_complexity_transform_constant}
\end{equation}

and using

\begin{equation}
\sum_{\ell=1}^{d}\log_{2}N_{\ell}
=
\log_{2}
\left(
\prod_{\ell=1}^{d}N_{\ell}
\right)
=
\log_{2}N,
\label{eq:supp_complexity_log_identity}
\end{equation}

the multidimensional transform cost satisfies

\begin{equation}
C_{\mathcal{T}}
\leq
c_{\mathcal{T},\max}N\log_{2}N.
\label{eq:supp_complexity_transform_bound}
\end{equation}

Consequently, both the forward and inverse transforms have complexity

\begin{equation}
C_{\mathcal{T}}
=
C_{\mathcal{T}^{-1}}
=
\mathcal{O}(N\log N).
\label{eq:supp_complexity_transform_order}
\end{equation}

Differences among Fourier, sine and cosine transforms modify only the constant prefactor. The logarithm is written with base two following the conventional FFT operation count, although its base does not affect the asymptotic complexity class.

A homogeneous linear propagation step has the form

\begin{equation}
\widehat{u}^{\,n}
=
\mathcal{T}\!\left[u^{n}\right],
\qquad
\widehat{u}^{\,n+1}
=
G\odot\widehat{u}^{\,n},
\qquad
u^{n+1}
=
\mathcal{T}^{-1}
\!\left[
\widehat{u}^{\,n+1}
\right],
\label{eq:supp_complexity_linear_step}
\end{equation}

where $G$ is a cached modal propagation factor and $\odot$ denotes elementwise multiplication. The forward and inverse transforms require $\mathcal{O}(N\log N)$ operations, whereas the modal multiplication requires $\mathcal{O}(N)$. Hence,

\begin{equation}
C_{\mathrm{linear}}
=
C_{\mathcal{T}}
+
C_{\mathcal{T}^{-1}}
+
\mathcal{O}(N)
=
\mathcal{O}(N\log N).
\label{eq:supp_complexity_linear_order}
\end{equation}

The same estimate applies to direct modal inversion of a transform-diagonal steady equation. For

\begin{equation}
\lambda(\mathbf{k})
\widehat{u}(\mathbf{k})
=
\widehat{f}(\mathbf{k}),
\label{eq:supp_complexity_modal_inversion_equation}
\end{equation}

the nonzero modal coefficients are obtained from

\begin{equation}
\widehat{u}(\mathbf{k})
=
\frac{\widehat{f}(\mathbf{k})}
{\lambda(\mathbf{k})}.
\label{eq:supp_complexity_modal_inversion}
\end{equation}

The modal division requires $\mathcal{O}(N)$ operations, so the forward and inverse transforms dominate the total cost:

\begin{equation}
C_{\mathrm{steady,direct}}
=
\mathcal{O}(N\log N).
\label{eq:supp_complexity_direct_steady}
\end{equation}

A zero-frequency mode, when present, is treated separately with constant additional work. Equation~\eqref{eq:supp_complexity_direct_steady} applies to constant-coefficient problems diagonalized by the selected transform and excludes outer iterations required by embedded geometries, variable coefficients or nonseparable constraints.

For a linear problem with a source term, the update may be written as

\begin{equation}
u^{n+1}
=
\mathcal{T}^{-1}
\left[
G\odot\mathcal{T}(u^{n})
+
H\odot\mathcal{T}(s^{n})
\right],
\label{eq:supp_complexity_source_step}
\end{equation}

where

\begin{equation}
H
=
\Delta t\,
\phi_{1}
\!\left(
\lambda\Delta t
\right)
\label{eq:supp_complexity_source_factor}
\end{equation}

is the cached source-integration factor. For a time-independent source, its transformed representation is evaluated once during preprocessing and reused. A time-dependent source supplied in physical space requires one additional forward transform per step. Since the number of transforms remains independent of $N$, the per-step complexity remains $\mathcal{O}(N\log N)$.

For a general numerical update, let $q_{\mathcal{T}}$ denote the total number of forward and inverse transforms performed during one complete step, and let $q_{\mathrm{p}}$ denote the number of pointwise array operations. These operations may include modal multiplication or division, nonlinear reaction updates, nonlinear products, mask multiplication, penalty relaxation, spectral filtering and diagnostic accumulation. The operation count can then be expressed as

\begin{equation}
C_{\mathrm{step}}
=
q_{\mathcal{T}}C_{\mathcal{T}}
+
q_{\mathrm{p}}c_{\mathrm{p}}N,
\label{eq:supp_complexity_general_step}
\end{equation}

where $c_{\mathrm{p}}$ is a representative constant for a pointwise operation. For a fixed numerical scheme, $q_{\mathcal{T}}$, $q_{\mathrm{p}}$ and $c_{\mathrm{p}}$ are independent of the spatial resolution. Therefore,

\begin{equation}
C_{\mathrm{step}}
=
\mathcal{O}
\!\left(
q_{\mathcal{T}}N\log N
+
q_{\mathrm{p}}N
\right)
=
\mathcal{O}(N\log N).
\label{eq:supp_complexity_general_order}
\end{equation}

Different equations and splitting schemes may require different numbers of transforms and therefore exhibit different practical execution times, even when their asymptotic complexity classes are identical.

For a Strang-splitting update,

\begin{equation}
u^{n+1}
=
\mathcal{P}_{\Delta t/2}
\circ
\mathcal{Q}_{\Delta t}
\circ
\mathcal{P}_{\Delta t/2}
(u^{n}),
\label{eq:supp_complexity_strang_step}
\end{equation}

the two linear half-steps are transform dominated. If the nonlinear flow $\mathcal{Q}_{\Delta t}$ is evaluated analytically and pointwise, as in the nonlinear Schr\"odinger and Allen--Cahn benchmarks, it requires $\mathcal{O}(N)$ operations. The complete splitting step therefore satisfies

\begin{equation}
C_{\mathrm{Strang}}
=
\mathcal{O}(N\log N).
\label{eq:supp_complexity_strang_order}
\end{equation}

When the nonlinear subproblem contains spectral derivatives, as in the viscous Burgers and incompressible Navier--Stokes equations, each derivative evaluation introduces additional transforms and modal multiplications. A fixed-stage Runge--Kutta or splitting scheme increases the computational prefactor but does not change the $\mathcal{O}(N\log N)$ scaling.

The two-thirds truncation rule introduces one spectral masking operation and therefore requires $\mathcal{O}(N)$ additional work. Under the three-halves padding rule, nonlinear products are evaluated on an enlarged array containing

\begin{equation}
N_{\mathrm{pad}}
=
\prod_{\ell=1}^{d}
\left\lceil
\frac{3N_{\ell}}{2}
\right\rceil
\label{eq:supp_complexity_padded_size_exact}
\end{equation}

spatial values. For sufficiently large $N_{\ell}$,

\begin{equation}
N_{\mathrm{pad}}
\approx
\left(
\frac{3}{2}
\right)^{d}N.
\label{eq:supp_complexity_padded_size}
\end{equation}

For a fixed spatial dimension $d$, the padded-transform cost is

\begin{equation}
C_{\mathrm{pad}}
=
\mathcal{O}
\!\left(
N_{\mathrm{pad}}
\log N_{\mathrm{pad}}
\right)
=
\mathcal{O}(N\log N).
\label{eq:supp_complexity_padded_order}
\end{equation}

Three-halves padding therefore preserves the asymptotic complexity class while increasing the arithmetic prefactor and temporary-memory requirement by a dimension-dependent factor.

For embedded geometries, signed-distance values, regularized masks, penalty factors and boundary labels are evaluated pointwise. A Dirichlet penalty update of the form

\begin{equation}
u
\leftarrow
g_{\mathrm{ext}}
+
\exp
\!\left[
-\frac{\Delta\tau}{\eta}\chi
\right]
\left(
u-g_{\mathrm{ext}}
\right)
\label{eq:supp_complexity_penalty_update}
\end{equation}

requires $\mathcal{O}(N)$ operations. Mask multiplication, restriction to the physical domain and evaluation of pointwise boundary residuals have the same linear cost.

Neumann and Robin corrections require evaluation of the normal derivative

\begin{equation}
\frac{\partial u}{\partial n}
=
\mathbf{n}\cdot\nabla u.
\label{eq:supp_complexity_normal_derivative}
\end{equation}

When the gradient is evaluated spectrally, a fixed number of additional transforms is required. One complete boundary-correction sweep therefore has complexity

\begin{equation}
C_{\mathrm{boundary}}
=
\mathcal{O}(N\log N).
\label{eq:supp_complexity_boundary_order}
\end{equation}

A fixed number of boundary-correction sweeps preserves the $\mathcal{O}(N\log N)$ per-step scaling.

For an embedded steady problem requiring $n_{\mathrm{it}}$ pseudo-time or outer-correction iterations, the total solver cost is

\begin{equation}
C_{\mathrm{embedded}}
=
C_{\mathrm{pre}}
+
n_{\mathrm{it}}
C_{\mathrm{iteration}}.
\label{eq:supp_complexity_embedded_decomposition}
\end{equation}

Since one iteration contains a fixed number of transforms and pointwise corrections,

\begin{equation}
C_{\mathrm{iteration}}
=
\mathcal{O}(N\log N),
\label{eq:supp_complexity_embedded_iteration}
\end{equation}

and hence

\begin{equation}
C_{\mathrm{embedded}}
=
\mathcal{O}(N)
+
\mathcal{O}
\!\left(
n_{\mathrm{it}}N\log N
\right).
\label{eq:supp_complexity_embedded_total}
\end{equation}

Thus, $\mathcal{O}(N\log N)$ describes the cost of one embedded-boundary iteration rather than necessarily the cost of the complete converged solution. The iteration count may depend on the spatial resolution, penalty parameter, mask-smoothing width, convergence tolerance and conditioning of the propagation--correction procedure.

The incompressible-flow projection algorithm contains a fixed number of velocity transforms, spectral derivatives, pressure inversions and pointwise corrections during each time step. Since the pressure Poisson equation is diagonal in the selected transform basis, no iterative algebraic pressure solver is required for the transform-compatible bulk problem. Consequently,

\begin{equation}
C_{\mathrm{projection}}
=
\mathcal{O}(N\log N)
\label{eq:supp_complexity_projection_order}
\end{equation}

per time step. Embedded-solid, inlet, outlet and wall corrections affect the constant prefactor but not the asymptotic scaling when the number of correction sweeps is fixed.

The preprocessing stage includes coordinate construction, wave-number generation, multidimensional modal-eigenvalue assembly, propagation-factor evaluation, signed-distance or level-set evaluation, mask construction, boundary labeling and reusable-array allocation. The one-dimensional coordinate and wave-number vectors require

\begin{equation}
\mathcal{O}
\!\left(
\sum_{\ell=1}^{d}N_{\ell}
\right)
\label{eq:supp_complexity_coordinate_preprocessing}
\end{equation}

operations and storage, whereas the multidimensional modal factors, masks and cached operators contain $N$ entries. For a prescribed static geometry evaluated pointwise, the overall preprocessing cost is

\begin{equation}
C_{\mathrm{pre}}
=
\mathcal{O}(N).
\label{eq:supp_complexity_preprocessing}
\end{equation}

These quantities are constructed once and cached, so their relative contribution decreases in long-time simulations, repeated-query calculations and parameter studies that reuse the same spatial discretization and geometry.

For a transient calculation containing $N_{t}$ time steps, the total arithmetic cost is

\begin{equation}
C_{\mathrm{total}}
=
C_{\mathrm{pre}}
+
N_{t}C_{\mathrm{step}}
+
C_{\mathrm{diag}},
\label{eq:supp_complexity_total_cost}
\end{equation}

where $C_{\mathrm{diag}}$ denotes the total cost of the requested error measures and physical diagnostics. A pointwise diagnostic over the complete domain requires $\mathcal{O}(N)$ operations, whereas a diagnostic involving spectral derivatives or additional transforms requires $\mathcal{O}(N\log N)$. If the number of diagnostic evaluations per time step is bounded independently of $N$, then

\begin{equation}
C_{\mathrm{total}}
=
\mathcal{O}(N)
+
\mathcal{O}
\!\left(
N_{t}N\log N
\right),
\label{eq:supp_complexity_transient_decomposition}
\end{equation}

and the dominant transient complexity is

\begin{equation}
C_{\mathrm{transient}}
=
\mathcal{O}
\!\left(
N_{t}N\log N
\right).
\label{eq:supp_complexity_transient_total}
\end{equation}

Equation~\eqref{eq:supp_complexity_transient_total} assumes that the number of transforms and pointwise stages per time step remains fixed. If the temporal increment is reduced as the spatial resolution increases, the resulting dependence of $N_{t}$ on $N$ must be included when interpreting the complete space--time complexity.

The storage requirement is determined by the simultaneously retained real- and complex-valued arrays. Let $n_{\mathrm{r}}$ and $n_{\mathrm{c}}$ denote the numbers of real and complex arrays containing $N$ values, respectively, and let $b_{\mathrm{r}}$ and $b_{\mathrm{c}}$ denote the corresponding numbers of bytes per value. The principal array storage is

\begin{equation}
M_{\mathrm{array}}
=
\left(
n_{\mathrm{r}}b_{\mathrm{r}}
+
n_{\mathrm{c}}b_{\mathrm{c}}
\right)N.
\label{eq:supp_complexity_array_memory}
\end{equation}

Including transform plans, backend workspaces and lower-order metadata gives

\begin{equation}
M
=
M_{\mathrm{array}}
+
M_{\mathrm{workspace}}
+
M_{\mathrm{metadata}}.
\label{eq:supp_complexity_total_memory}
\end{equation}

For a fixed algorithm and transform backend, the number of persistent and temporary arrays is independent of $N$. The dominant storage requirement therefore satisfies

\begin{equation}
M
=
\mathcal{O}(N).
\label{eq:supp_complexity_memory_order}
\end{equation}

In double-precision calculations, one real-valued array requires approximately

\begin{equation}
M_{\mathrm{real64}}
=
8N
\quad
\text{bytes},
\label{eq:supp_complexity_real64_memory}
\end{equation}

whereas one double-precision complex-valued array requires approximately

\begin{equation}
M_{\mathrm{complex128}}
=
16N
\quad
\text{bytes}.
\label{eq:supp_complexity_complex128_memory}
\end{equation}

These estimates exclude alignment, transform plans and backend-specific workspaces. Vector-valued equations, multistage nonlinear schemes, additional diagnostic fields and three-halves padding increase the number or size of the retained arrays but preserve linear asymptotic storage growth.

Unlike matrix-based discretizations, \OptiXDE{} does not assemble a global stiffness matrix, differentiation matrix or time-dependent system matrix. Its dominant persistent storage consists of the physical fields, transformed fields, cached modal operators, geometric masks, boundary labels and reusable work arrays. The principal asymptotic arithmetic and storage requirements are summarized in Table~\ref{tab:supp_complexity_summary}.

\begin{table}[t]
\centering
\caption{Asymptotic computational complexity and storage requirements of the principal \OptiXDE{} operations. Here, $N$ denotes the total number of spatial degrees of freedom, $N_t$ is the number of time steps and $n_{\mathrm{it}}$ is the number of embedded-boundary or correction iterations.}
\label{tab:supp_complexity_summary}
\begin{tabular*}{\textwidth}{@{\extracolsep{\fill}}lll@{}}
\toprule
Quantity & Asymptotic complexity & Description \\
\midrule
$C_{\mathrm{pre}}$ & $\mathcal{O}(N)$ & Preprocessing and modal-operator construction \\
$C_{\mathrm{step}}$ & $\mathcal{O}(N\log N)$ & One FFT-based propagation step \\
$C_{\mathrm{steady,direct}}$ & $\mathcal{O}(N\log N)$ & Direct solution of a constant-coefficient steady problem \\
$C_{\mathrm{transient}}$ & $\mathcal{O}(N_tN\log N)$ & Transient propagation over $N_t$ time steps \\
$C_{\mathrm{embedded}}$ & $\mathcal{O}(n_{\mathrm{it}}N\log N)$ & Iterative embedded-boundary or penalty correction \\
$M$ & $\mathcal{O}(N)$ & Persistent and reusable array storage \\
\botrule
\end{tabular*}
\end{table}

These estimates describe asymptotic arithmetic and storage growth rather than exact execution times. Practical performance additionally depends on transform planning, memory bandwidth, cache utilization, data layout, thread scheduling, accelerator-launch overhead, host--device transfers, device synchronization and the number of transforms required by each numerical scheme. The subsequent computational-performance and hardware-scaling studies examine whether the measured execution time, throughput and memory consumption are consistent with the theoretical estimates summarized in Table~\ref{tab:supp_complexity_summary}.

\subsection{Error measures and physical diagnostics}
\label{supp:error_measures}

Let $u_h$ denote the numerical solution and let $u_{\mathrm{ref}}$ denote either an analytical solution or a sufficiently resolved reference solution. Unless otherwise stated, the symbol $|\cdot|$ denotes the absolute value for real or complex scalar fields and the Euclidean magnitude for vector-valued fields. The set of sampled locations included in an error evaluation is denoted by $\Omega_h$, and its cardinality is denoted by $N_{\Omega}$. A small positive denominator safeguard $\varepsilon_{\mathrm{den}}$ is introduced whenever a normalized quantity may contain a vanishing reference value.

The pointwise field error at the sampled location $\mathbf{x}_j$ and time $t_n$ is defined as
\begin{equation}
e_j^n
=
\left|
u_h(\mathbf{x}_j,t_n)
-
u_{\mathrm{ref}}(\mathbf{x}_j,t_n)
\right|.
\label{eq:supp_pointwise_error}
\end{equation}

The normalized discrete absolute $L^2$ error, equivalently the root-mean-square error for uniform spatial sampling, is
\begin{equation}
e_{L^2}(t_n)
=
\left[
\frac{1}{N_{\Omega}}
\sum_{\mathbf{x}_j\in\Omega_h}
\left|
u_h(\mathbf{x}_j,t_n)
-
u_{\mathrm{ref}}(\mathbf{x}_j,t_n)
\right|^2
\right]^{1/2}.
\label{eq:supp_absolute_l2_error}
\end{equation}

The relative discrete $L^2$ error is
\begin{equation}
e_{\mathrm{rel}}(t_n)
=
\frac{
\left[
\displaystyle
\sum_{\mathbf{x}_j\in\Omega_h}
\left|
u_h(\mathbf{x}_j,t_n)
-
u_{\mathrm{ref}}(\mathbf{x}_j,t_n)
\right|^2
\right]^{1/2}
}{
\max\left\{
\left[
\displaystyle
\sum_{\mathbf{x}_j\in\Omega_h}
\left|
u_{\mathrm{ref}}(\mathbf{x}_j,t_n)
\right|^2
\right]^{1/2},
\varepsilon_{\mathrm{den}}
\right\}
}.
\label{eq:supp_relative_l2_error}
\end{equation}

The maximum pointwise error at time $t_n$ is
\begin{equation}
e_{\infty}(t_n)
=
\max_{\mathbf{x}_j\in\Omega_h}
\left|
u_h(\mathbf{x}_j,t_n)
-
u_{\mathrm{ref}}(\mathbf{x}_j,t_n)
\right|.
\label{eq:supp_maximum_error}
\end{equation}

For time-dependent problems, the terminal absolute error is denoted by
\begin{equation}
e(T)
=
e_{L^2}(T),
\label{eq:supp_terminal_error}
\end{equation}
and the maximum space--time error is
\begin{equation}
e_{\max}
=
\max_{0\leq n\leq N_t}
e_{\infty}(t_n)
=
\max_{0\leq n\leq N_t}
\max_{\mathbf{x}_j\in\Omega_h}
e_j^n.
\label{eq:supp_maximum_spacetime_error}
\end{equation}

Localized errors are evaluated by restricting the sampled locations to a diagnostic subdomain $\mathcal{D}\subseteq\Omega$. Let $\mathcal{D}_h$ denote the sampled locations belonging to $\mathcal{D}$ and let $N_{\mathcal{D}}$ denote their number. The regional root-mean-square and maximum errors are
\begin{equation}
e_{\mathcal{D}}(t_n)
=
\left[
\frac{1}{N_{\mathcal{D}}}
\sum_{\mathbf{x}_j\in\mathcal{D}_h}
\left|
u_h(\mathbf{x}_j,t_n)
-
u_{\mathrm{ref}}(\mathbf{x}_j,t_n)
\right|^2
\right]^{1/2},
\label{eq:supp_regional_l2_error}
\end{equation}
and
\begin{equation}
e_{\mathcal{D}}^{\infty}(t_n)
=
\max_{\mathbf{x}_j\in\mathcal{D}_h}
\left|
u_h(\mathbf{x}_j,t_n)
-
u_{\mathrm{ref}}(\mathbf{x}_j,t_n)
\right|.
\label{eq:supp_regional_maximum_error}
\end{equation}
The diagnostic region $\mathcal{D}$ is specified separately for each benchmark and may represent an embedded-interface neighbourhood, a singular-corner neighbourhood, a steep-gradient region or a smooth bulk region.

For a boundary condition written as $\mathcal{B}(u)=g$ on a boundary segment $\Gamma$, the maximum boundary residual is defined as
\begin{equation}
r_{\Gamma}^{\infty}(t_n)
=
\max_{\mathbf{x}_{\ell}^{\Gamma}\in\Gamma}
\left|
\mathcal{B}
\left[
\mathcal{I}_h u_h
\right]
\left(
\mathbf{x}_{\ell}^{\Gamma},t_n
\right)
-
g
\left(
\mathbf{x}_{\ell}^{\Gamma},t_n
\right)
\right|,
\label{eq:supp_boundary_residual_norm}
\end{equation}
where $\mathcal{I}_h$ denotes interpolation to the diagnostic boundary locations when the physical boundary does not coincide with the sampled locations. Separate values of $r_{\Gamma}^{\infty}$ are reported for different boundary segments when mixed conditions are imposed.

When a governing-equation residual is reported, let $\mathcal{R}_h(u_h)$ denote the discrete residual of the corresponding steady or time-dependent equation. Its normalized magnitude is written as
\begin{equation}
r_{\mathrm{PDE}}(t_n)
=
\frac{
\left\|
\mathcal{R}_h
\left(
u_h(\cdot,t_n)
\right)
\right\|_{2,\Omega}
}{
\mathcal{S}_{\mathcal{R}}(t_n)
+
\varepsilon_{\mathrm{den}}
},
\label{eq:supp_normalized_pde_residual}
\end{equation}
where $\mathcal{S}_{\mathcal{R}}$ is a problem-dependent residual scale specified in the corresponding benchmark. This normalization permits residual magnitudes to be compared across spatial resolutions without conflating the residual with the physical solution error.

For a sequence of calculations characterized by a refinement parameter $q$, the observed order between two successive values $q_1$ and $q_2$ is
\begin{equation}
p
=
\frac{
\log\left(e_1/e_2\right)
}{
\log\left(q_1/q_2\right)
},
\label{eq:supp_observed_convergence_order}
\end{equation}
where $q=h$ for a spatial-resolution study and $q=\Delta t$ for a temporal study. An observed order is not interpreted when the reported error is dominated by floating-point round-off, reference-solution uncertainty or a fixed-resolution error plateau.

For a scalar physical diagnostic $Q_h(t)$ that is theoretically conserved, the relative drift and its maximum value are defined as
\begin{equation}
\delta_Q(t_n)
=
\frac{
\left|
Q_h(t_n)-Q_h(0)
\right|
}{
\left|
Q_h(0)
\right|
+
\varepsilon_{\mathrm{den}}
},
\qquad
\varepsilon_Q
=
\max_{0\leq n\leq N_t}
\delta_Q(t_n).
\label{eq:supp_invariant_drift}
\end{equation}
The physical definition of $Q_h$ is problem dependent; examples include mass, Hamiltonian energy and total wave energy.

For a diagnostic quantity $Q_h(t)$ that is theoretically nonincreasing, the maximum normalized monotonicity violation is defined as
\begin{equation}
\delta_Q^{\uparrow}
=
\max_{0\leq n<N_t}
\frac{
\max\left[
0,
Q_h(t_{n+1})-Q_h(t_n)
\right]
}{
\left|
Q_h(0)
\right|
+
\varepsilon_{\mathrm{den}}
}.
\label{eq:supp_monotonicity_violation}
\end{equation}
This quantity is used only as a numerical consistency indicator; the corresponding energy, norm or dissipation functional is defined within the benchmark to which it applies.

\subsection{Software, hardware and reproducibility settings}
\label{supp:computational_environment}

Unless otherwise stated, all numerical experiments were performed using the same \OptiXDE{} implementation and IEEE 754 double-precision arithmetic. The reported performance measurements were obtained on a Google Colab runtime. The virtualized CPU was reported by the operating system as an Intel(R) Xeon(R) CPU @ 2.20GHz, with 12 logical CPUs exposed to the runtime and 83.5 GiB of system memory. Six CPU threads were used for the controlled CPU performance benchmarks. GPU-accelerated calculations were performed using an NVIDIA A100-SXM4-40GB with 39.5 GiB of device memory and CUDA runtime 12.9. Because the CPU hardware is exposed through a virtualized cloud environment, the operating-system-reported processor description and the number of logical CPUs available to the runtime are given rather than an inferred physical processor SKU.

The reference implementation used Python 3.13.15 and \OptiXDE{} 0.2.6. The controlled CPU--GPU transform benchmarks used the PyTorch 2.11.0+cu128 FFT interface on both devices, with CUDA transforms executed through cuFFT. The production cylinder calculation used CuPy 14.0.1 with cuFFT. NumPy 2.1.3 was used for auxiliary array operations and post-processing. Real-valued physical fields were stored in 64-bit floating-point format and complex spectral fields in 128-bit complex format. No mixed-precision arithmetic or tensor-core acceleration was used in the reported accuracy or performance studies.

The CPU thread count was fixed at six for the controlled CPU performance comparisons. All timed calculations included preliminary warm-up runs and repeated measurements, and GPU timings used explicit device synchronization before and after the measured region. The number of warm-up runs and recorded repetitions was benchmark-specific and is stated with the corresponding performance result. Unless otherwise noted, timing statistics are reported as medians together with interquartile ranges. The inclusion or exclusion of initialization, operator construction, host--device transfer, file input/output and post-processing is specified for each benchmark. The principal software and hardware settings are summarized in Table~\ref{tab:supp_computational_environment}. Source code, benchmark parameters and post-processing scripts will be released through the \OptiXDE{} repository upon publication.

\begin{table}[t]
\centering
\caption{Software and hardware environment used for the \OptiXDE{} numerical experiments and performance measurements.}
\label{tab:supp_computational_environment}
\begin{tabular}{ll}
\toprule
Setting & Specification \\
\midrule
Execution platform & Google Colab \\
Operating system & Linux-6.6.122+-x86\_64-with-glibc2.35 \\
CPU & Intel(R) Xeon(R) CPU @ 2.20GHz \\
Logical CPUs exposed & 12 \\
CPU threads used & 6 \\
System memory & 83.5 GiB \\
GPU & NVIDIA A100-SXM4-40GB \\
GPU memory & 39.5 GiB \\
CUDA runtime & 12.9 \\
GPU driver & 580.82.07 \\
Python version & 3.13.15 \\
\OptiXDE{} version & 0.2.6 \\
NumPy version & 2.1.3 \\
PyTorch version & 2.11.0+cu128 \\
CuPy version & 14.0.1 \\
cuFFT version & 11303 \\
CPU transform backend & PyTorch FFT for controlled CPU--GPU benchmarks \\
GPU transform backend & PyTorch/cuFFT for controlled benchmarks; CuPy/cuFFT for the cylinder calculation \\
Arithmetic precision & IEEE 754 double precision \\
Warm-up runs & Benchmark-specific; stated with each performance result \\
Recorded repetitions & Benchmark-specific; stated with each performance result \\
Reported statistic & Median with interquartile range \\
GPU synchronization & Explicit synchronization before and after timing \\
\bottomrule
\end{tabular}
\end{table}

\section{Additional benchmark results}
\label{supp:additional_benchmarks}

The common transform conventions, propagation operators, boundary treatments, geometry representation, nonlinear splitting rules, complexity estimates, error norms and computational environment are defined in Section~\ref{supp:numerical_implementation}. To avoid duplicating that material, the present section reports only the benchmark-specific governing data, selected implementation branch, reference solution, diagnostic quantities and numerical results. Whenever a benchmark is a direct specialization of a general operator introduced previously, the corresponding equation is cited rather than rederived.

\subsection{Transient diffusion equation}
\label{supp:diffusion}

This section provides the complete definition, numerical implementation and additional results for the transient diffusion benchmark summarized in the main text. The benchmark is designed to verify whether the closed-form propagation operator of OptiXDE reproduces the analytical exponential decay of individual Laplacian eigenmodes. The effects of spatial resolution, modal frequency and repeated propagation are examined separately to distinguish spatial-representation effects from accumulated finite-precision errors.

\subsubsection{Benchmark definition}
\label{supp:diffusion_definition}

The source-free transient diffusion equation is
\begin{equation}
\frac{\partial u}{\partial t}=\nu\nabla^2u,
\qquad
(\mathbf{x},t)\in\Omega\times(0,T],
\label{eq:supp_diffusion_governing}
\end{equation}
where $\Omega=[0,\pi]^2$, homogeneous Dirichlet conditions $u=0$ are imposed on $\partial\Omega$, and the diffusion coefficient is set to $\nu=1$.

The initial field is selected as a Laplacian eigenmode,
\begin{equation}
u(x,y,0)=10\sin(k_x x)\sin(k_y y),
\label{eq:supp_diffusion_initial}
\end{equation}
with the two modal combinations
\begin{equation}
(k_x,k_y)=(1,1)
\qquad\text{and}\qquad
(k_x,k_y)=(2,3).
\label{eq:supp_diffusion_modes}
\end{equation}

The corresponding analytical solution is
\begin{equation}
u_{\mathrm{ex}}(x,y,t)=10\exp\left[-\nu\left(k_x^2+k_y^2\right)t\right]\sin(k_x x)\sin(k_y y).
\label{eq:supp_diffusion_exact}
\end{equation}

All terminal quantities are evaluated at $T=1.5$. The fundamental mode $(1,1)$ represents a smooth and relatively slowly decaying field, whereas the mode $(2,3)$ contains more rapid spatial variations and undergoes substantially stronger attenuation. The corresponding analytical decay rates are
\begin{equation}
\lambda_{\mathrm{ex}}=\nu\left(k_x^2+k_y^2\right)=
\begin{cases}
2, & (k_x,k_y)=(1,1),\\
13, & (k_x,k_y)=(2,3).
\end{cases}
\label{eq:supp_diffusion_decay_rates}
\end{equation}

\subsubsection{Benchmark-specific spectral update}
\label{supp:diffusion_implementation}

The transform normalization, homogeneous-Dirichlet sampling convention and sine-basis wave numbers follow Section~\ref{supp:discrete_transforms}. For the present domain with $L_x=L_y=\pi$, the $N_x\times N_y$ interior samples are
\begin{equation}
x_i=\frac{i\pi}{N_x+1},
\qquad
y_j=\frac{j\pi}{N_y+1},
\qquad
i=1,\ldots,N_x,
\qquad
j=1,\ldots,N_y.
\label{eq:supp_diffusion_benchmark_sampling}
\end{equation}

The general diffusion propagator in Eq.~\eqref{eq:supp_diffusion_propagator} reduces, for the sine mode $(p,q)$, to the benchmark-specific modal update
\begin{equation}
\widetilde{u}_{p,q}^{\,n+1}
=
\exp\left[-\nu\left(p^2+q^2\right)\Delta t\right]
\widetilde{u}_{p,q}^{\,n}.
\label{eq:supp_diffusion_benchmark_update}
\end{equation}
Only the coefficient associated with the prescribed mode is nonzero in exact arithmetic. The zero mode is absent under the homogeneous Dirichlet representation.

The multiplier in Eq.~\eqref{eq:supp_diffusion_benchmark_update} is constructed once for each tested $\Delta t$ and reused. To examine accumulated finite-precision effects under the common OptiXDE workflow, every reported step performs a forward DST-I, the cached modal multiplication and an inverse DST-I, even though a single eigenmode could alternatively be retained entirely in transform space. Because the selected modes belong exactly to every tested sine basis, the resolution study measures transform-size and round-off sensitivity rather than a conventional spatial-discretization error. The stability and complexity properties follow Sections~\ref{supp:spectral_propagators} and~\ref{supp:computational_complexity} and are not rederived here.

\subsubsection{Eigenmode decay-rate verification}
\label{supp:diffusion_decay_measure}

Because the initial field contains a single Laplacian eigenmode, its spatial pattern remains unchanged and its amplitude decays as $\exp(-\lambda_{\mathrm{ex}}t)$, where $\lambda_{\mathrm{ex}}$ is given by Eq.~\eqref{eq:supp_diffusion_decay_rates}. To verify this mode-specific property independently of the general field-error measures defined in Section~\ref{supp:error_measures}, an effective numerical decay rate is inferred from the reduction in the discrete field norm:
\begin{equation}
\lambda_h=-\frac{1}{T}\ln\left[\frac{\left\|u_h(\cdot,T)\right\|_{2,h}}{\left\|u_h(\cdot,0)\right\|_{2,h}}\right].
\label{eq:supp_diffusion_decay_rate}
\end{equation}

The corresponding decay-rate discrepancy is
\begin{equation}
\varepsilon_{\lambda}=\left|\lambda_h-\lambda_{\mathrm{ex}}\right|.
\label{eq:supp_diffusion_decay_error}
\end{equation}

This diagnostic is specific to the present single-eigenmode benchmark and directly tests whether the closed-form OptiXDE propagator reproduces the prescribed exponential attenuation.

\subsubsection{Spatial-resolution study}
\label{supp:diffusion_resolution_study}

The influence of spatial resolution is examined using uniform tensor-product sampling with a fixed time step of $\Delta t=10^{-3}$ and a terminal time of $T=1.5$. Here, $N_x$ and $N_y$ denote the numbers of interior sampling points in the $x$- and $y$-directions, respectively, while the boundary values are imposed separately as zero. For the fundamental mode $(k_x,k_y)=(1,1)$, the resolution is varied as $N_x=N_y\in\{64,128,256,512\}$. The higher-frequency mode $(k_x,k_y)=(2,3)$ is additionally evaluated at spatial resolutions of $128^2$ and $256^2$ to verify that the propagation operator remains accurate for a more rapidly varying and more strongly attenuated eigenmode. Because both prescribed modes are represented exactly by the corresponding discrete sine bases, this study evaluates sensitivity to transform size and finite-precision operations rather than a conventional spatial-approximation convergence rate.

The complete numerical results are summarized in Table~\ref{tab:supp_diffusion_resolution}. For the $(1,1)$ mode, the terminal field error ranges from $3.320\times10^{-14}$ at a resolution of $64^2$ to $6.306\times10^{-14}$ at a resolution of $512^2$, while the decay-rate error ranges from $8.704\times10^{-14}$ to $1.679\times10^{-13}$. For the $(2,3)$ mode, the terminal field errors are $1.722\times10^{-18}$ and $1.253\times10^{-17}$ at resolutions of $128^2$ and $256^2$, respectively, with corresponding decay-rate errors of $1.474\times10^{-13}$ and $1.705\times10^{-13}$.

\begin{table}[t]
\centering
\caption{Spatial-resolution sensitivity of the transient diffusion eigenmode benchmark at $T=1.5$ using $\Delta t=10^{-3}$. Here, $N_x$ and $N_y$ denote the numbers of interior sampling points.}
\label{tab:supp_diffusion_resolution}
\begin{tabular}{cccc}
\toprule
$N_x\times N_y$ & Mode $(k_x,k_y)$ & $e(T)$ & $\varepsilon_\lambda$ \\
\midrule
$64^2$  & $(1,1)$ & $3.320\times10^{-14}$ & $8.704\times10^{-14}$ \\
$128^2$ & $(1,1)$ & $4.005\times10^{-14}$ & $1.064\times10^{-13}$ \\
$256^2$ & $(1,1)$ & $4.382\times10^{-14}$ & $1.168\times10^{-13}$ \\
$512^2$ & $(1,1)$ & $6.306\times10^{-14}$ & $1.679\times10^{-13}$ \\
$128^2$ & $(2,3)$ & $1.722\times10^{-18}$ & $1.474\times10^{-13}$ \\
$256^2$ & $(2,3)$ & $1.253\times10^{-17}$ & $1.705\times10^{-13}$ \\
\bottomrule
\end{tabular}
\end{table}

All errors remain within a round-off-dominated regime and therefore do not exhibit a conventional monotonic convergence trend. Once a prescribed eigenmode is represented by the discrete sine basis, increasing the spatial resolution does not reduce a leading spatial approximation error. Instead, it changes the transform size, normalization operations and accumulated finite-precision effects, which explains the small non-monotonic variation in the reported errors.

The substantially smaller absolute terminal errors of the $(2,3)$ mode should be interpreted in relation to its stronger analytical attenuation. At $T=1.5$, its amplitude is proportional to $\exp(-13T)=\exp(-19.5)$, whereas that of the $(1,1)$ mode is proportional to $\exp(-2T)=\exp(-3)$. The smaller value of $e(T)$ for the higher-frequency mode therefore does not imply intrinsically greater numerical accuracy. In contrast, the decay-rate errors of the two modes remain comparable, confirming that their mode-dependent exponential attenuation is reproduced consistently.

The dependence of $e(T)$ and $\varepsilon_{\lambda}$ on spatial resolution is shown in Fig.~\ref{fig:supp_diffusion_sensitivity}a,b. No empirical convergence order is reported because the measured errors are already governed by round-off effects rather than by an asymptotic spatial-approximation error.

\subsubsection{Modal field comparisons}
\label{supp:diffusion_fields}

\begin{figure}[t]
\centering
\includegraphics[width=\textwidth]{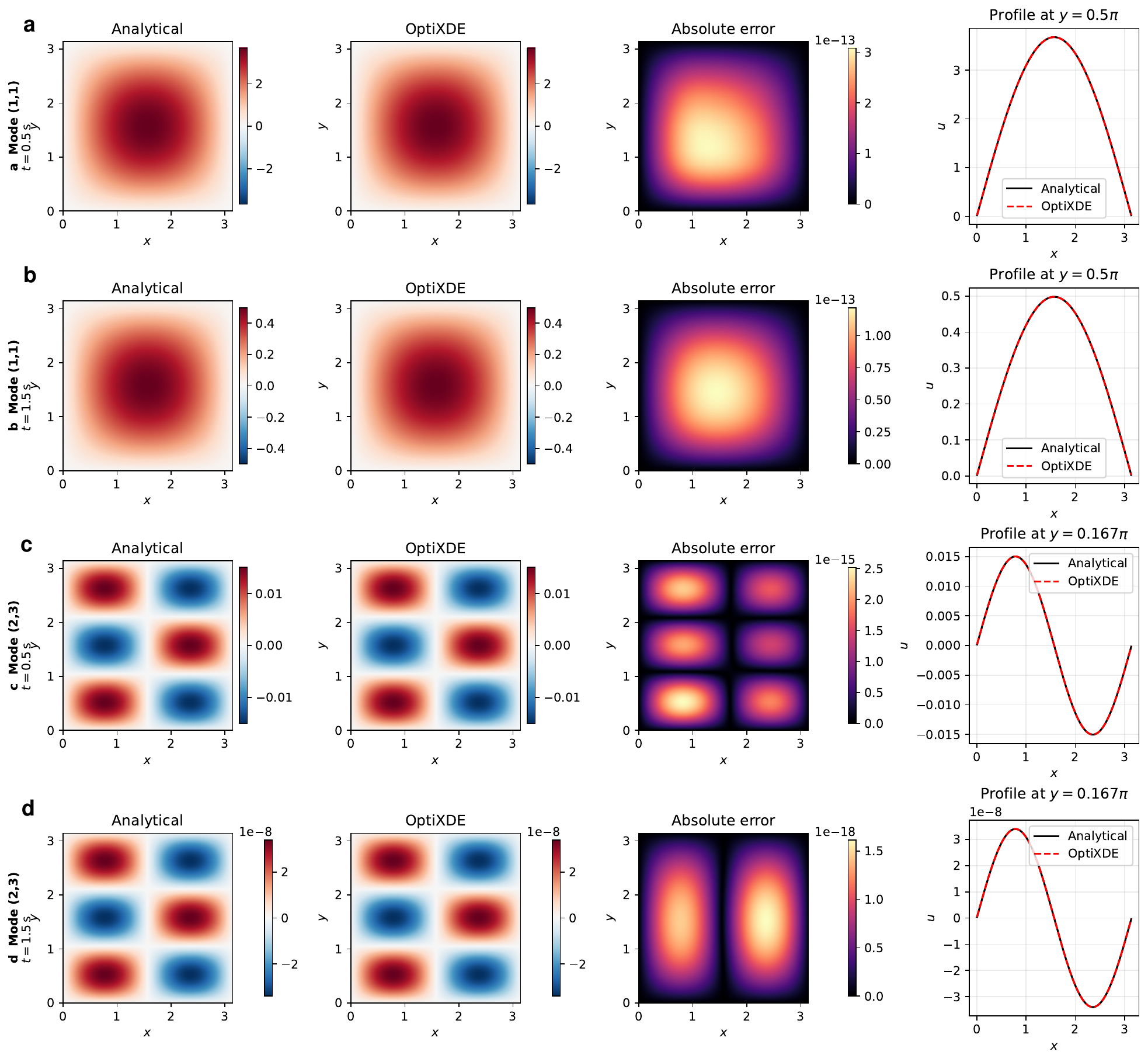}
\caption{\textbf{Transient diffusion eigenmodes at different times.} \textbf{a,b}, Analytical solution, \OptiXDE{} solution, absolute-error distribution and line profile for the fundamental mode $(k_x,k_y)=(1,1)$ at a spatial resolution of $512^2$ and at $t=0.5$ and $1.5$, respectively; the profiles are extracted at the sampled ordinate nearest to $y=\pi/2$. \textbf{c,d}, Corresponding comparisons for the higher-frequency mode $(k_x,k_y)=(2,3)$ at $t=0.5$ and $1.5$, respectively; the profiles are extracted at the sampled ordinate nearest to $y=\pi/6$. The numerical fields reproduce the analytical modal patterns without visible phase displacement or systematic amplitude distortion at either time.}
\label{fig:supp_diffusion_modes}
\end{figure}

Figure~\ref{fig:supp_diffusion_modes} compares the analytical and numerical fields for both tested modes at a spatial resolution of $512^2$, evaluated at the intermediate time $t=0.5$ and the terminal time $t=1.5$.

For the fundamental mode $(1,1)$, the solution retains a single smooth positive lobe and vanishes along the complete boundary. The \OptiXDE{} field is visually indistinguishable from the analytical solution at both times. The absolute-error distributions remain on the order of $10^{-13}$, many orders of magnitude below the corresponding solution amplitudes, and exhibit no coherent phase displacement or spatial bias.

For the higher-frequency mode $(2,3)$, two oscillations in the $x$-direction and three oscillations in the $y$-direction generate six alternating lobes separated by the expected nodal lines. These spatial features are reproduced without visible phase displacement or systematic amplitude distortion. Because the analytical decay rate is $\lambda_{\mathrm{ex}}=13$, the solution amplitude decreases strongly between $t=0.5$ and $t=1.5$; the associated absolute-error scale correspondingly decreases from approximately $10^{-15}$ to $10^{-18}$.

Representative profiles are extracted at the sampled ordinates nearest to $y=\pi/2$ for the $(1,1)$ mode and $y=\pi/6$ for the $(2,3)$ mode. In all four comparisons, the analytical and \OptiXDE{} curves overlap over the complete sampled interval.

\subsubsection{Time-step sensitivity}
\label{supp:diffusion_timestep}

The effect of repeated propagation is investigated at a fixed spatial resolution of $256^2$ by varying the time step from $\Delta t=10^{-1}$ to $\Delta t=10^{-4}$. Since the terminal time is fixed at $T=1.5$, the corresponding number of transform--propagate--inverse-transform cycles, $N_t=T/\Delta t$, increases from $15$ to $15\,000$.

The complete results are reported in Table~\ref{tab:supp_diffusion_timestep}. Reducing the time step does not improve the terminal solution because the modal evolution is already integrated analytically. For the fundamental mode, the terminal field error increases from $7.150\times10^{-16}$ at $\Delta t=10^{-1}$ to $4.416\times10^{-13}$ at $\Delta t=10^{-4}$. Over the same range, the decay-rate error increases from $1.776\times10^{-15}$ to $1.178\times10^{-12}$.

The higher-frequency mode exhibits the same overall behavior. Its terminal field error increases from $4.496\times10^{-19}$ to $1.872\times10^{-17}$, while the decay-rate error increases from $1.776\times10^{-15}$ to $2.203\times10^{-12}$. Its smaller absolute field error again reflects the substantially smaller terminal amplitude produced by the analytical factor $\exp(-13T)$ and should not be interpreted as a lower relative accumulation of round-off error.

\begin{table}[t]
\centering
\caption{Time-step sensitivity of the transient diffusion eigenmode benchmark at a spatial resolution of $256^2$ and $T=1.5$.}
\label{tab:supp_diffusion_timestep}
\begin{tabular}{ccccc}
\toprule
$\Delta t$ & $N_t$ & Mode $(k_x,k_y)$ & $e(T)$ & $\varepsilon_\lambda$ \\
\midrule
$1.000\times10^{-1}$ & 15    & $(1,1)$ & $7.150\times10^{-16}$ & $1.776\times10^{-15}$ \\
$1.000\times10^{-2}$ & 150   & $(1,1)$ & $6.764\times10^{-15}$ & $1.776\times10^{-14}$ \\
$1.000\times10^{-3}$ & 1500  & $(1,1)$ & $4.382\times10^{-14}$ & $1.168\times10^{-13}$ \\
$1.000\times10^{-4}$ & 15000 & $(1,1)$ & $4.416\times10^{-13}$ & $1.178\times10^{-12}$ \\
$1.000\times10^{-1}$ & 15    & $(2,3)$ & $4.496\times10^{-19}$ & $1.776\times10^{-15}$ \\
$1.000\times10^{-2}$ & 150   & $(2,3)$ & $7.695\times10^{-19}$ & $1.954\times10^{-14}$ \\
$1.000\times10^{-3}$ & 1500  & $(2,3)$ & $1.253\times10^{-17}$ & $1.705\times10^{-13}$ \\
$1.000\times10^{-4}$ & 15000 & $(2,3)$ & $1.872\times10^{-17}$ & $2.203\times10^{-12}$ \\
\bottomrule
\end{tabular}
\end{table}

\begin{figure}[t]
\centering
\includegraphics[width=\textwidth]{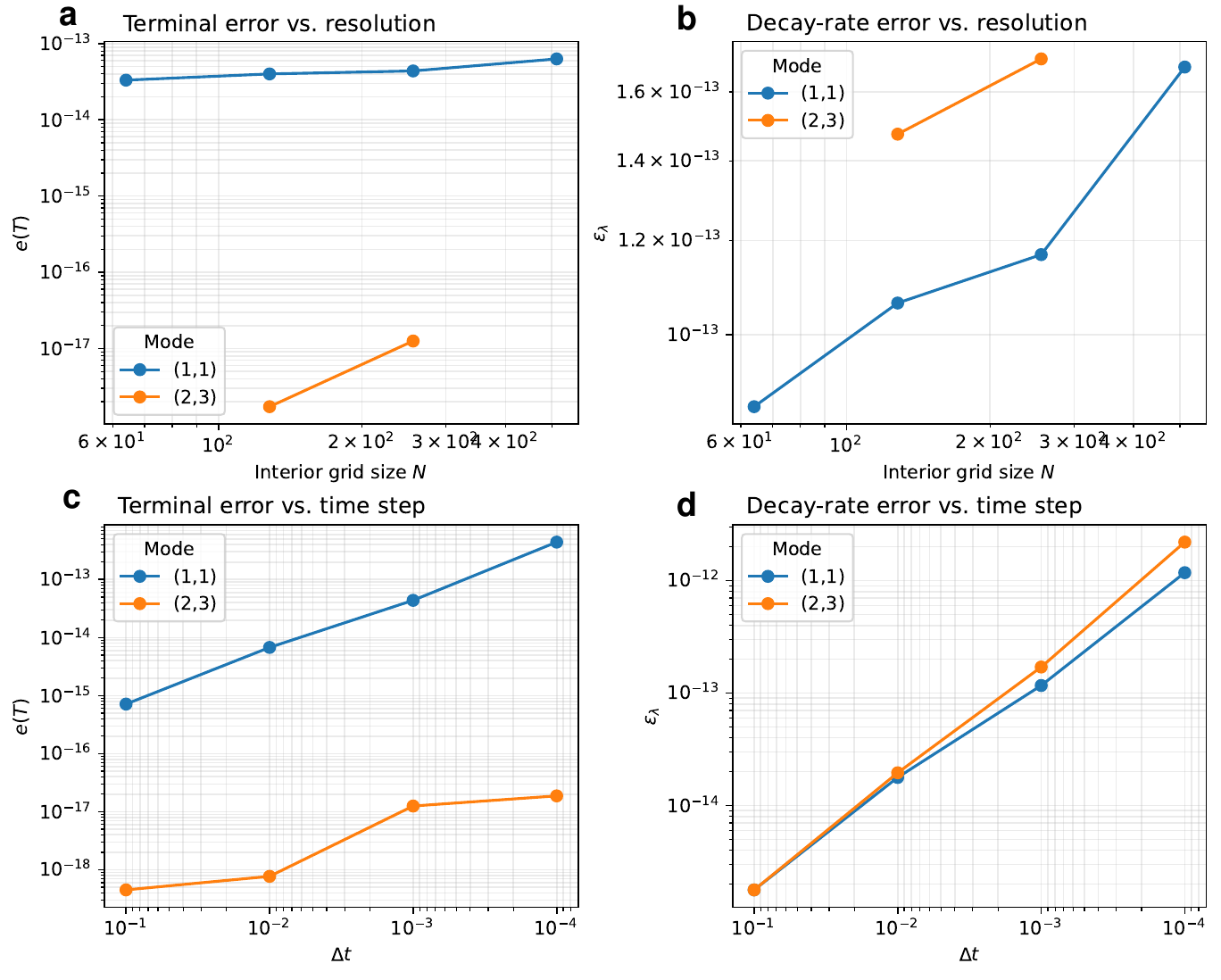}
\caption{\textbf{Spatial-resolution and time-step sensitivity of the transient diffusion benchmark.} \textbf{a}, Terminal error $e(T)$ as a function of spatial resolution for the $(1,1)$ mode and the available $(2,3)$-mode resolutions. \textbf{b}, Corresponding modal decay-rate error $\varepsilon_{\lambda}$. \textbf{c}, Terminal error as a function of the time step at a spatial resolution of $256^2$. \textbf{d}, Corresponding decay-rate error. Decreasing the time step increases the number of repeated transform cycles and consequently increases the accumulated round-off error.}
\label{fig:supp_diffusion_sensitivity}
\end{figure}

Figure~\ref{fig:supp_diffusion_sensitivity}c,d shows that both error measures generally increase as $\Delta t$ decreases. In particular, the decay-rate error grows by approximately one order of magnitude for each tenfold reduction in $\Delta t$, consistent with the corresponding tenfold increase in the number of repeated transforms rather than with a conventional temporal-convergence law.

\subsubsection{Finite-precision interpretation}
\label{supp:diffusion_roundoff}

The observed time-step dependence follows directly from the analytical form of the propagation operator. For a fixed modal decay rate $\lambda=\nu(k_x^2+k_y^2)$, exact propagation over the complete interval $T$ gives
\begin{equation}
\widetilde{u}(T)=\exp(-\lambda T)\widetilde{u}(0).
\label{eq:supp_diffusion_direct_propagation}
\end{equation}

When the same interval is divided into $N_t=T/\Delta t$ increments, the repeated update is formally
\begin{equation}
\widetilde{u}^{N_t}=\left[\exp(-\lambda\Delta t)\right]^{N_t}\widetilde{u}^{0}.
\label{eq:supp_diffusion_repeated_propagation}
\end{equation}

In exact arithmetic, Eqs.~\eqref{eq:supp_diffusion_direct_propagation} and \eqref{eq:supp_diffusion_repeated_propagation} are identical. In finite-precision arithmetic, however, each forward transform, spectral multiplication and inverse transform introduces a small round-off perturbation. Repeating the propagation cycle more frequently therefore increases the accumulated discrepancy.

This behavior differs fundamentally from conventional finite-difference or finite-element time integration, for which reducing $\Delta t$ typically decreases the temporal truncation error until a round-off plateau is reached. In the present benchmark, the modal time evolution is already integrated analytically. Consequently, no temporal truncation error is available to be reduced, and the principal effect of a smaller time step is to increase the number of finite-precision operations.

Taken together, the spatial-resolution and time-step studies show that increasing the spatial resolution does not produce a conventional convergence trend once the prescribed eigenmode is represented by the transform basis, whereas reducing the time step increases the number of finite-precision transform cycles without improving the analytically integrated modal evolution. These observations explain why the measured errors remain within, or close to, the round-off-dominated regime for the present linear constant-coefficient diffusion problem.

\subsubsection{Reproducibility information}
\label{supp:diffusion_reproducibility}

The principal numerical settings used in this benchmark are summarized in Table~\ref{tab:supp_diffusion_settings}. The software versions, hardware specifications, arithmetic precision and timing protocol are reported in Section~\ref{supp:computational_environment}. The source code, input parameters and post-processing scripts required to reproduce this benchmark will be made publicly available upon publication at \href{https://github.com/USERNAME/OptiXDE/tree/main/examples/transient_diffusion}{\texttt{examples/transient\_diffusion}}.

\begin{table}[t]
\centering
\caption{Numerical settings for the transient diffusion eigenmode benchmark.}
\label{tab:supp_diffusion_settings}
\begin{tabular}{ll}
\toprule
Setting & Value \\
\midrule
Spatial domain
&
$\Omega=[0,\pi]^2$
\\

Boundary condition
&
Homogeneous Dirichlet
\\

Diffusion coefficient
&
$\nu=1$
\\

Initial amplitude
&
$10$
\\

Tested modes
&
$(1,1)$ and $(2,3)$
\\

Terminal time
&
$T=1.5$
\\

Sampling convention
&
$N_x\times N_y$ interior samples; zero boundary values imposed separately
\\

Spatial-resolution study
&
$(1,1): 64^2,128^2,256^2,512^2$; $(2,3): 128^2,256^2$
\\

Time step in resolution study
&
$\Delta t=10^{-3}$
\\

Modal-field comparison
&
$512^2$ at $t=0.5$ and $1.5$
\\

Time-step study
&
Spatial resolution $256^2$; $\Delta t=10^{-1}$--$10^{-4}$
\\

Transform basis
&
Tensor-product sine basis
\\

Discrete transform
&
Two-dimensional DST-I
\\

Transform normalization
&
Orthonormal
\\

Update workflow
&
Transform--propagate--inverse transform at every time step
\\

Arithmetic precision
&
Double precision
\\

Propagation operator
&
$\exp[-\nu(k_x^2+k_y^2)\Delta t]$
\\
\bottomrule
\end{tabular}
\end{table}

\FloatBarrier

\subsection{Periodic Poisson equation}
\label{supp:poisson_periodic}

This section provides the complete formulation, spectral implementation and additional numerical results for the periodic Poisson benchmark summarized in the main text. The problem is selected to isolate the direct elliptic inversion kernel of \OptiXDE{} in an FFT-compatible setting, without geometric embedding, boundary penalization or iterative linear algebra.

\subsubsection{Benchmark definition}
\label{supp:poisson_periodic_definition}

The periodic Poisson equation is defined as
\begin{equation}
-\nabla^2u=f
\qquad
\text{in }\Omega,
\label{eq:supp_poisson_periodic_governing}
\end{equation}
where $\Omega=[0,2\pi)^2$ and periodic boundary conditions are imposed in both spatial directions.

A smooth manufactured solution is prescribed as
\begin{equation}
u_{\mathrm{ex}}(x,y)=\sin(x)\sin(y),
\label{eq:supp_poisson_periodic_exact}
\end{equation}
which gives the corresponding source term
\begin{equation}
f(x,y)=-\nabla^2u_{\mathrm{ex}}(x,y)=2\sin(x)\sin(y).
\label{eq:supp_poisson_periodic_source}
\end{equation}

For a periodic Poisson problem, existence requires the source term to have zero spatial mean, while the solution is determined only up to an arbitrary additive constant. In the present benchmark,
\begin{equation}
\int_{\Omega}f(x,y)\,\mathrm{d}\Omega=0,
\qquad
\int_{\Omega}u_{\mathrm{ex}}(x,y)\,\mathrm{d}\Omega=0,
\label{eq:supp_poisson_periodic_zero_mean}
\end{equation}
and uniqueness is imposed through the zero-mean constraint
\begin{equation}
\int_{\Omega}u(x,y)\,\mathrm{d}\Omega=0.
\label{eq:supp_poisson_periodic_constraint}
\end{equation}

The analytical solution consists of a single tensor-product Fourier mode that is represented directly by the periodic spectral basis. This benchmark therefore provides a baseline verification of the direct Fourier-domain inversion used by \OptiXDE{} for linear constant-coefficient elliptic equations.

\subsubsection{Benchmark-specific modal inversion}
\label{supp:poisson_periodic_inversion}

The periodic sampling, wave-number ordering, compatibility condition and treatment of the null mode follow Sections~\ref{supp:discrete_transforms} and~\ref{supp:spectral_propagators}. With $L_x=L_y=2\pi$, the physical wave numbers coincide with the signed integer FFT indices. Defining
\begin{equation}
K^2=k_x^2+k_y^2,
\label{eq:supp_poisson_periodic_benchmark_wavenumber}
\end{equation}
the general inverse-Laplacian relation in Eq.~\eqref{eq:supp_negative_laplacian_inversion} is applied to every nonzero mode, while the zero-frequency coefficient is fixed by Eq.~\eqref{eq:supp_poisson_zero_mean_selection}.

For the source in Eq.~\eqref{eq:supp_poisson_periodic_source}, only the four modes $(k_x,k_y)=(\pm1,\pm1)$ are active and satisfy $K^2=2$. The numerical solution is therefore obtained by one forward transform, one application of the cached inverse-Laplacian multiplier and one inverse transform. No embedded-boundary correction or outer iteration is used.

In addition to the field error defined in Section~\ref{supp:error_measures}, algebraic consistency is evaluated after reconstructing the physical-space solution. The reconstructed spectrum is
\begin{equation}
\widehat{u}_h^{\,\mathrm{rec}}
=
\mathcal{F}
\left[
u_h
\right],
\label{eq:supp_poisson_periodic_reconstructed_spectrum}
\end{equation}
and the physical-space residual and its relative norm are
\begin{equation}
r_h
=
\mathcal{F}^{-1}
\left[
K^2\widehat{u}_h^{\,\mathrm{rec}}-\widehat{f}
\right],
\qquad
\varepsilon_r
=
\frac{
\left\|r_h\right\|_{2,h}
}{
\left\|f\right\|_{2,h}
}.
\label{eq:supp_poisson_periodic_benchmark_residual}
\end{equation}
The zero-frequency coefficient is fixed to zero consistently in the inversion and reconstruction. Because the residual is evaluated after one inverse transform and a second forward transform, it includes the finite-precision effects associated with physical-space reconstruction in addition to the spectral inversion itself.

\subsubsection{Spatial-resolution study}
\label{supp:poisson_periodic_resolution_study}

Uniform periodic sampling is considered at spatial resolutions of $N_x=N_y\in\{64,128,256,512\}$. Because the manufactured solution is represented exactly by the Fourier basis at the coarsest tested resolution, this study does not measure a conventional spatial-convergence rate. Instead, it examines whether changes in transform size introduce systematic deterioration or resolution-dependent bias.

The complete results are reported in Table~\ref{tab:supp_poisson_periodic_resolution}. The discrete $L^2$ error remains between $2.556425\times10^{-15}$ and $2.585305\times10^{-15}$ over all tested resolutions, corresponding to a relative spread of approximately $1.1\%$. The relative algebraic residual $\varepsilon_r$ is evaluated from Eq.~\eqref{eq:supp_poisson_periodic_benchmark_residual} using the reconstructed physical field and the same discrete wave-number array used in the modal inversion.

\begin{table}[t]
\centering
\caption{Spatial-resolution sensitivity and algebraic consistency of the periodic Poisson benchmark for $u_{\mathrm{ex}}(x,y)=\sin(x)\sin(y)$.}
\label{tab:supp_poisson_periodic_resolution}
\begin{tabular*}{\linewidth}{@{\extracolsep{\fill}}ccc@{}}
\toprule
$N_x\times N_y$ & $e$ & $\varepsilon_r$ \\
\midrule
$64^2$ & $2.562192\times10^{-15}$ & $4.585209\times10^{-14}$ \\
$128^2$ & $2.556425\times10^{-15}$ & $2.207580\times10^{-13}$ \\
$256^2$ & $2.576418\times10^{-15}$ & $9.503076\times10^{-13}$ \\
$512^2$ & $2.585305\times10^{-15}$ & $4.403616\times10^{-12}$ \\
\bottomrule
\end{tabular*}
\end{table}

The nearly resolution-independent field error confirms that the active Fourier modes are already represented at the coarsest tested resolution and that the remaining discrepancy is dominated by finite-precision transform operations. Increasing the spatial resolution introduces additional Fourier modes that are absent from the manufactured solution and therefore does not improve the representation of the nonzero $(\pm1,\pm1)$ coefficients. The residual increases from $4.59\times10^{-14}$ to $4.40\times10^{-12}$ as the transform size increases, but remains within a round-off-dominated regime. This moderate growth is consistent with reconstruction noise in weak high-frequency coefficients being amplified by the factor $K^2$ during residual evaluation, rather than with a spatial-discretization error in the resolved analytical mode.

This behavior also distinguishes the periodic Poisson benchmark from the subsequent embedded-domain problem. In the present case, the geometry, boundary conditions and source term are fully compatible with the global Fourier basis. For nonperiodic or non-smooth domains, boundary enforcement, mask regularization and solution singularities introduce additional error sources that may dominate the finite-precision contribution.

\subsubsection{Field and error distributions}
\label{supp:poisson_periodic_fields}

\begin{figure}[t]
\centering
\suppfigure{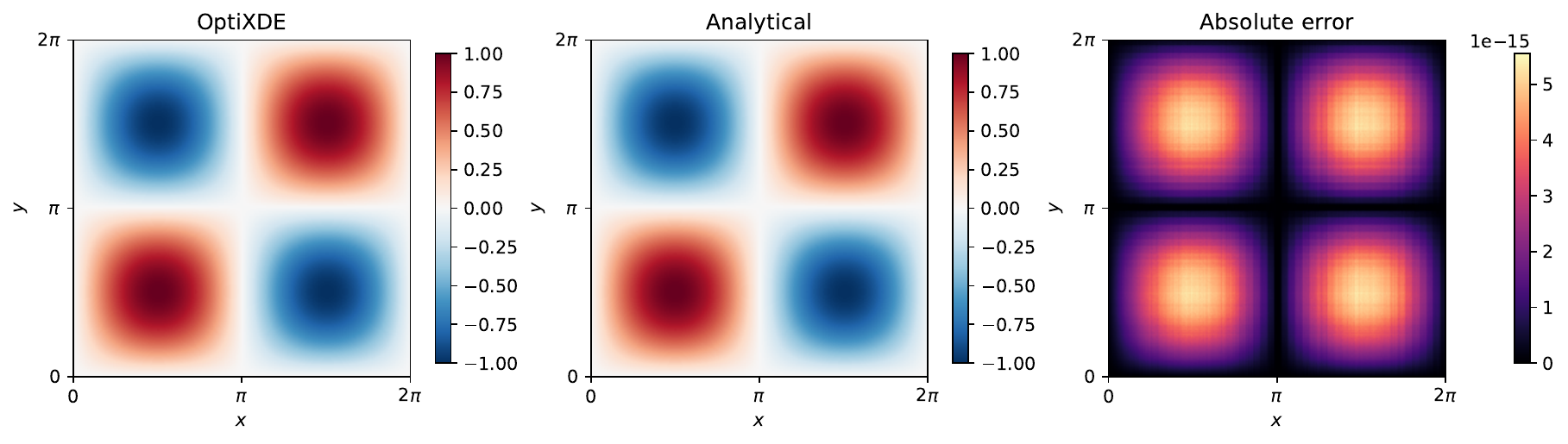}{\textwidth}
\caption{\textbf{Periodic Poisson inversion.} \OptiXDE{} solution, analytical solution and pointwise absolute error for $u_{\mathrm{ex}}(x,y)=\sin(x)\sin(y)$ on the periodic domain $\Omega=[0,2\pi)^2$ at a spatial resolution of $512^2$. The numerical and analytical fields are visually indistinguishable, while the pointwise absolute error remains below approximately $6\times10^{-15}$ and exhibits no boundary-localized or phase-related artifact.}
\label{fig:supp_poisson_periodic_fields}
\end{figure}

Figure~\ref{fig:supp_poisson_periodic_fields} compares the \OptiXDE{} solution, analytical solution and pointwise absolute error at a spatial resolution of $512^2$. The numerical and analytical fields are visually indistinguishable and reproduce the same four-lobe structure, nodal lines and extrema over the complete periodic domain.

The pointwise absolute error remains below approximately $6\times10^{-15}$ throughout the domain. The small structured variations in the error distribution are consistent with finite-precision transformation and reconstruction of the active Fourier modes. Importantly, the error is not concentrated near the domain boundaries and exhibits no boundary layer, phase displacement, spurious oscillation or systematic amplitude attenuation. This behavior is consistent with direct mode-wise inversion of the periodic Laplacian and the absence of embedded-boundary corrections.

Together, Table~\ref{tab:supp_poisson_periodic_resolution} and Fig.~\ref{fig:supp_poisson_periodic_fields} assess both physical-field accuracy and algebraic consistency. The regenerated values of $\varepsilon_r$ remain within a round-off-dominated regime and show the expected transform-size dependence of a residual reconstructed from the physical-space field. This benchmark therefore provides a steady elliptic counterpart to the transient diffusion eigenmode test and verifies the linear spectral kernel of \OptiXDE{} for both direct inversion and analytical time propagation.

\subsubsection{Reproducibility information}
\label{supp:poisson_periodic_reproducibility}

The benchmark-specific numerical settings are summarized in Table~\ref{tab:supp_poisson_periodic_settings}. The common software environment, hardware specifications, arithmetic precision and timing protocol are reported in Section~\ref{supp:computational_environment}. The source code, input parameters, generated CSV data and post-processing scripts used for this benchmark are organized in \texttt{examples/paper/periodic\_poisson}.

\begin{table}[t]
\centering
\caption{Numerical settings for the periodic Poisson benchmark.}
\label{tab:supp_poisson_periodic_settings}
\begin{tabular*}{\linewidth}{@{\extracolsep{\fill}}lp{0.62\linewidth}@{}}
\toprule
Setting & Value \\
\midrule
Spatial domain & $\Omega=[0,2\pi)^2$ \\
Boundary condition & Periodic in both spatial directions \\
Analytical solution & $u_{\mathrm{ex}}(x,y)=\sin(x)\sin(y)$ \\
Source term & $f(x,y)=2\sin(x)\sin(y)$ \\
Solution normalization & Zero spatial mean \\
Spatial resolutions & $64^2$, $128^2$, $256^2$ and $512^2$ \\
Sampling convention & Uniform periodic sampling without duplicated endpoints \\
Transform basis & Two-dimensional Fourier basis \\
Discrete transform & Two-dimensional FFT \\
Zero-frequency treatment & $\widehat{u}(0,0)=0$ \\
Inversion operator & $\mathcal{P}_{\mathrm{P}}=1/K^2$ for $K^2>0$ \\
Residual evaluation & Forward FFT of the reconstructed physical field \\
Residual definition & Eq.~\eqref{eq:supp_poisson_periodic_benchmark_residual} \\
Benchmark directory & \texttt{examples/paper/periodic\_poisson} \\
Result-data file & \texttt{poisson\_periodic\_results.csv} \\
Arithmetic precision & Double precision \\
\bottomrule
\end{tabular*}
\end{table}

\FloatBarrier

\subsection{Poisson equation on an L-shaped domain}
\label{supp:poisson_lshape}

This section provides the complete formulation and additional numerical results for the Poisson problem on an L-shaped domain. In contrast to the periodic Poisson benchmark, the present problem simultaneously examines geometric embedding, nonperiodic Dirichlet enforcement and a re-entrant-corner singularity. The numerical error is therefore governed by the interaction among reduced solution regularity, diffuse-interface representation and penalty enforcement rather than by finite-precision effects alone.

\subsubsection{Benchmark definition}
\label{supp:poisson_lshape_definition}

The physical domain is defined as
\begin{equation}
\Omega=(-1,1)^2\setminus\left([0,1]\times[-1,0]\right),
\label{eq:supp_poisson_lshape_domain}
\end{equation}
which forms an L-shaped region with a re-entrant corner of interior angle $3\pi/2$ at the origin. The governing boundary-value problem is
\begin{equation}
-\nabla^2u=f
\qquad
\text{in }\Omega,
\qquad
u=g
\qquad
\text{on }\partial\Omega.
\label{eq:supp_poisson_lshape_governing}
\end{equation}

Let $(r,\theta)$ denote polar coordinates centred at the re-entrant corner, with $r=\sqrt{x^2+y^2}$ and $0\leq\theta\leq3\pi/2$, where $\theta$ is measured counterclockwise from the positive $x$-axis. The exact reference solution is chosen as the classical re-entrant-corner singular function \cite{mitchell2013collection},
\begin{equation}
u_{\mathrm{ex}}(r,\theta)=r^{2/3}\sin\left(\frac{2\theta}{3}\right).
\label{eq:supp_poisson_lshape_exact}
\end{equation}

Because $\nabla^2u_{\mathrm{ex}}=0$ for $r>0$, the source and boundary data are prescribed as
\begin{equation}
f=0
\qquad
\text{in }\Omega,
\qquad
g=u_{\mathrm{ex}}\big|_{\partial\Omega}.
\label{eq:supp_poisson_lshape_data}
\end{equation}

Although $u_{\mathrm{ex}}$ remains bounded and continuous at the origin, its gradient is singular:
\begin{equation}
\frac{\partial u_{\mathrm{ex}}}{\partial r}
=
\frac{2}{3}r^{-1/3}\sin\left(\frac{2\theta}{3}\right),
\qquad
\frac{1}{r}\frac{\partial u_{\mathrm{ex}}}{\partial\theta}
=
\frac{2}{3}r^{-1/3}\cos\left(\frac{2\theta}{3}\right).
\label{eq:supp_poisson_lshape_gradient}
\end{equation}
Consequently, $|\nabla u_{\mathrm{ex}}|\sim r^{-1/3}$ as $r\rightarrow0$. The solution belongs to $H^{1+2/3-\delta}(\Omega)$ for any $\delta>0$, but not to $H^2(\Omega)$. This reduced regularity prevents the global exponential convergence normally associated with smooth spectral representations and makes the benchmark particularly sensitive to the representation of the embedded boundary and re-entrant corner.

\begin{figure}[t]
\centering
\suppfigure{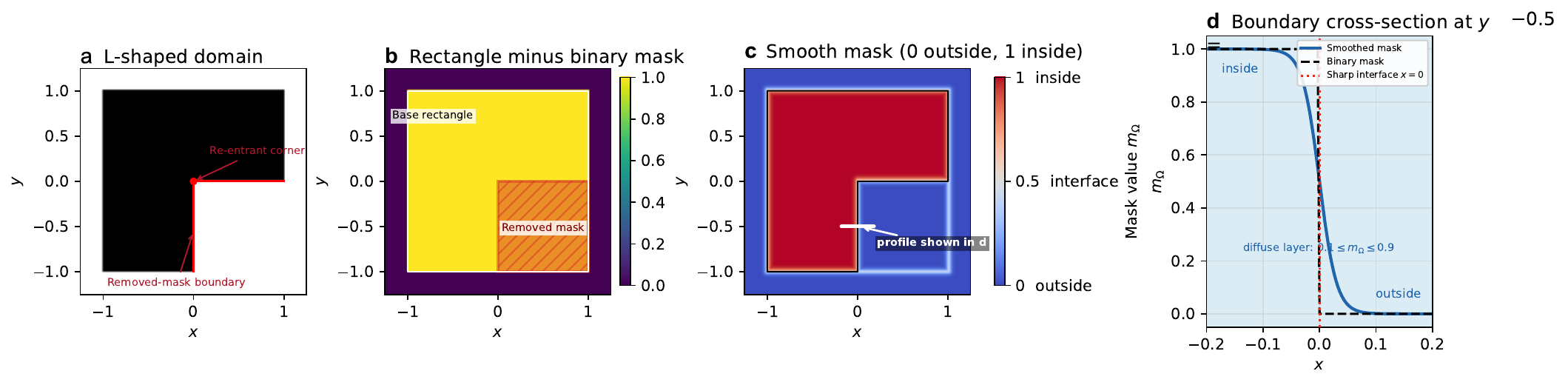}{\textwidth}
\caption{\textbf{Geometry construction and embedded representation of the L-shaped Poisson problem.} \textbf{a}, Physical domain $\Omega=(-1,1)^2\setminus([0,1]\times[-1,0])$, obtained by removing the lower-right rectangular region from the embedding box, with the re-entrant corner at $(0,0)$ and the removed-region boundary highlighted. \textbf{b}, Binary set-difference construction of the L-shaped domain within the embedding box $\widetilde{\Omega}=[-1,1]^2$, showing the base rectangle and the rectangular region to be removed. \textbf{c}, Regularized physical-domain mask $m_{\Omega}$, with $m_{\Omega}=1$ inside $\Omega$ and $m_{\Omega}=0$ outside; the white marker indicates the cross-section examined in panel \textbf{d}. \textbf{d}, Mask profile along $y=-0.5$ across the embedded boundary at $x=0$, comparing the binary and regularized masks. The vertical dashed line marks the sharp interface, and the shaded band denotes the diffuse transition region defined by $0.1\leq m_{\Omega}\leq0.9$.}
\label{fig:supp_poisson_lshape_geometry}
\end{figure}

\subsubsection{Benchmark-specific embedded representation}
\label{supp:poisson_lshape_embedding}

The L-shaped physical domain is embedded in $\widetilde{\Omega}=[-1,1]^2$ and sampled on an endpoint-inclusive Cartesian array,
\begin{equation}
x_i=-1+i\Delta x,
\qquad
y_j=-1+j\Delta y,
\qquad
\Delta x=\frac{2}{N_x-1},
\qquad
\Delta y=\frac{2}{N_y-1},
\qquad
h=\max(\Delta x,\Delta y).
\label{eq:supp_poisson_lshape_benchmark_sampling}
\end{equation}
Here, $N_x$ and $N_y$ denote the total numbers of stored samples including the outer boundary nodes. After the nonhomogeneous outer-boundary data are removed by lifting, the homogeneous sine transform is applied to the $(N_y-2)\times(N_x-2)$ interior-node array.

The signed-distance convention, binary interior mask $m_\Omega^0$, regularized interior mask $m_\Omega$ and exterior mask $m_{\mathrm{ext}}$ are those defined in Section~\ref{supp:geometry_masks}. For this benchmark, the penalty mask is identified as
\begin{equation}
\chi_\varepsilon
\equiv
m_{\mathrm{ext}}
=
1-m_\Omega,
\label{eq:supp_poisson_lshape_benchmark_penalty_mask}
\end{equation}
so that $\chi_\varepsilon\approx0$ in the physical domain and $\chi_\varepsilon\approx1$ in the removed lower-right region.

The mask thickness is prescribed as
\begin{equation}
\varepsilon
=
\max\left(0.03,2.5h\right),
\label{eq:supp_poisson_lshape_benchmark_smoothing}
\end{equation}
which is the specialization of Eq.~\eqref{eq:supp_mask_width}. All geometric arrays are constructed once during preprocessing and subsequently cached. Figure~\ref{fig:supp_poisson_lshape_geometry} illustrates the resulting set-difference geometry and diffuse transition; the general Heaviside regularization and signed-distance convention are not repeated here.

\subsubsection{Benchmark-specific pseudo-time correction}
\label{supp:poisson_lshape_penalized_solution}

The general embedded Dirichlet framework is given in Section~\ref{supp:boundary_enforcement}. For the present benchmark, the nonhomogeneous data on the outer transform box are first removed through the lifting decomposition
\begin{equation}
u=v+T_b,
\qquad
\widetilde{f}
=
f_{\mathrm{ext}}+\nabla^2T_b,
\qquad
\widetilde{g}_{\mathrm{ext}}
=
g_{\mathrm{ext}}-T_b,
\label{eq:supp_poisson_lshape_benchmark_lifting}
\end{equation}
where $T_b$ is the prescribed lifting field. For the present benchmark, $f_{\mathrm{ext}}=0$. The lifting field is constructed from the prescribed values on the four outer edges of the embedding box so that the transformed variable $v$ satisfies homogeneous Dirichlet conditions on $\partial\widetilde{\Omega}$.

The penalty coefficient is
\begin{equation}
\eta=0.003h^2.
\label{eq:supp_poisson_lshape_benchmark_penalty}
\end{equation}
The reference pseudo-time increment is chosen as
\begin{equation}
\Delta\tau=10h^2.
\label{eq:supp_poisson_lshape_pseudotime_step}
\end{equation}

At the $m$-th correction, the constant-coefficient bulk contribution is advanced using a backward-Euler Helmholtz step,
\begin{equation}
v^\star
=
\mathcal{T}^{-1}
\left\{
\frac{
\mathcal{T}
\left[
v^{(m)}+\Delta\tau\,\widetilde{f}
\right]
}{
1+\Delta\tau\,\Lambda
}
\right\},
\label{eq:supp_poisson_lshape_benchmark_helmholtz}
\end{equation}
where $\Lambda$ is the nonnegative spectrum of $-\nabla^2$ for the homogeneous transform on $\widetilde{\Omega}$. In the present implementation, $\mathcal{T}$ is the tensor-product orthonormal type-I discrete sine transform acting on the interior-node array.

The masked penalty term is then advanced by a pointwise backward-Euler correction,
\begin{equation}
v^{(m+1)}
=
\frac{
v^\star
+
\Delta\tau
\left(
\chi_\varepsilon/\eta
\right)
\widetilde{g}_{\mathrm{ext}}
}{
1+
\Delta\tau
\left(
\chi_\varepsilon/\eta
\right)
},
\qquad
u^{(m+1)}
=
v^{(m+1)}+T_b.
\label{eq:supp_poisson_lshape_benchmark_penalty_update}
\end{equation}
Equation~\eqref{eq:supp_poisson_lshape_benchmark_penalty_update} is the benchmark-specific steady-iteration choice and should be distinguished from the exact exponential relaxation map in Eq.~\eqref{eq:supp_exact_dirichlet_penalty}, which is available when the local penalty subproblem is integrated exactly. Both treatments retain a transform-diagonal bulk operator and a pointwise physical-space geometric correction.

Convergence is monitored only over the physical-domain samples:
\begin{equation}
r_\tau^{(m+1)}
=
\frac{
\left\|
 m_\Omega^0
\odot
\left(
u^{(m+1)}-u^{(m)}\right)
\right\|_2
}{
\left\|
 m_\Omega^0
\odot
u^{(m+1)}
\right\|_2
+
\epsilon_{\mathrm{mach}}
}.
\label{eq:supp_poisson_lshape_benchmark_convergence}
\end{equation}
The iteration is terminated when
\begin{equation}
r_\tau^{(m+1)}<10^{-8},
\label{eq:supp_poisson_lshape_stopping_tolerance}
\end{equation}
or when the maximum number of $12000$ outer corrections is reached. The independently initialized calculation uses a zero transformed field, $v^{(0)}=0$.

Thus, the transform-diagonal Helmholtz operator remains independent of the L-shaped geometry, whereas the irregular geometry enters through the cached physical-space mask, lifting field and target field.

\subsubsection{Benchmark-specific diagnostic regions}
\label{supp:poisson_lshape_error_measures}

The global root-mean-square error and regional errors follow Eqs.~\eqref{eq:supp_absolute_l2_error} and~\eqref{eq:supp_regional_l2_error}, respectively. For the present embedded problem, the sample set $\Omega_h$ is selected using the binary physical-domain mask $m_\Omega^0$. The global error is denoted by $e_\Omega$.

To separate the dominant localized error mechanisms, the physical domain is partitioned into three mutually exclusive regions:
\begin{equation}
\Omega_{\mathrm{c}}
=
\left\{
\mathbf{x}\in\Omega:
r\leq r_{\mathrm{c}}
\right\},
\qquad
\Omega_\Gamma
=
\left[
\left\{
\mathbf{x}\in\Omega:
|\phi(\mathbf{x})|\leq w_\Gamma
\right\}
\setminus
\Omega_{\mathrm{c}}
\right],
\qquad
\Omega_{\mathrm{b}}
=
\Omega
\setminus
\left(
\Omega_{\mathrm{c}}\cup\Omega_\Gamma
\right),
\label{eq:supp_poisson_lshape_benchmark_regions}
\end{equation}
with
\begin{equation}
r_{\mathrm{c}}=0.08,
\qquad
w_\Gamma=2\varepsilon.
\label{eq:supp_poisson_lshape_region_parameters}
\end{equation}
Here, $\Omega_{\mathrm{c}}$ isolates the re-entrant-corner neighbourhood, $\Omega_\Gamma$ isolates the diffuse-interface neighbourhood after removal of the corner region, and $\Omega_{\mathrm{b}}$ contains the remaining smooth bulk. Their root-mean-square errors are denoted by $e_{\Omega_{\mathrm{c}}}$, $e_{\Omega_\Gamma}$ and $e_{\Omega_{\mathrm{b}}}$, respectively.

\subsubsection{Spatial-resolution study}
\label{supp:poisson_lshape_grid}

The spatial-resolution study is performed using uniform Cartesian resolutions
\begin{equation}
N_x=N_y\in\left\{256,384,512,768\right\}.
\end{equation}
For every resolution, the reference embedded-boundary parameters are
\begin{equation}
\eta=0.003h^2,
\qquad
\varepsilon=\max(0.03,2.5h).
\label{eq:supp_poisson_lshape_resolution_parameters}
\end{equation}

For all four resolutions considered here, $2.5h<0.03$, and therefore $\varepsilon=0.03$ throughout the spatial-resolution study. Consequently, the results quantify sensitivity to the Cartesian spatial resolution at a fixed physical diffuse-interface thickness. They should not be interpreted as a sharp-interface asymptotic convergence study in which $\varepsilon\rightarrow0$ simultaneously with $h\rightarrow0$.

The global error $e_\Omega$ is evaluated over the complete physical L-shaped domain, whereas $e_{\Omega_{\mathrm{c}}}$, $e_{\Omega_\Gamma}$ and $e_{\Omega_{\mathrm{b}}}$ quantify the errors within the re-entrant-corner, diffuse-interface and smooth-bulk regions defined in Eq.~\eqref{eq:supp_poisson_lshape_benchmark_regions}.

\begin{table}[h]
\centering
\caption{Spatial-resolution study for the L-shaped Poisson benchmark using $\eta=0.003h^2$ and $\varepsilon=\max(0.03,2.5h)$. For all reported resolutions, $\varepsilon=0.03$.}
\label{tab:supp_poisson_lshape_grid}
\begin{tabular*}{\linewidth}{@{\extracolsep{\fill}}ccccc@{}}
\toprule
$N_x\times N_y$ & $e_\Omega$ & $e_{\Omega_{\mathrm{c}}}$ & $e_{\Omega_\Gamma}$ & $e_{\Omega_{\mathrm{b}}}$ \\
\midrule
$256^2$ & $1.5710\times10^{-3}$ & $2.6032\times10^{-5}$ & $1.5051\times10^{-4}$ & $1.7056\times10^{-3}$ \\
$384^2$ & $1.9895\times10^{-4}$ & $1.2482\times10^{-5}$ & $7.4071\times10^{-5}$ & $2.1394\times10^{-4}$ \\
$512^2$ & $3.1452\times10^{-5}$ & $8.5092\times10^{-6}$ & $4.4454\times10^{-5}$ & $2.8671\times10^{-5}$ \\
$768^2$ & $1.0007\times10^{-5}$ & $5.7914\times10^{-6}$ & $2.4958\times10^{-5}$ & $2.2466\times10^{-6}$ \\
\bottomrule
\end{tabular*}
\end{table}

As summarized in Table~\ref{tab:supp_poisson_lshape_grid}, the global error decreases systematically from $1.5710\times10^{-3}$ at $256^2$ to $1.0007\times10^{-5}$ at $768^2$, corresponding to a reduction of more than two orders of magnitude over the tested range. Because the solution contains a re-entrant-corner singularity, the penalty parameter scales with $h^2$, and the diffuse-interface thickness remains fixed at $\varepsilon=0.03$, the error reduction is not expected to follow a single asymptotic power law. The purpose of this study is therefore to establish systematic improvement with increasing spatial resolution and to identify the regions that limit the attainable accuracy.

The regional errors distinguish the different accuracy-limiting mechanisms. At $256^2$, the smooth-bulk error is the largest of the three regional measures, with $e_{\Omega_{\mathrm{b}}}=1.7056\times10^{-3}$. It subsequently decreases rapidly and reaches $2.2466\times10^{-6}$ at $768^2$, indicating that the smooth interior becomes increasingly well resolved. The corner-region error decreases consistently from $2.6032\times10^{-5}$ to $5.7914\times10^{-6}$, although its reduction is slower than that of the smooth bulk because the analytical gradient is singular at the re-entrant corner. The interface-region error also decreases monotonically, from $1.5051\times10^{-4}$ to $2.4958\times10^{-5}$, but remains larger than the corner- and bulk-region errors at the highest tested resolution. The diffuse embedded interface therefore becomes the principal localized accuracy limitation once the smooth interior has been sufficiently resolved.

These region-wise trends motivate the more detailed examination of the field and pointwise-error distributions at the highest tested spatial resolution of $768^2$.

\subsubsection{Field and localized-error distributions}
\label{supp:poisson_lshape_errors}

Two principal localized error mechanisms occur in this benchmark. The first is associated with the re-entrant corner, where the analytical gradient becomes unbounded as $r^{-1/3}$. The second is associated with the diffuse embedded interface, where the sharp physical geometry is represented by a regularized transition of finite thickness $\varepsilon$.

\begin{figure}[h]
\centering
\suppfigure{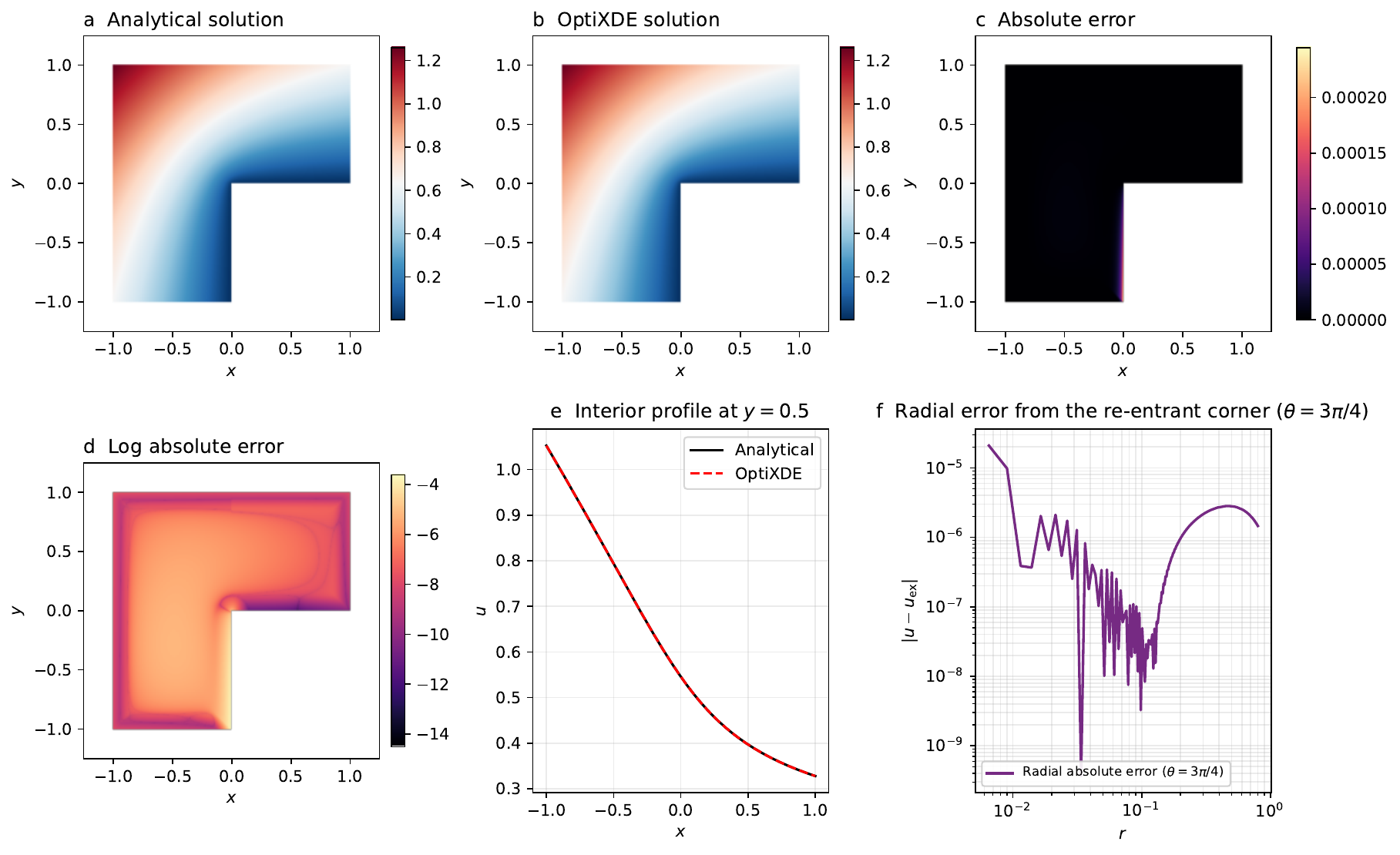}{\textwidth}
\caption{\textbf{Field and localized-error distributions for the L-shaped Poisson problem.} \textbf{a}, Analytical singular solution on the physical L-shaped domain. \textbf{b}, \OptiXDE{} solution at a spatial resolution of $768^2$. \textbf{c}, Pointwise absolute error $\lvert u_h-u_{\mathrm{ex}}\rvert$. \textbf{d}, Base-10 logarithm of the pointwise absolute error, $\log_{10}\lvert u_h-u_{\mathrm{ex}}\rvert$, showing error localization near the embedded boundary and the re-entrant corner. White regions in \textbf{a--d} lie outside the physical domain and are excluded from the error diagnostics. \textbf{e}, Interior profile along $y=0.5$; the analytical and numerical curves overlap visually at the plotted scale. \textbf{f}, Pointwise absolute error along the radial path $(x,y)=r(\cos\theta,\sin\theta)$, measured from the re-entrant corner with $\theta=3\pi/4$ and displayed as a function of $r$ on logarithmic axes.}
\label{fig:supp_poisson_lshape_error_fields}
\end{figure}

Figure~\ref{fig:supp_poisson_lshape_error_fields} presents the detailed field and error distributions at a spatial resolution of $768^2$. The analytical and \OptiXDE{} solutions are visually indistinguishable at the plotting scale. The pointwise absolute error is strongly localized and reaches approximately $2.4\times10^{-4}$ near the inner vertical embedded boundary, whereas the smooth interior is substantially more accurate. The logarithmic error distribution further reveals lower-amplitude discrepancies along the remaining embedded boundary and within the re-entrant-corner neighbourhood.

The interior profile along $y=0.5$ shows complete visual overlap between the analytical and \OptiXDE{} solutions. The radial error profile along $\theta=3\pi/4$ exhibits a distinctly nonmonotonic and multiscale structure. Starting from the order of $10^{-5}$ near the re-entrant corner, the pointwise error decreases rapidly to below $10^{-7}$. This initial reduction is followed by an oscillatory intermediate region, in which local minima fall below $10^{-9}$ while local peaks remain approximately between $10^{-7}$ and $10^{-6}$. At larger radial distances, the error increases again to a broad maximum of a few times $10^{-6}$ before decreasing slightly near the end of the sampled path. The radial discrepancy therefore cannot be represented by a single corner-distance power law over the complete plotted interval, indicating that reduced regularity near the re-entrant corner, diffuse-interface enforcement and global spectral reconstruction contribute over different spatial scales.

Away from the re-entrant corner and the diffuse interface, the analytical solution is harmonic and smooth, and the spectral representation achieves substantially greater local accuracy than suggested by the global error norm. This spatial separation explains why the analytical and numerical fields remain visually indistinguishable over most of the physical domain despite the localized error peaks observed near the geometrically challenging regions.

\subsubsection{Reproducibility information}
\label{supp:poisson_lshape_reproducibility}

The benchmark-specific numerical settings are summarized in Table~\ref{tab:supp_poisson_lshape_settings}. The common software environment, hardware specifications, arithmetic precision and timing protocol are reported in Section~\ref{supp:computational_environment}. Penalty-parameter, mask-smoothing and boundary-localization effects are examined separately in the accuracy and robustness studies rather than repeated within the present benchmark section.

The source code, input parameters and post-processing scripts required to reproduce this benchmark will be made publicly available upon publication at
\href{https://github.com/USERNAME/OptiXDE/tree/main/examples/l_shaped_poisson}{\texttt{examples/l\_shaped\_poisson}}.

\begin{table}[h]
\centering
\caption{Numerical settings for the L-shaped Poisson benchmark.}
\label{tab:supp_poisson_lshape_settings}
\begin{tabular}{ll}
\toprule
Setting & Value \\
\midrule
Physical domain & $\Omega=(-1,1)^2\setminus([0,1]\times[-1,0])$ \\
Embedding domain & $\widetilde{\Omega}=[-1,1]^2$ \\
Analytical solution & $u_{\mathrm{ex}}=r^{2/3}\sin(2\theta/3)$ \\
Source term & $f=0$ \\
Boundary condition & Nonhomogeneous Dirichlet data from $u_{\mathrm{ex}}$ \\
Spatial resolutions & $256^2$, $384^2$, $512^2$ and $768^2$ \\
Stored sampling convention & Endpoint-inclusive uniform Cartesian array \\
Homogeneous transform & Two-dimensional orthonormal DST-I on interior nodes \\
Reference penalty parameter & $\eta_{\mathrm{cells}}=0.003$ \\
Penalty scaling & $\eta=\eta_{\mathrm{cells}}h^2$ \\
Reference mask thickness & $\varepsilon_0=0.03$ \\
Mask regularization & $\varepsilon=\max(\varepsilon_0,2.5h)$ \\
Reference pseudo-time increment & $\Delta\tau=10h^2$ \\
Stopping tolerance & $r_\tau<10^{-8}$ \\
Maximum outer corrections & $12000$ \\
Corner diagnostic radius & $r_{\mathrm{c}}=0.08$ \\
Interface diagnostic width & $w_\Gamma=2\varepsilon$ \\
Independent initialization & $v^{(0)}=0$ \\
Arithmetic precision & Double precision \\
\bottomrule
\end{tabular}
\end{table}

\FloatBarrier

\subsection{Nonlinear Schr\"odinger equation}
\label{supp:schrodinger}

This section examines the nonlinear extension of \OptiXDE{} using the one-dimensional focusing nonlinear Schr\"odinger equation. In contrast to the preceding linear benchmarks, the present problem combines dispersive spectral propagation with an amplitude-dependent nonlinear phase rotation. It therefore provides a controlled test of nonlinear operator splitting, spatial spectral resolution, temporal convergence, nonlinear focusing and recurrence.

\subsubsection{Benchmark definition}
\label{supp:schrodinger_definition}

The one-dimensional focusing nonlinear Schr\"odinger equation is
\begin{equation}
\mathrm{i}
\frac{\partial h}{\partial t}
+
\frac{1}{2}
\frac{\partial^2 h}{\partial x^2}
+
|h|^2 h
=
0,
\qquad
(x,t)\in\Omega\times(0,T],
\label{eq:supp_schrodinger_governing}
\end{equation}
where $h(x,t)\in\mathbb{C}$ is the complex-valued wave field and $\mathrm{i}^2=-1$.

The computational domain and terminal time are
\begin{equation}
\Omega=[-15,15),
\qquad
T=\frac{\pi}{2}.
\label{eq:supp_schrodinger_domain}
\end{equation}

Periodic boundary conditions are imposed,
\begin{equation}
h(-15,t)=h(15,t),
\qquad
\frac{\partial h}{\partial x}(-15,t)
=
\frac{\partial h}{\partial x}(15,t).
\label{eq:supp_schrodinger_boundary}
\end{equation}

The initial condition is
\begin{equation}
h(x,0)
=
2\operatorname{sech}(x).
\label{eq:supp_schrodinger_initial}
\end{equation}

This initial field generates the classical second-order focusing soliton and is widely used as a nonlinear Schr\"odinger benchmark \citep{raissi2019physics}. The interval used here is deliberately larger than the shorter domain commonly adopted in learning-based benchmarks. Because the initial field decays exponentially away from the origin, the enlarged interval places its tails close to zero at the periodic boundaries and thereby reduces finite-domain contamination of the Fourier representation.

On the whole line, the corresponding analytical second-order soliton is
\begin{equation}
h_{\infty}(x,t)
=
4
\exp
\left(
\frac{\mathrm{i}t}{2}
\right)
\frac{
\cosh(3x)
+
3\exp(4\mathrm{i}t)\cosh(x)
}{
\cosh(4x)
+
4\cosh(2x)
+
3\cos(4t)
}.
\label{eq:supp_schrodinger_whole_line_solution}
\end{equation}

At $t=0$, Eq.~\eqref{eq:supp_schrodinger_whole_line_solution} reduces to Eq.~\eqref{eq:supp_schrodinger_initial}. The strongest focusing occurs at
\begin{equation}
t_{\mathrm{f}}
=
\frac{\pi}{4},
\label{eq:supp_schrodinger_focus_time}
\end{equation}
where the analytical peak amplitude is four. At the terminal time,
\begin{equation}
h_{\infty}
\left(
x,\frac{\pi}{2}
\right)
=
\exp
\left(
\frac{\mathrm{i}\pi}{4}
\right)
h_{\infty}(x,0),
\label{eq:supp_schrodinger_recurrence_exact}
\end{equation}
so that the initial amplitude profile is recovered with a spatially uniform phase rotation.

The whole-line expression is used as an additional physical diagnostic rather than as the principal numerical reference. The primary numerical errors are evaluated against an independently integrated periodic reference solution posed on the same finite computational interval.

The periodic spatial discretization is
\begin{equation}
x_j
=
-15+j\Delta x,
\qquad
\Delta x
=
\frac{30}{N_x},
\qquad
j=0,\ldots,N_x-1.
\label{eq:supp_schrodinger_grid}
\end{equation}

The principal temporal-convergence calculations use
\begin{equation}
N_x
=
2048,
\label{eq:supp_schrodinger_test_resolution}
\end{equation}
whereas the spatial-resolution study considers
\begin{equation}
N_x
\in
\left\{
256,\,
512,\,
1024,\,
2048
\right\}.
\label{eq:supp_schrodinger_spatial_resolutions}
\end{equation}

All convergence and reference calculations use strictly fixed time increments. Diagnostic output is recorded at the common times
\begin{equation}
t_q
=
\frac{qT}{256},
\qquad
q=0,\ldots,256,
\label{eq:supp_schrodinger_history_times}
\end{equation}
so that output sampling does not modify the underlying time-integration sequence.

\subsubsection{Benchmark-specific split-step settings}
\label{supp:schrodinger_splitting}

The nonlinear Schr\"odinger calculation uses the second-order Strang composition defined in Eq.~\eqref{eq:supp_strang_splitting}. The corresponding linear half-step and exact nonlinear phase rotation are given in Eqs.~\eqref{eq:supp_nls_first_linear_half_step}--\eqref{eq:supp_nls_second_linear_half_step} and are not repeated here.

For the present periodic interval, the domain length is
\begin{equation}
L
=
30,
\end{equation}
and the Fourier wave numbers are
\begin{equation}
k_m
=
\frac{2\pi}{30}\mu_m,
\label{eq:supp_schrodinger_wavenumbers}
\end{equation}
where $\mu_m$ follows the signed FFT ordering defined in Section~\ref{supp:discrete_transforms}.

For each prescribed fixed time increment $\Delta t$, the linear half-step multiplier implied by Eqs.~\eqref{eq:supp_nls_first_linear_half_step} and~\eqref{eq:supp_nls_second_linear_half_step} is constructed once and reused throughout the simulation. Both component subflows are evaluated in closed form and no nonlinear iteration is required. The baseline calculations are performed without spectral filtering or de-aliasing. Filtering sensitivity and invariant-preservation behaviour are considered separately in the accuracy and robustness study and are not used in the accuracy results reported below.

\subsubsection{Periodic reference solution and verification}
\label{supp:schrodinger_reference}

An independent Fourier pseudo-spectral fourth-order Runge--Kutta calculation is used to construct the periodic reference solution. The spatial and temporal resolutions are refined simultaneously through the hierarchy
\begin{equation}
\left(
N_x^{\mathrm{ref}},
N_t^{\mathrm{ref}}
\right)
=
\left\{
(1024,65536),
(2048,131072),
(4096,262144)
\right\}.
\label{eq:supp_schrodinger_reference_hierarchy}
\end{equation}

Because
\begin{equation}
\Delta t_{\mathrm{ref}}
=
\frac{T}{N_t^{\mathrm{ref}}},
\label{eq:supp_schrodinger_reference_dt}
\end{equation}
the corresponding time increments are approximately
\begin{equation}
2.396845\times10^{-5},
\qquad
1.198422\times10^{-5},
\qquad
5.992112\times10^{-6},
\end{equation}
respectively.

The finest calculation,
\begin{equation}
N_x^{\mathrm{ref}}
=
4096,
\qquad
N_t^{\mathrm{ref}}
=
262144,
\qquad
\Delta t_{\mathrm{ref}}
=
5.992112\times10^{-6},
\label{eq:supp_schrodinger_reference_finest}
\end{equation}
is used as the principal periodic reference.

The independent reference integrator advances the Fourier pseudo-spectral semi-discretization
\begin{equation}
\frac{\dd h}{\dd t}
=
\mathcal{F}^{-1}
\left[
-\frac{\mathrm{i}}{2}
k^2\widehat{h}
\right]
+
\mathrm{i}|h|^2h
\label{eq:supp_schrodinger_reference_rhs}
\end{equation}
using classical fourth-order Runge--Kutta time integration.

When two reference levels have different spatial resolutions, the finer solution is restricted to the coarser Fourier space by retaining the centered resolvable modes. Let $\mathcal{R}_{N}$ denote this periodic spectral restriction. The successive-reference difference is defined as
\begin{equation}
\delta_{\mathrm{ref}}^{(r)}
=
\max_q
\frac{
\left\|
h_{R_r}(\cdot,t_q)
-
\mathcal{R}_{N_r}
h_{R_{r+1}}(\cdot,t_q)
\right\|_2
}{
\left\|
\mathcal{R}_{N_r}
h_{R_{r+1}}(\cdot,t_q)
\right\|_2
},
\label{eq:supp_schrodinger_reference_difference}
\end{equation}
where $R_r$ and $R_{r+1}$ denote successive reference levels.

\begin{table}[h]
\centering
\caption{Independent-reference refinement for the nonlinear Schr\"odinger benchmark.}
\label{tab:supp_schrodinger_reference_verification}
\begin{tabular*}{\linewidth}{@{\extracolsep{\fill}}cccccc@{}}
\toprule
Level & $N_x$ & $N_t$ & $\Delta t_{\mathrm{ref}}$ & $\delta_{\mathrm{ref}}$ & $\varepsilon_E$ \\
\midrule
R1 & $1024$ & $65536$ & $2.397\times10^{-5}$ & $3.568614\times10^{-10}$ & $2.017434\times10^{-14}$ \\
R2 & $2048$ & $131072$ & $1.198\times10^{-5}$ & $1.417813\times10^{-10}$ & $1.408397\times10^{-14}$ \\
R3 & $4096$ & $262144$ & $5.992\times10^{-6}$ & -- & $9.325873\times10^{-15}$ \\
\bottomrule
\end{tabular*}
\end{table}

As summarized in Table~\ref{tab:supp_schrodinger_reference_verification}, the successive-reference differences are already below $4\times10^{-10}$ and decrease to $1.42\times10^{-10}$ between the two finest levels. The corresponding energy drifts remain at approximately $10^{-14}$. These results show that the independently generated periodic reference is substantially more accurate than the errors measured in the subsequent \OptiXDE{} calculations and therefore does not control the reported convergence trends.

When the reference and test resolutions differ, the principal $4096$-point field is transferred to the test resolution through the same spectral restriction procedure rather than through physical-space interpolation. This preserves the Fourier coefficients resolvable on the target discretization and avoids interpolation-induced amplitude or phase errors.

\begin{figure}[t]
\centering
\suppfigure{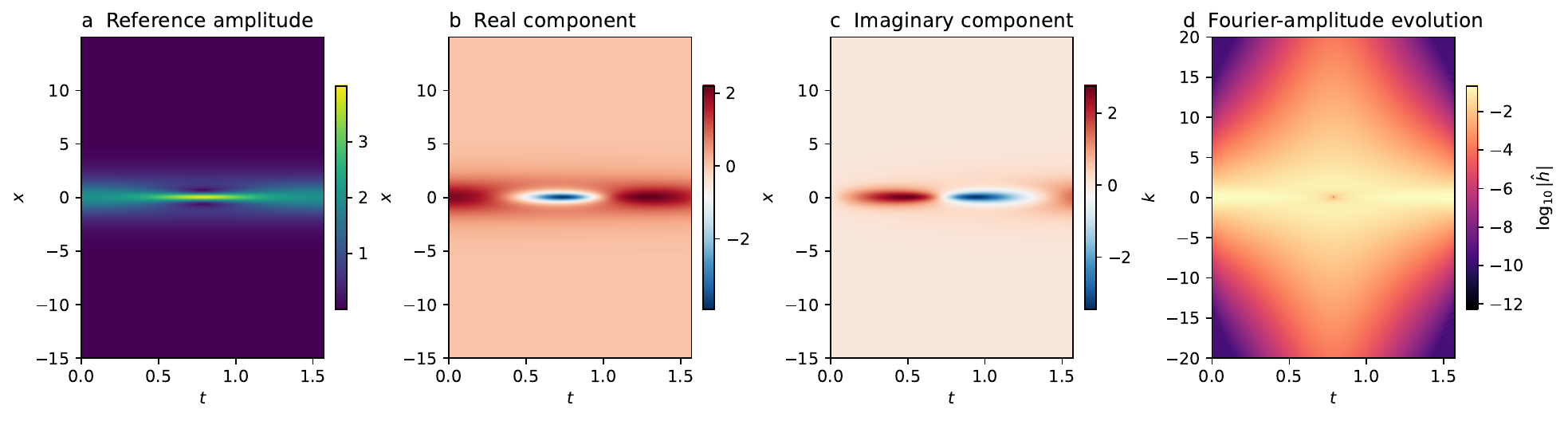}{\textwidth}
\caption{\textbf{Periodic reference evolution of the nonlinear Schr\"odinger benchmark.} \textbf{a}, Space--time distribution of the reference amplitude $\lvert h_{\mathrm{ref}}(x,t)\rvert$ on $\Omega=[-15,15)$. \textbf{b}, Real component of the reference field. \textbf{c}, Imaginary component. \textbf{d}, Logarithmic Fourier-amplitude evolution over the displayed wave-number interval $-20\leq k\leq20$. The spectrum broadens during nonlinear focusing near $t=\pi/4$ and contracts as the field approaches recurrence.}
\label{fig:supp_schrodinger_reference}
\end{figure}

Figure~\ref{fig:supp_schrodinger_reference} shows the complete principal reference evolution. The initially localized pulse undergoes strong nonlinear focusing near $t=\pi/4$ and subsequently broadens toward recurrence at $T=\pi/2$. The Fourier spectrum broadens substantially during focusing and contracts afterward, demonstrating reversible transfer toward higher resolved wave numbers.

\subsubsection{Spatial-resolution study}
\label{supp:schrodinger_spatial_resolution}

Spatial resolution is examined using
\begin{equation}
N_x
\in
\left\{
256,\,
512,\,
1024,\,
2048
\right\},
\end{equation}
with a common fixed step count
\begin{equation}
N_t
=
12800.
\end{equation}
The corresponding time increment is
\begin{equation}
\Delta t
=
\frac{\pi/2}{12800}
=
1.227185\times10^{-4}.
\label{eq:supp_schrodinger_spatial_dt}
\end{equation}

To measure successive spatial convergence independently of the periodic reference, the difference between a resolution $N$ and the next finer resolution $2N$ is defined as
\begin{equation}
\delta_N
=
\max_q
\frac{
\left\|
h_N(\cdot,t_q)
-
\mathcal{R}_{N}
h_{2N}(\cdot,t_q)
\right\|_2
}{
\left\|
\mathcal{R}_{N}
h_{2N}(\cdot,t_q)
\right\|_2
}.
\label{eq:supp_schrodinger_successive_spatial_difference}
\end{equation}

The terminal relative error with respect to the periodic reference is
\begin{equation}
e_{2,\mathrm{per}}(T)
=
\frac{
\left\|
h_h(\cdot,T)
-
h_{\mathrm{ref}}^{(N_x)}(\cdot,T)
\right\|_2
}{
\left\|
h_{\mathrm{ref}}^{(N_x)}(\cdot,T)
\right\|_2
}.
\label{eq:supp_schrodinger_periodic_terminal_error}
\end{equation}

To compare with the analytical whole-line solution while excluding the distant periodic boundaries, the interior region
\begin{equation}
\Omega_{\mathrm{int}}
=
\left\{
x:
|x|\leq10
\right\}
\label{eq:supp_schrodinger_interior_region}
\end{equation}
is introduced. The corresponding terminal relative error is
\begin{equation}
e_{2,\infty}(T)
=
\frac{
\left\|
h_h(\cdot,T)
-
h_{\infty}(\cdot,T)
\right\|_{2,\Omega_{\mathrm{int}}}
}{
\left\|
h_{\infty}(\cdot,T)
\right\|_{2,\Omega_{\mathrm{int}}}
}.
\label{eq:supp_schrodinger_whole_line_terminal_error}
\end{equation}

The high-frequency spectral-tail fraction is defined as
\begin{equation}
\rho_{\mathrm{tail}}(t)
=
\frac{
\displaystyle
\sum_{|\mu_m|\geq N_x/3}
\left|
\widehat{h}_m(t)
\right|^2
}{
\displaystyle
\sum_m
\left|
\widehat{h}_m(t)
\right|^2
},
\label{eq:supp_schrodinger_tail_ratio}
\end{equation}
with
\begin{equation}
\rho_{\mathrm{tail}}^{\max}
=
\max_q
\rho_{\mathrm{tail}}(t_q).
\label{eq:supp_schrodinger_tail_max}
\end{equation}

The focusing-amplitude error is
\begin{equation}
e_{\mathrm{peak}}
=
\left|
\max_x
\left|
h_h
\left(
x,\frac{\pi}{4}
\right)
\right|
-
4
\right|,
\label{eq:supp_schrodinger_peak_error}
\end{equation}
and the recurrence error is
\begin{equation}
e_{\mathrm{rec}}
=
\frac{
\left\|
h_h(\cdot,T)
-
\exp
\left(
\mathrm{i}\pi/4
\right)
h_h(\cdot,0)
\right\|_{2,\Omega_{\mathrm{int}}}
}{
\left\|
h_h(\cdot,0)
\right\|_{2,\Omega_{\mathrm{int}}}
}.
\label{eq:supp_schrodinger_recurrence_error}
\end{equation}

\begin{table}[h]
\centering
\scriptsize
\caption{Controlled spatial-resolution results at the fixed time step
$\Delta t=1.227185\times10^{-4}$.}
\label{tab:supp_nls_spatial_resolution}

\begin{tabular*}{\linewidth}{@{\extracolsep{\fill}}cccc@{}}
\toprule
$N_x$ & $\delta_N$ & $e_{2,\mathrm{per}}(T)$ & $e_{2,\infty}(T)$ \\
\midrule
$256$  & $8.988691\times10^{-5}$ & $5.234771\times10^{-5}$ & $5.234774\times10^{-5}$ \\
$512$  & $3.232984\times10^{-9}$ & $1.147892\times10^{-6}$ & $1.148834\times10^{-6}$ \\
$1024$ & $3.568402\times10^{-10}$ & $1.147891\times10^{-6}$ & $1.148828\times10^{-6}$ \\
$2048$ & -- & $1.147895\times10^{-6}$ & $1.148826\times10^{-6}$ \\
\bottomrule
\end{tabular*}

\vspace{0.6em}

\begin{tabular*}{\linewidth}{@{\extracolsep{\fill}}cccc@{}}
\toprule
$N_x$ & $\rho_{\mathrm{tail}}^{\max}$ & $e_{\mathrm{peak}}$ & $e_{\mathrm{rec}}$ \\
\midrule
$256$  & $3.235999\times10^{-6}$  & $2.545535\times10^{-4}$ & $5.234774\times10^{-5}$ \\
$512$  & $4.989137\times10^{-12}$ & $3.425421\times10^{-7}$ & $1.148834\times10^{-6}$ \\
$1024$ & $4.822414\times10^{-19}$ & $3.440816\times10^{-7}$ & $1.148828\times10^{-6}$ \\
$2048$ & $6.017475\times10^{-20}$ & $3.442807\times10^{-7}$ & $1.148826\times10^{-6}$ \\
\bottomrule
\end{tabular*}

\end{table}

The spatial study shows a clear transition between the $256$-point calculation and the higher resolutions. The successive difference decreases from $8.99\times10^{-5}$ between $256$ and $512$ points to $3.23\times10^{-9}$ between $512$ and $1024$ points, and then to $3.57\times10^{-10}$ between $1024$ and $2048$ points. The maximum spectral-tail fraction simultaneously decreases from $3.24\times10^{-6}$ at $N_x=256$ to $4.99\times10^{-12}$ at $N_x=512$ and below $10^{-18}$ at the two finest resolutions.

The periodic-reference and whole-line terminal errors decrease to approximately $1.15\times10^{-6}$ at $N_x=512$ and remain essentially unchanged under further spatial refinement. Because the successive spatial differences are several orders of magnitude smaller at the higher resolutions, this plateau is not caused by inadequate Fourier resolution. Instead, it reflects the remaining temporal splitting error associated with the fixed time step used in the spatial study. The focusing-amplitude error shows the same behaviour, decreasing from $2.55\times10^{-4}$ at $N_x=256$ to approximately $3.4\times10^{-7}$ for $N_x\geq512$.

These results establish that the $N_x=2048$ discretization used in the temporal-convergence study is well within the spatially resolved regime.

\begin{figure}[h]
\centering
\suppfigure{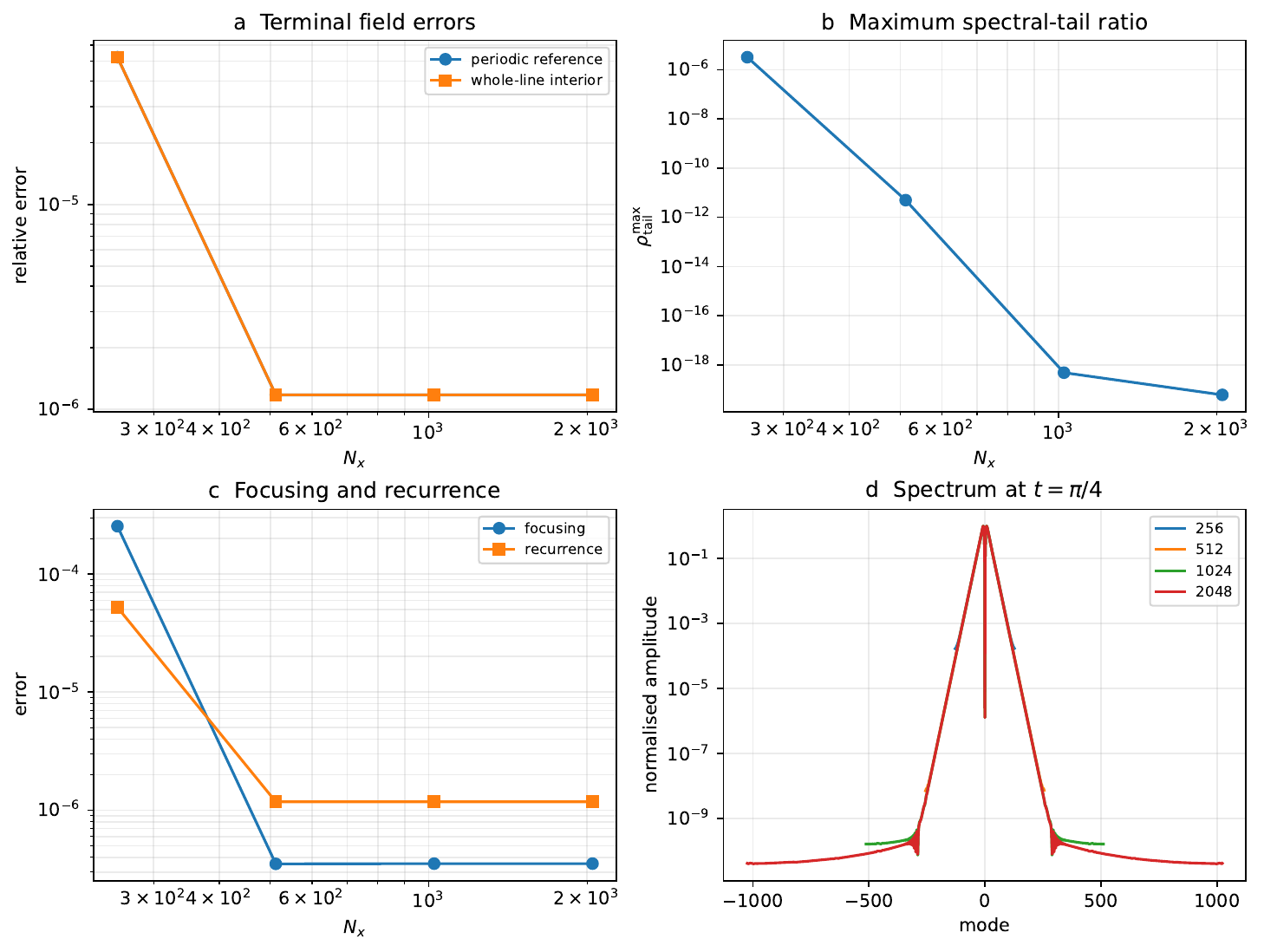}{\textwidth}
\caption{\textbf{Spatial-resolution assessment of the nonlinear Schr\"odinger benchmark.} \textbf{a}, Terminal relative errors with respect to the periodic numerical reference and the analytical whole-line solution in the interior region $|x|\leq10$. \textbf{b}, Maximum high-frequency spectral-tail ratio. \textbf{c}, Focusing-amplitude and recurrence errors. \textbf{d}, Normalized Fourier spectra at $t=\pi/4$ for the tested spatial resolutions.}
\label{fig:supp_schrodinger_spatial_resolution}
\end{figure}

\subsubsection{Temporal convergence}
\label{supp:schrodinger_convergence}

Temporal convergence is evaluated at the fixed spatial resolution
\begin{equation}
N_x
=
2048.
\end{equation}

The total numbers of fixed Strang steps are
\begin{equation}
N_t
\in
\left\{
256,\,
512,\,
1024,\,
2048,\,
4096
\right\},
\label{eq:supp_schrodinger_temporal_step_counts}
\end{equation}
giving
\begin{equation}
\Delta t
=
\frac{T}{N_t}
\in
\left\{
6.135923\times10^{-3},
3.067962\times10^{-3},
1.533981\times10^{-3},
7.669904\times10^{-4},
3.834952\times10^{-4}
\right\}.
\label{eq:supp_schrodinger_benchmark_time_steps}
\end{equation}

Thus, every refinement halves the time increment exactly while preserving the same terminal time and common diagnostic-output schedule.

At each stored comparison time, the pointwise error, maximum recorded space--time error and relative discrete $L^2$ error are evaluated using Eqs.~\eqref{eq:supp_pointwise_error}, \eqref{eq:supp_maximum_spacetime_error} and~\eqref{eq:supp_relative_l2_error}, respectively, with $u_h=h_h$ and $u_{\mathrm{ref}}=h_{\mathrm{ref}}$. The observed temporal order is evaluated from Eq.~\eqref{eq:supp_observed_convergence_order} with $q=\Delta t$ and $e=e_{\max}$.

\begin{table}[h]
\centering
\small
\caption{Temporal convergence of the nonlinear Schr\"odinger benchmark at $N_x=2048$.}
\label{tab:supp_schrodinger_convergence}
\begin{tabular*}{\linewidth}{@{\extracolsep{\fill}}cccc@{}}
\toprule
$N_t$ & $\Delta t$ & $e_{\max}$ & Order $p$ \\
\midrule
$256$ & $6.135923\times10^{-3}$ & $1.100375\times10^{-2}$ & -- \\
$512$ & $3.067962\times10^{-3}$ & $2.761080\times10^{-3}$ & $1.995$ \\
$1024$ & $1.533981\times10^{-3}$ & $6.909122\times10^{-4}$ & $1.999$ \\
$2048$ & $7.669904\times10^{-4}$ & $1.727683\times10^{-4}$ & $2.000$ \\
$4096$ & $3.834952\times10^{-4}$ & $4.319465\times10^{-5}$ & $2.000$ \\
\bottomrule
\end{tabular*}
\end{table}

The temporal results exhibit a uniform second-order convergence trend over the complete refinement sequence. Successive time-step halvings reduce $e_{\max}$ from $1.100375\times10^{-2}$ to $2.761080\times10^{-3}$, $6.909122\times10^{-4}$, $1.727683\times10^{-4}$ and finally $4.319465\times10^{-5}$. The corresponding observed orders are
\begin{equation}
1.995,
\qquad
1.999,
\qquad
2.000,
\qquad
2.000,
\end{equation}
in close agreement with the theoretical second-order accuracy of the Strang composition.

Unlike the previous fixed-resolution calculation in which spatial and reference errors obscured the asymptotic temporal regime, the enlarged domain, independently verified reference and spatially resolved $N_x=2048$ discretization provide a clean temporal-convergence sequence over the complete tested range.

\begin{figure}[h]
\centering
\suppfigure{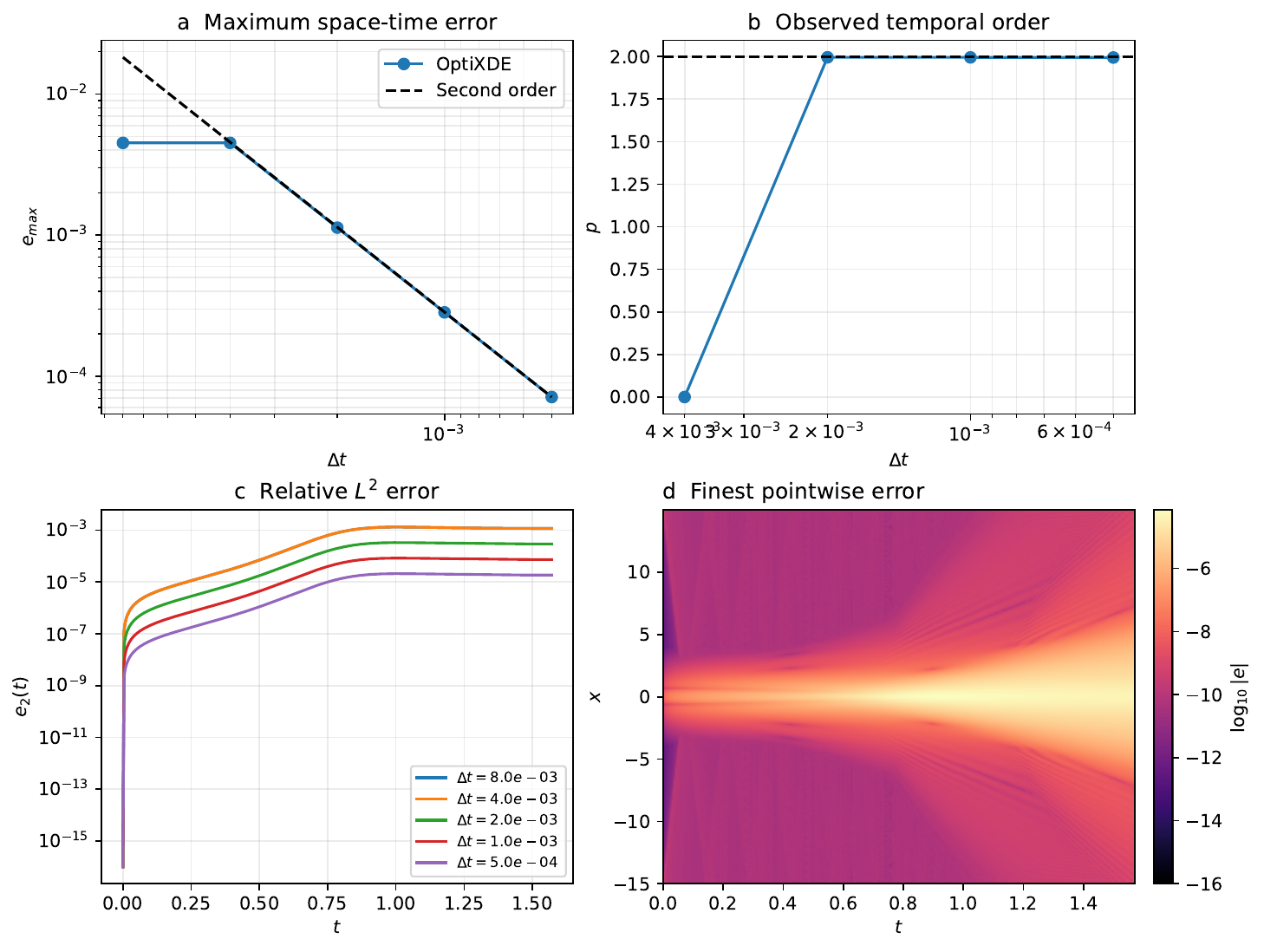}{\textwidth}
\caption{\textbf{Temporal convergence of the nonlinear Schr\"odinger benchmark at $N_x=2048$.} \textbf{a}, Maximum recorded space--time error $e_{\max}$ as a function of the fixed time increment, together with a second-order reference slope. \textbf{b}, Observed temporal order between successive refinements. \textbf{c}, Relative $L^2$ error histories for the tested time increments. \textbf{d}, Logarithmic pointwise complex-error distribution for the finest calculation. The complete refinement sequence exhibits second-order temporal convergence.}
\label{fig:supp_schrodinger_convergence}
\end{figure}

\subsubsection{Focusing, recurrence and space--time fields}
\label{supp:schrodinger_snapshots}

The finest temporal calculation,
\begin{equation}
N_x
=
2048,
\qquad
N_t
=
4096,
\qquad
\Delta t
=
3.834952\times10^{-4},
\label{eq:supp_schrodinger_snapshot_case}
\end{equation}
is used for the detailed field comparisons.

Three diagnostic times aligned exactly with the common output schedule are selected:
\begin{equation}
t
\in
\left\{
\frac{3\pi}{16},
\frac{\pi}{4},
\frac{5\pi}{16}
\right\}.
\label{eq:supp_schrodinger_snapshot_times}
\end{equation}

These correspond approximately to
\begin{equation}
t
\approx
0.5890,
\qquad
0.7854,
\qquad
0.9817,
\end{equation}
and represent the approach to focusing, the maximum-focusing state and the subsequent broadening stage.

At each selected time, the amplitude
\begin{equation}
A(x,t)
=
|h(x,t)|,
\label{eq:supp_schrodinger_amplitude}
\end{equation}
the real component
\begin{equation}
u(x,t)
=
\operatorname{Re}
\left[
h(x,t)
\right],
\label{eq:supp_schrodinger_real_component}
\end{equation}
and the imaginary component
\begin{equation}
v(x,t)
=
\operatorname{Im}
\left[
h(x,t)
\right]
\label{eq:supp_schrodinger_imaginary_component}
\end{equation}
are compared directly with the spectrally restricted periodic reference.

The real- and imaginary-component comparisons provide a phase-sensitive verification that is more stringent than amplitude agreement alone. For the finest calculation, the maximum complex error over all $257$ stored comparison times and all spatial points is
\begin{equation}
e_{\max}
=
4.319465\times10^{-5}.
\label{eq:supp_schrodinger_finest_emax}
\end{equation}

\begin{figure}[t]
\centering
\suppfigure{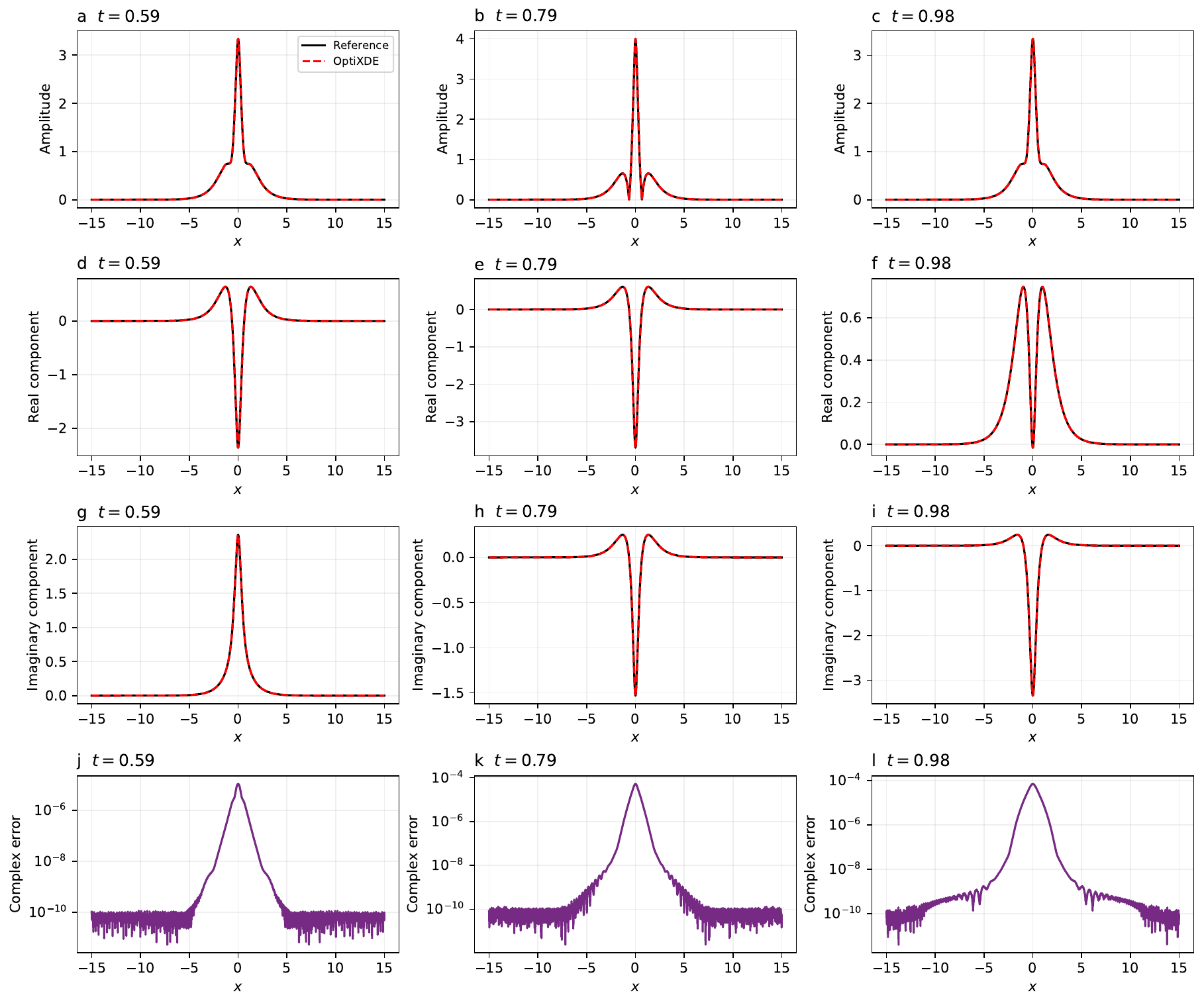}{\textwidth}
\caption{\textbf{Temporal snapshots of the nonlinear Schr\"odinger solution for $N_x=2048$ and $\Delta t=3.834952\times10^{-4}$.} \textbf{a--c}, Reference and \OptiXDE{} amplitude profiles at $t=3\pi/16$, $\pi/4$ and $5\pi/16$, respectively. \textbf{d--f}, Corresponding real components. \textbf{g--i}, Corresponding imaginary components. \textbf{j--l}, Pointwise complex errors. The comparisons verify both amplitude and phase through the nonlinear focusing cycle.}
\label{fig:supp_schrodinger_snapshots}
\end{figure}

Figure~\ref{fig:supp_schrodinger_snapshots} shows that the numerical solution reproduces the strong central focusing at $t=\pi/4$ together with the associated phase evolution before and after the focusing event.

\begin{figure}[t]
\centering
\suppfigure{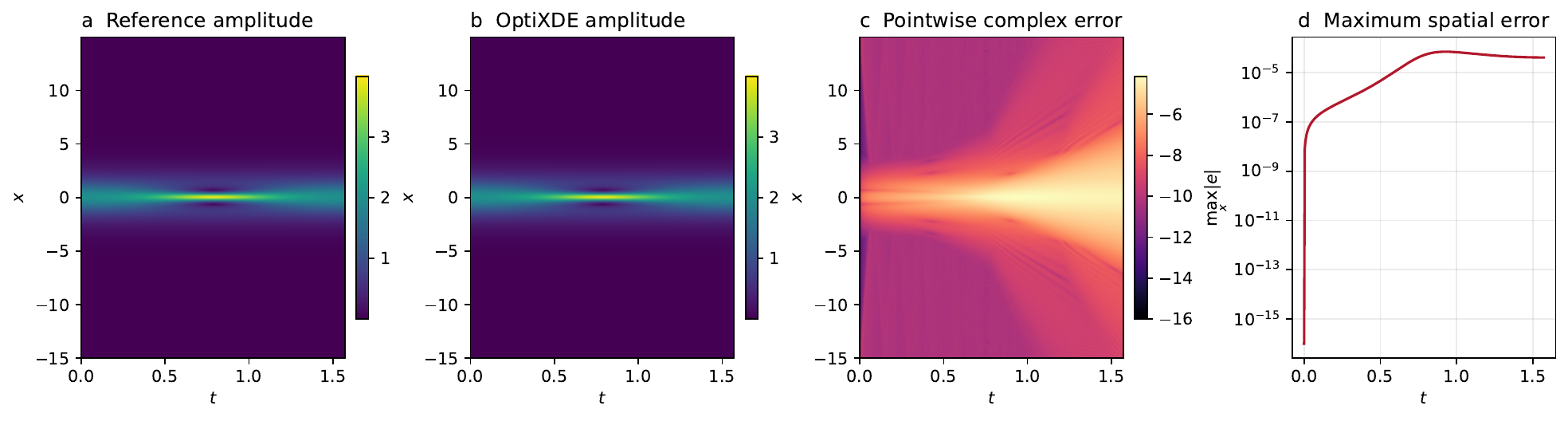}{\textwidth}
\caption{\textbf{Space--time comparison for the nonlinear Schr\"odinger benchmark at $N_x=2048$ and $\Delta t=3.834952\times10^{-4}$.} \textbf{a}, Periodic reference amplitude $\lvert h_{\mathrm{ref}}(x,t)\rvert$. \textbf{b}, \OptiXDE{} amplitude. \textbf{c}, Logarithmic pointwise complex error. \textbf{d}, Maximum spatial error at the stored diagnostic times. The error remains bounded through the focusing and recurrence process, with a maximum recorded value of $4.319465\times10^{-5}$.}
\label{fig:supp_schrodinger_spacetime}
\end{figure}

The complete space--time comparison for the same calculation is shown in Fig.~\ref{fig:supp_schrodinger_spacetime}. The reference and \OptiXDE{} amplitude fields reproduce the same focusing and recurrence pattern, while the complex field error remains bounded throughout the complete interval.

Taken together, the reference-refinement, spatial-resolution and temporal-convergence studies provide mutually independent checks of the nonlinear calculation. The independent RK4 reference is converged to approximately the $10^{-10}$ level in successive relative field differences, spatial differences fall below $10^{-9}$ once $N_x\geq512$, and the complete fixed-step refinement sequence exhibits second-order temporal convergence. The field comparisons further reproduce the characteristic nonlinear focusing, phase evolution and recurrence of the second-order soliton.

\FloatBarrier

\subsection{Burgers equation}
\label{supp:burgers}

This section provides the complete formulation, conservative split-step implementation, reference-solution verification and accuracy assessment for the one-dimensional viscous Burgers benchmark summarized in the main text. The problem combines nonlinear convection with linear diffusion and develops a narrow viscous internal layer around the centre of the domain. It therefore provides a stringent nonlinear test of the analytical diffusion propagator, conservative pseudo-spectral flux evaluation and temporal splitting used in \OptiXDE{}.

\subsubsection{Benchmark definition}
\label{supp:burgers_definition}

The one-dimensional viscous Burgers equation is
\begin{equation}
\frac{\partial u}{\partial t}
+
u\frac{\partial u}{\partial x}
=
\nu\frac{\partial^2u}{\partial x^2},
\qquad
(x,t)\in\Omega\times(0,T],
\label{eq:supp_burgers_governing}
\end{equation}
where $u(x,t)$ denotes the scalar velocity field and $\nu$ is the kinematic viscosity. The physical domain, terminal time and viscosity are
\begin{equation}
\Omega=[-1,1],
\qquad
T=1,
\qquad
\nu=\frac{0.01}{\pi}.
\label{eq:supp_burgers_domain_parameters}
\end{equation}
Homogeneous Dirichlet boundary conditions are imposed at both endpoints,
\begin{equation}
u(-1,t)=u(1,t)=0,
\label{eq:supp_burgers_boundary}
\end{equation}
and the initial condition is
\begin{equation}
u(x,0)=-\sin(\pi x).
\label{eq:supp_burgers_initial}
\end{equation}
The initially smooth profile develops a steep negative gradient around $x=0$ while remaining continuous because $\nu>0$.

The initial field is odd with respect to the origin, and this symmetry is preserved by the governing equation,
\begin{equation}
u(-x,t)=-u(x,t).
\label{eq:supp_burgers_odd_symmetry}
\end{equation}
Consequently, the physical Dirichlet problem can be represented by a period-two odd continuation. The numerical field is stored on the FFT-compatible sampling interval
\begin{equation}
\widetilde{\Omega}=[-1,1),
\label{eq:supp_burgers_fft_domain}
\end{equation}
with uniform samples
\begin{equation}
x_j=-1+j\Delta x,
\qquad
\Delta x=\frac{2}{N_x},
\qquad
j=0,\ldots,N_x-1.
\label{eq:supp_burgers_sampling}
\end{equation}
For even $N_x$, the FFT-compatible modal wave numbers are
\begin{equation}
\kappa
=
\pi
\left\{
0,1,\ldots,\frac{N_x}{2}-1,
-\frac{N_x}{2},\ldots,-1
\right\}.
\label{eq:supp_burgers_wavenumbers}
\end{equation}

The conservative form of Eq.~\eqref{eq:supp_burgers_governing} is
\begin{equation}
\frac{\partial u}{\partial t}
=
\nu\frac{\partial^2u}{\partial x^2}
-
\frac{1}{2}
\frac{\partial}{\partial x}\left(u^2\right).
\label{eq:supp_burgers_conservative_form}
\end{equation}
This form is used in the numerical implementation because the quadratic nonlinear contribution can be evaluated as a de-aliased physical-space product followed by a single spectral derivative.

The continuous kinetic energy is
\begin{equation}
E(t)=\frac{1}{2}\int_{-1}^{1}u^2\,\mathrm{d}x,
\label{eq:supp_burgers_energy}
\end{equation}
and satisfies
\begin{equation}
\frac{\mathrm{d}E}{\mathrm{d}t}
=
-\nu\int_{-1}^{1}
\left(\frac{\partial u}{\partial x}\right)^2
\mathrm{d}x
\leq0,
\label{eq:supp_burgers_energy_dissipation}
\end{equation}
where the boundary contribution vanishes under Eq.~\eqref{eq:supp_burgers_boundary}. The associated discrete consistency behaviour is examined separately in the accuracy and robustness study.

\subsubsection{Conservative split-step implementation and Nyquist treatment}
\label{supp:burgers_splitting}

The Burgers calculation uses the symmetric Strang composition introduced in Section~\ref{supp:nonlinear_splitting}. The diffusion subproblem is propagated analytically in Fourier space, whereas the nonlinear conservative subproblem is advanced by classical fourth-order Runge--Kutta.

For an even number of spatial samples, the Nyquist entry is treated differently in odd- and even-order spectral operators. Let $\kappa_m$ denote the FFT wave-number array in Eq.~\eqref{eq:supp_burgers_wavenumbers}. The wave-number array used for a first derivative is
\begin{equation}
\kappa_m^{(1)}
=
\begin{cases}
\kappa_m, & m\neq N_x/2,\\
0, & m=N_x/2,
\end{cases}
\label{eq:supp_burgers_first_derivative_wavenumber}
\end{equation}
whereas the diffusion operator retains the full Nyquist magnitude,
\begin{equation}
\kappa_m^{(2)}=\kappa_m.
\label{eq:supp_burgers_diffusion_wavenumber}
\end{equation}
Accordingly,
\begin{equation}
\widehat{D_xu}_m
=
\mathrm{i}\kappa_m^{(1)}\widehat{u}_m,
\qquad
\widehat{D_{xx}u}_m
=
-\left(\kappa_m^{(2)}\right)^2\widehat{u}_m.
\label{eq:supp_burgers_derivative_operators}
\end{equation}
This separation removes the ambiguous Nyquist contribution from the odd first derivative while retaining the dissipative action of the even second derivative on the Nyquist mode.

The exact diffusion half-step is therefore
\begin{equation}
\widehat{u}^{\star}_m
=
\exp\left[
-\frac{1}{2}\nu
\left(\kappa_m^{(2)}\right)^2
\Delta t
\right]
\widehat{u}^{n}_m.
\label{eq:supp_burgers_diffusion_half_step}
\end{equation}
For the nonlinear subproblem,
\begin{equation}
\mathcal{N}(u)
=
-\frac{1}{2}
\frac{\partial}{\partial x}
\left(u^2\right),
\qquad
\widehat{\mathcal{N}}_m(u)
=
-\frac{\mathrm{i}\kappa_m^{(1)}}{2}
\widehat{u^2}_{\mathrm{da},m},
\label{eq:supp_burgers_benchmark_nonlinearity}
\end{equation}
where $\widehat{u^2}_{\mathrm{da}}$ denotes the Fourier representation of the de-aliased quadratic product.

All principal Burgers calculations use three-halves zero padding for the nonlinear product. The influence of this choice is assessed independently against no de-aliasing and classical two-thirds truncation in the Burgers robustness study, so that the present benchmark focuses on reference accuracy and spatial--temporal convergence.

The explicit nonlinear stage is monitored using the advective Courant number
\begin{equation}
C_{\mathrm{a}}(t)
=
\frac{\Delta t\max_x|u(x,t)|}{\Delta x}.
\label{eq:supp_burgers_benchmark_courant}
\end{equation}
The spatial-resolution study uses
\begin{equation}
N_x\in\{128,256,512,1024\},
\qquad
\Delta t=6.25\times10^{-5},
\label{eq:supp_burgers_spatial_settings}
\end{equation}
whereas the temporal-convergence study fixes $N_x=1024$ and uses
\begin{equation}
\Delta t
\in
\left\{
1.0\times10^{-3},
5.0\times10^{-4},
2.5\times10^{-4},
1.25\times10^{-4},
6.25\times10^{-5}
\right\}.
\label{eq:supp_burgers_time_steps}
\end{equation}
Unless otherwise stated, the principal field and viscous-layer diagnostics use $N_x=1024$, $\Delta t=6.25\times10^{-5}$ and three-halves zero padding.

\subsubsection{Cole--Hopf reference solution}
\label{supp:burgers_reference}

A semi-analytical reference solution is constructed using the Cole--Hopf transformation
\begin{equation}
u(x,t)
=
-2\nu
\frac{\partial_x\psi(x,t)}{\psi(x,t)},
\label{eq:supp_burgers_cole_hopf}
\end{equation}
which converts Eq.~\eqref{eq:supp_burgers_governing} into the linear heat equation
\begin{equation}
\frac{\partial\psi}{\partial t}
=
\nu\frac{\partial^2\psi}{\partial x^2}.
\label{eq:supp_burgers_cole_hopf_heat}
\end{equation}
For Eq.~\eqref{eq:supp_burgers_initial}, the transformed initial field is
\begin{equation}
\psi(x,0)
=
\exp\left[
-\frac{\cos(\pi x)}{2\pi\nu}
\right].
\label{eq:supp_burgers_cole_hopf_initial}
\end{equation}
The derivative $\partial_x\psi$ vanishes at $x=\pm1$, consistently producing the homogeneous boundary values of $u$.

Defining
\begin{equation}
a=\frac{1}{2\pi\nu},
\label{eq:supp_burgers_bessel_parameter}
\end{equation}
the transformed initial field has the cosine expansion
\begin{equation}
\psi(x,0)
=
I_0(a)
+
2\sum_{m=1}^{\infty}
(-1)^m I_m(a)\cos(m\pi x),
\label{eq:supp_burgers_initial_cosine_series}
\end{equation}
where $I_m$ is the modified Bessel function of the first kind. The heat equation gives
\begin{equation}
\psi(x,t)
=
I_0(a)
+
2\sum_{m=1}^{\infty}
(-1)^m I_m(a)
\exp\left(-\nu m^2\pi^2t\right)
\cos(m\pi x).
\label{eq:supp_burgers_cole_hopf_series}
\end{equation}
Substitution into Eq.~\eqref{eq:supp_burgers_cole_hopf} yields
\begin{equation}
\begin{split}
u_{\mathrm{ref}}(x,t)
=
\frac{
4\pi\nu
\displaystyle\sum_{m=1}^{\infty}
m(-1)^m I_m(a)
\exp\left(-\nu m^2\pi^2t\right)
\sin(m\pi x)
}{
\displaystyle I_0(a)
+
2\displaystyle\sum_{m=1}^{\infty}
(-1)^m I_m(a)
\exp\left(-\nu m^2\pi^2t\right)
\cos(m\pi x)
}.
\end{split}
\label{eq:supp_burgers_reference_solution}
\end{equation}

For the selected viscosity, $a=50$. Direct evaluation of the unscaled Bessel representation is poorly conditioned because the coefficients span a wide dynamic range. The reference is therefore evaluated using an equivalent log-scaled Cole--Hopf heat-convolution procedure. The stabilized convolution preserves the quotient in Eq.~\eqref{eq:supp_burgers_cole_hopf} while avoiding unnecessary loss of significance.

\begin{table}[t]
\centering
\caption{Settings and refinement verification for the Cole--Hopf reference solution of the Burgers benchmark.}
\label{tab:supp_burgers_reference_settings}
\begin{tabular*}{\linewidth}{@{\extracolsep{\fill}}ll@{}}
\toprule
Setting & Value \\
\midrule
Reference sampling interval & $\widetilde{\Omega}=[-1,1)$ \\
Terminal time & $T=1$ \\
Viscosity & $\nu=0.01/\pi$ \\
Reference formulation & Log-scaled Cole--Hopf heat convolution \\
Reference sampling resolution & $1024$ \\
Primary periodic quadrature points & $8192$ \\
Refinement-check quadrature points & $4096$ \\
Series-equivalent Bessel parameter & $a=50$ \\
Arithmetic precision & IEEE 754 binary64 \\
Quadrature refinement difference & $6.098041\times10^{-13}$ \\
\bottomrule
\end{tabular*}
\end{table}

The reported difference between the $8192$- and $4096$-point quadrature evaluations is $6.098041\times10^{-13}$, more than five orders of magnitude below the smallest terminal field error reported below. The reference evaluation therefore does not control the measured \OptiXDE{} errors.

\begin{figure}[t]
\centering
\suppfigure{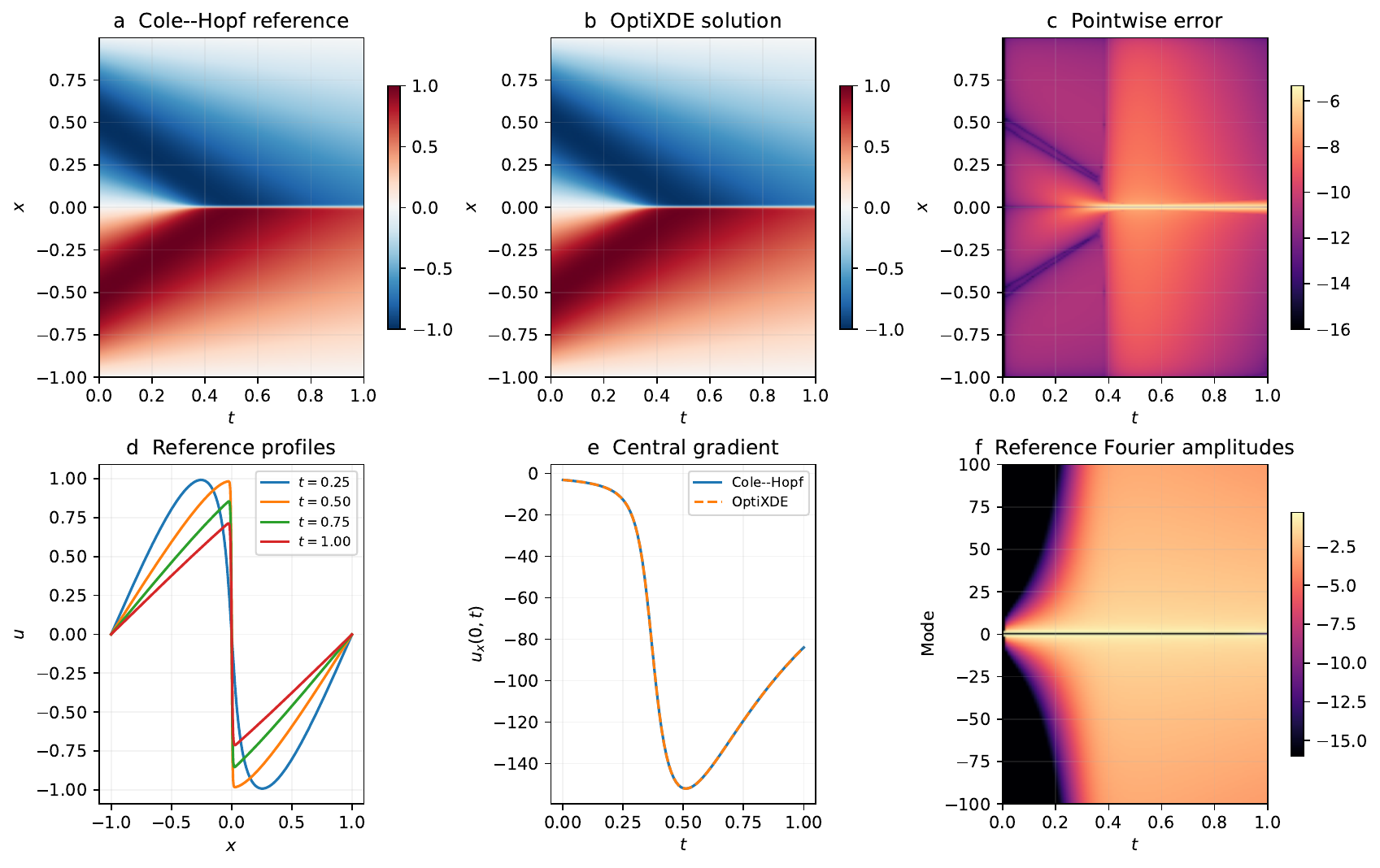}{\textwidth}
\caption{\textbf{Reference and principal numerical evolution of the viscous Burgers benchmark.} \textbf{a}, Space--time distribution of the Cole--Hopf reference solution $u_{\mathrm{ref}}(x,t)$. \textbf{b}, Corresponding principal \OptiXDE{} solution obtained with $N_x=1024$, $\Delta t=6.25\times10^{-5}$ and three-halves zero padding. \textbf{c}, Pointwise absolute error represented by $\log_{10}|u_h-u_{\mathrm{ref}}|$. \textbf{d}, Reference profiles at $t=0.25$, $0.50$, $0.75$ and $1.00$. \textbf{e}, Central gradient $\partial_xu(0,t)$ for the Cole--Hopf reference and \OptiXDE{} solution. \textbf{f}, Evolution of the reference Fourier amplitudes represented by $\log_{10}|\widehat{u}_{\mathrm{ref}}|$.}
\label{fig:supp_burgers_reference}
\end{figure}

Figure~\ref{fig:supp_burgers_reference} shows the progressive steepening of the initially smooth odd profile and the formation of a narrow viscous layer around $x=0$. The strongest central gradient occurs near $t=0.5$, after which viscous dissipation gradually relaxes the layer. At the scale of the solution field, the principal \OptiXDE{} result is visually indistinguishable from the Cole--Hopf reference, including through the maximum-steepening stage. The pointwise-error map shows that the largest discrepancies remain localized around the central layer and the period during which high-wave-number content is generated most rapidly.

\subsubsection{Spatial-resolution study}
\label{supp:burgers_spatial_resolution}

The spatial-resolution study uses Eq.~\eqref{eq:supp_burgers_spatial_settings} and three-halves padding. The common field errors $e_{\max}$ and $e(T)$ follow Section~\ref{supp:error_measures}. To quantify successive-resolution consistency independently of the analytical reference, we define
\begin{equation}
\delta_N
=
\max_{t\in\mathcal{T}_{\mathrm{out}}}
\frac{
\left\|
u_N(t)-\mathcal{R}_{2N\rightarrow N}u_{2N}(t)
\right\|_{2,h}
}{
\left\|
\mathcal{R}_{2N\rightarrow N}u_{2N}(t)
\right\|_{2,h}
},
\label{eq:supp_burgers_successive_resolution_difference}
\end{equation}
where $\mathcal{R}_{2N\rightarrow N}$ denotes restriction by exact selection of coincident nodes from the finer nested periodic sampling grid, and $\mathcal{T}_{\mathrm{out}}$ is the set of stored comparison times.

The high-wave-number occupation is measured by
\begin{equation}
\rho_{\mathrm{high}}(t)
=
\frac{
\displaystyle\sum_{|\kappa|>2\kappa_{\max}/3}
|\widehat{u}(\kappa,t)|^2
}{
\displaystyle\sum_{\kappa}
|\widehat{u}(\kappa,t)|^2
},
\label{eq:supp_burgers_high_wavenumber_ratio}
\end{equation}
with
\begin{equation}
\rho_{\mathrm{high}}^{\max}
=
\max_{0\leq t\leq T}
\rho_{\mathrm{high}}(t).
\label{eq:supp_burgers_maximum_high_wavenumber_ratio}
\end{equation}
This quantity is used here only as a spectral-resolution diagnostic; its sensitivity to the nonlinear-product treatment is considered separately in the robustness study.

Because the strongest gradient remains centred at $x=0$, the absolute central-gradient error is
\begin{equation}
\varepsilon_{\mathrm{g}}(t)
=
\left|
\left.\frac{\partial u_h}{\partial x}\right|_{x=0}
-
\left.\frac{\partial u_{\mathrm{ref}}}{\partial x}\right|_{x=0}
\right|.
\label{eq:supp_burgers_gradient_error}
\end{equation}
A characteristic viscous-layer thickness is defined by
\begin{equation}
\delta_{\mathrm{s}}(t)
=
\frac{2\max_x|u(x,t)|}{|\partial_xu(0,t)|},
\label{eq:supp_burgers_layer_thickness}
\end{equation}
and the corresponding relative thickness error is
\begin{equation}
\varepsilon_{\delta}(t)
=
\frac{
|\delta_{\mathrm{s},h}(t)-\delta_{\mathrm{s},\mathrm{ref}}(t)|
}{
\delta_{\mathrm{s},\mathrm{ref}}(t)
}.
\label{eq:supp_burgers_layer_thickness_error}
\end{equation}

\begin{table}[t]
\centering\scriptsize
\caption{Spatial-resolution study of the Burgers benchmark using three-halves padding at fixed $\Delta t=6.25\times10^{-5}$.}
\label{tab:supp_burgers_spatial_resolution}
\begin{tabular*}{\linewidth}{@{\extracolsep{\fill}}ccccccc@{}}
\toprule
$N_x$ & $\delta_N$ & $e_{\max}$ & $e(T)$ & $\max_t\varepsilon_{\mathrm{g}}$ & $\max_t\varepsilon_\delta$ & $\rho_{\mathrm{high}}^{\max}$ \\
\midrule
$128$ & $3.717373\times10^{-2}$ & $8.929470\times10^{-2}$ & $9.095883\times10^{-3}$ & $2.540457\times10^{1}$ & $1.801626\times10^{-1}$ & $2.790068\times10^{-3}$ \\
$256$ & $4.034870\times10^{-3}$ & $1.158858\times10^{-2}$ & $5.088851\times10^{-4}$ & $7.976599\times10^{0}$ & $4.754095\times10^{-2}$ & $1.419871\times10^{-4}$ \\
$512$ & $5.918097\times10^{-5}$ & $1.944313\times10^{-4}$ & $1.898306\times10^{-6}$ & $2.802315\times10^{-1}$ & $1.959462\times10^{-3}$ & $5.966732\times10^{-7}$ \\
$1024$ & -- & $4.590928\times10^{-6}$ & $8.508945\times10^{-8}$ & $2.123683\times10^{-3}$ & $1.400289\times10^{-5}$ & $1.066850\times10^{-11}$ \\
\bottomrule
\end{tabular*}
\end{table}

\begin{figure}[t]
\centering
\suppfigure{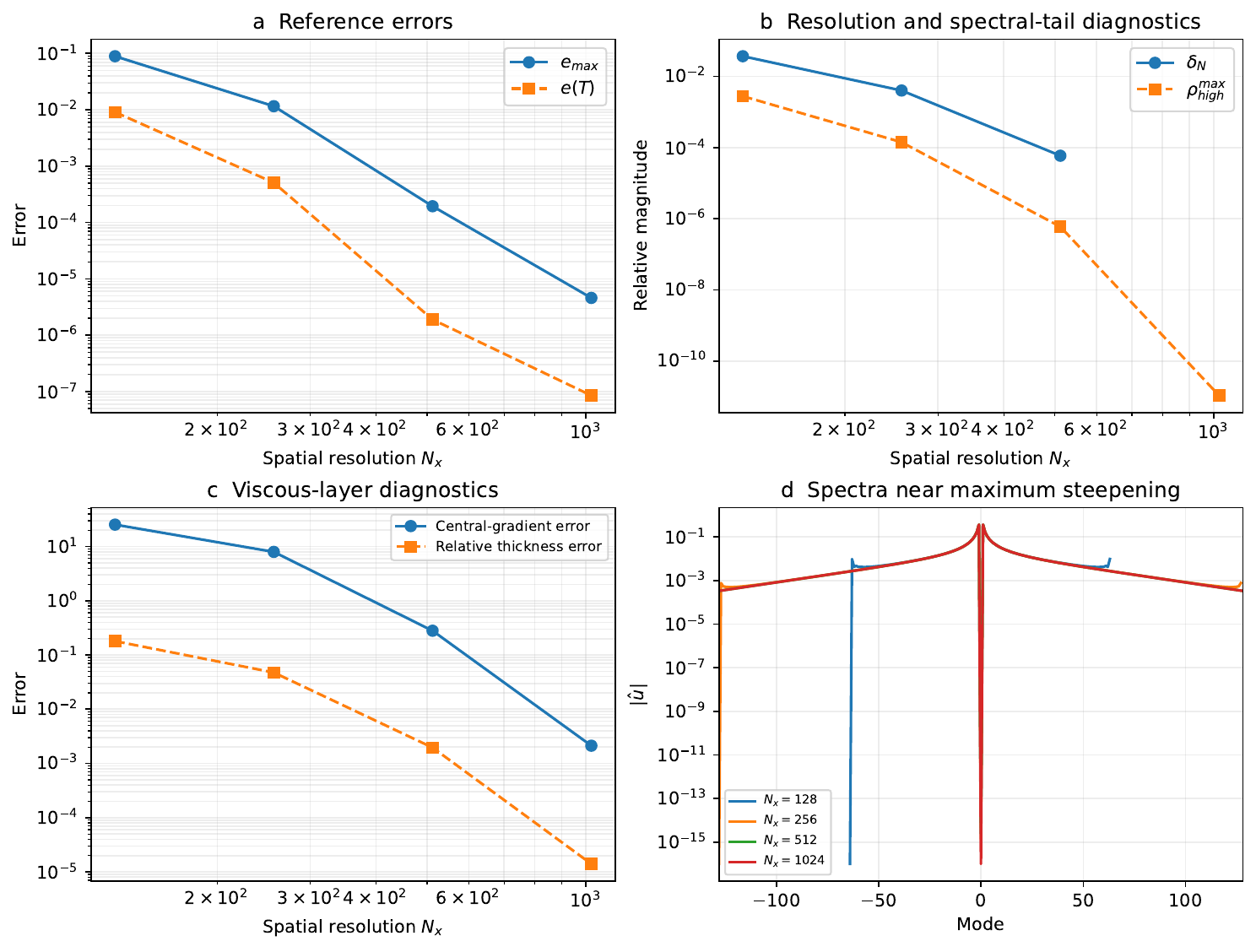}{\textwidth}
\caption{\textbf{Spatial-resolution study of the viscous Burgers benchmark.} \textbf{a}, Maximum space--time error $e_{\max}$ and terminal field error $e(T)$ as functions of $N_x$. \textbf{b}, Successive-resolution difference $\delta_N$ and maximum high-wave-number ratio $\rho_{\mathrm{high}}^{\max}$. \textbf{c}, Maximum central-gradient error and maximum relative viscous-layer-thickness error. \textbf{d}, Fourier-amplitude spectra near the time of maximum steepening for $N_x=128$, $256$, $512$ and $1024$.}
\label{fig:supp_burgers_spatial_resolution}
\end{figure}

Table~\ref{tab:supp_burgers_spatial_resolution} and Fig.~\ref{fig:supp_burgers_spatial_resolution} demonstrate rapid convergence under spatial refinement. Increasing $N_x$ from $128$ to $1024$ reduces $e_{\max}$ from $8.929470\times10^{-2}$ to $4.590928\times10^{-6}$ and reduces $e(T)$ from $9.095883\times10^{-3}$ to $8.508945\times10^{-8}$. The successive-resolution difference decreases from $3.717373\times10^{-2}$ at $N_x=128$ to $5.918097\times10^{-5}$ at $N_x=512$, while the central-gradient and layer-thickness diagnostics exhibit the same systematic refinement trend.

The finest $N_x=1024$ spatial result reaches the error level of the finest temporal calculation at the same $\Delta t=6.25\times10^{-5}$. The final point should therefore be interpreted as approaching the temporal-splitting floor rather than as a pure spatial asymptote. This separation is examined directly in the following time-step study.

\subsubsection{Temporal convergence}
\label{supp:burgers_convergence}

The temporal-convergence study fixes $N_x=1024$ and uses the five time increments in Eq.~\eqref{eq:supp_burgers_time_steps}. All cases use three-halves padding so that the spatial representation and nonlinear-product treatment remain unchanged. The observed order based on $e_{\max}$ is evaluated using the common definition in Eq.~\eqref{eq:supp_observed_convergence_order} with $q=\Delta t$.

\begin{table}[t]
\centering\scriptsize
\caption{Fixed-step temporal convergence of the Burgers benchmark at $N_x=1024$ using three-halves padding.}
\label{tab:supp_burgers_convergence}
\begin{tabular*}{\linewidth}{@{\extracolsep{\fill}}cccccc@{}}
\toprule
$N_t$ & $\Delta t$ & $e_{\max}$ & Order $p$ & $e(T)$ & $\max_t C_{\mathrm{a}}$ \\
\midrule
$1000$ & $1.000\times10^{-3}$ & $1.161069\times10^{-3}$ & -- & $2.174178\times10^{-5}$ & $0.5120$ \\
$2000$ & $5.000\times10^{-4}$ & $2.913281\times10^{-4}$ & $1.995$ & $5.442919\times10^{-6}$ & $0.2560$ \\
$4000$ & $2.500\times10^{-4}$ & $7.291245\times10^{-5}$ & $1.998$ & $1.361210\times10^{-6}$ & $0.1280$ \\
$8000$ & $1.250\times10^{-4}$ & $1.825756\times10^{-5}$ & $1.998$ & $3.403361\times10^{-7}$ & $0.0640$ \\
$16000$ & $6.250\times10^{-5}$ & $4.590928\times10^{-6}$ & $1.992$ & $8.508945\times10^{-8}$ & $0.0320$ \\
\bottomrule
\end{tabular*}
\end{table}

\begin{figure}[t]
\centering
\suppfigure{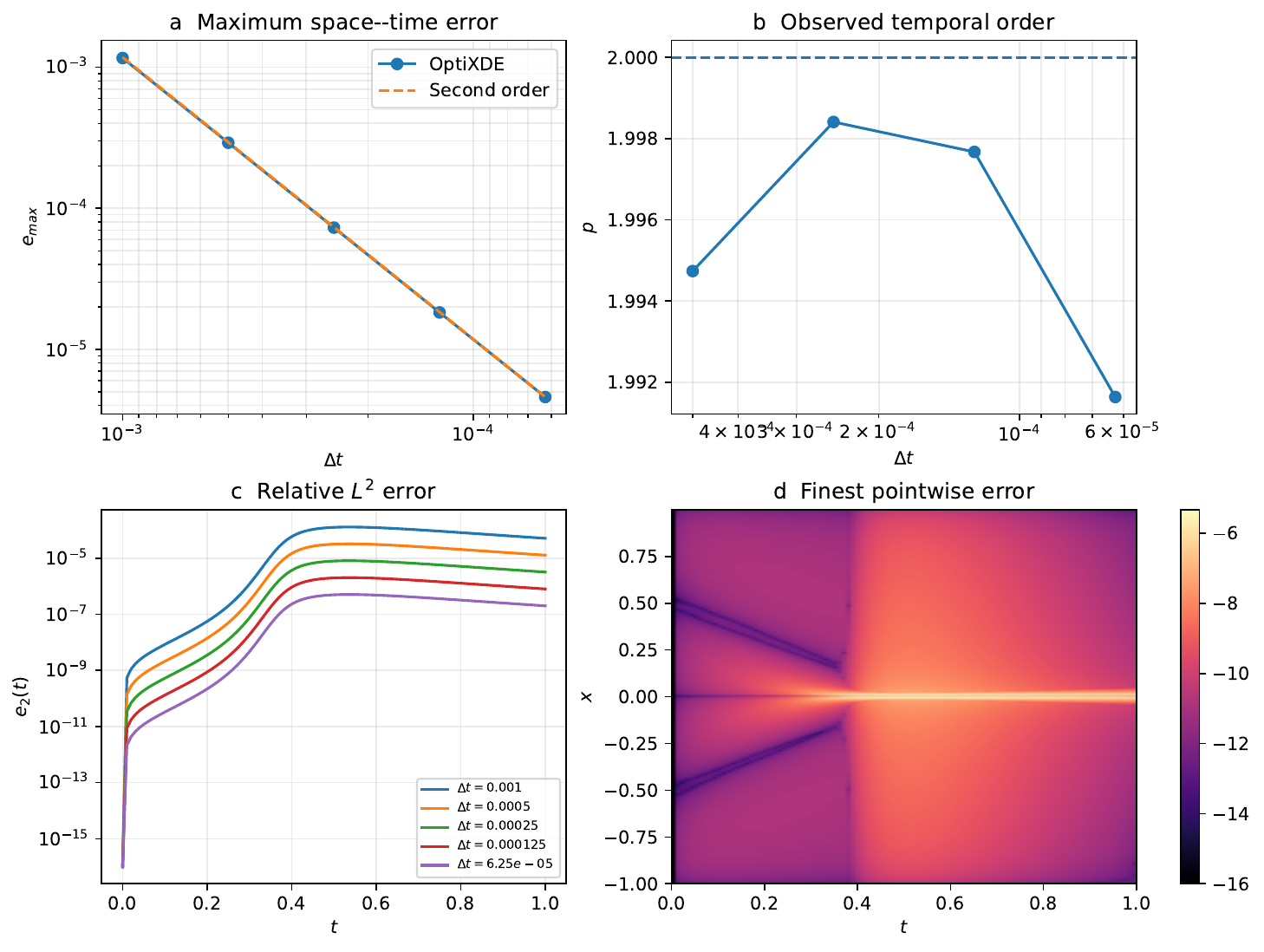}{\textwidth}
\caption{\textbf{Temporal convergence of the viscous Burgers benchmark.} \textbf{a}, Maximum space--time error $e_{\max}$ as a function of $\Delta t$, together with a second-order reference slope. \textbf{b}, Observed temporal order between successive time-step refinements. \textbf{c}, Relative $L^2$-error histories for the five tested time increments. \textbf{d}, Pointwise absolute-error distribution, represented by $\log_{10}|u_h-u_{\mathrm{ref}}|$, for the finest time-step calculation.}
\label{fig:supp_burgers_convergence}
\end{figure}

Each halving of $\Delta t$ reduces $e_{\max}$ by approximately a factor of four. The observed orders are $1.995$, $1.998$, $1.998$ and $1.992$, providing direct numerical verification of the formal second-order accuracy of the Strang diffusion--convection composition. The terminal error shows the same approximately fourfold reduction, decreasing from $2.174178\times10^{-5}$ to $8.508945\times10^{-8}$. The maximum advective Courant number decreases from $0.5120$ to $0.0320$ over the same refinement sequence.

The relative-error histories in Fig.~\ref{fig:supp_burgers_convergence}(c) remain extremely small during the early smooth stage and increase most strongly as the viscous layer forms. Their nearly parallel separation with time-step refinement, together with the slope in Fig.~\ref{fig:supp_burgers_convergence}(a), confirms that the measured sequence is governed predominantly by temporal splitting error rather than by the spatial or reference-solution floor.

\subsubsection{Viscous-layer error analysis}
\label{supp:burgers_shock_error}

Although the viscous Burgers solution remains smooth, nonlinear convection generates a narrow internal layer around $x=0$. For localized diagnostics, the viscous-layer and outer regions are defined as
\begin{equation}
\Omega_{\mathrm{s}}
=
\left\{x\in[-1,1]:|x|\leq0.15\right\},
\qquad
\Omega_{\mathrm{b}}
=
\Omega\setminus\Omega_{\mathrm{s}}.
\label{eq:supp_burgers_regions}
\end{equation}
The regional root-mean-square and maximum errors follow the common definitions in Section~\ref{supp:error_measures} and are denoted by $e_{\mathrm{s}}$, $e_{\mathrm{b}}$ and $e_{\mathrm{s}}^{\infty}$. The gradient and thickness errors are given by Eqs.~\eqref{eq:supp_burgers_gradient_error} and~\eqref{eq:supp_burgers_layer_thickness_error}.

\begin{table}[t]
\centering
\caption{Viscous-layer accuracy of the Burgers calculation on $|x|\leq0.15$ using the finest principal three-halves calculation.}
\label{tab:supp_burgers_shock_error}
\begin{tabular*}{\linewidth}{@{\extracolsep{\fill}}ccccc@{}}
\toprule
$t$ & $e_{\mathrm{s}}$ & $e_{\mathrm{s}}^\infty$ & $\varepsilon_{\mathrm{g}}$ & $\varepsilon_\delta$ \\
\midrule
$0.25$ & $1.411001\times10^{-9}$ & $2.797755\times10^{-9}$ & $1.854281\times10^{-7}$ & $1.453257\times10^{-8}$ \\
$0.50$ & $8.176659\times10^{-7}$ & $4.569925\times10^{-6}$ & $2.113719\times10^{-3}$ & $1.388259\times10^{-5}$ \\
$0.75$ & $4.814902\times10^{-7}$ & $2.555347\times10^{-6}$ & $1.100815\times10^{-3}$ & $9.210525\times10^{-6}$ \\
$1.00$ & $2.201304\times10^{-7}$ & $1.044744\times10^{-6}$ & $3.933379\times10^{-4}$ & $4.671863\times10^{-6}$ \\
\bottomrule
\end{tabular*}
\end{table}

\begin{figure}[t]
\centering
\suppfigure{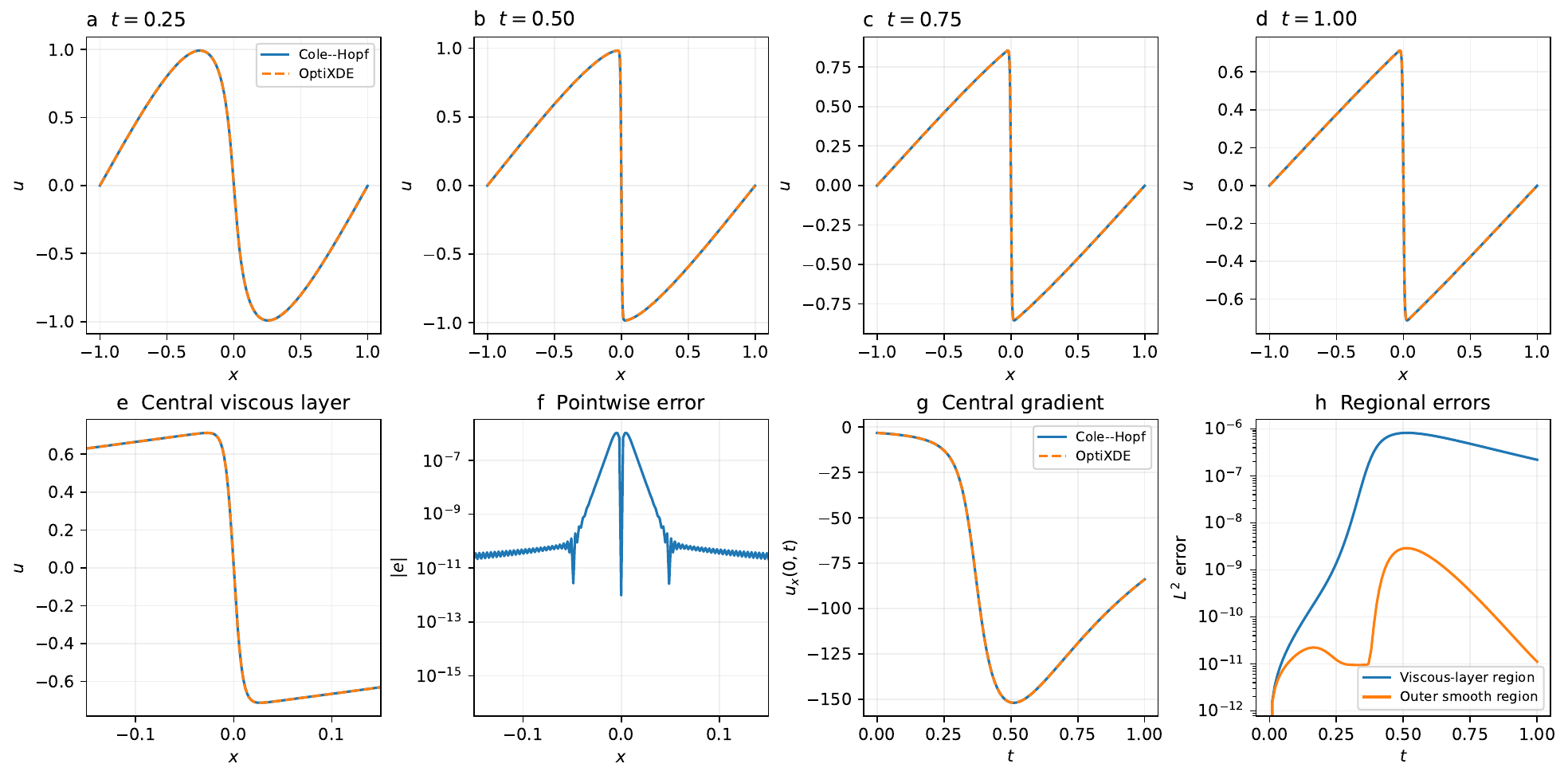}{\textwidth}
\caption{\textbf{Viscous-layer accuracy of the Burgers benchmark.} \textbf{a--d}, Cole--Hopf reference and principal \OptiXDE{} profiles at $t=0.25$, $0.50$, $0.75$ and $1.00$, respectively. \textbf{e}, Enlarged terminal comparison around the central viscous layer. \textbf{f}, Terminal pointwise absolute error in the neighbourhood of $x=0$. \textbf{g}, Reference and numerical central gradients as functions of time. \textbf{h}, Viscous-layer and outer-region $L^2$ errors throughout the simulation.}
\label{fig:supp_burgers_shock_error}
\end{figure}

The principal calculation remains highly accurate throughout the complete steepening process. At $t=0.25$, the viscous-layer root-mean-square and maximum errors are $1.411001\times10^{-9}$ and $2.797755\times10^{-9}$, respectively. The largest localized discrepancy occurs near the maximum-steepening stage at $t=0.50$, where $e_{\mathrm{s}}^{\infty}=4.569925\times10^{-6}$ and the absolute central-gradient error is $2.113719\times10^{-3}$. Even at this time, the relative viscous-layer-thickness error is only $1.388259\times10^{-5}$. By $t=1$, $e_{\mathrm{s}}^{\infty}$ has decreased to $1.044744\times10^{-6}$ and the relative thickness error to $4.671863\times10^{-6}$.

Figure~\ref{fig:supp_burgers_shock_error} confirms that the remaining error is strongly localized around the narrow central layer. The reference and numerical profiles are visually coincident at the full-domain scale, the enlarged terminal comparison shows only a small localized discrepancy, and the central-gradient histories remain nearly indistinguishable. The outer smooth-region error remains substantially below the viscous-layer error after the onset of rapid steepening.

Taken together, the independent reference verification, spatial-resolution study, temporal-convergence study and localized viscous-layer diagnostics establish the accuracy of the Burgers implementation. The reference quadrature error is negligible relative to the measured numerical error, spatial refinement progressively resolves the narrow internal layer, and the complete diffusion--convection split exhibits the expected second-order temporal convergence.

\subsubsection{Reproducibility information}
\label{supp:burgers_reproducibility}

All reported calculations use uniform periodic sampling without duplication of the right endpoint and IEEE 754 binary64 arithmetic. The physical Dirichlet condition is represented through the odd period-two continuation, the nonlinear term is evaluated in conservative flux form, and the principal nonlinear-product treatment is three-halves zero padding. The common software environment, hardware configuration and timing protocol are reported in Section~\ref{supp:computational_environment}. The benchmark-specific reference settings are listed in Table~\ref{tab:supp_burgers_reference_settings}, and the complete spatial and temporal parameter sequences are given in Eqs.~\eqref{eq:supp_burgers_spatial_settings} and~\eqref{eq:supp_burgers_time_steps}. Sensitivity to nonlinear de-aliasing and additional physical-consistency diagnostics are reported separately in the accuracy and robustness study.

\FloatBarrier

\subsection{Allen--Cahn equation}
\label{supp:allen_cahn}

This section provides the benchmark-specific formulation, numerical implementation and physical and numerical verification for the two-dimensional Allen--Cahn equation. In contrast to the conservative nonlinear Schr\"odinger equation and the convection-dominated Burgers equation, the Allen--Cahn equation describes a dissipative reaction--diffusion process in which a nonconserved phase field relaxes toward the stable states \(u=\pm1\) while a diffuse interface evolves under curvature. The present benchmark therefore examines whether \OptiXDE{} reproduces the expected interface motion, free-energy dissipation and phase-field bounds while retaining the second-order temporal accuracy of the Strang composition.

\subsubsection{Benchmark definition and independent reference}
\label{supp:allen_cahn_definition}

The two-dimensional Allen--Cahn equation is defined as
\begin{equation}
\frac{\partial u}{\partial t}
=
\epsilon^2\nabla^2u
+
u-u^3,
\qquad
(x,y)\in\Omega=[0,2\pi)^2,
\qquad
t\in(0,T],
\label{eq:supp_allen_cahn_governing}
\end{equation}
where \(u(x,y,t)\) is a nonconserved phase-field variable and \(\epsilon>0\) controls the characteristic diffuse-interface thickness. Periodic boundary conditions are imposed in both spatial directions:
\begin{equation}
u(x+2\pi,y,t)=u(x,y,t),
\qquad
u(x,y+2\pi,t)=u(x,y,t).
\label{eq:supp_allen_cahn_periodic}
\end{equation}

The periodic square is represented by the uniform Cartesian sampling
\begin{equation}
x_i=\frac{2\pi i}{N_x},
\qquad
y_j=\frac{2\pi j}{N_y},
\qquad
i=0,\ldots,N_x-1,
\qquad
j=0,\ldots,N_y-1.
\label{eq:supp_allen_cahn_sampling}
\end{equation}
For even \(N_x\) and \(N_y\), the corresponding FFT-compatible wave numbers are
\begin{equation}
k_x=
\left\{
0,1,\ldots,\frac{N_x}{2}-1,
-\frac{N_x}{2},\ldots,-1
\right\},
\qquad
k_y=
\left\{
0,1,\ldots,\frac{N_y}{2}-1,
-\frac{N_y}{2},\ldots,-1
\right\},
\label{eq:supp_allen_cahn_wavenumbers}
\end{equation}
with
\begin{equation}
K^2=k_x^2+k_y^2.
\label{eq:supp_allen_cahn_wavenumber}
\end{equation}

The initial field consists of a circular region of the \(u\approx1\) phase embedded in the \(u\approx-1\) phase:
\begin{equation}
u(x,y,0)
=
\tanh
\left[
\frac{R_0-r(x,y)}{\sqrt{2}\epsilon}
\right],
\label{eq:supp_allen_cahn_initial}
\end{equation}
where
\begin{equation}
r(x,y)
=
\sqrt{(x-\pi)^2+(y-\pi)^2},
\qquad
R_0=\frac{\pi}{2}.
\label{eq:supp_allen_cahn_radius_definition}
\end{equation}
The benchmark parameters are
\begin{equation}
\epsilon=0.04,
\qquad
T=100.
\label{eq:supp_allen_cahn_parameters}
\end{equation}

The hyperbolic-tangent profile is the one-dimensional stationary transition associated with the double-well potential. It provides a smooth interface connecting the two stable bulk phases, with \(u=0\) identifying the centre of the diffuse layer. The initial configuration is radially symmetric and subsequently contracts under curvature.

The principal calculation used for the phase-field evolution and physical diagnostics is performed using
\begin{equation}
N_x=N_y=256,
\qquad
\Delta t=5.0\times10^{-2}.
\label{eq:supp_allen_cahn_primary_discretization}
\end{equation}

An independent reference solution is generated using a Fourier pseudo-spectral spatial discretization combined with the fourth-order exponential time-differencing Runge--Kutta method (ETDRK4). The primary reference uses a spatial resolution of \(512^2\) and a time step of \(5.0\times10^{-2}\). Temporal reference consistency is assessed by comparison with \(\Delta t=1.0\times10^{-1}\) at the same spatial resolution, while spatial reference consistency is assessed by comparing the \(256^2\) and \(512^2\) ETDRK4 solutions at \(\Delta t=5.0\times10^{-2}\).

\begin{table}[t]
\centering
\caption{Independent-reference verification for the two-dimensional Allen--Cahn benchmark.}
\label{tab:supp_allen_cahn_reference}
\begin{tabular*}{\linewidth}{@{\extracolsep{\fill}}ll@{}}
\toprule
Setting or diagnostic & Value \\
\midrule
Reference method & Fourier pseudo-spectral ETDRK4 \\
Primary reference resolution & \(512^2\) \\
Primary reference time step & \(5.0\times10^{-2}\) \\
Temporal refinement check & \(512^2\), \(\Delta t=1.0\times10^{-1}\) \\
Temporal reference difference & \(2.432253\times10^{-7}\) \\
Spatial refinement check & \(256^2\) versus \(512^2\) at \(\Delta t=5.0\times10^{-2}\) \\
Spatial reference difference & \(5.200814\times10^{-7}\) \\
\bottomrule
\end{tabular*}
\end{table}

Both reference differences remain below the smallest OptiXDE error reported in the temporal-convergence study and are substantially below the error of the principal calculation. The independently generated ETDRK4 solution therefore provides a sufficiently resolved reference for the present convergence assessment.

For the spatial-resolution study, the \(512^2\) ETDRK4 field is periodically restricted to each tested resolution before the field error is evaluated. For the temporal-convergence study, a same-resolution \(256^2\) ETDRK4 reference is used so that temporal splitting errors are compared without introducing differences associated with spatial restriction.

The numerical accuracy is measured using the common relative discrete $L^2$ error $e_{\mathrm{rel}}(t)$ and maximum pointwise error $e_{\infty}(t)$ defined in Eqs.~\eqref{eq:supp_relative_l2_error} and~\eqref{eq:supp_maximum_error}, respectively.

\subsubsection{Benchmark-specific splitting settings}
\label{supp:allen_cahn_splitting}

The Allen--Cahn calculation is the periodic two-dimensional specialization of the nonlinear framework introduced in Section~\ref{supp:nonlinear_splitting}. The diffusion and reaction contributions are advanced separately using the second-order Strang composition. In the present parameterization, the diffusion coefficient is
\begin{equation}
\kappa_{\mathrm{AC}}=\epsilon^2.
\label{eq:supp_allen_cahn_diffusion_coefficient}
\end{equation}

For each tested \(\Delta t\), the cached diffusion half-step is obtained from Eq.~\eqref{eq:supp_allen_cahn_linear_flow} by setting \(\kappa_{\mathrm{AC}}=\epsilon^2\) and \(\tau=\Delta t/2\). The complete time step applies this analytical diffusion propagator before and after the exact pointwise reaction flow defined in Section~\ref{supp:nonlinear_splitting}. No nonlinear iteration is required.

Because the diffusion and reaction subflows are individually evaluated in closed form but do not commute, their symmetric composition is second-order accurate in time. The numerical temporal order is verified independently below rather than inferred solely from the formal splitting construction.

\subsubsection{Phase-field evolution and curvature-flow verification}
\label{supp:allen_cahn_evolution}

The diffuse interface is identified by the zero level set
\begin{equation}
\Gamma(t)
=
\left\{
(x,y)\in\Omega:
u(x,y,t)=0
\right\}.
\label{eq:supp_allen_cahn_interface}
\end{equation}
To quantify the positive-phase geometry without explicitly reconstructing the complete interface, a regularized area is defined as
\begin{equation}
A_+(t)
=
\int_{\Omega}
H_{\delta}\!\left(u(x,y,t)\right)
\,\dd\Omega,
\label{eq:supp_allen_cahn_positive_area}
\end{equation}
where
\begin{equation}
H_{\delta}(u)
=
\frac{1}{2}
\left[
1+
\tanh\left(\frac{u}{\delta}\right)
\right].
\label{eq:supp_allen_cahn_smooth_heaviside}
\end{equation}
The principal diagnostic uses
\begin{equation}
\delta=0.05,
\label{eq:supp_allen_cahn_heaviside_delta}
\end{equation}
and the corresponding equivalent interface radius is
\begin{equation}
R_h(t)
=
\sqrt{
\frac{A_+(t)}{\pi}
}.
\label{eq:supp_allen_cahn_equivalent_radius}
\end{equation}

For a well-resolved circular interface, the leading-order sharp-interface motion associated with Eq.~\eqref{eq:supp_allen_cahn_governing} satisfies
\begin{equation}
\frac{\dd R_{\kappa}}{\dd t}
=
-\frac{\epsilon^2}{R_{\kappa}},
\label{eq:supp_allen_cahn_curvature_velocity}
\end{equation}
and therefore
\begin{equation}
R_{\kappa}(t)
=
\sqrt{
R_0^2-2\epsilon^2t
}.
\label{eq:supp_allen_cahn_curvature_radius}
\end{equation}
This expression is used as a leading-order physical reference rather than as an exact finite-\(\epsilon\) solution of the complete phase-field equation.

\begin{figure}[t]
\centering
\includegraphics[width=\textwidth]{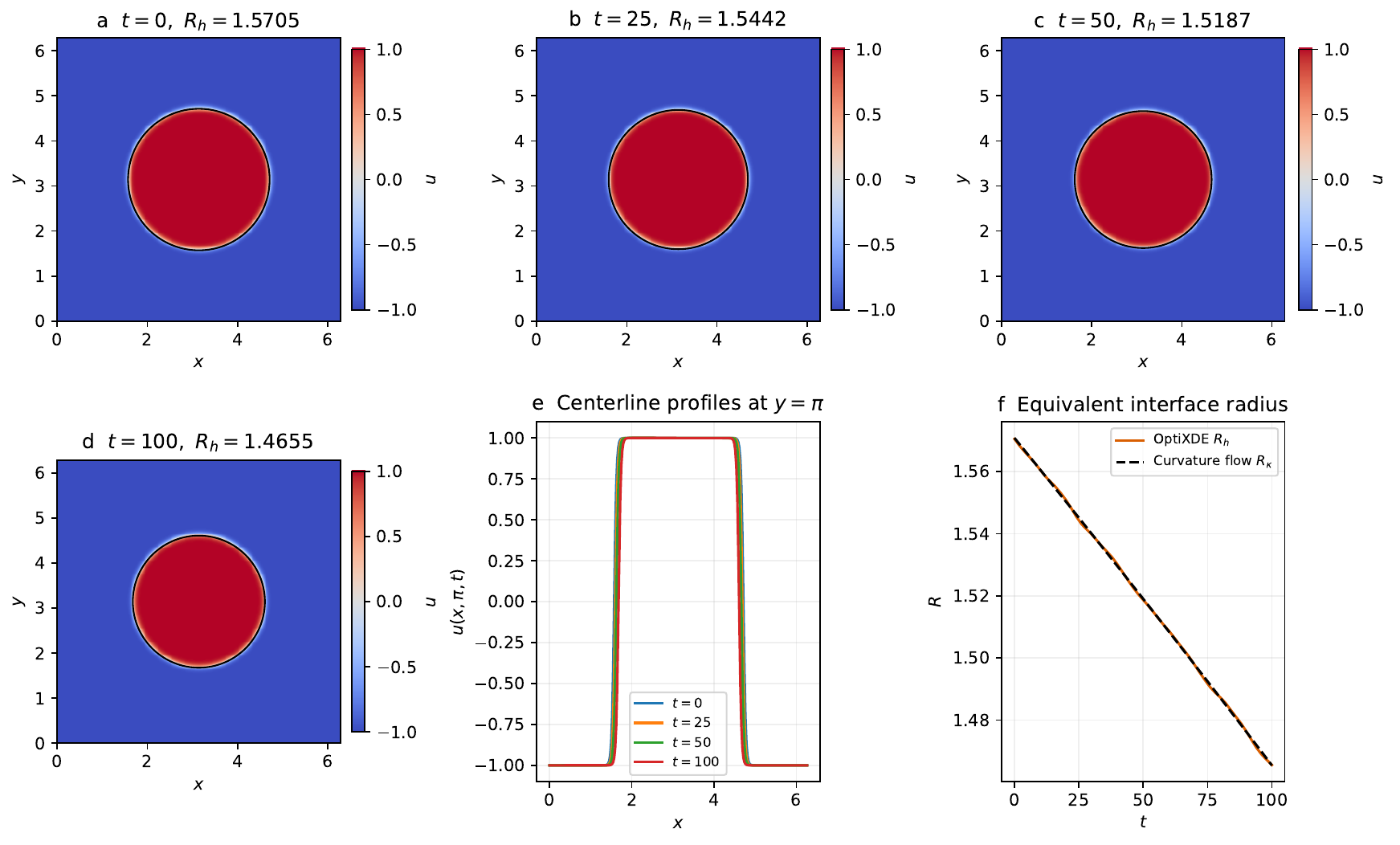}
\caption{\textbf{Phase-field evolution of the two-dimensional Allen--Cahn benchmark.}
\textbf{a--d}, Phase-field distributions at \(t=0\), \(25\), \(50\) and \(100\), respectively, using a common color range of \([-1,1]\); the black \(u=0\) contour identifies the centre of the diffuse interface.
\textbf{e}, Centerline profiles along \(y=\pi\) at the selected times.
\textbf{f}, Equivalent interface radius \(R_h(t)\) compared with the leading-order curvature-flow prediction \(R_{\kappa}(t)=\sqrt{R_0^2-2\epsilon^2t}\).}
\label{fig:supp_allen_cahn_evolution}
\end{figure}

Figure~\ref{fig:supp_allen_cahn_evolution}a--d shows continuous contraction of the initially circular positive-phase region while the diffuse interface remains smooth. The centerline profiles in Fig.~\ref{fig:supp_allen_cahn_evolution}e retain almost unchanged transition slopes as the two interfaces move toward the centre, indicating that the interface thickness is maintained during the evolution.

The equivalent radius decreases from \(R_h(0)=1.57049045\) to
\begin{equation}
R_h(100)=1.46550643.
\label{eq:supp_allen_cahn_terminal_radius}
\end{equation}
The numerical radius closely follows the leading-order curvature-flow prediction throughout the calculation. The maximum relative discrepancy,
\begin{equation}
\varepsilon_{R,\max}
=
\max_t
\frac{
\left|
R_h(t)-R_{\kappa}(t)
\right|
}{
R_{\kappa}(t)
},
\label{eq:supp_allen_cahn_radius_error}
\end{equation}
is
\begin{equation}
\varepsilon_{R,\max}
=
6.245849\times10^{-4}.
\label{eq:supp_allen_cahn_radius_error_value}
\end{equation}
The discrepancy therefore remains below \(0.063\%\) over the complete simulated interval.

The dependence of the equivalent-radius diagnostic on the regularization parameter is also small. At \(T=100\), the choices \(\delta=0.025\), \(0.05\) and \(0.10\) give \(R_h=1.46573999\), \(1.46550643\) and \(1.46534308\), respectively. Relative to the principal choice \(\delta=0.05\), the maximum variation is \(1.593675\times10^{-4}\). The reported interface-radius conclusions are therefore insensitive to the selected regularization within the tested range.

\subsubsection{Free-energy and phase-field diagnostics}
\label{supp:allen_cahn_energy}

Equation~\eqref{eq:supp_allen_cahn_governing} is the \(L^2\)-gradient flow associated with the free-energy functional
\begin{equation}
\mathcal{E}[u]
=
\int_{\Omega}
\left[
\frac{\epsilon^2}{2}
\left|
\nabla u
\right|^2
+
\frac{1}{4}
\left(
u^2-1
\right)^2
\right]
\,\dd\Omega.
\label{eq:supp_allen_cahn_free_energy}
\end{equation}
Its variational derivative is
\begin{equation}
\frac{\delta\mathcal{E}}{\delta u}
=
-\epsilon^2\nabla^2u
+
u^3-u,
\label{eq:supp_allen_cahn_energy_variation}
\end{equation}
so that
\begin{equation}
\frac{\partial u}{\partial t}
=
-
\frac{\delta\mathcal{E}}{\delta u}.
\label{eq:supp_allen_cahn_gradient_flow}
\end{equation}
Multiplying Eq.~\eqref{eq:supp_allen_cahn_gradient_flow} by \(\partial_tu\), integrating over the periodic domain and applying integration by parts gives
\begin{equation}
\frac{\dd\mathcal{E}}{\dd t}
=
-
\int_{\Omega}
\left|
\frac{\partial u}{\partial t}
\right|^2
\,\dd\Omega
\leq0.
\label{eq:supp_allen_cahn_energy_dissipation}
\end{equation}
The continuous free energy therefore decreases monotonically.

In the discrete calculation, the spatial gradients entering the energy are evaluated spectrally:
\begin{equation}
D_xu_h
=
\mathcal{F}^{-1}
\left[
\mathrm{i}k_x\widehat{u}_h
\right],
\qquad
D_yu_h
=
\mathcal{F}^{-1}
\left[
\mathrm{i}k_y\widehat{u}_h
\right].
\label{eq:supp_allen_cahn_spectral_gradients}
\end{equation}
The discrete energy is evaluated using the periodic trapezoidal rule,
\begin{equation}
\mathcal{E}_h^n
=
\Delta x\Delta y
\sum_{i=0}^{N_x-1}
\sum_{j=0}^{N_y-1}
\left[
\frac{\epsilon^2}{2}
\left(
\left|D_xu_{i,j}^n\right|^2
+
\left|D_yu_{i,j}^n\right|^2
\right)
+
\frac{1}{4}
\left[
\left(u_{i,j}^n\right)^2-1
\right]^2
\right],
\label{eq:supp_allen_cahn_discrete_energy}
\end{equation}
and the discrete energy-dissipation rate is
\begin{equation}
\mathcal{D}_{\mathcal{E}}^{n+1/2}
=
-
\frac{
\mathcal{E}_h^{n+1}-\mathcal{E}_h^n
}{
\Delta t
}.
\label{eq:supp_allen_cahn_discrete_dissipation}
\end{equation}

For the one-dimensional equilibrium transition, the interfacial energy per unit length is
\begin{equation}
\sigma_{\epsilon}
=
\int_{-\infty}^{\infty}
\left[
\frac{\epsilon^2}{2}
\left(
\frac{\dd u}{\dd s}
\right)^2
+
\frac{1}{4}
\left(
u^2-1
\right)^2
\right]
\,\dd s
=
\frac{2\sqrt{2}}{3}\epsilon.
\label{eq:supp_allen_cahn_surface_energy}
\end{equation}
The corresponding leading-order energy of a circular interface is
\begin{equation}
\mathcal{E}_{\Gamma}(t)
\approx
2\pi\sigma_{\epsilon}R_{\kappa}(t).
\label{eq:supp_allen_cahn_interface_energy}
\end{equation}

The principal calculation gives
\begin{equation}
\mathcal{E}_h(0)=0.372206094,
\qquad
\mathcal{E}_h(100)=0.347220812.
\label{eq:supp_allen_cahn_energy_values}
\end{equation}
No sampled increase in the discrete free energy is detected:
\begin{equation}
\max_n
\max
\left(
\mathcal{E}_h^{n+1}-\mathcal{E}_h^n,
0
\right)
=
0.
\label{eq:supp_allen_cahn_energy_violation}
\end{equation}

\begin{figure}[t]
\centering
\includegraphics[width=\textwidth]{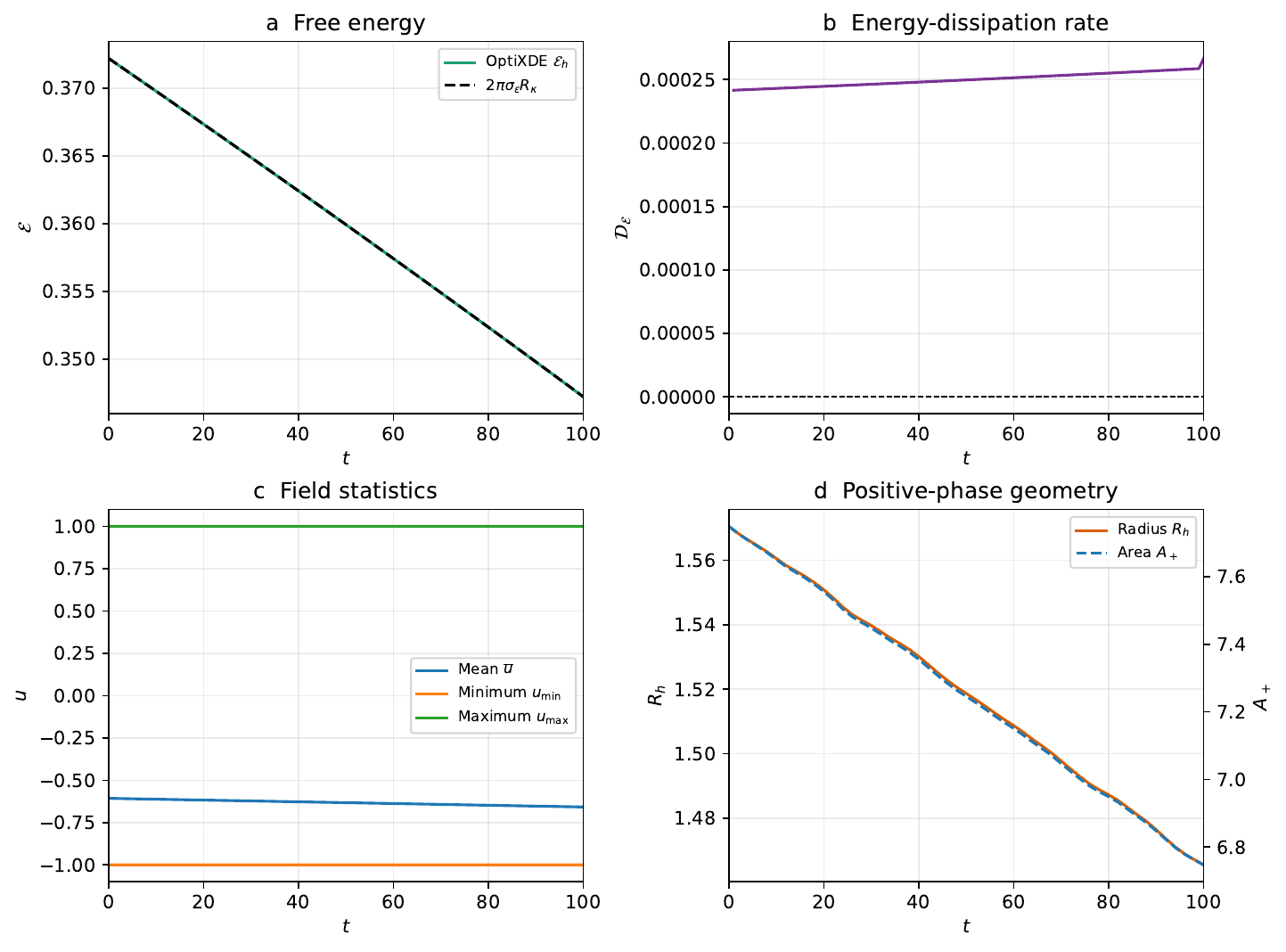}
\caption{\textbf{Physical and geometric diagnostics for the Allen--Cahn benchmark.}
\textbf{a}, Discrete free-energy evolution \(\mathcal{E}_h(t)\) together with the leading-order circular-interface estimate \(2\pi\sigma_{\epsilon}R_{\kappa}(t)\).
\textbf{b}, Discrete energy-dissipation rate \(\mathcal{D}_{\mathcal{E}}\).
\textbf{c}, Spatial mean and minimum and maximum phase-field values.
\textbf{d}, Equivalent interface radius \(R_h(t)\) and positive-phase area \(A_+(t)\).
\textbf{e}, Interface-isotropy diagnostic \(\varepsilon_{\mathrm{iso}}(t)\).
\textbf{f}, Relative differences between the curvature-flow prediction and the area-based radius \(R_h\) and the mean radius extracted directly from the \(u=0\) contour.}
\label{fig:supp_allen_cahn_diagnostics}
\end{figure}

Figure~\ref{fig:supp_allen_cahn_diagnostics}a,b confirms the gradient-flow structure numerically. The free energy decreases smoothly throughout the calculation, while the discrete dissipation rate remains nonnegative. The energy curve closely follows the leading-order interfacial estimate, showing that the decrease in energy is primarily associated with the reduction in circumference of the contracting circular interface.

Although the spatial integral of the phase field is sometimes referred to as its mass, the Allen--Cahn equation is not mass conservative. Defining
\begin{equation}
M(t)
=
\int_{\Omega}
u(x,y,t)
\,\dd\Omega,
\qquad
\overline{u}(t)
=
\frac{M(t)}{|\Omega|},
\qquad
|\Omega|=4\pi^2,
\label{eq:supp_allen_cahn_mass}
\end{equation}
integration of Eq.~\eqref{eq:supp_allen_cahn_governing} over the periodic domain gives
\begin{equation}
\frac{\dd M}{\dd t}
=
\int_{\Omega}
\left(
u-u^3
\right)
\,\dd\Omega,
\label{eq:supp_allen_cahn_mass_evolution}
\end{equation}
because the integral of the Laplacian vanishes under periodic boundary conditions.

For a sharp circular configuration with \(u\approx1\) inside the interface and \(u\approx-1\) outside,
\begin{equation}
\overline{u}_{\Gamma}(t)
\approx
\frac{R_{\kappa}^2(t)}{2\pi}-1.
\label{eq:supp_allen_cahn_sharp_mean}
\end{equation}
The numerical decrease in \(\overline{u}\) shown in Fig.~\ref{fig:supp_allen_cahn_diagnostics}c is therefore a physical consequence of the shrinking positive-phase region rather than spurious loss of a conserved quantity.

The phase-field bounds are monitored using
\begin{equation}
u_{\min}(t)
=
\min_{i,j}u_h(x_i,y_j,t),
\qquad
u_{\max}(t)
=
\max_{i,j}u_h(x_i,y_j,t).
\label{eq:supp_allen_cahn_extrema}
\end{equation}
The maximum sampled violation of the nominal interval \([-1,1]\),
\begin{equation}
\varepsilon_{\mathrm{bound}}
=
\max_t
\left\{
\max\left[-1-u_{\min}(t),0\right],
\max\left[u_{\max}(t)-1,0\right]
\right\},
\label{eq:supp_allen_cahn_bound_violation}
\end{equation}
is
\begin{equation}
\varepsilon_{\mathrm{bound}}
=
1.467695\times10^{-7}.
\label{eq:supp_allen_cahn_bound_violation_value}
\end{equation}
At the terminal time,
\begin{equation}
u_{\min}(T)=-1.00000013,
\qquad
u_{\max}(T)=1.00000014,
\label{eq:supp_allen_cahn_terminal_bounds}
\end{equation}
showing that the equilibrium bulk values are preserved to approximately seven decimal places.

To quantify directional distortion of the nominally circular interface, the \(u=0\) radius is extracted along equally spaced polar directions. If \(R_{\min}(t)\), \(R_{\max}(t)\) and \(\overline{R}(t)\) denote the minimum, maximum and mean extracted radii, the interface-isotropy error is defined as
\begin{equation}
\varepsilon_{\mathrm{iso}}(t)
=
\frac{
R_{\max}(t)-R_{\min}(t)
}{
\overline{R}(t)
}.
\label{eq:supp_allen_cahn_isotropy}
\end{equation}
For the principal calculation,
\begin{equation}
\max_t
\varepsilon_{\mathrm{iso}}(t)
=
2.144430\times10^{-4}.
\label{eq:supp_allen_cahn_isotropy_value}
\end{equation}
The small magnitude of this quantity is consistent with the nearly circular contours in Fig.~\ref{fig:supp_allen_cahn_evolution} and shows that the resolved evolution is not materially affected by Cartesian directional bias.

The leading-order curvature-flow approximation predicts extinction when \(R_{\kappa}=0\), corresponding to
\begin{equation}
t_{\mathrm{ext}}
=
\frac{R_0^2}{2\epsilon^2}
=
771.06.
\label{eq:supp_allen_cahn_extinction_time}
\end{equation}
Because \(T=100\ll t_{\mathrm{ext}}\), the present benchmark examines sustained curvature-driven contraction rather than final phase extinction.

\subsubsection{Spatial-resolution and temporal-convergence studies}
\label{supp:allen_cahn_sensitivity}

Spatial resolution is controlled primarily by the diffuse-interface thickness. For the initial hyperbolic-tangent profile, the distance between the levels \(u=-0.9\) and \(u=0.9\) is
\begin{equation}
w_{0.9}
=
2\sqrt{2}\epsilon
\operatorname{arctanh}(0.9)
=
0.166563.
\label{eq:supp_allen_cahn_interface_width}
\end{equation}
The number of sampling intervals across this transition is
\begin{equation}
n_{\Gamma}
=
\frac{w_{0.9}}{\Delta x},
\qquad
\Delta x=\frac{2\pi}{N_x}.
\label{eq:supp_allen_cahn_interface_points}
\end{equation}

\begin{figure}[t]
\centering
\suppfigure{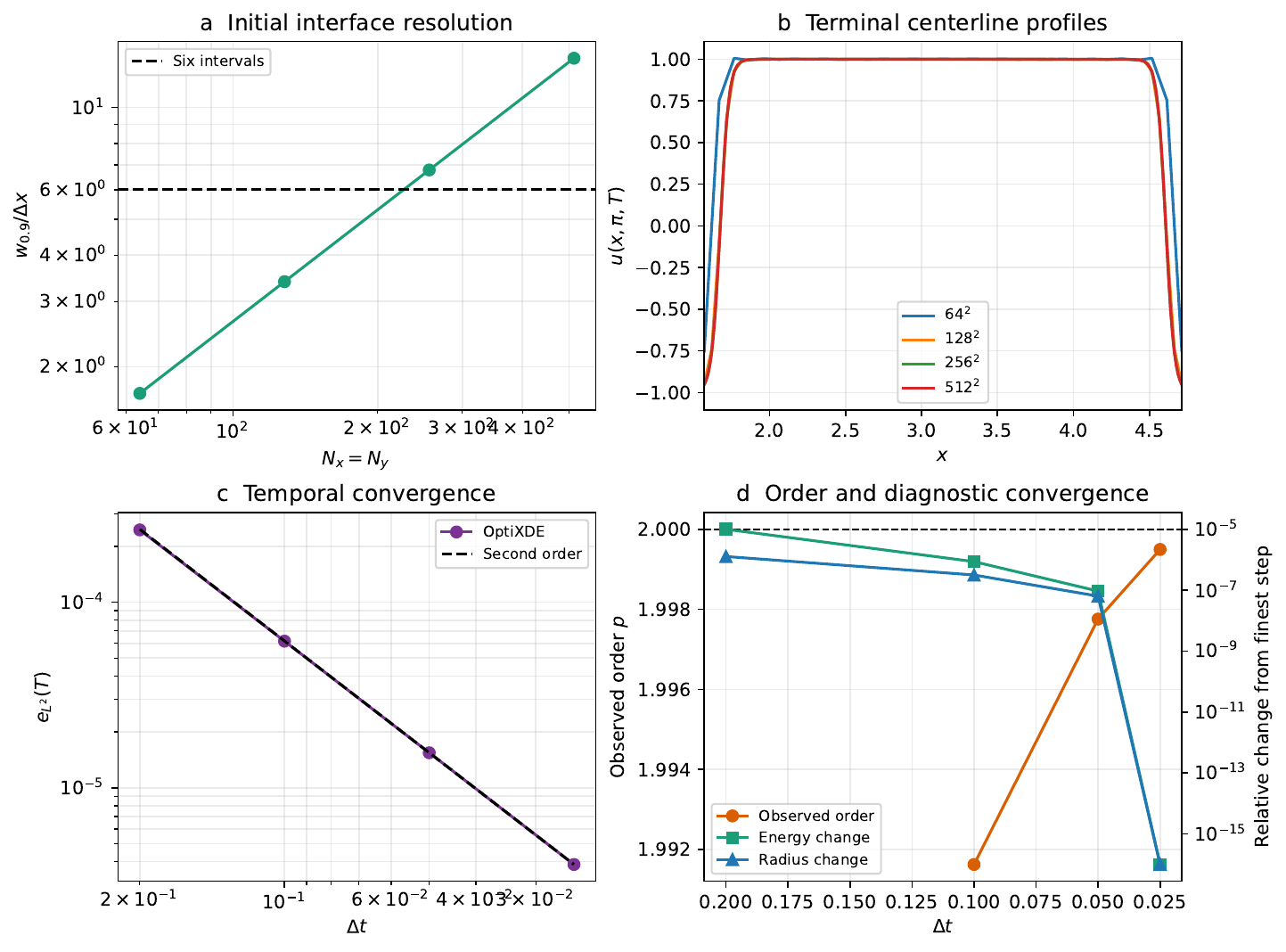}{\textwidth}
\caption{\textbf{Spatial-resolution and temporal-convergence assessment of the Allen--Cahn benchmark.}
\textbf{a}, Number of sampling intervals across the initial \(u=-0.9\) to \(u=0.9\) transition as a function of spatial resolution; the dashed line denotes a six-interval practical resolution indicator.
\textbf{b}, Terminal relative \(L^2\) error and successive-resolution difference \(\delta_N\).
\textbf{c}, Terminal centerline profiles for the tested spatial resolutions.
\textbf{d}, Terminal spectral-tail ratio and interface-isotropy error as functions of spatial resolution.
\textbf{e}, Terminal relative \(L^2\) error as a function of \(\Delta t\), together with a second-order reference slope.
\textbf{f}, Observed temporal order and relative changes in terminal free energy and equivalent radius with respect to the finest tested time increment.}
\label{fig:supp_allen_cahn_sensitivity}
\end{figure}

The spatial-resolution study uses
\begin{equation}
N_x=N_y
\in
\left\{
64,128,256,512
\right\},
\qquad
\Delta t=5.0\times10^{-2}.
\label{eq:supp_allen_cahn_spatial_resolutions}
\end{equation}
The six-interval line shown in Fig.~\ref{fig:supp_allen_cahn_sensitivity}a is used only as a practical interface-resolution indicator; sufficiency of a given resolution is assessed independently using quantitative field and physical diagnostics.

For two successive spatial resolutions, the relative field difference is defined as
\begin{equation}
\delta_N
=
\frac{
\left\|
u_N-\mathcal{R}_N u_{2N}
\right\|_2
}{
\left\|
\mathcal{R}_N u_{2N}
\right\|_2
},
\label{eq:supp_allen_cahn_successive_difference}
\end{equation}
where \(\mathcal{R}_N\) denotes periodic restriction of the finer solution to the coarser sampling points.

The relative energy diagnostic used in the spatial study is
\begin{equation}
\varepsilon_{\mathcal{E}}(T)
=
\frac{
\left|
\mathcal{E}_h(T)-\mathcal{E}_{\mathrm{ref}}(T)
\right|
}{
\left|
\mathcal{E}_{\mathrm{ref}}(T)
\right|
},
\label{eq:supp_allen_cahn_energy_error}
\end{equation}
whereas the radius error is evaluated against the curvature-flow prediction:
\begin{equation}
\varepsilon_R(T)
=
\frac{
\left|
R_h(T)-R_{\kappa}(T)
\right|
}{
R_{\kappa}(T)
}.
\label{eq:supp_allen_cahn_spatial_radius_error}
\end{equation}

Spectral resolution is monitored through the high-wave-number power fraction
\begin{equation}
\rho_{\mathrm{tail}}(T)
=
\frac{
\displaystyle
\sum_{|\mathbf{k}|\geq(2/3)k_{\mathrm{Ny}}}
\left|
\widehat{u}_h(\mathbf{k},T)
\right|^2
}{
\displaystyle
\sum_{\mathbf{k}}
\left|
\widehat{u}_h(\mathbf{k},T)
\right|^2
},
\label{eq:supp_allen_cahn_spectral_tail}
\end{equation}
where \(k_{\mathrm{Ny}}\) is the one-dimensional Nyquist wave number.

\begin{table}[t]
\centering
\small
\caption{Controlled spatial-resolution results for the two-dimensional Allen--Cahn benchmark at \(\Delta t=5.0\times10^{-2}\).}
\label{tab:supp_allen_cahn_spatial}

\begin{tabular*}{\linewidth}{@{\extracolsep{\fill}}ccccc@{}}
\toprule
\(N_x=N_y\) & \(w_{0.9}/\Delta x\) & \(\delta_N\) & \(e_{\mathrm{rel}}(T)\) & \(e_{\infty}(T)\) \\
\midrule
\(64^2\)  & \(1.70\)  & \(2.112408\times10^{-1}\) & \(2.124906\times10^{-1}\) & \(1.501291\times10^{0}\) \\
\(128^2\) & \(3.39\)  & \(2.198349\times10^{-3}\) & \(2.197974\times10^{-3}\) & \(3.425871\times10^{-2}\) \\
\(256^2\) & \(6.79\)  & \(5.722078\times10^{-7}\) & \(1.542557\times10^{-5}\) & \(1.256167\times10^{-4}\) \\
\(512^2\) & \(13.57\) & -- & \(1.544428\times10^{-5}\) & \(1.191990\times10^{-4}\) \\
\bottomrule
\end{tabular*}

\vspace{0.6em}

\begin{tabular*}{\linewidth}{@{\extracolsep{\fill}}ccccc@{}}
\toprule
\(N_x=N_y\) & \(\varepsilon_{\mathcal{E}}(T)\) & \(\varepsilon_R(T)\) & \(\varepsilon_{\mathrm{iso}}(T)\) & \(\rho_{\mathrm{tail}}(T)\) \\
\midrule
\(64^2\)  & \(1.584602\times10^{-3}\) & \(6.440024\times10^{-2}\) & \(3.964092\times10^{-2}\) & \(2.331053\times10^{-3}\) \\
\(128^2\) & \(4.865903\times10^{-4}\) & \(4.354846\times10^{-4}\) & \(2.194202\times10^{-3}\) & \(2.649222\times10^{-5}\) \\
\(256^2\) & \(8.942018\times10^{-8}\) & \(7.171250\times10^{-5}\) & \(1.784393\times10^{-4}\) & \(1.334742\times10^{-8}\) \\
\(512^2\) & \(1.118458\times10^{-7}\) & \(7.281794\times10^{-6}\) & \(2.660069\times10^{-5}\) & \(2.685689\times10^{-15}\) \\
\bottomrule
\end{tabular*}
\end{table}

At \(64^2\), only \(1.70\) intervals span the initial diffuse transition, and both the field and geometric errors are large. Increasing the resolution to \(128^2\) reduces the relative \(L^2\) error by nearly two orders of magnitude, while the \(256^2\) calculation provides \(6.79\) intervals across the initial transition and reduces the successive-resolution difference to
\begin{equation}
\delta_{256}
=
5.722078\times10^{-7}.
\label{eq:supp_allen_cahn_spatial_difference_value}
\end{equation}

The field errors at \(256^2\) and \(512^2\) remain at approximately \(1.54\times10^{-5}\) despite the six-order-of-magnitude reduction in \(\delta_N\) relative to the \(64^2\) calculation. This plateau is therefore not caused by inadequate spatial representation. It reflects the fixed temporal splitting error associated with \(\Delta t=5.0\times10^{-2}\). The same conclusion is supported by the terminal spectral-tail ratio, which decreases from \(2.331053\times10^{-3}\) at \(64^2\) to \(1.334742\times10^{-8}\) at \(256^2\) and \(2.685689\times10^{-15}\) at \(512^2\).

The interface-isotropy diagnostic also converges strongly with spatial refinement, decreasing from \(3.964092\times10^{-2}\) at \(64^2\) to \(1.784393\times10^{-4}\) at \(256^2\) and \(2.660069\times10^{-5}\) at \(512^2\). Together, these results establish that the principal \(256^2\) discretization is within the spatially resolved regime for the reported calculations.

Temporal convergence is evaluated at fixed resolution
\begin{equation}
N_x=N_y=256
\label{eq:supp_allen_cahn_temporal_resolution}
\end{equation}
using
\begin{equation}
\Delta t
\in
\left\{
2.0\times10^{-1},
1.0\times10^{-1},
5.0\times10^{-2},
2.5\times10^{-2}
\right\}.
\label{eq:supp_allen_cahn_time_steps}
\end{equation}
The observed temporal order between successive time increments is evaluated using Eq.~\eqref{eq:supp_observed_convergence_order} with $q=\Delta t$ and $e=e_{\mathrm{rel}}(T)$.

\begin{table}[t]
\centering
\caption{Temporal-convergence and terminal physical diagnostics for the two-dimensional Allen--Cahn benchmark with \(N_x=N_y=256\).}
\label{tab:supp_allen_cahn_timestep}
\begin{tabular*}{\linewidth}{@{\extracolsep{\fill}}cccccc@{}}
\toprule
\(\Delta t\) & \(e_{\mathrm{rel}}(T)\) & \(e_{\infty}(T)\) & Order \(p\) & \(\mathcal{E}_h(T)\) & \(R_h(T)\) \\
\midrule
\(2.0\times10^{-1}\) & \(2.452780\times10^{-4}\) & \(1.879522\times10^{-3}\) & -- & \(0.34722423\) & \(1.46550466\) \\
\(1.0\times10^{-1}\) & \(6.167643\times10^{-5}\) & \(4.766109\times10^{-4}\) & \(1.992\) & \(0.34722108\) & \(1.46550606\) \\
\(5.0\times10^{-2}\) & \(1.544310\times10^{-5}\) & \(1.196627\times10^{-4}\) & \(1.998\) & \(0.34722081\) & \(1.46550643\) \\
\(2.5\times10^{-2}\) & \(3.862123\times10^{-6}\) & \(3.002313\times10^{-5}\) & \(1.999\) & \(0.34722078\) & \(1.46550653\) \\
\bottomrule
\end{tabular*}
\end{table}

As \(\Delta t\) is reduced from \(2.0\times10^{-1}\) to \(2.5\times10^{-2}\), the terminal relative \(L^2\) error decreases from \(2.452780\times10^{-4}\) to \(3.862123\times10^{-6}\). Each time-step halving reduces the error by approximately a factor of four, and the observed orders,
\begin{equation}
p
=
1.992,\quad
1.998,\quad
1.999,
\label{eq:supp_allen_cahn_observed_orders}
\end{equation}
confirm the theoretical second-order temporal accuracy of the Strang composition.

The terminal physical diagnostics converge more rapidly than the complete phase field. The terminal energy changes only from \(0.34722423\) at \(\Delta t=0.2\) to \(0.34722078\) at the finest tested increment, while the equivalent radius changes from \(1.46550466\) to \(1.46550653\). The principal choice \(\Delta t=0.05\) therefore lies within a regime in which the principal physical diagnostics are already insensitive to further time-step refinement, while the field error continues to display the expected second-order reduction.

Figure~\ref{fig:supp_allen_cahn_sensitivity} separates the spatial and temporal error mechanisms. Spatially, the \(64^2\) and \(128^2\) calculations underresolve the diffuse transition, whereas the \(256^2\) and \(512^2\) terminal profiles are nearly indistinguishable and the successive-resolution difference has fallen below \(10^{-6}\). The simultaneous decay of the spectral-tail and isotropy diagnostics confirms that the resolved interface is both spectrally well represented and nearly direction independent.

Temporally, the terminal field error follows a uniform second-order trend over the complete tested range. Together with the independently verified reference solution, monotonic free-energy evolution, phase-bound error of only \(1.467695\times10^{-7}\), maximum interface-isotropy error of \(2.144430\times10^{-4}\), maximum curvature-radius discrepancy of \(6.245849\times10^{-4}\), and weak dependence of \(R_h\) on the regularization parameter \(\delta\), these results demonstrate that the principal \(256^2\), \(\Delta t=5.0\times10^{-2}\) calculation is spatially resolved, temporally convergent and physically consistent.

\FloatBarrier

\subsection{Incompressible Navier--Stokes equations: Taylor--Green vortex}
\label{supp:navier_stokes_taylor_green}

\subsubsection{Benchmark definition}
\label{supp:navier_stokes_taylor_green_definition}

The two-dimensional Taylor--Green vortex is considered as a strict analytical verification problem for the periodic incompressible Navier--Stokes solver. In contrast to the isolated-cylinder flow considered subsequently, this benchmark admits closed-form velocity, pressure and vorticity fields and therefore permits direct assessment of the spectral representation, vorticity--velocity recovery, incompressibility, pressure reconstruction and viscous dissipation.

The incompressible Navier--Stokes equations are
\begin{equation}
\frac{\partial \mathbf{u}}{\partial t}
+
(\mathbf{u}\cdot\nabla)\mathbf{u}
=
-\nabla p
+
\nu\nabla^2\mathbf{u},
\qquad
\nabla\cdot\mathbf{u}=0,
\label{eq:supp_tg_navier_stokes}
\end{equation}
where $\mathbf{u}=(u,v)^{\mathrm{T}}$ is the velocity, $p$ is the pressure normalized by the constant density, and $\nu$ is the kinematic viscosity.

The problem is defined on the doubly periodic domain
\begin{equation}
\Omega=[0,2\pi)\times[0,2\pi),
\label{eq:supp_tg_domain}
\end{equation}
with $\nu=0.01$. The analytical Taylor--Green solution is
\begin{equation}
u_{\mathrm{ex}}(x,y,t)
=
\sin x\,\cos y\,\exp(-2\nu t),
\label{eq:supp_tg_u_exact}
\end{equation}
\begin{equation}
v_{\mathrm{ex}}(x,y,t)
=
-\cos x\,\sin y\,\exp(-2\nu t),
\label{eq:supp_tg_v_exact}
\end{equation}
and
\begin{equation}
p_{\mathrm{ex}}(x,y,t)
=
\frac{1}{4}
\left(
\cos 2x+\cos 2y
\right)
\exp(-4\nu t).
\label{eq:supp_tg_p_exact}
\end{equation}

The scalar vorticity is defined as
\begin{equation}
\omega
=
\frac{\partial v}{\partial x}
-
\frac{\partial u}{\partial y},
\label{eq:supp_tg_vorticity_definition}
\end{equation}
which gives the analytical vorticity
\begin{equation}
\omega_{\mathrm{ex}}(x,y,t)
=
2\sin x\,\sin y\,\exp(-2\nu t).
\label{eq:supp_tg_omega_exact}
\end{equation}

The analytical solution at $t=0$ is used as the initial condition. The reference calculation is advanced to $T=1$ using a spatial resolution of $128^2$ and a time step of $\Delta t=0.01$.

\subsubsection{Vorticity--streamfunction spectral formulation}
\label{supp:navier_stokes_taylor_green_formulation}

For this doubly periodic benchmark, \OptiXDE{} employs a vorticity--streamfunction formulation rather than a primitive-variable velocity--pressure projection scheme. Introducing the streamfunction $\psi$ according to
\begin{equation}
u=\frac{\partial\psi}{\partial y},
\qquad
v=-\frac{\partial\psi}{\partial x},
\label{eq:supp_tg_streamfunction_velocity}
\end{equation}
gives
\begin{equation}
-\nabla^2\psi=\omega.
\label{eq:supp_tg_streamfunction_poisson}
\end{equation}
The Navier--Stokes equations are then advanced through the scalar vorticity transport equation
\begin{equation}
\frac{\partial\omega}{\partial t}
+
u\frac{\partial\omega}{\partial x}
+
v\frac{\partial\omega}{\partial y}
=
\nu\nabla^2\omega.
\label{eq:supp_tg_vorticity_equation}
\end{equation}

For each non-zero Fourier mode, the streamfunction is recovered directly from
\begin{equation}
\widehat{\psi}
=
\frac{\widehat{\omega}}{K^2},
\qquad
K^2=k_x^2+k_y^2,
\label{eq:supp_tg_streamfunction_spectral}
\end{equation}
with the zero mode fixed to zero. The velocity components are subsequently obtained as
\begin{equation}
\widehat{u}
=
\mathrm{i}k_y\widehat{\psi},
\qquad
\widehat{v}
=
-\mathrm{i}k_x\widehat{\psi}.
\label{eq:supp_tg_velocity_spectral}
\end{equation}
This construction satisfies the incompressibility constraint analytically in spectral space because
\begin{equation}
\mathrm{i}k_x\widehat{u}
+
\mathrm{i}k_y\widehat{v}
=
0.
\label{eq:supp_tg_spectral_incompressibility}
\end{equation}

The viscous contribution is advanced using the general diffusion propagator in Eq.~\eqref{eq:supp_diffusion_propagator}, with $u=\omega$, $\kappa=\nu$ and $|\mathbf{k}|^2=K^2$. A complete time step uses a symmetric composition consisting of a viscous half-step, an explicit nonlinear full-step and a second viscous half-step. The nonlinear vorticity-advection term
\begin{equation}
\mathcal{N}(\omega)
=
-u\omega_x-v\omega_y
\label{eq:supp_tg_nonlinear_operator}
\end{equation}
is evaluated pseudospectrally using second-order Runge--Kutta integration, with the two-dimensional $2/3$ rule applied to the nonlinear products.

The present Taylor--Green mode has the additional analytical property
\begin{equation}
\mathbf{u}_{\mathrm{ex}}\cdot\nabla\omega_{\mathrm{ex}}
=
u_{\mathrm{ex}}
\frac{\partial\omega_{\mathrm{ex}}}{\partial x}
+
v_{\mathrm{ex}}
\frac{\partial\omega_{\mathrm{ex}}}{\partial y}
=
0.
\label{eq:supp_tg_nonlinear_cancellation}
\end{equation}
Consequently, the vorticity dynamics of this particular solution reduce exactly to viscous modal decay even though the complete nonlinear solver is executed. The numerical quantity $\|\mathbf{u}\cdot\nabla\omega\|_{\infty}$ therefore provides an additional internal consistency check for velocity recovery, spectral differentiation and pseudospectral nonlinear evaluation.

Pressure is reconstructed independently from the computed velocity field. Taking the divergence of Eq.~\eqref{eq:supp_tg_navier_stokes} and using $\nabla\cdot\mathbf{u}=0$ gives
\begin{equation}
-\nabla^2p
=
\nabla\cdot
\left[
(\mathbf{u}\cdot\nabla)\mathbf{u}
\right].
\label{eq:supp_tg_pressure_poisson}
\end{equation}
The periodic pressure Poisson equation is solved spectrally with the spatial mean of $p$ fixed to zero. The pressure used for verification is therefore reconstructed from the computed velocity rather than prescribed from the analytical Taylor--Green expression.

\subsubsection{Field accuracy and incompressibility}
\label{supp:navier_stokes_taylor_green_accuracy}

\begin{figure}[t]
\centering
\includegraphics[width=\textwidth]{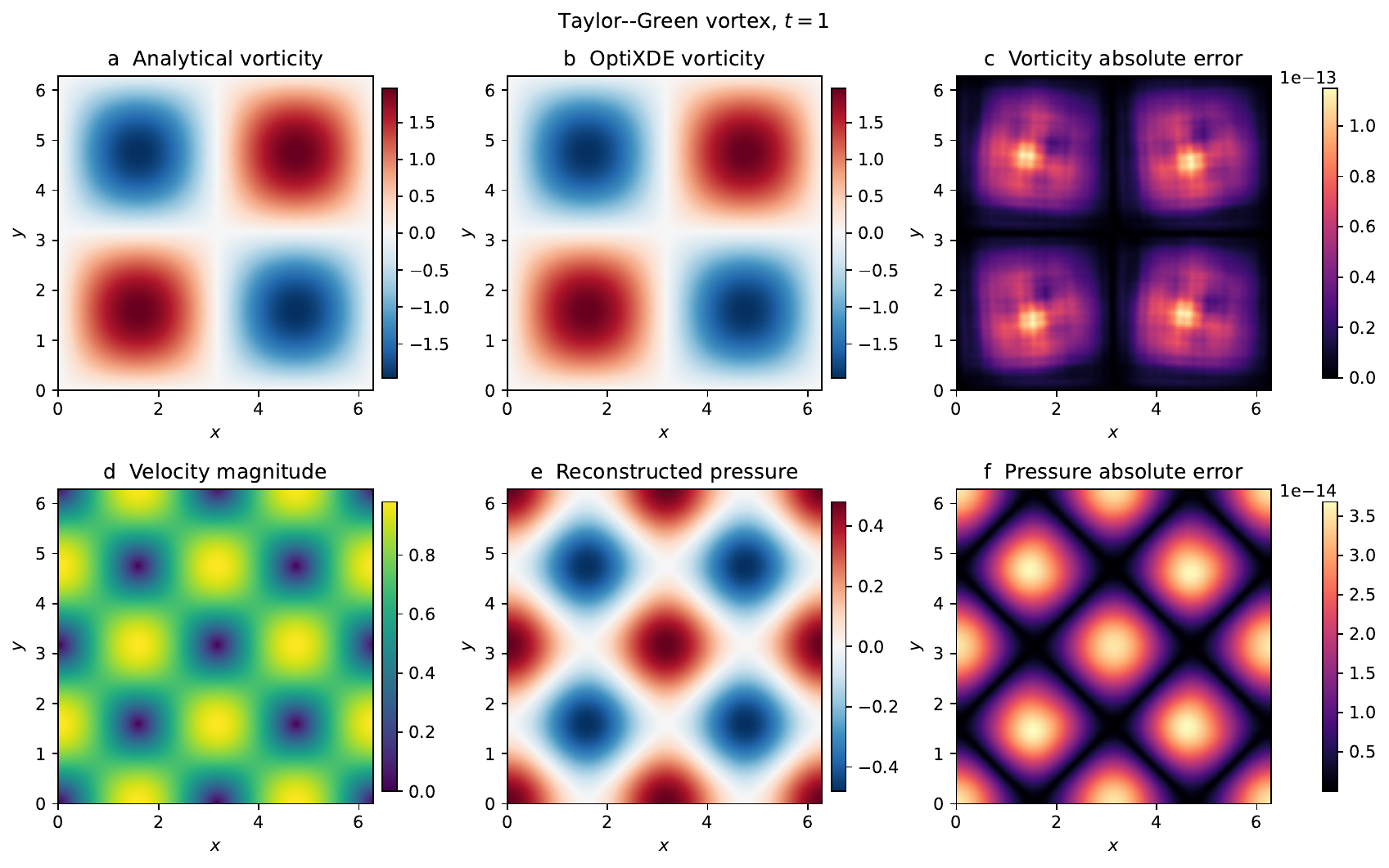}
\caption{\textbf{Field verification for the Taylor--Green vortex at $t=1$.} \textbf{a}, Analytical vorticity. \textbf{b}, Vorticity computed by \OptiXDE{}. \textbf{c}, Pointwise absolute vorticity error. \textbf{d}, Velocity magnitude recovered through the spectral streamfunction formulation. \textbf{e}, Pressure reconstructed from the numerical velocity through the periodic pressure Poisson equation. \textbf{f}, Pointwise absolute pressure error. The calculation uses a $128^2$ periodic resolution, $\nu=0.01$ and $\Delta t=0.01$.}
\label{fig:supp_tg_fields}
\end{figure}

Figure~\ref{fig:supp_tg_fields} compares the analytical and numerical fields at $T=1$. The analytical and \OptiXDE{} vorticity distributions are visually indistinguishable at the plotting scale, while the pointwise vorticity error remains of order $10^{-13}$. The velocity recovered through the streamfunction retains the expected Taylor--Green structure, and the independently reconstructed pressure reproduces the analytical pressure field with pointwise errors of approximately $10^{-14}$--$10^{-13}$.

At $T=1$, the relative field errors are
\begin{equation}
e_{\omega,L_2}
=
3.8340\times10^{-14},
\qquad
e_{u,L_2}
=
3.7865\times10^{-14},
\qquad
e_{v,L_2}
=
3.7853\times10^{-14},
\label{eq:supp_tg_terminal_velocity_errors}
\end{equation}
and
\begin{equation}
e_{p,L_2}
=
7.3426\times10^{-14}.
\label{eq:supp_tg_terminal_pressure_error}
\end{equation}
The corresponding incompressibility residual is
\begin{equation}
\left\|
\nabla\cdot\mathbf{u}
\right\|_{\infty}
=
2.4314\times10^{-14}.
\label{eq:supp_tg_terminal_divergence}
\end{equation}
The velocity field therefore remains divergence-free to within double-precision round-off, consistent with the spectral construction in Eq.~\eqref{eq:supp_tg_spectral_incompressibility}.

\subsubsection{Kinetic-energy and enstrophy decay}
\label{supp:navier_stokes_taylor_green_dissipation}

The kinetic energy and enstrophy are defined as
\begin{equation}
E(t)
=
\frac{1}{2}
\int_{\Omega}
\left(
u^2+v^2
\right)
\,\mathrm{d}\Omega
\label{eq:supp_tg_kinetic_energy}
\end{equation}
and
\begin{equation}
Z(t)
=
\frac{1}{2}
\int_{\Omega}
\omega^2
\,\mathrm{d}\Omega,
\label{eq:supp_tg_enstrophy}
\end{equation}
respectively. For the analytical Taylor--Green solution,
\begin{equation}
E_{\mathrm{ex}}(t)
=
\pi^2\exp(-4\nu t),
\qquad
Z_{\mathrm{ex}}(t)
=
2\pi^2\exp(-4\nu t).
\label{eq:supp_tg_energy_exact}
\end{equation}

\begin{figure}[t]
\centering
\includegraphics[width=\textwidth]{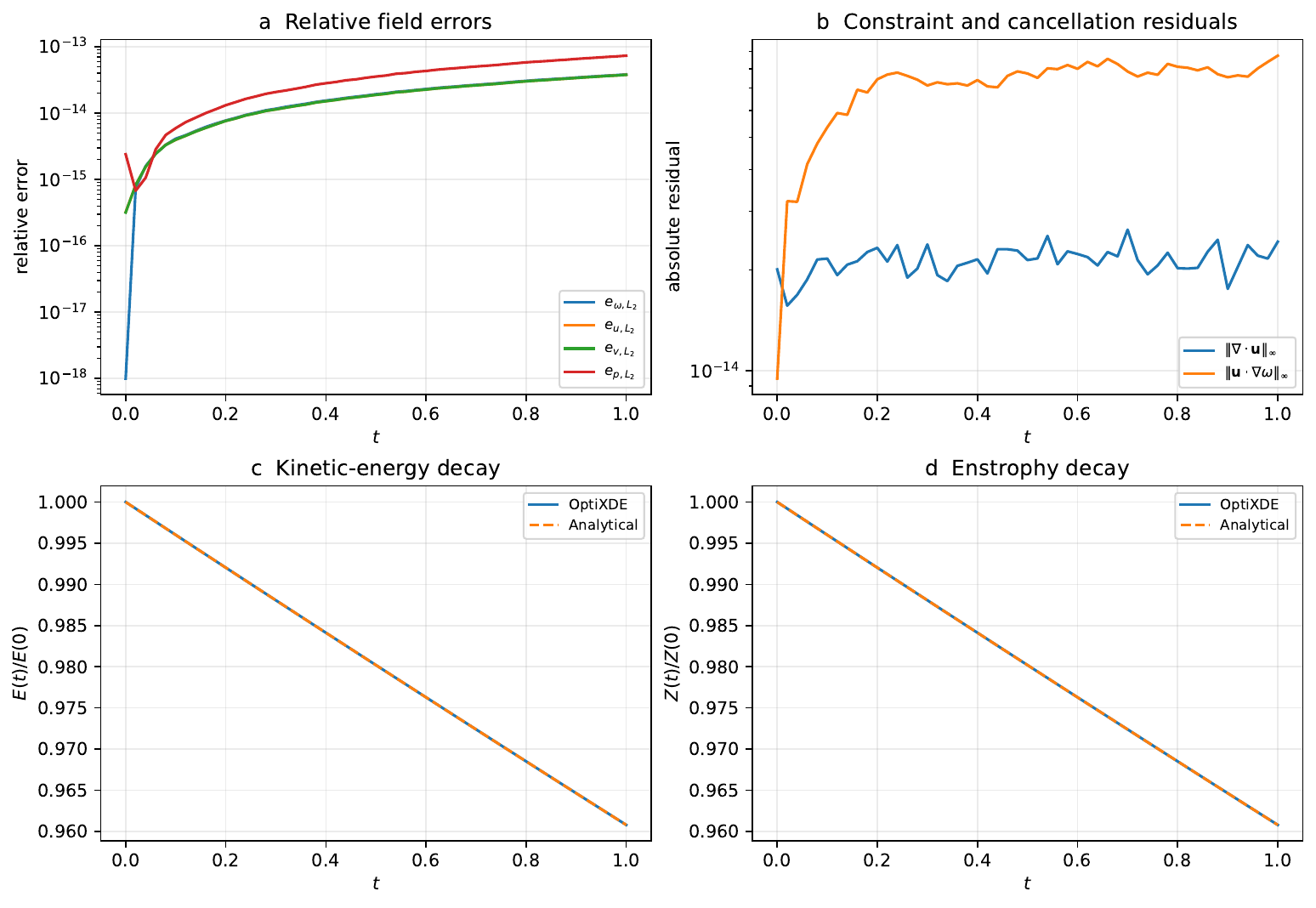}
\caption{\textbf{Accuracy and physical diagnostics for the Taylor--Green vortex.} \textbf{a}, Relative $L_2$ errors of the vorticity, two velocity components and reconstructed pressure. \textbf{b}, Maximum-norm incompressibility residual $\|\nabla\cdot\mathbf{u}\|_{\infty}$ and nonlinear-cancellation residual $\|\mathbf{u}\cdot\nabla\omega\|_{\infty}$. \textbf{c}, Numerical and analytical kinetic-energy decay normalized by their initial values. \textbf{d}, Numerical and analytical enstrophy decay normalized by their initial values.}
\label{fig:supp_tg_diagnostics}
\end{figure}

Figure~\ref{fig:supp_tg_diagnostics} shows the complete temporal verification. The relative errors in $\omega$, $u$, $v$ and $p$ remain between approximately $10^{-15}$ and $10^{-13}$ over the complete integration interval. The incompressibility residual and nonlinear-cancellation residual remain close to the floating-point-error regime, while the numerical kinetic-energy and enstrophy histories are visually indistinguishable from their analytical exponential-decay curves.

At the terminal time, the relative kinetic-energy and enstrophy errors are
\begin{equation}
e_E
=
7.5306\times10^{-14},
\qquad
e_Z
=
7.4744\times10^{-14},
\label{eq:supp_tg_terminal_energy_errors}
\end{equation}
while the nonlinear-cancellation diagnostic gives
\begin{equation}
\left\|
\mathbf{u}\cdot\nabla\omega
\right\|_{\infty}
=
8.7578\times10^{-14}.
\label{eq:supp_tg_terminal_nonlinear_residual}
\end{equation}
The simultaneous agreement of the field variables, divergence constraint, nonlinear-cancellation condition and integral dissipation measures confirms that the periodic vorticity--streamfunction implementation reproduces this analytical Navier--Stokes solution to near machine precision.

\subsubsection{Resolution and time-step sensitivity}
\label{supp:navier_stokes_taylor_green_sensitivity}

The dependence of the verification results on spatial resolution and time-step size is examined through two controlled studies. The spatial-resolution study uses
\begin{equation}
N=32,\ 64,\ 128,\ 256
\label{eq:supp_tg_resolution_set}
\end{equation}
with $\Delta t=0.01$, whereas the time-step study uses
\begin{equation}
\Delta t
=
0.04,\ 0.02,\ 0.01,\ 0.005
\label{eq:supp_tg_timestep_set}
\end{equation}
at a fixed spatial resolution of $128^2$.

\begin{figure}[t]
\centering
\includegraphics[width=\textwidth]{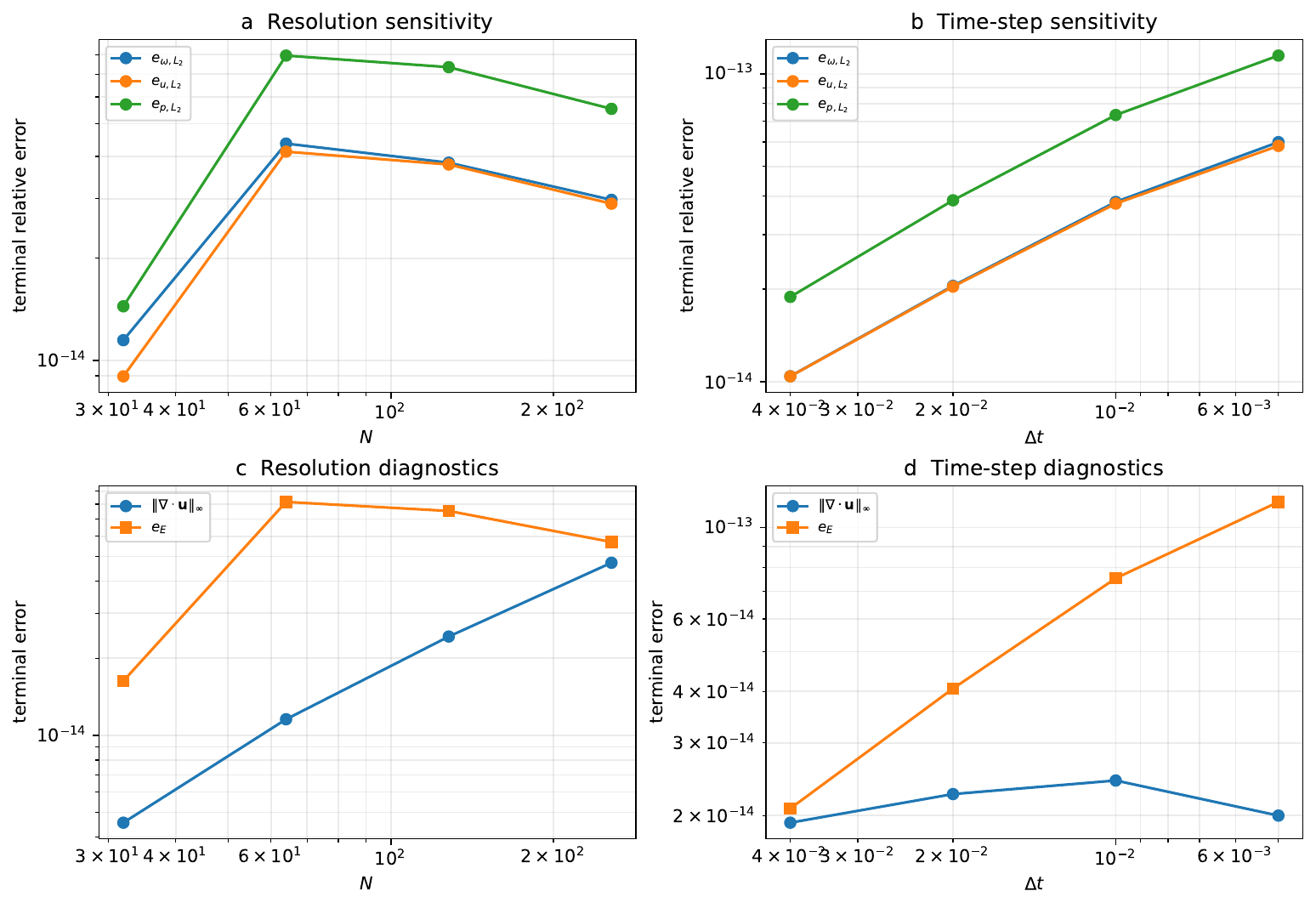}
\caption{\textbf{Resolution and time-step sensitivity of the Taylor--Green vortex verification.} \textbf{a}, Terminal relative errors in vorticity, velocity and reconstructed pressure as functions of spatial resolution. \textbf{b}, Corresponding terminal errors as functions of time-step size on a $128^2$ resolution. \textbf{c}, Terminal incompressibility and kinetic-energy diagnostics for the spatial-resolution study. \textbf{d}, Corresponding diagnostics for the time-step study. All tested configurations remain in or close to the double-precision round-off regime rather than exhibiting a conventional algebraic convergence trend.}
\label{fig:supp_tg_sensitivity}
\end{figure}

As shown in Fig.~\ref{fig:supp_tg_sensitivity}a,c, all tested resolutions already resolve the low-order Fourier content of the analytical Taylor--Green solution. The terminal field errors consequently remain in the $10^{-14}$--$10^{-13}$ range rather than exhibiting a conventional algebraic spatial-convergence trend. The terminal vorticity errors are $1.146\times10^{-14}$, $4.366\times10^{-14}$, $3.834\times10^{-14}$ and $2.978\times10^{-14}$ for $N=32$, $64$, $128$ and $256$, respectively. This non-monotonic variation is consistent with a finite-precision plateau: once the analytical Fourier mode is represented by the discrete transform basis, further spatial refinement does not reduce an appreciable spatial truncation error.

The time-step study shown in Fig.~\ref{fig:supp_tg_sensitivity}b,d exhibits a corresponding finite-precision behaviour. The terminal vorticity errors are $1.042\times10^{-14}$, $2.048\times10^{-14}$, $3.834\times10^{-14}$ and $5.995\times10^{-14}$ for $\Delta t=0.04$, $0.02$, $0.01$ and $0.005$, respectively. The pressure and kinetic-energy diagnostics exhibit the same mild increase as the time step is reduced. This trend is not interpreted as a loss of temporal consistency. For the present Taylor--Green mode, the nonlinear vorticity-advection term cancels analytically according to Eq.~\eqref{eq:supp_tg_nonlinear_cancellation}, while the remaining viscous Fourier mode is propagated by the exact exponential operator in Eq.~\eqref{eq:supp_diffusion_propagator}. Reducing $\Delta t$ therefore increases the number of repeated transform--propagate--inverse-transform cycles without reducing an appreciable temporal truncation error, and the remaining variation is dominated by accumulated floating-point round-off.

The sensitivity results therefore show that the near-machine-precision agreement of the reference calculation is not specific to a particular spatial resolution or time step. Once the Taylor--Green Fourier mode is resolved, the remaining discrepancy is governed predominantly by finite-precision effects associated with repeated Fourier transforms and spectral operations.

\subsubsection{Reproducibility information}
\label{supp:navier_stokes_taylor_green_reproducibility}

The benchmark-specific numerical settings are summarized in Table~\ref{tab:supp_taylor_green_settings}. The common software environment, hardware specifications, arithmetic precision and timing protocol are reported in Section~\ref{supp:computational_environment}. The resolution and time-step sensitivity studies are reported in Section~\ref{supp:navier_stokes_taylor_green_sensitivity} and are therefore not repeated here. Because the nonlinear vorticity-advection term vanishes analytically for the present Taylor--Green mode, this benchmark is used specifically as a strict verification of the periodic incompressible-flow formulation; non-trivial nonlinear vortex dynamics are examined separately using the $Re=200$ isolated-cylinder benchmark.

The source code, input parameters and post-processing scripts required to reproduce this benchmark will be made publicly available upon publication at
\href{https://github.com/USERNAME/OptiXDE/tree/main/examples/navier_stokes_taylor_green}{\texttt{examples/navier\_stokes\_taylor\_green}}.

\begin{table}[h]
\centering
\caption{Numerical settings for the Taylor--Green vortex benchmark.}
\label{tab:supp_taylor_green_settings}
\begin{tabular}{ll}
\toprule
Setting & Value \\
\midrule
Physical domain & $\Omega=[0,2\pi)\times[0,2\pi)$ \\
Boundary condition & Doubly periodic \\
Kinematic viscosity & $\nu=0.01$ \\
Initial velocity & $u_0=\sin x\cos y,\quad v_0=-\cos x\sin y$ \\
Initial vorticity & $\omega_0=2\sin x\sin y$ \\
Flow formulation & Vorticity--streamfunction \\
Spatial transform & Two-dimensional discrete Fourier transform \\
Streamfunction recovery & $\widehat{\psi}=\widehat{\omega}/K^2$ for $K^2>0$ \\
Velocity recovery & $\widehat{u}=\mathrm{i}k_y\widehat{\psi},\quad \widehat{v}=-\mathrm{i}k_x\widehat{\psi}$ \\
Viscous propagation & Closed-form exponential spectral propagator \\
Nonlinear integration & Second-order Runge--Kutta \\
Nonlinear dealiasing & Two-dimensional $2/3$ rule \\
Pressure reconstruction & Periodic spectral pressure Poisson equation \\
Pressure normalization & Zero spatial mean \\
Reference spatial resolution & $128^2$ \\
Reference time step & $\Delta t=0.01$ \\
Terminal time & $T=1$ \\
Resolution-sensitivity cases & $32^2$, $64^2$, $128^2$ and $256^2$ \\
Time-step-sensitivity cases & $0.04$, $0.02$, $0.01$ and $0.005$ \\
Arithmetic precision & Double precision \\
\bottomrule
\end{tabular}
\end{table}
\FloatBarrier

\subsection{Incompressible Navier--Stokes equations: flow past a circular cylinder}
\label{supp:navier_stokes_cylinder}

\subsubsection{Benchmark definition}
\label{supp:navier_stokes_cylinder_definition}

The two-dimensional flow past an isolated circular cylinder is considered as a nonlinear physical demonstration of the incompressible-flow capability of \OptiXDE{}. In contrast to the Taylor--Green vortex in Section~\ref{supp:navier_stokes_taylor_green}, this problem does not admit a closed-form solution and is therefore not used as the primary analytical-accuracy benchmark. Instead, it examines the ability of the transform-based formulation to represent an embedded solid boundary, maintain incompressibility and resolve sustained vortex shedding over a long-time nonlinear simulation.

The Reynolds number is defined as
\begin{equation}
Re
=
\frac{U_{\infty}D}{\nu}
=
200,
\label{eq:supp_cyl_reynolds}
\end{equation}
where $D$ is the cylinder diameter and $U_{\infty}$ is the free-stream velocity. The dimensionless parameters used throughout the benchmark are
\begin{equation}
D=1,
\qquad
U_{\infty}=1,
\qquad
\nu=0.005.
\label{eq:supp_cyl_physical_parameters}
\end{equation}

The cylinder is embedded in the doubly periodic Fourier domain
\begin{equation}
\widetilde{\Omega}
=
[0,40D)\times[0,20D),
\label{eq:supp_cyl_domain}
\end{equation}
with centre
\begin{equation}
(x_c,y_c)
=
(10D,10D).
\label{eq:supp_cyl_center}
\end{equation}
The transverse domain height of $20D$ provides a large separation between periodic cylinder images. In the streamwise direction, a smooth fringe region of width $5D$ is introduced around the periodic boundary to relax the downstream wake toward the uniform inflow before it re-enters the upstream side of the domain.

A small asymmetric transverse-velocity perturbation is introduced initially to trigger vortex shedding. The initial velocity is
\begin{equation}
u(x,y,0)=U_{\infty},
\label{eq:supp_cyl_initial_u}
\end{equation}
and
\begin{equation}
v(x,y,0)
=
0.02
\exp
\left[
-\frac{(x-x_c-D)^2}{(0.55D)^2}
-\frac{(y-y_c-0.30D)^2}{(0.50D)^2}
\right].
\label{eq:supp_cyl_initial_v}
\end{equation}
The initial field is projected onto the divergence-free Fourier subspace before time integration.

All final production calculations use
\begin{equation}
\Delta t=0.01,
\qquad
T=100,
\label{eq:supp_cyl_time_parameters}
\end{equation}
and the statistically stationary interval
\begin{equation}
50\leq t\leq100
\label{eq:supp_cyl_statistics_interval}
\end{equation}
is used to evaluate the force and shedding-frequency statistics.

\subsubsection{Vorticity--streamfunction formulation and embedded-cylinder treatment}
\label{supp:navier_stokes_cylinder_formulation}

As in the Taylor--Green calculation, the periodic incompressible flow is represented in terms of the scalar vorticity
\begin{equation}
\omega
=
\frac{\partial v}{\partial x}
-
\frac{\partial u}{\partial y}.
\label{eq:supp_cyl_vorticity_definition}
\end{equation}
Between the embedded-boundary and fringe constraint operations, the vorticity satisfies
\begin{equation}
\frac{\partial\omega}{\partial t}
+
u\frac{\partial\omega}{\partial x}
+
v\frac{\partial\omega}{\partial y}
=
\nu\nabla^2\omega.
\label{eq:supp_cyl_vorticity_equation}
\end{equation}

The velocity is decomposed into a uniform mean flow and a zero-mean fluctuating component,
\begin{equation}
u
=
U_{\infty}
+
\frac{\partial\psi}{\partial y},
\qquad
v
=
-\frac{\partial\psi}{\partial x},
\label{eq:supp_cyl_streamfunction_velocity}
\end{equation}
where
\begin{equation}
-\nabla^2\psi=\omega.
\label{eq:supp_cyl_streamfunction_poisson}
\end{equation}
For each non-zero Fourier mode,
\begin{equation}
\widehat{\psi}
=
\frac{\widehat{\omega}}{K^2},
\qquad
K^2=k_x^2+k_y^2,
\label{eq:supp_cyl_streamfunction_spectral}
\end{equation}
and the velocity fluctuations are recovered spectrally. The uniform streamwise mode $U_{\infty}$ is then restored explicitly.

The cylinder is represented by a compact-support smooth mask $\chi$. Defining
\begin{equation}
r
=
\sqrt{(x-x_c)^2+(y-y_c)^2},
\qquad
s
=
\frac{D}{2}-r,
\label{eq:supp_cyl_mask_distance}
\end{equation}
the compact mask is
\begin{equation}
\chi(s)
=
\begin{cases}
1,
& s\geq\varepsilon_{\chi},
\\[4pt]
\dfrac{1}{2}
\left[
1+\dfrac{s}{\varepsilon_{\chi}}
+
\dfrac{1}{\pi}
\sin
\left(
\dfrac{\pi s}{\varepsilon_{\chi}}
\right)
\right],
& |s|<\varepsilon_{\chi},
\\[8pt]
0,
& s\leq-\varepsilon_{\chi}.
\end{cases}
\label{eq:supp_cyl_compact_mask}
\end{equation}
The transition half-width is fixed at
\begin{equation}
\varepsilon_{\chi}=1.0h,
\label{eq:supp_cyl_mask_width}
\end{equation}
where
\begin{equation}
h=\max(\Delta x,\Delta y).
\end{equation}
The mask therefore has exactly finite support outside the nominal cylinder boundary, rather than the long weak tail associated with a hyperbolic-tangent mask.

The stationary no-slip cylinder is imposed through exact Brinkman relaxation in velocity space. Over an operator interval $\tau$,
\begin{equation}
u^{\star}
=
u
\exp
\left(
-\frac{\tau\chi}{\eta}
\right),
\qquad
v^{\star}
=
v
\exp
\left(
-\frac{\tau\chi}{\eta}
\right),
\label{eq:supp_cyl_brinkman_relaxation}
\end{equation}
with the final penalty parameter
\begin{equation}
\eta=0.0025.
\label{eq:supp_cyl_penalty}
\end{equation}
The values $\eta=0.0025$ and $\varepsilon_{\chi}=1.0h$ were selected from controlled penalty-parameter and mask-width sensitivity calculations and were subsequently held fixed for all production-resolution calculations.

A smooth streamwise fringe mask $\chi_f$ is defined using the periodic distance
\begin{equation}
d_f(x)=\min(x,L_x-x),
\label{eq:supp_cyl_fringe_distance}
\end{equation}
such that
\begin{equation}
\chi_f(x)
=
\begin{cases}
\dfrac{1}{2}
\left[
1+
\cos
\left(
\dfrac{\pi d_f}{w_f}
\right)
\right],
& d_f<w_f,
\\[6pt]
0,
& d_f\geq w_f,
\end{cases}
\label{eq:supp_cyl_fringe_mask}
\end{equation}
where
\begin{equation}
w_f=5D.
\end{equation}
Following the cylinder relaxation, the velocity is relaxed toward the uniform target state $(U_{\infty},0)$ according to
\begin{equation}
u^{\star\star}
=
U_{\infty}
+
\left(
u^{\star}-U_{\infty}
\right)
\exp
\left(
-\frac{\tau\chi_f}{\tau_f}
\right),
\label{eq:supp_cyl_fringe_u}
\end{equation}
and
\begin{equation}
v^{\star\star}
=
v^{\star}
\exp
\left(
-\frac{\tau\chi_f}{\tau_f}
\right),
\label{eq:supp_cyl_fringe_v}
\end{equation}
with
\begin{equation}
\tau_f=0.25.
\label{eq:supp_cyl_fringe_tau}
\end{equation}

Because the pointwise cylinder and fringe operations do not individually preserve incompressibility, each combined constraint operation is followed by a Fourier Helmholtz projection. For each non-zero Fourier mode,
\begin{equation}
\widehat{\mathbf{u}}^{\,\perp}
=
\left(
\mathbf{I}
-
\frac{\mathbf{k}\mathbf{k}^{\mathrm{T}}}{K^2}
\right)
\widehat{\mathbf{u}},
\label{eq:supp_cyl_helmholtz_projection}
\end{equation}
after which the projected velocity is converted back to the vorticity representation.

One complete time step is written schematically as
\begin{equation}
\mathcal{S}_{\Delta t}
=
\mathcal{C}_{\Delta t/2}
\mathcal{D}_{\Delta t/2}
\mathcal{A}_{\Delta t}
\mathcal{D}_{\Delta t/2}
\mathcal{C}_{\Delta t/2},
\label{eq:supp_cyl_time_splitting}
\end{equation}
where $\mathcal{C}$ denotes the combined cylinder/fringe relaxation and incompressibility projection, $\mathcal{D}$ denotes the closed-form viscous propagator given by Eq.~\eqref{eq:supp_diffusion_propagator} with $u=\omega$, $\kappa=\nu$ and $|\mathbf{k}|^2=K^2$, and $\mathcal{A}$ denotes the nonlinear vorticity-advection update. The nonlinear full-step is advanced using second-order Runge--Kutta integration. Nonlinear products are evaluated pseudospectrally and filtered using the two-dimensional $2/3$ de-aliasing rule.

The pressure shown in the field diagnostics is reconstructed a posteriori and is not used as the primary time-evolution variable. Defining the continuous Brinkman and fringe forcing as
\begin{equation}
\mathbf{f}_{\mathrm{cyl}}
=
-\frac{\chi}{\eta}\mathbf{u},
\label{eq:supp_cyl_brinkman_force}
\end{equation}
and
\begin{equation}
\mathbf{f}_{\mathrm{fringe}}
=
-\frac{\chi_f}{\tau_f}
\begin{pmatrix}
u-U_{\infty}\\
v
\end{pmatrix},
\label{eq:supp_cyl_fringe_force}
\end{equation}
the zero-mean pressure is obtained from
\begin{equation}
-\nabla^2p
=
\nabla\cdot
\left[
(\mathbf{u}\cdot\nabla)\mathbf{u}
\right]
-
\nabla\cdot
\left(
\mathbf{f}_{\mathrm{cyl}}
+
\mathbf{f}_{\mathrm{fringe}}
\right),
\label{eq:supp_cyl_pressure_poisson}
\end{equation}
which is inverted directly in Fourier space.

The hydrodynamic force is evaluated from the momentum actually removed by the exact cylinder-relaxation operation. For each constraint half-step,
\begin{equation}
\mathbf{J}_{\mathrm{cyl}}
=
\rho
\int_{\widetilde{\Omega}}
\left(
\mathbf{u}^{-}
-
\mathbf{u}^{+}
\right)
\,\mathrm{d}\Omega,
\label{eq:supp_cyl_impulse}
\end{equation}
where $\mathbf{u}^{-}$ and $\mathbf{u}^{+}$ denote the velocity immediately before and after the cylinder relaxation, respectively. The impulses from the two constraint half-steps are accumulated over the diagnostic interval and divided by its duration to obtain the cylinder force. With $\rho=1$, the drag and lift coefficients are
\begin{equation}
C_D
=
\frac{F_x}
{\tfrac{1}{2}\rho U_{\infty}^2D},
\qquad
C_L
=
\frac{F_y}
{\tfrac{1}{2}\rho U_{\infty}^2D}.
\label{eq:supp_cyl_force_coefficients}
\end{equation}

\subsubsection{Spatial-resolution assessment}
\label{supp:navier_stokes_cylinder_resolution}

A controlled spatial-resolution study is performed using
\begin{equation}
N_D=24,\ 32,\ 40,
\label{eq:supp_cyl_resolution_set}
\end{equation}
where $N_D$ denotes the number of points per cylinder diameter. Because the physical domain is fixed at $40D\times20D$, the corresponding Fourier resolutions are $960\times480$, $1280\times640$ and $1600\times800$, respectively. The physical parameters, time step, final time, compact-mask definition, penalty parameter, fringe treatment and statistical interval are identical in all three calculations.

\begin{table}[t]
\centering
\caption{Spatial-resolution assessment for the two-dimensional isolated-cylinder calculation at $Re=200$. The Strouhal number is evaluated independently from the lift coefficient and from a downstream wake probe.}
\label{tab:supp_cylinder_resolution}
\begin{tabular*}{\linewidth}{@{\extracolsep{\fill}}ccccccc@{}}
\toprule
$N_D$
& Resolution
& $\overline{C_D}$
& $C_{L,\mathrm{rms}}$
& $St_{C_L}$
& $St_{\mathrm{probe}}$
& $\max\|\nabla\cdot\mathbf{u}\|_{\infty}$ \\
\midrule
24 & $960\times480$  & 1.49722 & 0.56529 & 0.19753 & 0.19754 & $4.78\times10^{-14}$ \\
32 & $1280\times640$ & 1.47520 & 0.54654 & 0.19884 & 0.19885 & $6.84\times10^{-14}$ \\
40 & $1600\times800$ & 1.46716 & 0.54049 & 0.19955 & 0.19957 & $9.79\times10^{-14}$ \\
\bottomrule
\end{tabular*}
\end{table}

The principal wake statistics vary systematically with increasing spatial resolution. From $N_D=24$ to $N_D=32$, the mean drag coefficient changes by approximately $1.47\%$, whereas the change from $N_D=32$ to $N_D=40$ decreases to approximately $0.55\%$. The corresponding changes in $C_{L,\mathrm{rms}}$ decrease from approximately $3.32\%$ to $1.11\%$. The Strouhal number changes by approximately $0.66\%$ between $N_D=24$ and $N_D=32$ and by only approximately $0.36\%$ between $N_D=32$ and $N_D=40$. The decreasing successive changes indicate that the principal integral-force and shedding-frequency diagnostics are approaching a resolution-independent regime.

The statistical stationarity of the finest calculation was additionally examined by dividing the sampling window into
\begin{equation}
50\leq t\leq75
\end{equation}
and
\begin{equation}
75\leq t\leq100.
\end{equation}
The corresponding mean drag coefficients are
\begin{equation}
\overline{C_D}
=
1.46696
\qquad\text{and}\qquad
1.46741,
\label{eq:supp_cyl_split_drag}
\end{equation}
respectively, a difference of approximately $0.03\%$. The corresponding lift root-mean-square values are $0.54146$ and $0.54000$. These results indicate that the interval $50\leq t\leq100$ is sufficiently stationary for the reported force and frequency statistics. The $N_D=40$ calculation is therefore adopted for the field and diagnostic results below.

\subsubsection{Wake structure and flow fields}
\label{supp:navier_stokes_cylinder_fields}

\begin{figure}[t]
\centering
\includegraphics[width=\textwidth]{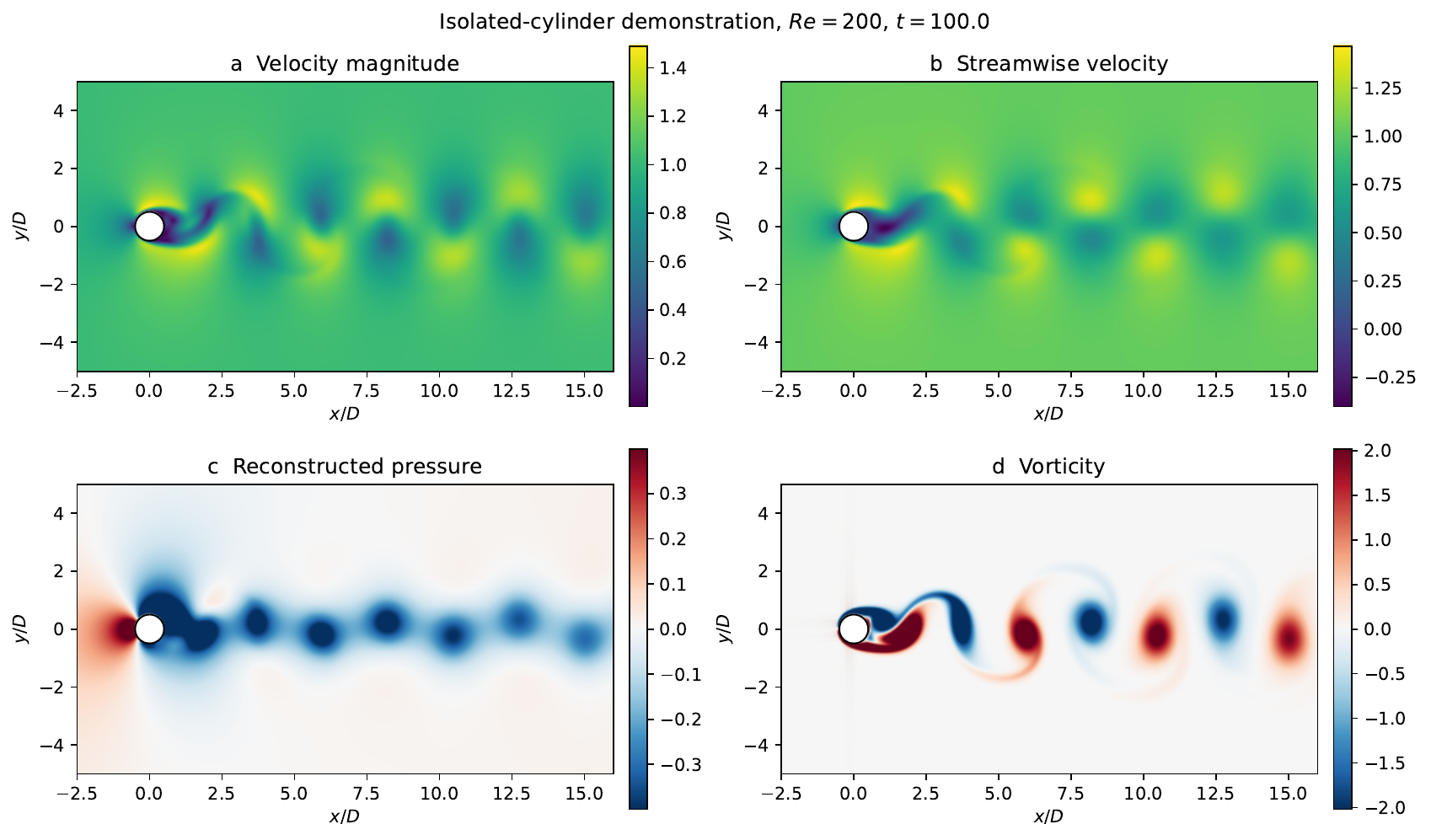}
\caption{\textbf{Instantaneous flow fields for the two-dimensional isolated-cylinder calculation at $Re=200$.} \textbf{a}, Velocity magnitude. \textbf{b}, Streamwise velocity. \textbf{c}, Reconstructed pressure. \textbf{d}, Vorticity. The calculation uses $N_D=40$, corresponding to a Fourier resolution of $1600\times800$, and the snapshot is taken at $t=100$.}
\label{fig:supp_cyl_fields}
\end{figure}

Figure~\ref{fig:supp_cyl_fields} shows the instantaneous $N_D=40$ solution at $t=100$. The velocity-magnitude field contains the downstream velocity deficit and its alternating transverse displacement. The streamwise-velocity distribution further resolves the periodically displaced wake behind the cylinder. The reconstructed pressure field exhibits the high-pressure upstream stagnation region and the alternating pressure structures associated with the unsteady wake.

The vorticity field provides the clearest representation of the nonlinear wake dynamics. Alternating positive and negative vortices are released from the two sides of the cylinder and convected downstream, forming a sustained K\'arm\'an-type vortex street. The wake remains well separated from the transverse periodic boundaries over the plotted region, while the downstream fringe treatment prevents the developed wake from being recirculated directly through the streamwise periodic boundary.

\subsubsection{Force, frequency and constraint diagnostics}
\label{supp:navier_stokes_cylinder_diagnostics}

\begin{figure}[t]
\centering
\includegraphics[width=\textwidth]{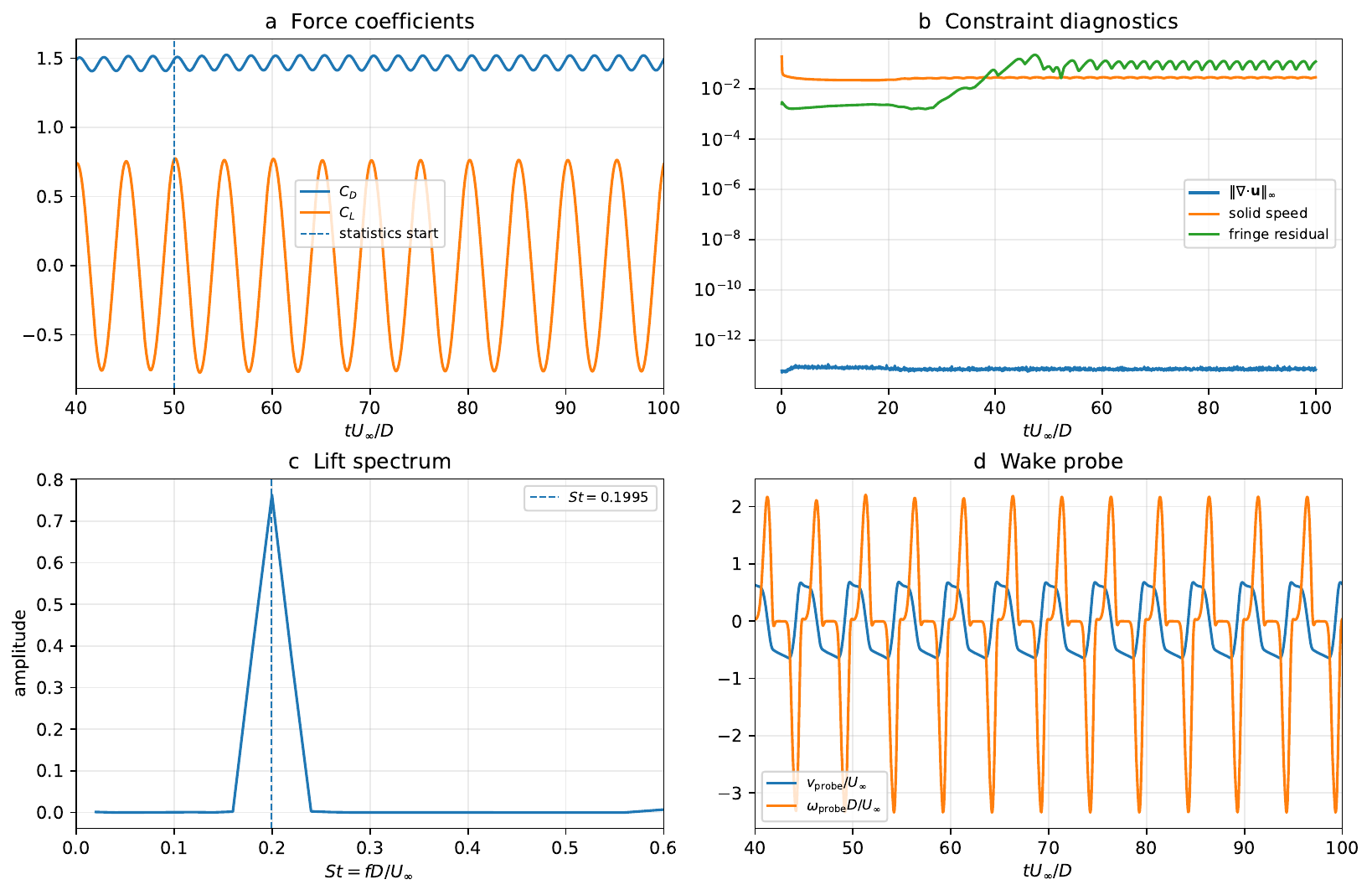}
\caption{\textbf{Temporal and spectral diagnostics for the two-dimensional isolated-cylinder calculation at $Re=200$.} \textbf{a}, Drag and lift coefficients. \textbf{b}, Incompressibility, solid-region and fringe diagnostics. \textbf{c}, Lift-coefficient spectrum. \textbf{d}, Downstream wake-probe signals. The calculation uses $N_D=40$, and statistical quantities are evaluated over $50\leq t\leq100$.}
\label{fig:supp_cyl_diagnostics}
\end{figure}

Figure~\ref{fig:supp_cyl_diagnostics} summarizes the temporal and spectral diagnostics for the $N_D=40$ calculation. After the initial wake-development stage, the drag and lift coefficients enter a sustained periodic state. Over the statistical interval $50\leq t\leq100$,
\begin{equation}
\overline{C_D}
=
1.46716,
\qquad
C_{L,\mathrm{rms}}
=
0.54049.
\label{eq:supp_cyl_final_force_statistics}
\end{equation}

The vortex-shedding frequency is first extracted from the lift-coefficient signal. The corresponding Strouhal number,
\begin{equation}
St
=
\frac{f_sD}{U_{\infty}},
\label{eq:supp_cyl_strouhal_definition}
\end{equation}
is
\begin{equation}
St_{C_L}
=
0.19955.
\label{eq:supp_cyl_strouhal_lift}
\end{equation}

An independent local measurement is obtained from the transverse velocity at the downstream probe
\begin{equation}
(x_p,y_p)
=
(x_c+3D,\ y_c+0.5D).
\label{eq:supp_cyl_probe_location}
\end{equation}
The dominant frequency of this signal gives
\begin{equation}
St_{\mathrm{probe}}
=
0.19957.
\label{eq:supp_cyl_strouhal_probe}
\end{equation}
The close agreement between the two independently determined values indicates that the dominant oscillation of the integrated lift force corresponds directly to the resolved periodic wake-shedding mode. Approximately ten shedding cycles are contained within the statistical interval.

The Fourier projection maintains incompressibility close to double-precision round-off throughout the calculation. For the highest-resolution case,
\begin{equation}
\max_t
\left\|
\nabla\cdot\mathbf{u}
\right\|_{\infty}
=
9.79\times10^{-14}.
\label{eq:supp_cyl_final_divergence}
\end{equation}
The maximum Courant number is
\begin{equation}
CFL_{\max}
=
0.771,
\label{eq:supp_cyl_final_cfl}
\end{equation}
and therefore remains below unity over the complete calculation. The solid-region velocity and fringe-relaxation diagnostics also remain bounded throughout the statistically stationary interval.

The combined field, force and spectral diagnostics show that the embedded Fourier formulation develops and sustains the expected alternating wake dynamics while retaining the divergence-free constraint to approximately machine precision. Together with the strict Taylor--Green verification in Section~\ref{supp:navier_stokes_taylor_green}, this benchmark therefore provides a complementary physical test of the nonlinear, embedded-boundary and long-time components of the incompressible-flow implementation.

\subsubsection{Reproducibility information}
\label{supp:navier_stokes_cylinder_reproducibility}

The benchmark-specific numerical settings are summarized in Table~\ref{tab:supp_cylinder_settings}. The common software environment, hardware specifications, arithmetic precision and timing protocol are reported in Section~\ref{supp:computational_environment}. The spatial-resolution calculations reported above use the same frozen numerical configuration, with only $N_D$ and the resulting Fourier resolution changed between cases.

The source code, input parameters and post-processing scripts required to reproduce this benchmark will be made publicly available upon publication at
\href{https://github.com/USERNAME/OptiXDE/tree/main/examples/navier_stokes_cylinder_re200}{\texttt{examples/navier\_stokes\_cylinder\_re200}}.

\begin{table}[h]
\centering
\caption{Numerical settings for the two-dimensional isolated-cylinder benchmark at $Re=200$.}
\label{tab:supp_cylinder_settings}
\begin{tabular}{ll}
\toprule
Setting & Value \\
\midrule
Reynolds number & $Re=200$ \\
Cylinder diameter & $D=1$ \\
Free-stream velocity & $U_{\infty}=1$ \\
Kinematic viscosity & $\nu=0.005$ \\
Fourier domain & $40D\times20D$ \\
Cylinder centre & $(10D,10D)$ \\
Boundary representation & Compact smooth Brinkman mask \\
Mask transition half-width & $\varepsilon_{\chi}=1.0h$ \\
Penalty parameter & $\eta=0.0025$ \\
Fringe width & $5D$ \\
Fringe relaxation time & $\tau_f=0.25$ \\
Flow formulation & Vorticity--streamfunction with uniform mean flow \\
Spatial transform & Two-dimensional discrete Fourier transform \\
Viscous propagation & Closed-form exponential spectral propagator \\
Nonlinear integration & Second-order Runge--Kutta \\
Nonlinear dealiasing & Two-dimensional $2/3$ rule \\
Constraint treatment & Exact Brinkman/fringe relaxation with Fourier projection \\
Pressure reconstruction & Periodic spectral pressure Poisson equation \\
Pressure normalization & Zero spatial mean \\
Force evaluation & Brinkman-relaxation momentum impulse \\
Time step & $\Delta t=0.01$ \\
Terminal time & $T=100$ \\
Statistical interval & $50\leq t\leq100$ \\
Diagnostic sampling interval & $0.05$ \\
Wake probe location & $(x_c+3D,\ y_c+0.5D)$ \\
Initial perturbation amplitude & $0.02$ \\
Spatial-resolution cases & $N_D=24,\ 32,\ 40$ \\
Fourier resolutions & $960\times480$, $1280\times640$, $1600\times800$ \\
Reference production resolution & $N_D=40$ \\
Arithmetic precision & Double precision \\
\bottomrule
\end{tabular}
\end{table}

\FloatBarrier

\subsection{Scalar ODE: exponential decay}
\label{subsec:ode}

To demonstrate that the OptiXDE framework also accommodates ordinary differential equations in addition to PDEs, we consider the scalar linear decay problem
\begin{equation}
\frac{\mathrm{d}u}{\mathrm{d}t}=-\lambda u,
\qquad
u(0)=1,
\end{equation}
with the analytical solution
\begin{equation}
u_{\mathrm{ex}}(t)=\exp(-\lambda t).
\end{equation}
This spatially independent problem bypasses the Fourier transform and directly exercises the exponential time-propagation mechanism used in OptiXDE. We take $\lambda=5$ and integrate the solution to $T=1$.

The numerical solution is evaluated using the terminal error $|u_h(T)-u_{\mathrm{ex}}(T)|$. In addition, the decay rate is inferred from the numerical solution as
\begin{equation}
\widehat{\lambda}_h
=
-\frac{1}{T}\log\left(\frac{u_h(T)}{u_h(0)}\right),
\qquad
\varepsilon_\lambda
=
|\widehat{\lambda}_h-\lambda|.
\end{equation}
For this homogeneous linear ODE, the exponential propagator advances the solution analytically over each time step. Consequently, the remaining discrepancies are dominated by floating-point round-off rather than time-discretization error. As shown in Table~\ref{tab:ode}, both error measures remain close to machine precision for all tested time-step sizes.

\begin{figure}[H]
\centering
\includegraphics[width=0.7\linewidth]{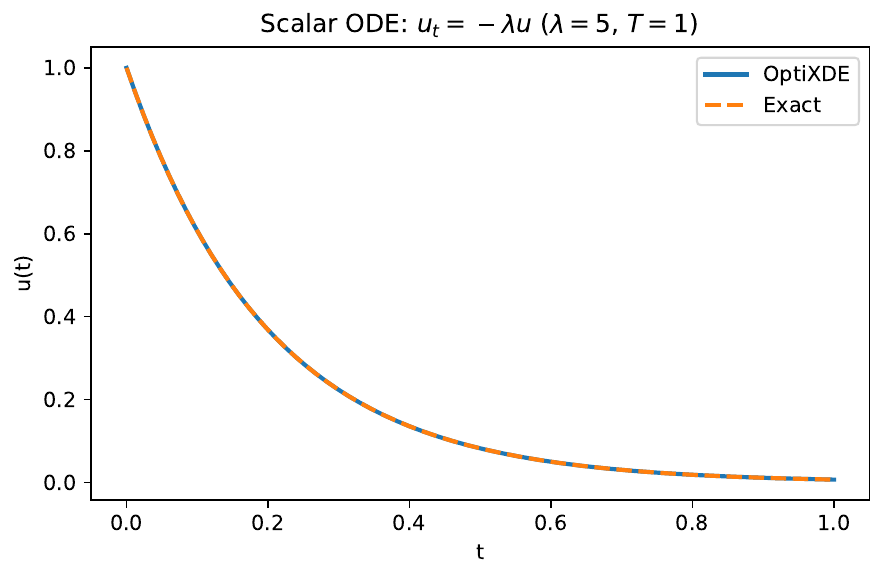}
\caption{Scalar ODE exponential decay computed by OptiXDE and compared with the analytical solution for $\lambda=5$.}
\label{fig:ode}
\end{figure}

Figure~\ref{fig:ode} compares the OptiXDE solution with the analytical exponential decay. The two curves overlap throughout the simulated interval, confirming the expected exponential evolution.

\begin{table}[htbp]
\centering
\caption{\textbf{Accuracy of the scalar ODE benchmark at $T=1$.}}
\label{tab:ode}
\begin{tabular}{ccc}
\toprule
$\Delta t$ &
$\left|u_h(T)-u_{\mathrm{ex}}(T)\right|$ &
$\varepsilon_{\lambda}$ \\
\midrule
$1.0\times10^{-2}$ & $6.939\times10^{-18}$ & $8.882\times10^{-16}$ \\
$5.0\times10^{-3}$ & $7.026\times10^{-17}$ & $1.066\times10^{-14}$ \\
$1.0\times10^{-3}$ & $4.250\times10^{-17}$ & $6.217\times10^{-15}$ \\
\bottomrule
\end{tabular}
\end{table}

\section{Accuracy and robustness studies}
\label{supp:accuracy_robustness}

This section examines the numerical accuracy, arithmetic sensitivity and robustness of \OptiXDE{} beyond the baseline benchmark configurations presented in Section~\ref{supp:additional_benchmarks}. The studies are designed to distinguish errors associated with spectral representation, repeated floating-point operations, embedded-boundary enforcement, mask regularization and nonlinear aliasing. Unless otherwise stated, the spatial discretization, propagation operators and error definitions follow Sections~\ref{supp:discrete_transforms}, \ref{supp:spectral_propagators} and \ref{supp:error_measures}, respectively. All calculations are performed using the same software environment and transform backends described in Section~\ref{supp:computational_environment}.

\subsection{Spectral propagation and finite-precision robustness}
\label{supp:spectral_precision_robustness}

The transient-diffusion benchmark in Section~\ref{supp:diffusion} already establishes the accuracy of the closed-form propagation operator with respect to spatial resolution and repeated time stepping. The present study therefore focuses only on three complementary robustness questions: whether the same accuracy is retained across a broad portion of the resolved spectrum, how finite-precision errors accumulate under repeated propagation, and how the attainable accuracy changes between single- and double-precision arithmetic. All calculations use the same homogeneous-Dirichlet formulation and propagation kernel as the reference benchmark.

To characterize the position of an eigenmode $\boldsymbol{\kappa}=(k_x,k_y)$ within the two-dimensional sine spectrum, we define the normalized modal radius
\begin{equation}
\rho_{\boldsymbol{\kappa}}
=
\frac{\sqrt{k_x^2+k_y^2}}{\sqrt{2}\,N},
\label{eq:supp_robustness_modal_radius}
\end{equation}
where $N=256$ is the number of interior samples in each spatial direction and $\rho_{\boldsymbol{\kappa}}=1$ corresponds to the spectral corner. Modes from $(1,1)$ to $(224,224)$ are examined, spanning $0.004\leq\rho_{\boldsymbol{\kappa}}\leq0.875$. To prevent the stronger attenuation of high-wave-number modes from biasing the comparison, each mode is propagated to the same analytical attenuation factor $\exp(-2)$, corresponding to the mode-dependent terminal time $T_{\boldsymbol{\kappa}}=2/\lambda_{\boldsymbol{\kappa}}$. Relative field and decay-rate errors are evaluated using the common discrete norms and modal diagnostics introduced previously.

\begin{figure}[h]
\centering
\includegraphics[width=\textwidth]{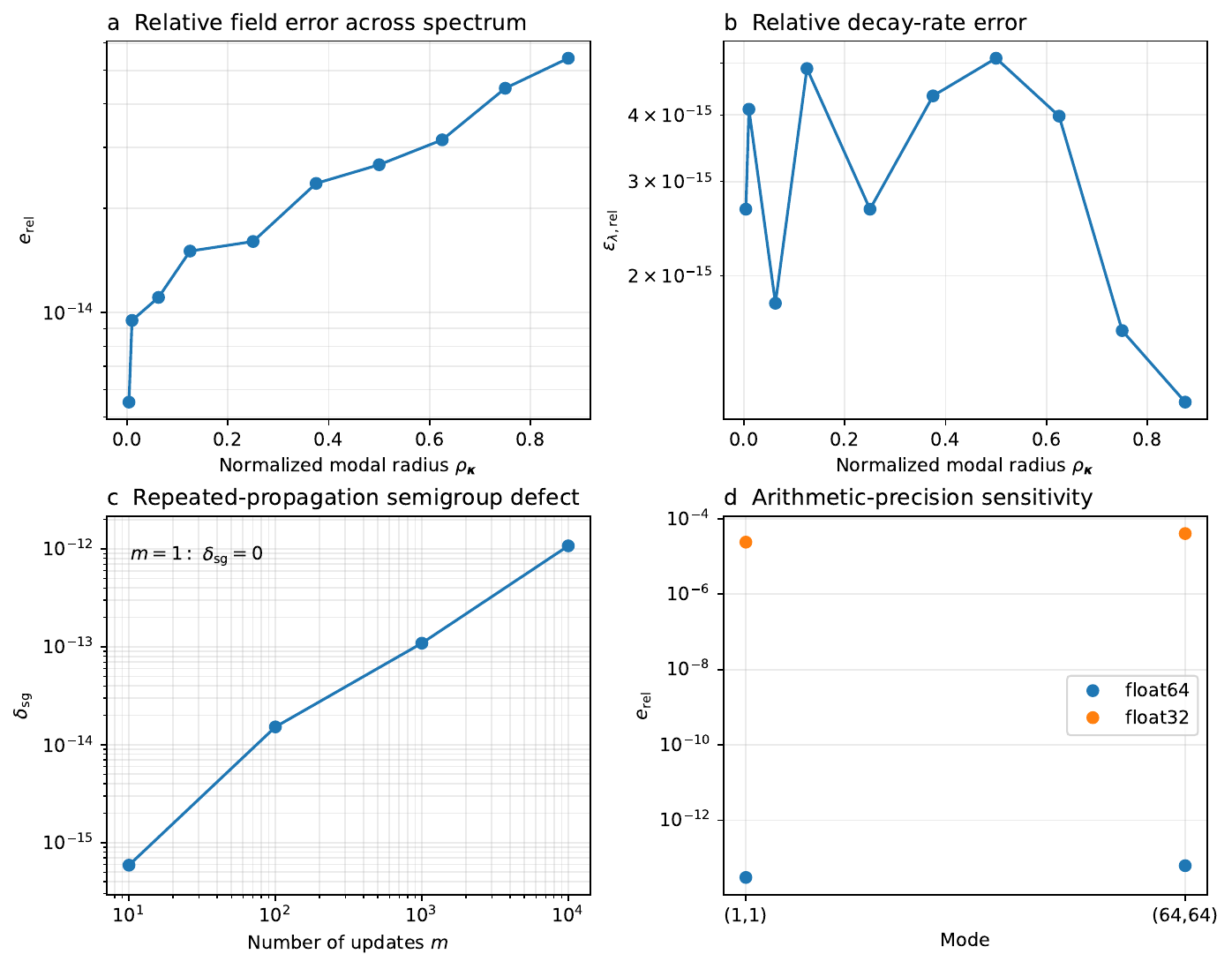}
\caption{\textbf{Spectral propagation and finite-precision robustness of the transient-diffusion benchmark.}
\textbf{a}, Relative field error as a function of normalized modal radius for resolved sine eigenmodes from $(1,1)$ to $(224,224)$, with each mode propagated to the same analytical attenuation factor $\exp(-2)$.
\textbf{b}, Corresponding relative decay-rate error.
\textbf{c}, Semigroup defect for repeated propagation of the $(1,1)$ mode at $T=1.5$; the single-step value is zero and is indicated separately because it cannot be displayed on the logarithmic ordinate.
\textbf{d}, Relative field errors obtained using single- and double-precision arithmetic for the representative modes $(1,1)$ and $(64,64)$.}
\label{fig:supp_diffusion_robustness}
\end{figure}

As shown in Fig.~\ref{fig:supp_diffusion_robustness}a,b, the relative field error remains between $5.52\times10^{-15}$ and $5.42\times10^{-14}$ over the complete tested spectral range. A moderate increase is observed toward the highest retained modes, but the error remains within the double-precision round-off regime. The relative decay-rate error remains between approximately $1.16\times10^{-15}$ and $5.11\times10^{-15}$ and exhibits no systematic deterioration with modal frequency. The analytical propagation multiplier therefore remains accurate over a broad fraction of the resolvable sine spectrum, including modes approaching the spectral cutoff.

Repeated-propagation robustness is examined using the fundamental mode $(1,1)$ at the fixed terminal time $T=1.5$. The same terminal state is computed either by one propagation over $T$ or by $m$ successive updates of duration $T/m$. The resulting finite-precision semigroup defect is measured as
\begin{equation}
\delta_{\mathrm{sg}}
=
\frac{
\left\|
u_h^{\mathrm{repeat}}(T)
-
u_h^{\mathrm{single}}(T)
\right\|_{2,h}
}{
\left\|
u_h^{\mathrm{single}}(T)
\right\|_{2,h}
}.
\label{eq:supp_robustness_semigroup_defect}
\end{equation}

Figure~\ref{fig:supp_diffusion_robustness}c shows that $\delta_{\mathrm{sg}}$ is zero for the single-step reference, $5.88\times10^{-16}$ for $m=10$, $1.52\times10^{-14}$ for $m=100$, $1.09\times10^{-13}$ for $m=1000$, and $1.08\times10^{-12}$ after $10^4$ repeated updates. This gradual increase is consistent with the finite-precision accumulation identified in the time-step study of Section~\ref{supp:diffusion_timestep}, rather than with a conventional temporal truncation error. Even after $10^4$ transform--propagate--inverse-transform cycles, the discrepancy remains approximately $10^{-12}$ in double precision.

Arithmetic sensitivity is evaluated using the representative modes $(1,1)$ and $(64,64)$ with otherwise identical numerical settings. The solver preserves the requested floating-point type, so the comparison represents native single- and double-precision calculations. As shown in Fig.~\ref{fig:supp_diffusion_robustness}d, the relative field errors for the $(1,1)$ mode are $3.06\times10^{-14}$ and $2.40\times10^{-5}$ in double and single precision, respectively. The corresponding values for the $(64,64)$ mode are $6.27\times10^{-14}$ and $4.05\times10^{-5}$. The direct single-to-double-precision discrepancies are $2.41\times10^{-5}$ and $4.06\times10^{-5}$ for the two modes.

These results identify two distinct arithmetic regimes. Double precision preserves the near-round-off accuracy of the analytical spectral propagator, whereas single precision introduces an error floor of order $10^{-5}$ for the present repeated-transform workflow. Single precision can therefore remain appropriate when the required physical accuracy is substantially above this level, while double precision is required for verification studies intended to resolve near-machine-precision behavior.

Together, the three tests show that the closed-form diffusion propagator remains accurate across the resolved spectrum, accumulates only small finite-precision discrepancies under repeated application, and exhibits a predictable accuracy floor determined by the selected arithmetic precision.

\subsection{Embedded-boundary robustness for the L-shaped Poisson problem}
\label{supp:accuracy_lshape_robustness}

The robustness of the embedded-boundary formulation was examined using the L-shaped Poisson benchmark at a fixed spatial resolution of \(512^2\). The study isolates three parameters associated with the physical-space boundary treatment and pseudo-time iteration: the penalty coefficient, the diffuse-interface width and the pseudo-time increment. Unless otherwise stated, the production values $\eta/h^2=0.003$, $\varepsilon=0.03$ and $\Delta\tau/h^2=10$ from Section~\ref{supp:poisson_lshape_penalized_solution} were used. Each parameter case was initialized independently and iterated using the same stopping criterion as the principal benchmark, Eq.~\eqref{eq:supp_poisson_lshape_stopping_tolerance}. The global physical-domain error $e_\Omega$ and the regional errors in the re-entrant-corner region, diffuse-interface region and remaining bulk region were evaluated using the same definitions as those employed for the principal L-shaped benchmark.

The penalty sensitivity was first examined using
\begin{equation}
\frac{\eta}{h^2}
\in
\left\{
0.0015,\,
0.003,\,
0.006
\right\}.
\label{eq:supp_lshape_penalty_values}
\end{equation}
The corresponding global errors were
\begin{equation}
e_\Omega
=
1.604\times10^{-5},
\quad
3.137\times10^{-5},
\quad
6.161\times10^{-5},
\label{eq:supp_lshape_penalty_errors}
\end{equation}
respectively. The interface error exhibited the same systematic trend, increasing from \(2.223\times10^{-5}\) to \(8.888\times10^{-5}\) over the tested interval. The approximately monotonic reduction in error with decreasing \(\eta\) is consistent with the finite-penalty contribution of the embedded-boundary formulation. All three cases nevertheless converged to the prescribed pseudo-time tolerance, indicating that the production value \(\eta/h^2=0.003\) lies within a stable operating range rather than at a narrowly tuned isolated point.

The influence of the regularized-mask width was assessed using
\begin{equation}
\varepsilon
\in
\left\{
0.020,\,
0.025,\,
0.030,\,
0.035,\,
0.040
\right\},
\label{eq:supp_lshape_mask_values}
\end{equation}
while retaining \(\eta/h^2=0.003\) and \(\Delta\tau/h^2=10\). The global errors were
\begin{equation}
e_\Omega
=
4.348\times10^{-4},
\;
1.052\times10^{-4},
\;
3.137\times10^{-5},
\;
1.647\times10^{-5},
\;
1.317\times10^{-5},
\label{eq:supp_lshape_mask_errors}
\end{equation}
for increasing \(\varepsilon\). The pronounced deterioration for the narrowest interface is accompanied by a bulk error of \(4.584\times10^{-4}\), substantially larger than the corresponding corner error. This behaviour indicates that insufficient smoothing does not merely increase the local interface discrepancy but also introduces high-frequency contamination into the global spectral representation.

Increasing the smoothing width beyond the production value continues to reduce the discrete RMS error, although the improvement becomes progressively weaker between \(\varepsilon=0.035\) and \(0.040\). The smallest global error among the tested cases therefore occurs at \(\varepsilon=0.040\); however, this does not imply that the widest diffuse interface is universally preferable. Increasing \(\varepsilon\) also broadens the regularized representation of the physical boundary. The production value \(\varepsilon=0.03\) is consequently retained as a compromise between suppressing spectral contamination and limiting geometric-interface broadening, rather than as the point that minimizes the RMS error alone.

Pseudo-time sensitivity was examined using
\begin{equation}
\frac{\Delta\tau}{h^2}
\in
\left\{
5,\,
10,\,
20
\right\}.
\label{eq:supp_lshape_pseudotime_values}
\end{equation}
The number of correction steps required to satisfy Eq.~\eqref{eq:supp_poisson_lshape_stopping_tolerance} decreased from \(4600\) to \(2600\) and then to \(1400\) as \(\Delta\tau/h^2\) increased. The corresponding global errors were
\begin{equation}
e_\Omega
=
2.196\times10^{-5},
\quad
3.137\times10^{-5},
\quad
1.638\times10^{-4}.
\label{eq:supp_lshape_pseudotime_errors}
\end{equation}
Thus, increasing \(\Delta\tau/h^2\) accelerates the fixed-point correction, but excessively large pseudo-time increments degrade the converged spatial solution. In particular, increasing \(\Delta\tau/h^2\) from \(5\) to \(10\) reduces the required correction count by approximately \(43\%\) while preserving a global error of order \(10^{-5}\). Increasing it further to \(20\) reduces the correction count to \(1400\), but raises the global error by more than a factor of five relative to the production case. The selected value
\begin{equation}
\frac{\Delta\tau}{h^2}=10
\label{eq:supp_lshape_pseudotime_selected}
\end{equation}
therefore provides a practical accuracy--efficiency compromise for the present embedded-domain problem.

\begin{figure}[h]
\centering
\suppfigure{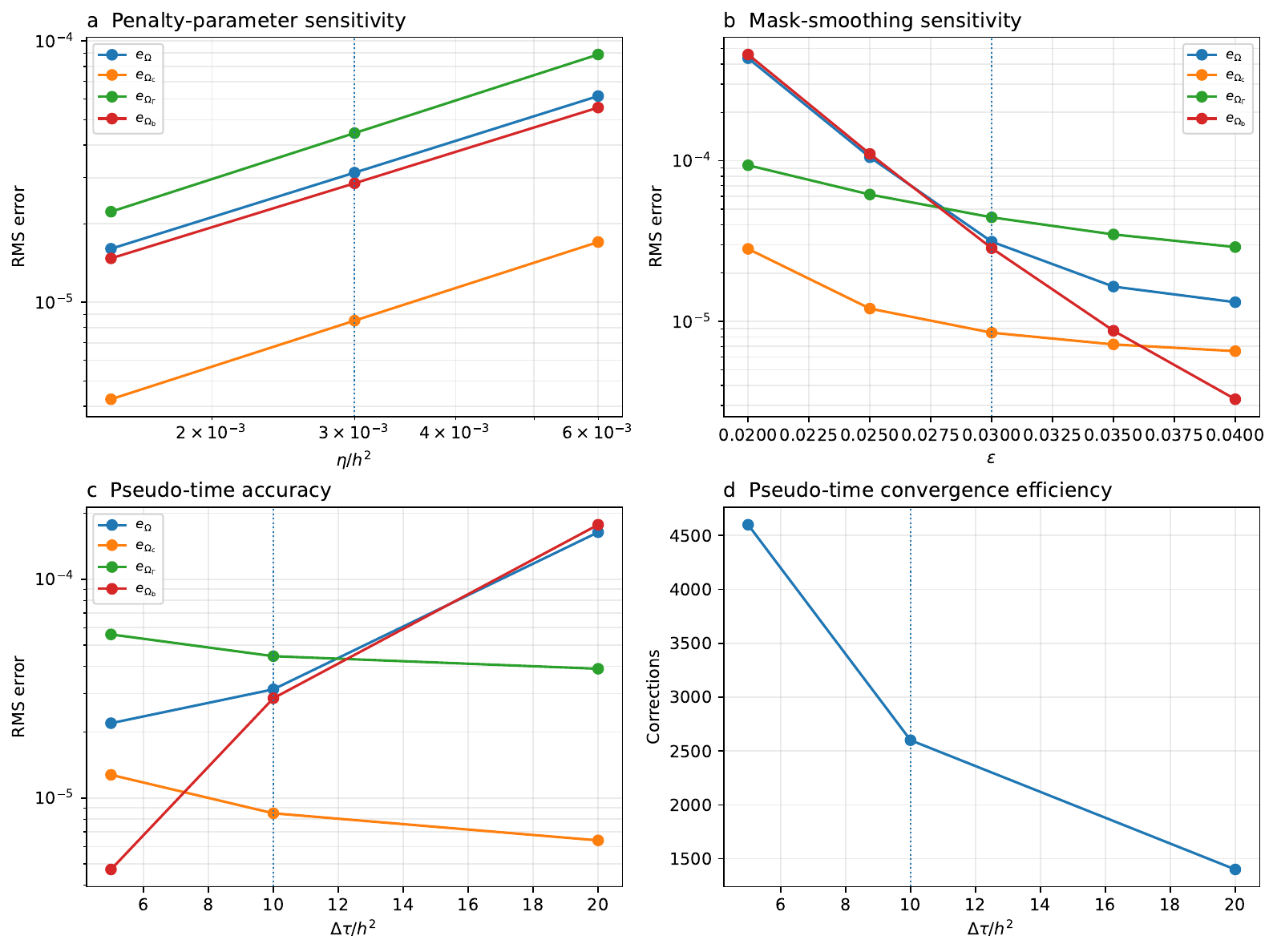}{\textwidth}
\caption{\textbf{Embedded-boundary robustness of the L-shaped Poisson benchmark at \(512^2\) resolution.}
\textbf{a}, Sensitivity of the global and regional errors to the scaled penalty coefficient \(\eta/h^2\), with the production value \(\eta/h^2=0.003\) indicated.
\textbf{b}, Sensitivity to the diffuse-interface width \(\varepsilon\). Narrow interfaces produce substantially larger global and bulk errors, whereas further smoothing progressively reduces the discrete error at the cost of a broader geometric transition; the production value is \(\varepsilon=0.03\).
\textbf{c}, Global and regional errors obtained for the tested pseudo-time increments \(\Delta\tau/h^2\).
\textbf{d}, Number of correction steps required to reach \(r_\tau<10^{-8}\), demonstrating the accuracy--efficiency trade-off associated with the pseudo-time increment. The production value \(\Delta\tau/h^2=10\) provides a compromise between solution accuracy and convergence speed.}
\label{fig:supp_lshape_robustness}
\end{figure}

Figure~\ref{fig:supp_lshape_robustness} summarizes these three controlled studies. Together with the localized-error analysis in the principal benchmark, the results distinguish the principal error mechanisms of the embedded formulation from floating-point transform errors. In contrast to the exactly representable periodic benchmarks, the L-shaped problem is governed by finite-penalty effects, diffuse-interface regularization and the re-entrant-corner singularity. The robustness study therefore does not seek machine-precision invariance with respect to the embedding parameters; instead, it establishes a reproducible parameter regime in which the solution remains converged and the accuracy--regularization--cost trade-offs can be identified explicitly.

\subsection{Nonlinear Schr\"odinger equation: invariant and spectral robustness}
\label{supp:robustness_schrodinger}

The nonlinear Schr\"odinger benchmark in Section~\ref{supp:schrodinger} establishes spatial resolution, second-order temporal convergence and accurate reproduction of the nonlinear focusing and recurrence cycle. Here, robustness is examined through two complementary properties: preservation of the continuous invariants during repeated split-step propagation and sensitivity of the resolved solution to explicit high-wave-number filtering.

The focusing nonlinear Schr\"odinger equation conserves the mass
\begin{equation}
M(t)
=
\int_{\Omega}
|h(x,t)|^2
\,\dd x
\label{eq:supp_schrodinger_mass}
\end{equation}
and the Hamiltonian energy
\begin{equation}
E(t)
=
\int_{\Omega}
\left[
\frac{1}{2}
\left|
\frac{\partial h}{\partial x}
\right|^2
-
\frac{1}{2}
|h|^4
\right]
\,\dd x.
\label{eq:supp_schrodinger_energy}
\end{equation}

For the continuous initial profile restricted to $[-15,15]$,
\begin{equation}
M(0)
=
8\tanh(15),
\label{eq:supp_schrodinger_initial_mass}
\end{equation}
and
\begin{equation}
E(0)
=
-16\tanh(15)
+
\frac{20}{3}
\tanh^3(15).
\label{eq:supp_schrodinger_initial_energy}
\end{equation}

The numerical conservation diagnostics are evaluated relative to the corresponding discrete initial values. The discrete mass is
\begin{equation}
M_h^n
=
\Delta x
\sum_{j=0}^{N_x-1}
|h_j^n|^2.
\label{eq:supp_schrodinger_discrete_mass}
\end{equation}

The derivative $D_xh^n$ is evaluated using the common spectral-differentiation rule in Eq.~\eqref{eq:supp_spectral_derivatives}, and the discrete energy is
\begin{equation}
E_h^n
=
\Delta x
\sum_{j=0}^{N_x-1}
\left[
\frac{1}{2}
\left|
\left(
D_xh^n
\right)_j
\right|^2
-
\frac{1}{2}
|h_j^n|^4
\right].
\label{eq:supp_schrodinger_discrete_energy}
\end{equation}

For $Q_h=M_h$ or $Q_h=E_h$, the maximum recorded relative drift $\varepsilon_Q$ is evaluated using the common invariant-drift definition in Eq.~\eqref{eq:supp_invariant_drift}. The denominator safeguard in that definition is inactive here because both initial invariant values are nonzero.

The five fixed-step calculations in Section~\ref{supp:schrodinger_convergence} provide a repeated-propagation test over
\begin{equation}
N_t
\in
\left\{
256,\,
512,\,
1024,\,
2048,\,
4096
\right\}.
\label{eq:supp_schrodinger_robustness_step_counts}
\end{equation}
Because the dispersive Fourier multiplier and the exact nonlinear phase rotation both have unit modulus, the discrete mass is preserved close to floating-point round-off. The maximum relative mass drift remains between
\begin{equation}
5.329071\times10^{-14}
\leq
\varepsilon_M
\leq
8.493206\times10^{-13}
\label{eq:supp_schrodinger_mass_drift_range}
\end{equation}
over the complete refinement sequence. The small increase in the accumulated mass drift as $\Delta t$ decreases is associated with the larger number of Fourier transforms and pointwise operations rather than a systematic loss of the invariant.

The Hamiltonian energy is not preserved exactly by the separated subflows, but its bounded excursion decreases systematically with temporal refinement:
\begin{equation}
1.676698\times10^{-3}
\rightarrow
4.208425\times10^{-4}
\rightarrow
1.053153\times10^{-4}
\rightarrow
2.633536\times10^{-5}
\rightarrow
6.584251\times10^{-6}.
\label{eq:supp_schrodinger_energy_drift_sequence}
\end{equation}
Each halving of the time increment reduces the maximum energy excursion by approximately a factor of four. The invariant behaviour therefore provides an independent physical-consistency check of the second-order symmetric splitting beyond the reference-field convergence reported in Section~\ref{supp:schrodinger_convergence}.

\begin{figure}[t]
\centering
\suppfigure{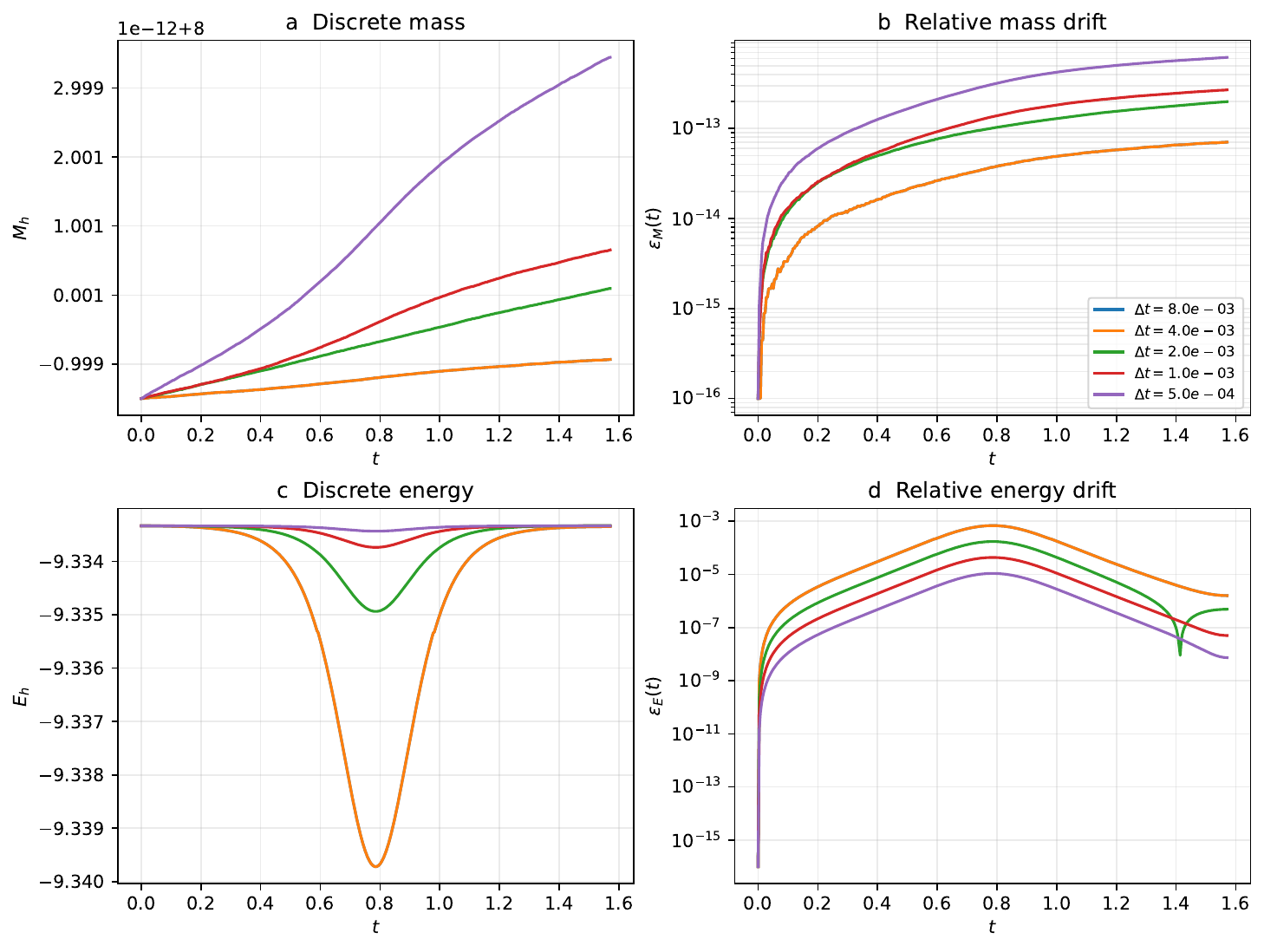}{\textwidth}
\caption{\textbf{Invariant preservation of the nonlinear Schr\"odinger benchmark under repeated split-step propagation.} \textbf{a}, Discrete mass $M_h(t)$ for the five tested fixed time increments. \textbf{b}, Relative mass drift $\varepsilon_M(t)$. \textbf{c}, Discrete Hamiltonian energy $E_h(t)$. \textbf{d}, Relative energy drift $\varepsilon_E(t)$. Mass remains close to floating-point accuracy throughout the complete propagation interval, whereas the bounded energy excursion decreases approximately quadratically under temporal refinement.}
\label{fig:supp_schrodinger_conservation}
\end{figure}

A second robustness test examines whether the well-resolved principal solution depends materially on unresolved high-wave-number content. The finest temporal calculation defined in Eq.~\eqref{eq:supp_schrodinger_snapshot_case} is repeated with an explicit spectral filter applied after every complete Strang step. In the filtered calculation, Fourier coefficients satisfying
\begin{equation}
|\mu_m|>\frac{N_x}{3}
\label{eq:supp_schrodinger_filter_cutoff}
\end{equation}
are set to zero, while the time step, spatial resolution, initial condition and all other solver settings remain unchanged.

Let $h_{\mathrm{f}}$ and $h_{\mathrm{u}}$ denote the filtered and unfiltered solutions, respectively. Their relative difference is measured by
\begin{equation}
\delta_{\mathrm{f}}(t_q)
=
\frac{
\left\|
h_{\mathrm{f}}(\cdot,t_q)
-
h_{\mathrm{u}}(\cdot,t_q)
\right\|_2
}{
\left\|
h_{\mathrm{u}}(\cdot,t_q)
\right\|_2
}.
\label{eq:supp_schrodinger_filter_difference}
\end{equation}

The maximum difference over the complete focusing--recurrence trajectory is
\begin{equation}
\max_q\delta_{\mathrm{f}}(t_q)
=
2.455593\times10^{-10},
\label{eq:supp_schrodinger_filter_max_difference}
\end{equation}
and the terminal difference is
\begin{equation}
\delta_{\mathrm{f}}(T)
=
2.453568\times10^{-10}.
\label{eq:supp_schrodinger_filter_terminal_difference}
\end{equation}
The change in the focusing-amplitude error is only
\begin{equation}
\Delta e_{\mathrm{peak}}
=
2.138223\times10^{-10},
\label{eq:supp_schrodinger_filter_peak_difference}
\end{equation}
while the difference between the maximum energy drifts of the filtered and unfiltered calculations is
\begin{equation}
\Delta\varepsilon_E
=
5.978075\times10^{-13}.
\label{eq:supp_schrodinger_filter_energy_difference}
\end{equation}

These differences are several orders of magnitude smaller than the numerical error of the principal calculation. In particular, the unfiltered finest temporal calculation has
\begin{equation}
e_{\max}
=
4.319465\times10^{-5},
\label{eq:supp_schrodinger_robustness_baseline_error}
\end{equation}
so that the filtered--unfiltered difference is negligible relative to the measured discretization error. This observation is consistent with the spatial-resolution study, in which the resolved spectral-tail fraction at $N_x=2048$ falls to approximately $6.0\times10^{-20}$.

The filtering experiment therefore shows that the principal nonlinear Schr\"odinger result is not sustained by unresolved spectral content or by an accidental accumulation of high-frequency modes. Once the focusing spectrum is adequately resolved, explicitly removing the upper spectral band produces no material change in the field evolution, focusing amplitude or invariant behaviour. Together with near-round-off mass preservation and the systematic reduction of the Hamiltonian-energy excursion, this provides a complementary robustness check for repeated nonlinear split-step propagation.

\subsection{Burgers equation: nonlinear spectral robustness}
\label{supp:robustness_burgers}

The Burgers benchmark in Section~\ref{supp:burgers} establishes spatial convergence, second-order temporal convergence and localized accuracy through the formation of the narrow viscous layer. The present study examines a complementary question: whether the nonlinear pseudo-spectral treatment remains physically consistent and how the solution depends on the treatment of quadratic aliasing. The analysis therefore focuses on de-aliasing sensitivity, high-wave-number behaviour, kinetic-energy monotonicity, boundary preservation and odd-symmetry preservation, without repeating the spatial- and temporal-refinement studies reported previously.

Three treatments of the nonlinear conservative flux are compared under otherwise identical conditions:
\begin{equation}
N_x=256,
\qquad
\Delta t=5.0\times10^{-4},
\qquad
T=1,
\label{eq:supp_burgers_robustness_settings}
\end{equation}
namely no de-aliasing, classical two-thirds truncation and three-halves zero padding. The same Cole--Hopf reference solution and the same fixed-step Strang/RK4 integration sequence are used in all cases. The principal Burgers calculations reported in Section~\ref{supp:burgers} use three-halves padding.

The maximum field error, terminal error and maximum viscous-layer error obtained without de-aliasing are
\begin{equation}
e_{\max}=2.065024\times10^{-2},
\qquad
e(T)=1.119278\times10^{-3},
\qquad
\max_t e_{\mathrm{s}}^{\infty}
=
2.065024\times10^{-2}.
\label{eq:supp_burgers_none_results}
\end{equation}
Two-thirds truncation suppresses the spectral content above the prescribed cutoff almost completely,
\begin{equation}
\rho_{\mathrm{high}}^{\max}
=
5.783084\times10^{-31},
\label{eq:supp_burgers_twothirds_highmode}
\end{equation}
but produces the largest field errors among the three treatments,
\begin{equation}
e_{\max}=3.840180\times10^{-2},
\qquad
e(T)=2.538407\times10^{-3}.
\label{eq:supp_burgers_twothirds_results}
\end{equation}
The near-zero value of $\rho_{\mathrm{high}}^{\max}$ is therefore not an indicator of superior accuracy; it results directly from removal of the modes above the two-thirds cutoff and also removes resolved high-wave-number content associated with the narrow viscous layer.

Three-halves padding gives the smallest errors,
\begin{equation}
e_{\max}=1.157029\times10^{-2},
\qquad
e(T)=5.068437\times10^{-4},
\qquad
\max_t e_{\mathrm{s}}^{\infty}
=
1.157029\times10^{-2},
\label{eq:supp_burgers_threehalves_results}
\end{equation}
while retaining a finite high-wave-number occupation,
\begin{equation}
\rho_{\mathrm{high}}^{\max}
=
1.416450\times10^{-4}.
\label{eq:supp_burgers_threehalves_highmode}
\end{equation}
Relative to the calculation without de-aliasing, three-halves padding reduces $e_{\max}$ by approximately $44.0\%$ and $e(T)$ by approximately $54.7\%$. Relative to two-thirds truncation, the corresponding reductions are approximately $69.9\%$ and $80.0\%$. These results show that suppressing high-wave-number content alone is not equivalent to suppressing aliasing error. For this benchmark, three-halves padding removes the unresolved quadratic convolution while retaining the resolved high-frequency content required to represent the viscous layer.

\begin{figure}[t]
\centering
\suppfigure{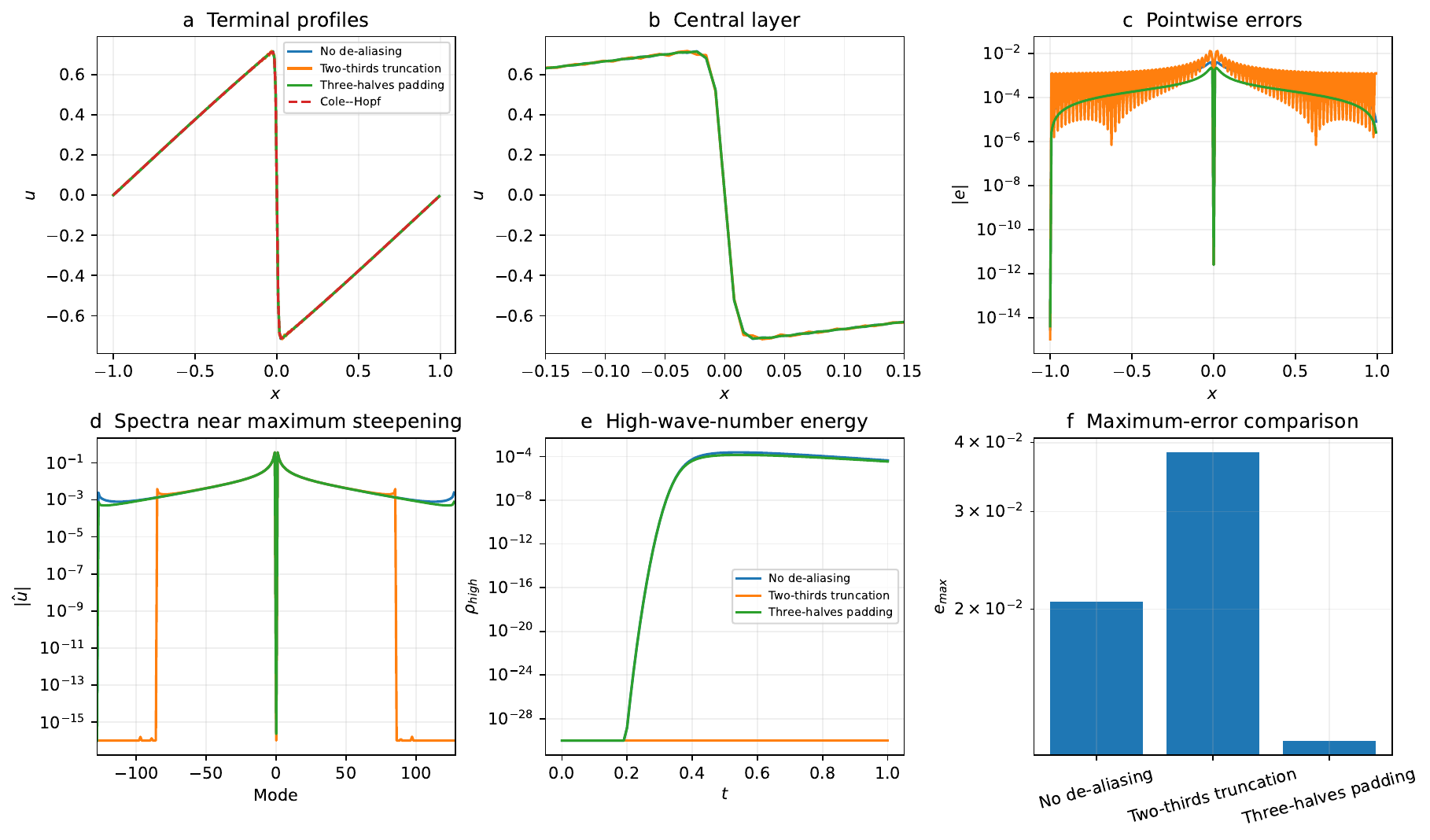}{\textwidth}
\caption{\textbf{Robustness of the Burgers solution with respect to nonlinear de-aliasing.} \textbf{a}, Terminal profiles obtained without de-aliasing, with two-thirds truncation and with three-halves zero padding, together with the Cole--Hopf reference. \textbf{b}, Enlarged comparison across the central viscous layer. \textbf{c}, Terminal pointwise absolute errors. \textbf{d}, Fourier-amplitude spectra near the maximum-steepening stage. \textbf{e}, High-wave-number ratio $\rho_{\mathrm{high}}(t)$ for the three nonlinear-product treatments. \textbf{f}, Maximum space--time error $e_{\max}$. Three-halves padding gives the smallest field error while preserving resolved high-wave-number content, whereas two-thirds truncation imposes a sharp spectral cutoff and produces the largest error.}
\label{fig:supp_burgers_dealiasing}
\end{figure}

Figure~\ref{fig:supp_burgers_dealiasing} localizes the effect of the nonlinear-product treatment. The full-domain terminal profiles remain visually similar, whereas the differences become apparent in the central viscous layer and in the pointwise-error distribution. The two-thirds spectrum exhibits the imposed sharp cutoff, while the calculation without de-aliasing retains a larger high-wave-number population that contains both physically resolved content and aliased contributions. Three-halves padding provides the best balance: the nonlinear product is evaluated on an enlarged spectral representation before restriction to the original Fourier band, avoiding direct quadratic aliasing without indiscriminately removing the resolved spectral tail.

Additional physical-consistency diagnostics are evaluated using the $N_x=1024$ temporal-refinement sequence of Section~\ref{supp:burgers_convergence}. The discrete kinetic energy is
\begin{equation}
E_h^n
=
\frac{\Delta x}{2}
\sum_{j=0}^{N_x-1}
\left(u_j^n\right)^2,
\label{eq:supp_burgers_robustness_energy}
\end{equation}
and the maximum normalized positive energy increment $\delta_E^{\uparrow}$ is evaluated using the common monotonicity-violation definition in Eq.~\eqref{eq:supp_monotonicity_violation} with $Q_h=E_h$.
No positive energy increment is detected for any of the five tested time increments,
\begin{equation}
\delta_E^{\uparrow}=0,
\qquad
6.25\times10^{-5}
\leq
\Delta t
\leq
1.0\times10^{-3},
\label{eq:supp_burgers_zero_energy_violation}
\end{equation}
consistent with the dissipative energy law in Eq.~\eqref{eq:supp_burgers_energy_dissipation}.

Boundary preservation is quantified by
\begin{equation}
r_{\mathrm{b}}
=
\max_{t\in\mathcal{T}_{\mathrm{out}}}
|u_h(-1,t)|,
\label{eq:supp_burgers_boundary_robustness}
\end{equation}
where $x=1$ is the periodic counterpart of the stored sample at $x=-1$. Across the complete temporal-refinement sequence,
\begin{equation}
1.29\times10^{-14}
\leq
r_{\mathrm{b}}
\leq
2.11\times10^{-13}.
\label{eq:supp_burgers_boundary_robustness_results}
\end{equation}
The homogeneous boundary value is therefore preserved to approximately round-off level over as many as $16000$ fixed time steps.

Because the exact solution remains odd, symmetry preservation provides an additional independent diagnostic. We define
\begin{equation}
r_{\mathrm{sym}}
=
\frac{
\displaystyle
\max_{t,x}
|u_h(x,t)+u_h(-x,t)|
}{
\displaystyle
\max_{t,x}|u_h(x,t)|
}.
\label{eq:supp_burgers_symmetry_robustness}
\end{equation}
The measured defect remains between
\begin{equation}
5.31\times10^{-13}
\leq
r_{\mathrm{sym}}
\leq
8.86\times10^{-12}
\label{eq:supp_burgers_symmetry_robustness_results}
\end{equation}
over the five temporal refinements. The slight accumulation with increasing step count remains at a negligible level and is consistent with floating-point accumulation rather than a systematic loss of symmetry.

Taken together, these tests show that the nonlinear Burgers implementation is robust with respect to both spectral treatment and repeated propagation. Three-halves padding gives the lowest error among the tested nonlinear-product treatments while retaining dynamically relevant resolved high-wave-number content. At the same time, the complete split-step evolution preserves monotone kinetic-energy decay, the homogeneous endpoint value and odd symmetry to near-machine precision. These results support the use of three-halves padding as the default nonlinear-product treatment for the principal Burgers calculations.

\subsection{Large-time-step robustness of the Allen--Cahn equation}
\label{supp:allen_cahn_robustness}

The robustness of the nonlinear splitting strategy was further examined by deliberately increasing the time step beyond the range used for the temporal-convergence assessment. The spatial resolution, interface parameter and final time were fixed at the principal benchmark values $N_x=N_y=256$, $\epsilon=0.04$ and $T=100$ from Eqs.~\eqref{eq:supp_allen_cahn_primary_discretization} and~\eqref{eq:supp_allen_cahn_parameters}, while the time step was varied over
\begin{equation}
\Delta t
\in
\left\{
0.05,\,
0.1,\,
0.2,\,
0.5,\,
1,\,
2,\,
5,\,
10
\right\}.
\end{equation}
The solution obtained with \(\Delta t=0.05\), which corresponds to the principal discretization used in the benchmark, was taken as the baseline for the field and physical-diagnostic comparisons.

\begin{figure}[h]
\centering
\includegraphics[width=\textwidth]{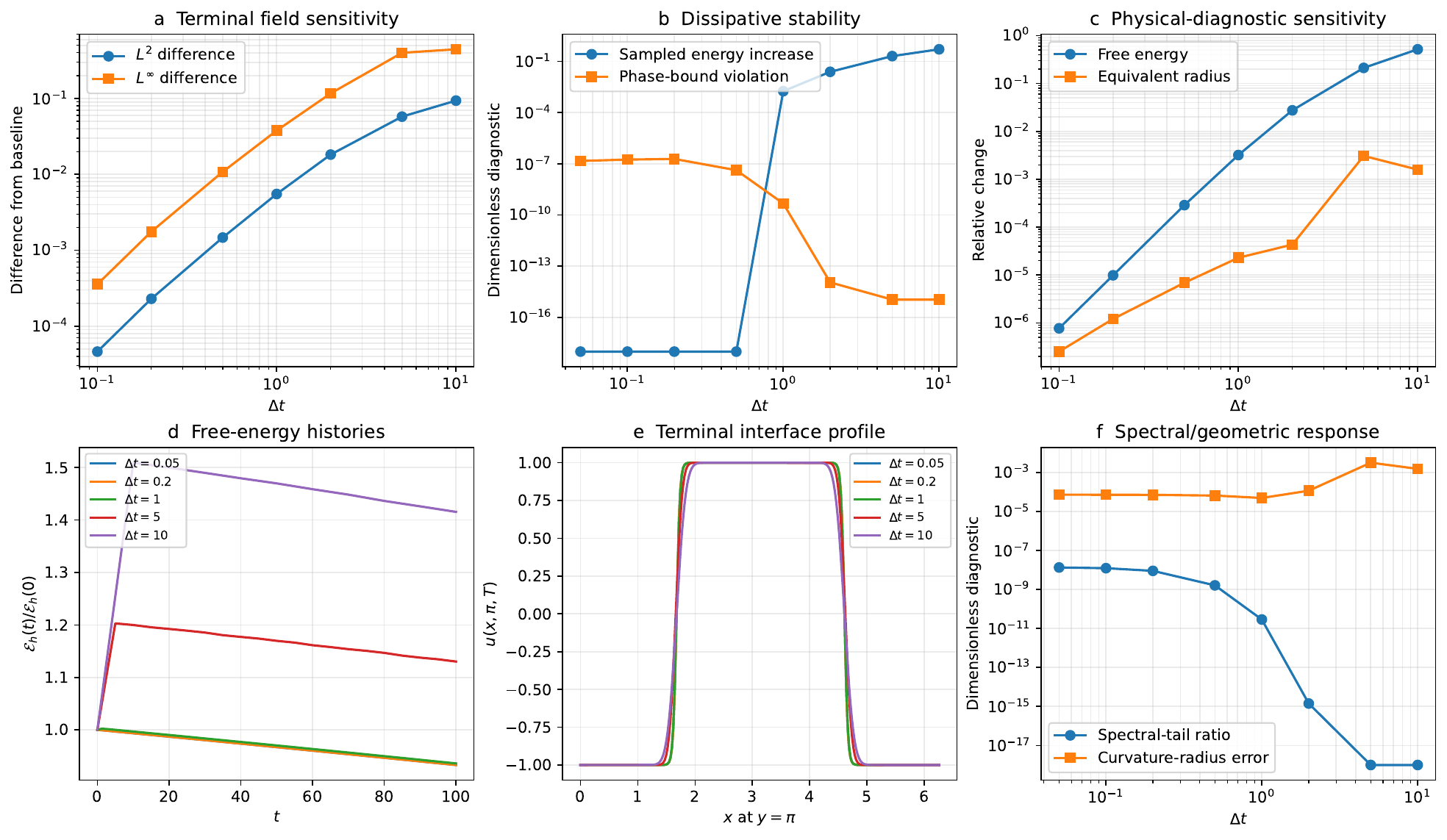}
\caption{\textbf{Large-time-step robustness of the two-dimensional Allen--Cahn calculation.}
\textbf{a}, Terminal relative \(L^2\) and maximum pointwise differences with respect to the \(\Delta t=0.05\) baseline.
\textbf{b}, Maximum sampled relative free-energy increase and phase-bound violation as functions of the time step.
\textbf{c}, Relative changes in the terminal free energy and equivalent interface radius with respect to the baseline calculation.
\textbf{d}, Normalized free-energy histories for representative time steps, showing the onset and progressive growth of non-monotonic energy evolution for excessively large time increments.
\textbf{e}, Terminal centerline phase-field profiles along \(y=\pi\).
\textbf{f}, Terminal spectral-tail ratio and relative curvature-radius error.
All calculations use \(N_x=N_y=256\), \(\epsilon=0.04\) and \(T=100\).}
\label{fig:supp_allen_cahn_robustness}
\end{figure}

Figure~\ref{fig:supp_allen_cahn_robustness}a shows that the terminal field remains close to the baseline over a substantially wider range of time steps than that required for high-accuracy calculations. At \(\Delta t=0.5\), the relative \(L^2\) difference from the baseline is only \(1.471\times10^{-3}\), while the maximum pointwise difference is \(1.081\times10^{-2}\). Increasing the time step further produces a progressive loss of field accuracy: the relative \(L^2\) difference increases to \(5.516\times10^{-3}\), \(1.824\times10^{-2}\), \(5.768\times10^{-2}\) and \(9.378\times10^{-2}\) for \(\Delta t=1\), \(2\), \(5\) and \(10\), respectively.

A more restrictive robustness criterion is provided by the free-energy evolution. For
\begin{equation}
\Delta t\leq0.5,
\end{equation}
no sampled increase in the discrete free energy is detected, indicating that the dissipative character of the Allen--Cahn dynamics is retained throughout the simulated interval. At \(\Delta t=1\), however, a nonzero relative energy increase of \(1.791\times10^{-3}\) first appears. This quantity subsequently grows to \(2.442\times10^{-2}\), \(2.029\times10^{-1}\) and \(5.094\times10^{-1}\) for \(\Delta t=2\), \(5\) and \(10\), respectively. The free-energy histories in Fig.~\ref{fig:supp_allen_cahn_robustness}d therefore identify the transition from a dissipative regime to an over-aggressive time-stepping regime more clearly than solution boundedness alone.

Importantly, loss of discrete energy dissipation is not accompanied by numerical blow-up. All tested calculations remain finite, and the phase-field values remain confined to the nominal interval \(u\in[-1,1]\) to high accuracy. Indeed, the measured phase-bound violation decreases from approximately \(10^{-7}\) for the smaller time steps to the level of round-off for the largest tested increments. Consequently, preservation of the phase bounds alone is not a sufficient indicator of physical fidelity for large time steps.

The principal geometric diagnostic is considerably less sensitive than the complete phase field. Relative to the \(\Delta t=0.05\) baseline, the terminal equivalent-radius variation remains below \(4.4\times10^{-5}\) for \(\Delta t\leq2\). Even for \(\Delta t=5\) and \(10\), the corresponding changes are only \(3.072\times10^{-3}\) and \(1.588\times10^{-3}\), respectively. In contrast, the terminal free-energy difference increases much more rapidly, reaching \(2.115\times10^{-1}\) at \(\Delta t=5\) and \(5.174\times10^{-1}\) at \(\Delta t=10\). This separation demonstrates that apparently reasonable interface locations can persist after the underlying energetic evolution has already become inaccurate.

The high-wave-number spectral-tail ratio also decreases as \(\Delta t\) is enlarged, reaching approximately machine precision for the largest time steps. This reduction should not be interpreted as improved accuracy. Instead, it reflects increasingly strong effective damping of the high-frequency content, while the field and free-energy differences simultaneously increase. Spectral-tail decay must therefore be interpreted together with physical and solution-based diagnostics rather than used as an isolated robustness measure.

Overall, the Allen--Cahn results reveal three distinct regimes. For \(\Delta t\leq0.5\), the solution remains close to the baseline while retaining monotonic free-energy dissipation. Near \(\Delta t=1\), the first measurable loss of discrete dissipativity appears despite continued phase boundedness and modest field error. For \(\Delta t\geq2\), the calculations remain numerically bounded but exhibit progressively larger energetic and field-level deviations. These results show that the exact subflows used in the Strang composition provide substantial large-step numerical robustness, while also defining a practical distinction between numerical boundedness and physically reliable dissipative evolution.

\subsection{Embedded-boundary and domain robustness of the cylinder flow}
\label{supp:cylinder_robustness}

The robustness of the embedded-boundary treatment was further assessed using the nonlinear $Re=200$ isolated-cylinder calculation. Because the spatial-resolution dependence has already been examined using $N_D=24$, $32$ and $40$, the present study instead isolates three parameters specific to the embedded Fourier formulation: the Brinkman penalty parameter $\eta$, the compact-mask transition width $\varepsilon_{\chi}$ and the transverse periodic-domain height $L_y$.

All robustness calculations use $N_D=32$ together with the production time interval $\Delta t=0.01$, $T=100$ and $50\leq t\leq100$ defined in Eqs.~\eqref{eq:supp_cyl_time_parameters} and~\eqref{eq:supp_cyl_statistics_interval}. The reference configuration uses the production values $\eta=0.0025$ and $\varepsilon_{\chi}=h$ from Eqs.~\eqref{eq:supp_cyl_penalty} and~\eqref{eq:supp_cyl_mask_width}, together with $L_y=20D$ from Eq.~\eqref{eq:supp_cyl_domain}. Only one parameter is varied at a time, with all remaining numerical and physical settings held fixed.

The penalty-parameter study considers
\begin{equation}
\eta
\in
\left\{
0.01,\,
0.005,\,
0.0025
\right\}.
\end{equation}
As the penalty is strengthened, the mean drag decreases systematically from
\begin{equation}
\overline{C_D}
=
1.53104
\end{equation}
at $\eta=0.01$ to $1.49480$ at $\eta=0.005$ and $1.47520$ at $\eta=0.0025$. Relative to the final setting, the corresponding differences are approximately $3.78\%$ and $1.33\%$ for $\eta=0.01$ and $0.005$, respectively. The solid-region velocity residual simultaneously decreases from $7.17\times10^{-2}$ to $2.44\times10^{-2}$, consistent with increasingly strong enforcement of the stationary-cylinder constraint.

The dominant shedding frequency is considerably less sensitive to the penalty parameter. The corresponding Strouhal numbers are
\begin{equation}
St
=
0.19972,\quad
0.19923,\quad
0.19884
\end{equation}
for $\eta=0.01$, $0.005$ and $0.0025$, respectively. The maximum variation relative to the final setting is therefore below $0.45\%$. This separation indicates that the local force magnitude is more sensitive to the strength of the diffuse no-slip enforcement than the large-scale vortex-shedding frequency.

The compact-mask sensitivity is examined using
\begin{equation}
\frac{\varepsilon_{\chi}}{h}
\in
\left\{
0.75,\,
1.0,\,
1.5
\right\}.
\end{equation}
Reducing the transition half-width from $1.0h$ to $0.75h$ changes the mean drag from $1.47520$ to $1.47205$, corresponding to a relative difference of only approximately $0.21\%$. The corresponding Strouhal number changes from $0.19884$ to $0.19959$.

In contrast, broadening the diffuse transition to $1.5h$ produces a more noticeable modification of the effective embedded boundary. The mean drag increases to
\begin{equation}
\overline{C_D}=1.53871,
\end{equation}
approximately $4.30\%$ above the reference value, while the Strouhal number decreases to
\begin{equation}
St=0.19611.
\end{equation}
The corresponding lift root-mean-square value also increases from $0.54654$ to $0.58996$. These results show that an excessively broad transition layer can alter the effective hydrodynamic geometry even though the nominal cylinder diameter is unchanged. The adopted choice $\varepsilon_{\chi}=h$ therefore provides a compact boundary representation without introducing the stronger geometric sensitivity observed for the wider mask.

Sensitivity to transverse periodic-image separation is examined using
\begin{equation}
\frac{L_y}{D}
\in
\left\{
16,\,
20,\,
24
\right\},
\end{equation}
while the streamwise length $L_x=40D$, cylinder streamwise position, fringe width and all other numerical parameters remain unchanged. The mean drag values are
\begin{equation}
\overline{C_D}
=
1.48096,\quad
1.47520,\quad
1.47279
\end{equation}
for $L_y/D=16$, $20$ and $24$, respectively. Relative to the $20D$ domain, the changes are approximately $0.39\%$ and $0.16\%$.

The corresponding Strouhal numbers,
\begin{equation}
St
=
0.19923,\quad
0.19884,\quad
0.19867,
\end{equation}
vary by less than $0.2\%$. The weak dependence of both the mean force and shedding frequency on transverse-domain height confirms that the adopted $20D$ separation between periodic cylinder images is sufficient for the reported wake statistics.

Throughout all tested configurations, the Fourier Helmholtz projection maintains the incompressibility residual at approximately $10^{-13}$ or below, and the maximum Courant number remains approximately $0.61$--$0.62$. The parameter variations therefore modify the embedded-boundary representation and associated physical observables without degrading the divergence-free character of the numerical solution.

\begin{figure}[h]
\centering
\includegraphics[width=\textwidth]{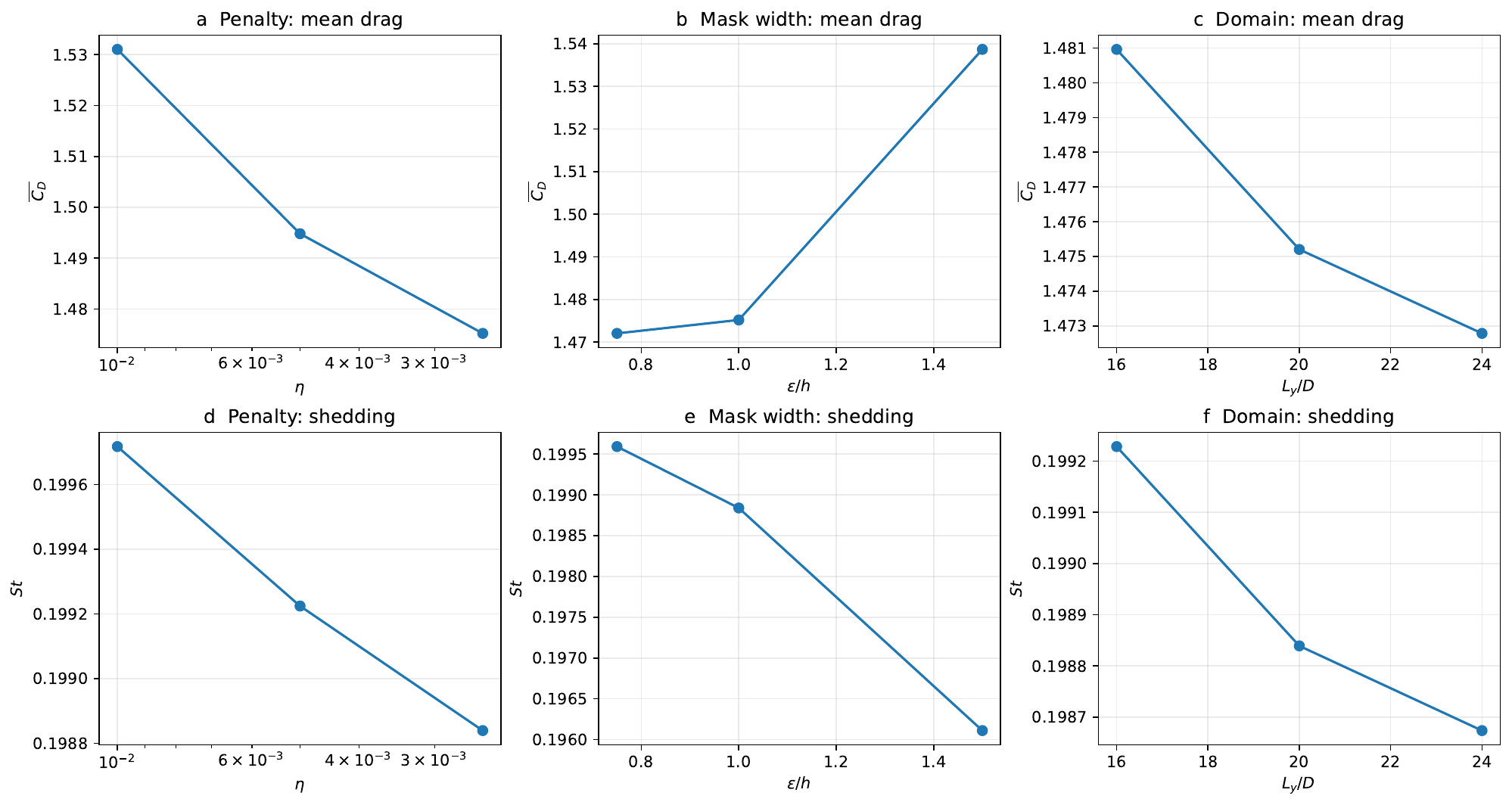}
\caption{\textbf{Embedded-boundary and computational-domain robustness of the $Re=200$ isolated-cylinder calculation.}
\textbf{a}, Mean drag coefficient as a function of the Brinkman penalty parameter $\eta$.
\textbf{b}, Mean drag coefficient as a function of the compact-mask transition half-width $\varepsilon_{\chi}/h$.
\textbf{c}, Mean drag coefficient as a function of the transverse periodic-domain height $L_y/D$.
\textbf{d--f}, Corresponding Strouhal numbers obtained from the dominant vortex-shedding frequency.
The reference configuration uses $N_D=32$, $\eta=0.0025$, $\varepsilon_{\chi}/h=1.0$ and $L_y/D=20$.}
\label{fig:supp_cylinder_robustness}
\end{figure}

Figure~\ref{fig:supp_cylinder_robustness} summarizes these results. The penalty study demonstrates progressive strengthening of the no-slip enforcement, whereas the mask-width study identifies the increased geometric influence associated with an excessively broad diffuse interface. In contrast, the weak transverse-domain dependence confirms that periodic-image interactions are already small for the production domain. Together, these controlled studies support the production choices $\eta=0.0025$, $\varepsilon_{\chi}=h$ and $L_y=20D$ specified in Section~\ref{supp:navier_stokes_cylinder}.

\section{Computational performance and hardware scaling}
\label{supp:performance}

This section evaluates the measured computational cost and hardware scaling of \OptiXDE{} after the numerical accuracy and robustness of the corresponding algorithms have been established in Sections~\ref{supp:additional_benchmarks} and~\ref{supp:accuracy_robustness}. The asymptotic operation counts and storage requirements are derived in Section~\ref{supp:computational_complexity} and are not repeated here. The present section therefore focuses on measured execution time, CPU--GPU scaling, application-level cost, de-aliasing overhead and peak memory usage.

\subsection{Timing protocol}
\label{supp:timing_protocol}

All performance measurements use the software, hardware and arithmetic precision reported in Section~\ref{supp:computational_environment}. File input, visualization and disk output are excluded from the timed regions. Reusable wave-number arrays, propagation operators, masks and other static quantities are constructed before the repeated propagation loop unless preprocessing time is explicitly reported. CPU calculations use the prescribed thread count for each comparison, whereas GPU calculations keep the evolving fields and reusable operators resident on the device whenever device-resident timing is reported.

Each configuration is preceded by unrecorded warm-up executions and is then repeated independently. The representative wall time is the median of the recorded measurements, with the interquartile range retained as a measure of timing variability. GPU timings are obtained with explicit device synchronization before and after the measured region.

For a transient calculation, the measured execution time is separated as
\begin{equation}
t_{\mathrm{total}}
=
t_{\mathrm{pre}}
+
t_{\mathrm{loop}}
+
t_{\mathrm{diag}},
\label{eq:supp_performance_time_decomposition}
\end{equation}
where \(t_{\mathrm{pre}}\) denotes one-time preprocessing, \(t_{\mathrm{loop}}\) contains all numerical time-stepping operations and \(t_{\mathrm{diag}}\) contains optional physical diagnostics. The mean per-step execution time is
\begin{equation}
t_{\mathrm{step}}
=
\frac{t_{\mathrm{loop}}}{N_t}.
\label{eq:supp_performance_step_time}
\end{equation}

For accelerator measurements, device-resident execution is used as the principal metric for repeated propagation because this corresponds to the normal execution mode of a time-dependent GPU simulation. A second transfer-inclusive metric explicitly measures one host-to-device transfer, one numerical propagation step and one device-to-host transfer, thereby identifying the problem size at which GPU acceleration remains beneficial when data movement cannot be amortized over repeated device-resident updates.

\subsection{Transform-kernel scaling and hardware acceleration}
\label{supp:kernel_scaling}

The two-dimensional diffusion problem is used to isolate the fundamental transform-based propagation kernel because its update consists only of a forward transform, multiplication by a cached spectral propagator and an inverse transform. The numerical solution itself has already been verified in Section~\ref{supp:diffusion}; the present calculation therefore measures the computational behaviour of the same validated \OptiXDE{} update over grids ranging from \(64^2\) to \(4096^2\).

For a spatial discretization containing
\begin{equation}
N=N_xN_y,
\label{eq:supp_performance_total_points}
\end{equation}
the expected per-step complexity is \(\mathcal{O}(N\log N)\), as derived in Section~\ref{supp:computational_complexity}. The measured quantity
\begin{equation}
\widetilde{t}
=
\frac{t_{\mathrm{step}}}{N\log_2N}
\label{eq:supp_performance_normalized_time}
\end{equation}
is therefore used to distinguish the transform-dominated regime from fixed execution overheads.

\begin{figure}[t]
\centering
\suppfigure{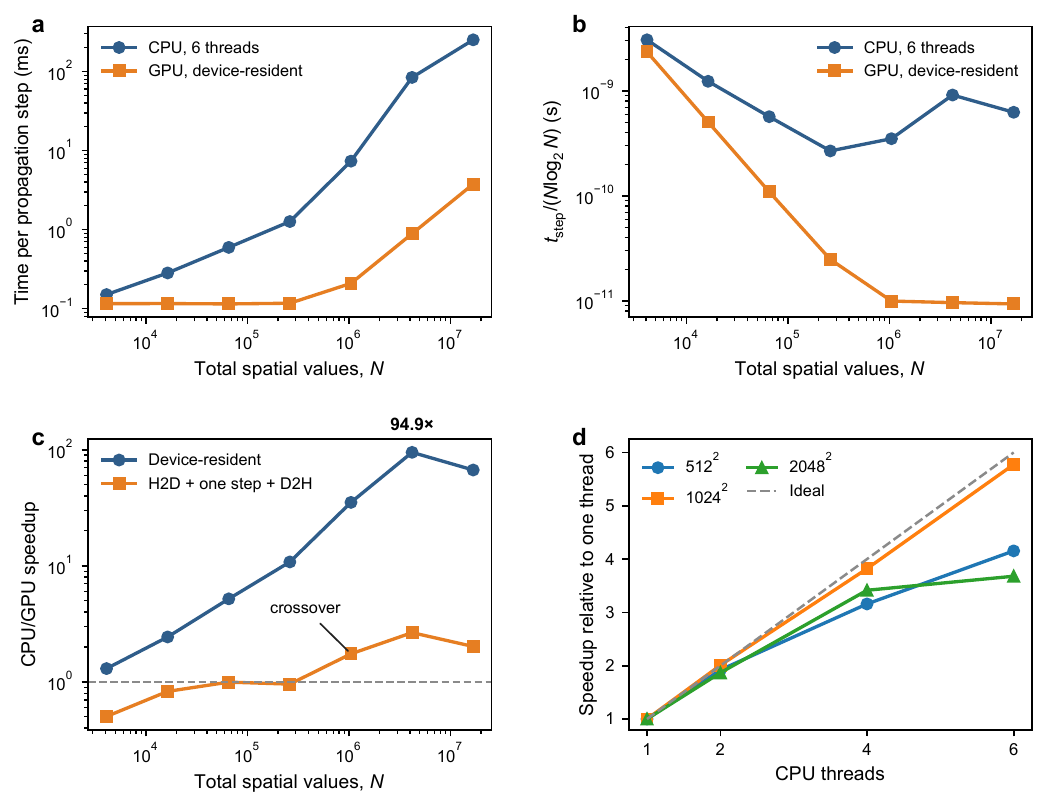}{\textwidth}
\caption{\textbf{Transform-kernel and hardware scaling of \OptiXDE{}.} \textbf{a}, CPU and device-resident GPU execution time per propagation step for the two-dimensional diffusion benchmark as a function of the total number of spatial values \(N\). \textbf{b}, Normalized execution time \(t_{\mathrm{step}}/(N\log_2N)\), showing the transition from fixed GPU execution overhead at small problem sizes toward the transform-dominated regime at larger resolutions. \textbf{c}, CPU/GPU speedup for device-resident execution and for a transfer-inclusive operation consisting of one host-to-device transfer, one \OptiXDE{} propagation step and one device-to-host transfer. The maximum measured device-resident acceleration is \(94.9\times\), whereas transfer costs delay the practical crossover to larger grids. \textbf{d}, Shared-memory CPU scaling for \(512^2\), \(1024^2\) and \(2048^2\) grids using one, two, four and six threads; the dashed line denotes ideal linear speedup.}
\label{fig:supp_performance_kernel_scaling}
\end{figure}

GPU acceleration is quantified by
\begin{equation}
S_{\mathrm{GPU}}
=
\frac{t_{\mathrm{CPU}}}{t_{\mathrm{GPU}}},
\label{eq:supp_performance_gpu_speedup}
\end{equation}
where CPU and GPU calculations use identical spatial resolutions and numerical operations. The measured device-resident GPU advantage increases strongly with problem size (Fig.~\ref{fig:supp_performance_kernel_scaling}a,c). At the smallest grids, the GPU time approaches an approximately resolution-independent floor because kernel-launch and dispatch costs are comparable to the transform cost. Once the arrays become sufficiently large, this fixed overhead is amortized and the acceleration rises rapidly, reaching a maximum measured speedup of \(94.9\times\) at \(2048^2\), where the median CPU and GPU propagation times are \(84.33\) and \(0.889\) ms per step, respectively. At \(4096^2\), the corresponding speedup decreases to \(66.9\times\), while remaining strongly favourable to GPU execution. This reduction shows that the precise hardware acceleration is governed by FFT workload, memory traffic and backend-specific execution characteristics rather than by the asymptotic operation count alone.

The normalized timing in Fig.~\ref{fig:supp_performance_kernel_scaling}b further separates these regimes. The pronounced decrease of the GPU-normalized cost at small and intermediate resolutions reflects progressive amortization of fixed accelerator overhead, whereas the large-grid calculations enter a transform-dominated regime consistent with the \(\mathcal{O}(N\log N)\) complexity derived analytically. The CPU curve exhibits substantially weaker fixed-overhead effects because the transforms are already large relative to function-dispatch costs over most of the tested range.

Including host--device transfers alters the practical crossover without changing the device-resident kernel performance. For the smallest grids, transferring one field to the GPU and returning the result costs more than the accelerated propagation itself, and the transfer-inclusive GPU calculation is therefore slower than the CPU calculation (Fig.~\ref{fig:supp_performance_kernel_scaling}c). At \(256^2\), CPU and transfer-inclusive GPU execution are essentially equal, whereas a clear and sustained transfer-inclusive advantage emerges at \(1024^2\) and above. At \(2048^2\), the transfer-inclusive speedup is \(2.66\times\), compared with \(94.9\times\) for device-resident execution. This distinction shows that data residency is critical for small and moderate workloads, whereas sufficiently large transform problems retain a net GPU advantage even when a complete host--device round trip is included.

Shared-memory CPU scaling is shown in Fig.~\ref{fig:supp_performance_kernel_scaling}d. Increasing the thread count from one to six reduces the propagation time at all three representative resolutions, although the efficiency depends on problem size and departs from ideal linear scaling as memory traffic, FFT threading overhead and shared-resource contention become important. The \(1024^2\) case approaches the ideal trend most closely, whereas the \(2048^2\) calculation shows stronger saturation at the highest thread count. These results confirm that the same transform-based kernel benefits from both shared-memory CPU parallelism and massively parallel GPU execution, with the relative benefit determined by workload size.

\subsection{Application-level performance}
\label{supp:application_performance}

\begin{figure}[t]
\centering
\suppfigure{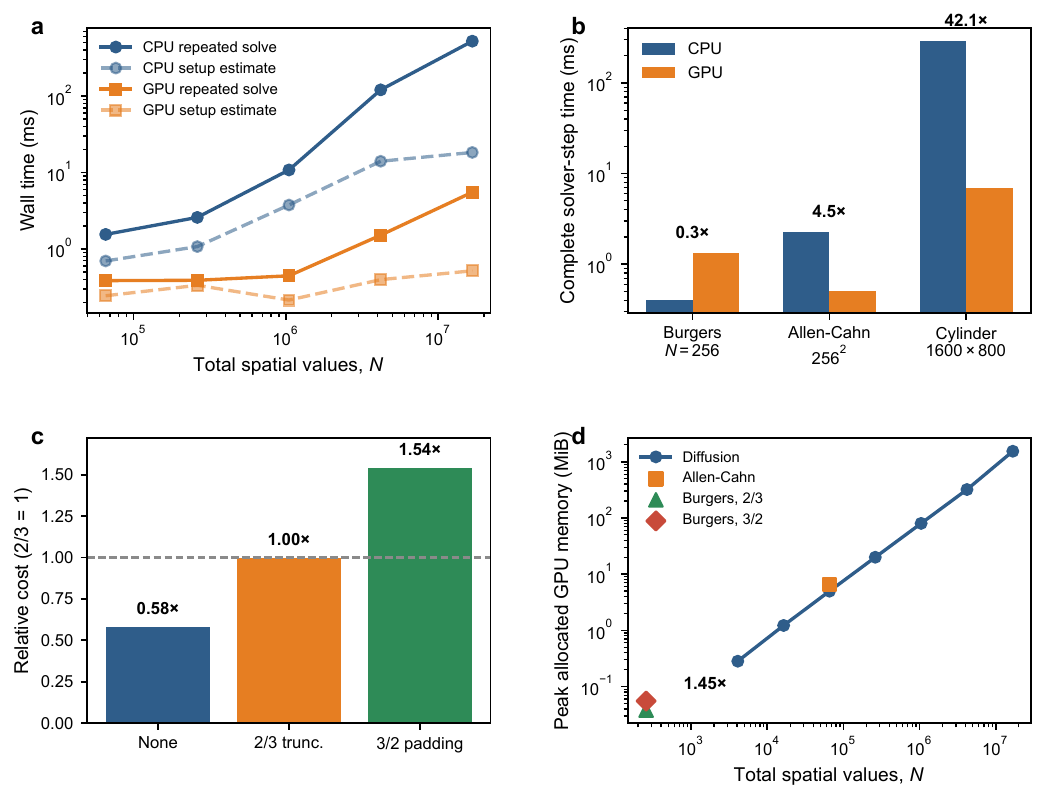}{\textwidth}
\caption{\textbf{Application-level computational performance of \OptiXDE{}.} \textbf{a}, Setup estimates and repeated-solve wall times for the periodic Poisson problem on CPU and GPU, separating one-time operator construction from repeated spectral inversion. \textbf{b}, Complete solver-step times for representative nonlinear calculations using matched numerical configurations. The small one-dimensional Burgers calculation remains dominated by accelerator overhead, whereas the \(256^2\) Allen--Cahn and \(1600\times800\) cylinder calculations obtain approximately \(4.5\times\) and \(42.1\times\) GPU acceleration, respectively. \textbf{c}, Relative computational cost of the Burgers nonlinear-product treatments, normalized by the two-thirds truncation result. Three-halves zero padding increases the measured cost to \(1.54\times\). \textbf{d}, Measured peak allocated GPU memory for representative linear and nonlinear calculations. Diffusion provides the baseline memory scaling, while the Burgers comparison shows that three-halves padding increases the peak allocated memory by approximately \(1.45\times\) relative to two-thirds truncation.}
\label{fig:supp_performance_application}
\end{figure}

Kernel timings do not include the additional operations required by steady inversion, nonlinear splitting, spectral differentiation, de-aliasing or embedded-boundary enforcement. Representative application-level measurements are therefore reported for the periodic Poisson, Burgers, Allen--Cahn and isolated-cylinder calculations whose numerical accuracy has already been established in the preceding sections (Fig.~\ref{fig:supp_performance_application}). The principal CPU--GPU measurements are also summarized numerically in Table~\ref{tab:supp_performance_summary}.

The periodic Poisson problem represents direct steady spectral inversion. Its setup contribution and repeated-solve time are separated in Fig.~\ref{fig:supp_performance_application}a to distinguish one-time construction and caching of the inverse spectral operator from repeated right-hand-side evaluation. The repeated CPU and GPU solve times both increase systematically with problem size, while the GPU curve remains substantially below the CPU curve at large \(N\). At \(4096^2\), the median cached solve time decreases from \(522.46\) ms on the CPU to \(5.522\) ms on the GPU, corresponding to a \(94.6\times\) acceleration. The estimated one-time setup costs at the same resolution are \(18.30\) ms on the CPU and \(0.517\) ms on the GPU, substantially smaller than the corresponding repeated-solve cost. These results demonstrate the practical benefit of caching transform-domain operators when multiple solves are performed with the same spatial discretization.

The nonlinear application measurements expose a clear workload-dependent CPU--GPU crossover (Fig.~\ref{fig:supp_performance_application}b). For the one-dimensional Burgers benchmark at \(N=256\) with two-thirds de-aliasing, the median complete solver-step time is \(0.402\) ms on the CPU and \(1.343\) ms on the GPU, giving a GPU speedup of only \(0.30\times\). The workload per update is therefore too small to amortize accelerator launch and dispatch overhead. For the two-dimensional Allen--Cahn calculation at \(256^2\), the corresponding times are \(2.267\) and \(0.504\) ms per step, producing a \(4.50\times\) GPU acceleration. The effect becomes much stronger for the matched \(1600\times800\) isolated-cylinder benchmark, for which CPU and GPU execute the same numerical operations through the same \OptiXDE{} Torch backend family. The median step time decreases from \(289.77\) ms on six CPU threads to \(6.886\) ms on the GPU, corresponding to a \(42.1\times\) acceleration. The CPU--GPU terminal relative \(L_2\) difference remains \(1.19\times10^{-13}\), confirming numerical consistency between the two hardware paths. The progression from Burgers to Allen--Cahn and then to the cylinder calculation therefore shows that accelerator efficiency is controlled primarily by the amount of parallel spectral work available per update rather than by the use of a GPU alone.

The Burgers equation additionally quantifies the computational cost of nonlinear de-aliasing under an otherwise identical Strang--Runge--Kutta integration procedure. Taking two-thirds truncation as the unit-cost reference, omitting de-aliasing reduces the measured time to \(0.58\times\), whereas three-halves zero padding increases it to \(1.54\times\) (Fig.~\ref{fig:supp_performance_application}c). Combined with the accuracy study in Section~\ref{supp:burgers_robustness}, this result makes the numerical trade-off explicit: three-halves padding incurs approximately \(54\%\) more computational cost than two-thirds truncation but avoids the loss of resolved nonlinear content associated with simple high-wave-number truncation.

\begin{table}[t]
\centering
\small
\caption{Representative computational-performance measurements for \OptiXDE{}. Diffusion timings correspond to device-resident propagation, Poisson timings to repeated solves with the inverse operator cached, and Burgers timings to the two-thirds de-aliased public-solver configuration. The Allen--Cahn and cylinder entries report complete device-resident solver-step times. CPU and GPU calculations use identical spatial resolutions and matched numerical configurations.}
\label{tab:supp_performance_summary}
\begin{tabular}{lcccc}
\toprule
Benchmark & Resolution & CPU time & GPU time & CPU/GPU speedup \\
\midrule
Diffusion propagation & \(2048^2\) & \(84.328\) ms/step & \(0.889\) ms/step & \(94.9\times\) \\
Periodic Poisson & \(4096^2\) & \(522.459\) ms/solve & \(5.522\) ms/solve & \(94.6\times\) \\
Burgers, two-thirds & \(N=256\) & \(0.402\) ms/step & \(1.343\) ms/step & \(0.30\times\) \\
Allen--Cahn & \(256^2\) & \(2.267\) ms/step & \(0.504\) ms/step & \(4.50\times\) \\
Cylinder, \(Re=200\) & \(1600\times800\) & \(289.765\) ms/step & \(6.886\) ms/step & \(42.1\times\) \\
\bottomrule
\end{tabular}
\end{table}

\subsection{Memory scaling}
\label{supp:memory_scaling}

The theoretical storage requirement of \OptiXDE{} is \(\mathcal{O}(N)\), as derived in Section~\ref{supp:computational_complexity}. The measured peak allocated GPU memory is used here to assess the practical scaling of the implemented algorithms and to quantify the additional prefactors introduced by nonlinear work arrays and de-aliasing. Allocated rather than backend-reserved memory is used for the principal comparison because the latter is affected by allocator caching and block-reservation policies.

For the linear diffusion kernel, the measured peak allocated memory increases systematically with the total number of spatial values over the complete resolution sequence (Fig.~\ref{fig:supp_performance_application}d). It grows from approximately \(0.28\) MiB at \(64^2\) to approximately \(1.50\) GiB at \(4096^2\), consistent with the linear storage requirement of the persistent physical and spectral fields, cached propagator and FFT work arrays. The Allen--Cahn calculation at \(256^2\) requires approximately \(6.5\) MiB of peak allocated GPU memory and lies on the same practical scale while additionally containing temporary arrays associated with nonlinear splitting. The Burgers cases are much smaller one-dimensional calculations and therefore require substantially less absolute memory, but they isolate the additional cost of nonlinear-product treatment.

At \(N=256\), the measured peak allocated memory is \(0.0381\) MiB for two-thirds truncation and \(0.0552\) MiB for three-halves padding. The latter is therefore \(1.45\times\) the two-thirds value (Fig.~\ref{fig:supp_performance_application}d). This increase results from the temporary padded spectral and physical arrays required to evaluate the quadratic convolution on the enlarged \(3N/2\) grid. It changes the constant prefactor of the storage requirement but not its asymptotic dependence on the number of resolved spatial values.

Taken together, Figs.~\ref{fig:supp_performance_kernel_scaling} and~\ref{fig:supp_performance_application} show three complementary aspects of computational performance. First, the large-grid propagation measurements enter a transform-dominated regime consistent with the \(\mathcal{O}(N\log N)\) arithmetic complexity derived in Section~\ref{supp:computational_complexity}. Second, GPU acceleration is strongly workload dependent, ranging from no benefit for the small one-dimensional Burgers calculation to \(42.1\times\) for the matched embedded-cylinder update and \(94.9\times\) for the large diffusion propagation kernel. Third, the measured memory behaviour is consistent with linear storage, while nonlinear splitting and three-halves padding introduce predictable constant-factor overheads.

\section{Algorithmic details}
\label{supp:algorithmic_details}

This section records the execution order used in the benchmark calculations. The transform conventions, wave-number construction, closed-form modal operators, boundary representations, geometry masks, penalty formulations, nonlinear subflows, complexity estimates and diagnostic definitions have already been established in Section~\ref{supp:numerical_implementation} and are not re-derived here. The purpose of the present section is therefore limited to the assembly of those previously defined components into reproducible computational procedures.

\subsection{Preprocessing and operator caching}
\label{supp:preprocessing}

Preprocessing is performed once for every fixed combination of spatial resolution, domain dimensions, transform family, physical parameters, propagation interval, arithmetic precision and computational backend. The coordinate arrays and transform-compatible wave numbers are constructed according to Section~\ref{supp:discrete_transforms}. The modal eigenvalues, propagation factors, source-integration factors and inverse multipliers are then generated from the definitions in Section~\ref{supp:spectral_propagators}. When an embedded geometry is present, the signed-distance field, binary physical-domain mask, smoothed interior and exterior masks, interface-localization weights, boundary partitions, normals and prescribed-data extensions are constructed according to Section~\ref{supp:geometry_masks}.

Only arrays that remain unchanged during the subsequent propagation or correction loop are retained in the persistent cache. These may include the transform metadata, componentwise wave numbers, squared wave-number magnitude, differentiation multipliers, full-step and half-step propagators, Helmholtz and Poisson inverse factors, source-integration factors, de-aliasing operators, geometry masks, boundary labels and prescribed solid or boundary fields. Time-dependent source terms, moving geometries and time-dependent boundary extensions are updated outside the persistent cache at the stage where they are required.

A cached object is invalidated whenever any quantity affecting its numerical value changes. In particular, changing the spatial resolution, domain length, transform type, physical coefficient, propagation interval, precision or device requires reconstruction of the affected operator. Changing only the evolving solution field does not trigger reconstruction.

\begin{algorithm}[t]
\caption{Preprocessing and operator caching}
\label{alg:supp_preprocessing}
\begin{algorithmic}[1]
\Require Transform domain $\widetilde{\Omega}$, physical domain $\Omega$, spatial resolutions, transform family, PDE parameters, propagation or pseudo-time increments, backend, precision and optional geometry description
\Ensure Persistent spatial, spectral, geometric and operator cache
\State Select the CPU or GPU backend and arithmetic precision
\State Construct the coordinates, spacings and transform metadata
\State Construct the wave-number arrays and spectral differentiation multipliers
\State Construct the PDE-specific modal eigenvalues, propagators and inverse factors
\If{a source-integration factor is independent of time}
    \State Construct and cache the corresponding modal source factor
\EndIf
\If{a nonlinear pseudospectral product is present}
    \State Construct and cache the selected truncation mask or padding map
\EndIf
\If{$\Omega\neq\widetilde{\Omega}$}
    \State Construct the signed-distance or level-set field
    \State Construct the physical-domain, exterior, interface and segment masks
    \State Construct the boundary normals and prescribed-data extensions
\EndIf
\State Transfer all persistent arrays to the selected computational device
\State Return the cache
\end{algorithmic}
\end{algorithm}

The preprocessing and storage estimates associated with these arrays are given in Section~\ref{supp:computational_complexity}; they are not repeated here.

\subsection{Linear propagation algorithm}
\label{supp:linear_algorithm}

The linear solver uses the modal update derived in Section~\ref{supp:spectral_propagators}. A transient homogeneous step applies the cached propagation factor to the transformed state. When a source term is present, the transformed source is combined with the cached source-integration factor from Eq.~\eqref{eq:supp_cached_source_factor}. A time-independent source is transformed once during preprocessing, whereas a time-dependent source is evaluated and transformed at the temporal location prescribed by the benchmark-specific source approximation.

For a steady transform-diagonal problem, the right-hand side is transformed once and multiplied by the cached inverse operator. Any compatibility condition and null-space normalization are imposed before the inverse transform. In particular, periodic and pure-Neumann Poisson problems use the compatibility and zero-mode treatments defined in Eqs.~\eqref{eq:supp_poisson_compatibility_spectral}--\eqref{eq:supp_poisson_zero_mean_selection} and Eqs.~\eqref{eq:supp_neumann_compatibility}--\eqref{eq:supp_neumann_zero_mode_constraint}, respectively.

\begin{algorithm}[t]
\caption{Linear spectral propagation or direct inversion}
\label{alg:supp_linear_propagation}
\begin{algorithmic}[1]
\Require Initial field $u^0$ or right-hand side $f$, cached transform pair $(\mathcal{T},\mathcal{T}^{-1})$, cached modal operator, optional source specification and requested output times
\Ensure Transient field history or steady solution
\If{the problem is steady and transform diagonal}
    \State $\widehat{f}\gets\mathcal{T}[f]$
    \State Enforce the required compatibility condition
    \State Apply the cached inverse factor to every admissible mode
    \State Set each null-space coefficient according to the prescribed normalization
    \State $u\gets\mathcal{T}^{-1}[\widehat{u}]$
    \State Apply the lifting recovery or physical-domain restriction, when required
\Else
    \State Initialize $u\gets u^0$
    \For{each propagation interval $[t_n,t_{n+1}]$}
        \State $\widehat{u}\gets\mathcal{T}[u]$
        \If{a source term is absent}
            \State $\widehat{u}\gets G\odot\widehat{u}$
        \Else
            \State Evaluate the benchmark-prescribed approximation of the source over $[t_n,t_{n+1}]$
            \State Retrieve or compute its transformed coefficients $\widehat{s}^{\,n}$
            \State $\widehat{u}\gets G\odot\widehat{u}+H\odot\widehat{s}^{\,n}$
        \EndIf
        \State $u\gets\mathcal{T}^{-1}[\widehat{u}]$
        \State Recover any nonhomogeneous lifting field or apply the prescribed physical-space constraint
        \State Discard only round-off-level imaginary components for real-valued problems
        \State Evaluate and store the requested diagnostics
    \EndFor
\EndIf
\State Return the solution on the physical domain
\end{algorithmic}
\end{algorithm}

Here, $G$ denotes the cached homogeneous propagation factor from Eq.~\eqref{eq:supp_general_propagator}, and $H$ denotes the cached source-integration factor from Eq.~\eqref{eq:supp_cached_source_factor}. Complex-valued equations retain the complete complex inverse-transform result.

\subsection{Embedded-boundary algorithm}
\label{supp:embedded_algorithm}

The general propagate--enforce composition and the associated masks are defined in Sections~\ref{supp:geometry_masks} and~\ref{supp:boundary_enforcement}. For the steady embedded Poisson benchmark, this composition is implemented as a pseudo-time fixed-point iteration. The constant-coefficient diffusion contribution is advanced by a backward-Euler Helmholtz correction, followed by a pointwise backward-Euler realization of the Dirichlet penalty subproblem. This subsection specifies that benchmark implementation without redefining the underlying geometry or penalty model.

Let $\chi_D$ denote the cached Dirichlet enforcement mask and let $\widetilde{g}_{\mathrm{ext}}$ denote the prescribed-data extension. One bulk correction is
\begin{equation}
v^\star
=
\mathcal{T}^{-1}
\left\{
\mathcal{P}_{H}
\odot
\mathcal{T}
\left[
v^{(m)}
+
\Delta\tau\,\widetilde{f}
\right]
\right\},
\qquad
\mathcal{P}_{H}
=
\frac{1}{1+\Delta\tau\,\Lambda},
\label{eq:supp_algorithm_embedded_diffusion}
\end{equation}
where $\Lambda$ is the nonnegative Laplacian eigenvalue array associated with the selected transform. The local penalty correction is
\begin{equation}
v^{(m+1)}
=
\frac{
v^\star
+
\beta_D\chi_D\widetilde{g}_{\mathrm{ext}}
}{
1+\beta_D\chi_D
},
\qquad
\beta_D
=
\frac{\Delta\tau}{\eta_D}.
\label{eq:supp_algorithm_embedded_penalty}
\end{equation}
The denominator is strictly positive for $\Delta\tau>0$, $\eta_D>0$ and $\chi_D\geq0$. The correction therefore leaves the field unchanged where $\chi_D=0$ and relaxes it toward $\widetilde{g}_{\mathrm{ext}}$ inside the enforcement region.

Convergence is measured only over the physical domain using the binary mask $m_{\Omega}^{0}$ defined in Eq.~\eqref{eq:supp_binary_mask}:
\begin{equation}
r_{\mathrm{upd}}^{(m+1)}
=
\frac{
\left\|
m_{\Omega}^{0}
\left(
v^{(m+1)}-v^{(m)}
\right)
\right\|_2
}{
\left\|
m_{\Omega}^{0}v^{(m+1)}
\right\|_2
+
\delta_{\mathrm{upd}}
},
\label{eq:supp_algorithm_embedded_convergence}
\end{equation}
where $\delta_{\mathrm{upd}}>0$ prevents division by zero. When a residual-based stopping criterion is also prescribed, both the relative update and the physical-domain residual must satisfy their tolerances.

\begin{algorithm}[t]
\caption{Embedded-boundary correction for a steady penalized elliptic problem}
\label{alg:supp_embedded_boundary}
\begin{algorithmic}[1]
\Require Extended source $\widetilde{f}$, boundary-data extension $\widetilde{g}_{\mathrm{ext}}$, Dirichlet mask $\chi_D$, binary physical-domain mask $m_{\Omega}^{0}$, penalty parameter $\eta_D$, pseudo-time increment $\Delta\tau$, tolerances and maximum correction count $M_{\max}$
\Ensure Converged penalized solution restricted to $\Omega$
\State Initialize $v^{(0)}$
\State Retrieve $\mathcal{P}_{H}=(1+\Delta\tau\Lambda)^{-1}$
\State Compute $\beta_D\gets\Delta\tau/\eta_D$
\For{$m=0,\ldots,M_{\max}-1$}
    \State $b\gets v^{(m)}+\Delta\tau\,\widetilde{f}$
    \State $\widehat{b}\gets\mathcal{T}[b]$
    \State $v^\star\gets\mathcal{T}^{-1}[\mathcal{P}_{H}\odot\widehat{b}]$
    \State $v^{(m+1)}\gets\left(v^\star+\beta_D\chi_D\widetilde{g}_{\mathrm{ext}}\right)/\left(1+\beta_D\chi_D\right)$
    \State Apply any exact condition prescribed on the outer boundary of $\widetilde{\Omega}$
    \State Evaluate $r_{\mathrm{upd}}^{(m+1)}$ and any prescribed physical-domain residual
    \If{all stopping criteria are satisfied}
        \State \textbf{break}
    \EndIf
\EndFor
\State Restrict $v^{(m+1)}$ using $m_{\Omega}^{0}$
\State Return the physical-domain solution and convergence history
\end{algorithmic}
\end{algorithm}

The finite values of $\eta_D$ and the mask-smoothing width introduce penalty and diffuse-interface errors. Accordingly, the embedded-domain benchmark is assessed using the physical-domain, interface-band and boundary diagnostics defined in Section~\ref{supp:error_measures}, rather than by comparison with the round-off-level behavior of an exactly transform-diagonal periodic problem.

\subsection{Nonlinear splitting algorithm}
\label{supp:nonlinear_algorithm}

The nonlinear equations, exact nonlinear subflows, spectral derivative evaluations and de-aliasing rules are defined in Section~\ref{supp:nonlinear_splitting}. The implementation uses the second-order Strang composition in Eq.~\eqref{eq:supp_strang_splitting}. The present subsection specifies only the order in which the cached linear half-step, nonlinear flow and optional physical-space constraints are applied.

When the nonlinear subflow is available analytically, it is evaluated pointwise without iteration. This applies to the local phase rotation in the nonlinear Schr\"odinger benchmark and the exact reaction flow in the Allen--Cahn benchmark. When no exact nonlinear flow is available, the benchmark-prescribed stage integrator is used. In particular, the Burgers nonlinear subproblem is advanced with the classical fourth-order Runge--Kutta procedure described in Section~\ref{supp:nonlinear_splitting}. At every stage involving derivative-dependent products, the derivatives are evaluated spectrally, products are formed in physical space and the selected de-aliasing operation is applied before the stage increment is assembled.

\begin{algorithm}[t]
\caption{Strang-split nonlinear propagation}
\label{alg:supp_nonlinear_splitting}
\begin{algorithmic}[1]
\Require Field $u^n$, cached linear half-step factor $G_{1/2}$, nonlinear flow specification, step size $\Delta t$, optional de-aliasing operator and optional boundary-enforcement map
\Ensure Field $u^{n+1}$
\State $\widehat{u}\gets\mathcal{T}[u^n]$
\State $u^a\gets\mathcal{T}^{-1}[G_{1/2}\odot\widehat{u}]$
\If{an exact nonlinear flow is available}
    \State $u^b\gets\mathcal{Q}_{\Delta t}[u^a]$
\Else
    \State Initialize the benchmark-prescribed nonlinear stage integrator with $u^a$
    \For{each nonlinear stage}
        \State Evaluate all required spectral derivatives
        \State Form the nonlinear products or couplings in physical space
        \State Apply the selected de-aliasing operation
        \State Update the stage value using the prescribed stage coefficients
    \EndFor
    \State Set $u^b$ to the completed nonlinear-stage solution
\EndIf
\State $\widehat{u}\gets\mathcal{T}[u^b]$
\State $u^{n+1}\gets\mathcal{T}^{-1}[G_{1/2}\odot\widehat{u}]$
\If{boundary or embedded-domain enforcement is required}
    \State Apply the benchmark-prescribed symmetric enforcement map
\EndIf
\State Evaluate the requested invariants, residuals and output quantities
\State Return $u^{n+1}$
\end{algorithmic}
\end{algorithm}

This procedure deliberately refers to the benchmark-specific nonlinear subflow instead of introducing a second generic Runge--Kutta formula. It therefore preserves the exact nonlinear updates used for the nonlinear Schr\"odinger and Allen--Cahn equations and the fourth-order stage update used for the Burgers equation.

\subsection{Incompressible-flow projection algorithm}
\label{supp:flow_algorithm}

For incompressible flow, the velocity $\mathbf{u}$ and pressure $p$ satisfy
\begin{equation}
\frac{\partial\mathbf{u}}{\partial t}
+
\mathcal{N}(\mathbf{u})
=
-\nabla p
+
\nu\nabla^2\mathbf{u}
+
\mathbf{f}
-
\frac{\chi_{\mathrm{s}}}{\eta_{\mathrm{f}}}
\left(
\mathbf{u}-\mathbf{u}_{\mathrm{s}}
\right),
\qquad
\nabla\cdot\mathbf{u}
=
0,
\label{eq:supp_algorithm_incompressible_equations}
\end{equation}
where $\nu$ is the kinematic viscosity, $\mathbf{f}$ is a prescribed forcing term, $\chi_{\mathrm{s}}$ is the regularized solid-region mask, $\eta_{\mathrm{f}}>0$ is the flow-penalty parameter and $\mathbf{u}_{\mathrm{s}}$ is the prescribed solid velocity. A stationary obstacle corresponds to $\mathbf{u}_{\mathrm{s}}=\mathbf{0}$.

The incremental projection step consists of an explicit nonlinear and forcing update, an implicit spectral diffusion solve, a pointwise Brinkman correction and a pressure-increment projection. The predictor right-hand side is
\begin{equation}
\mathbf{r}^{n}
=
\mathbf{u}^{n}
+
\Delta t
\left[
-\mathcal{N}^{n+1/2}
+
\mathbf{f}^{n+1/2}
-
\nabla p^{n}
\right],
\label{eq:supp_algorithm_flow_predictor_rhs}
\end{equation}
where the temporal approximation of $\mathcal{N}^{n+1/2}$ is specified with the corresponding flow benchmark. Viscous diffusion is advanced through
\begin{equation}
\mathbf{u}^{\dagger}
=
\mathcal{T}^{-1}
\left\{
\mathcal{P}_{\nu}
\odot
\mathcal{T}
\left[
\mathbf{r}^{n}
\right]
\right\},
\qquad
\mathcal{P}_{\nu}
=
\frac{1}{1+\nu\Delta t\,\Lambda}.
\label{eq:supp_algorithm_flow_diffusion}
\end{equation}
The solid-region velocity is then imposed pointwise:
\begin{equation}
\mathbf{u}^{\star}
=
\frac{
\mathbf{u}^{\dagger}
+
\beta_{\mathrm{f}}\chi_{\mathrm{s}}\mathbf{u}_{\mathrm{s}}
}{
1+\beta_{\mathrm{f}}\chi_{\mathrm{s}}
},
\qquad
\beta_{\mathrm{f}}
=
\frac{\Delta t}{\eta_{\mathrm{f}}}.
\label{eq:supp_algorithm_flow_penalty}
\end{equation}

The pressure increment $\phi_p$ satisfies
\begin{equation}
\nabla^2\phi_p
=
\frac{1}{\Delta t}
\nabla\cdot\mathbf{u}^{\star}.
\label{eq:supp_algorithm_pressure_poisson}
\end{equation}
For a periodic transform, its nonzero Fourier modes are
\begin{equation}
\widehat{\phi}_p(\mathbf{k})
=
-
\frac{
\mathrm{i}\mathbf{k}\cdot
\widehat{\mathbf{u}}^{\star}(\mathbf{k})
}{
\Delta t\,\Lambda(\mathbf{k})
},
\qquad
\Lambda(\mathbf{k})>0,
\label{eq:supp_algorithm_pressure_increment_spectral}
\end{equation}
and the zero-frequency pressure-increment mode is set to zero. The corrected velocity and pressure are
\begin{equation}
\mathbf{u}^{n+1}
=
\mathbf{u}^{\star}
-
\Delta t\,\nabla\phi_p,
\qquad
p^{n+1}
=
p^{n}
+
\phi_p.
\label{eq:supp_algorithm_projection_correction}
\end{equation}

\begin{algorithm}[t]
\caption{Embedded incompressible-flow projection step}
\label{alg:supp_flow_projection}
\begin{algorithmic}[1]
\Require Velocity $\mathbf{u}^{n}$, pressure $p^{n}$, viscosity $\nu$, forcing $\mathbf{f}$, solid mask $\chi_{\mathrm{s}}$, solid velocity $\mathbf{u}_{\mathrm{s}}$, penalty parameter $\eta_{\mathrm{f}}$, step size $\Delta t$ and cached spectral operators
\Ensure Velocity $\mathbf{u}^{n+1}$ and pressure $p^{n+1}$
\State Apply the prescribed inlet, wall and outer-boundary extensions
\State Evaluate $\mathcal{N}^{n+1/2}$ using spectral derivatives and physical-space products
\State Apply the selected de-aliasing operation to the nonlinear term
\State Evaluate $\nabla p^n$ using the cached differentiation multipliers
\State $\mathbf{r}^{n}\gets\mathbf{u}^{n}+\Delta t[-\mathcal{N}^{n+1/2}+\mathbf{f}^{n+1/2}-\nabla p^{n}]$
\State $\widehat{\mathbf{r}}\gets\mathcal{T}[\mathbf{r}^{n}]$
\State $\mathbf{u}^{\dagger}\gets\mathcal{T}^{-1}[\mathcal{P}_{\nu}\odot\widehat{\mathbf{r}}]$
\State $\beta_{\mathrm{f}}\gets\Delta t/\eta_{\mathrm{f}}$
\State $\mathbf{u}^{\star}\gets(\mathbf{u}^{\dagger}+\beta_{\mathrm{f}}\chi_{\mathrm{s}}\mathbf{u}_{\mathrm{s}})/(1+\beta_{\mathrm{f}}\chi_{\mathrm{s}})$
\State $d^{\star}\gets\nabla\cdot\mathbf{u}^{\star}$
\State Solve $\nabla^2\phi_p=d^{\star}/\Delta t$ using the cached Poisson inverse
\State $\mathbf{u}^{n+1}\gets\mathbf{u}^{\star}-\Delta t\,\nabla\phi_p$
\State $p^{n+1}\gets p^{n}+\phi_p$
\State Apply the prescribed velocity and pressure-correction boundary treatments
\State Evaluate the fluid-region divergence, boundary mismatch and update measures
\State Return $\mathbf{u}^{n+1}$ and $p^{n+1}$
\end{algorithmic}
\end{algorithm}

The projection removes the irrotational component of the intermediate velocity, whereas the Brinkman correction suppresses relative motion in the embedded solid. Divergence and update measures are evaluated over the physical fluid region using the masks and diagnostic conventions of Sections~\ref{supp:geometry_masks} and~\ref{supp:error_measures}. The transform and pointwise-operation count of the projection step is included in the complexity analysis of Section~\ref{supp:computational_complexity}.

\bibliography{}